\documentclass{article}
\usepackage{graphicx} 
\usepackage[margin=1in]{geometry}
\usepackage{amsfonts,amsmath,amssymb,amsthm,bbm}
\usepackage{float,mathtools,mathrsfs}
\usepackage[shortlabels]{enumitem}
\usepackage{array}
\usepackage{graphicx}
\usepackage{algorithmicx}
\usepackage[ruled]{algorithm}
\usepackage{algpseudocode}
\usepackage[numbers]{natbib}
\usepackage[colorlinks]{hyperref}

\newenvironment{proof*}{\paragraph{Proof:}}{\hfill$\blacksquare$}

\newtheorem{theorem}{Theorem}
\newtheorem{remark}{Remark}
\newtheorem*{theorem*}{Theorem}
\newtheorem{lemma}[theorem]{Lemma}
\newtheorem{proposition}{Proposition}
\newtheorem{corollary}[theorem]{Corollary}
\newtheorem{definition}{Definition}

\def\E{\mathbb{E}}

\def\Var{\mathrm{Var}}

\makeatletter
\newcommand{\mylabel}[1]{%
    \protected@edef\@currentlabel{#1}
    }
\makeatother

\usepackage{amsfonts,amsmath,amssymb,amsthm,bbm}
\usepackage{float,mathtools}
\usepackage{changepage}
\usepackage[shortlabels]{enumitem}
\usepackage{graphicx}
\usepackage{algorithmicx}
\usepackage[ruled]{algorithm}
\usepackage{algpseudocode}
\usepackage[numbers]{natbib}
\usepackage[colorlinks]{hyperref}
\title{Uniform Central Limit Theorems for Degree Distributions of Sparse Graphon Models}
\author{Robert Lunde and Soumendra N. Lahiri}
\date{}

\begin{document}
\maketitle

\begin{abstract}
The degree distribution is among the most basic summaries of a network, yet
no functional limit theorem for the empirical degree distribution of graphon
models has been available, even for dense, bounded graphons, for which only
finite-dimensional results were known.  We study the empirical degree
distribution of sparse $L^p$ graphon models, allowing unbounded kernels and
hence heavy-tailed degrees.  For graphons whose latent expected degrees have a
density, we prove a finite-dimensional central limit theorem whenever
$n\rho_n/\log n\to\infty$ and $n\rho_n^{2+\beta}\to0$, where $\rho_n$ is the
sparsity level and the kernel has $2+\beta$ moments for some $\beta>0$; the
limiting covariance coincides with that of the dense case.  Under additional
local regularity, we establish a functional central limit theorem on compact
intervals in the interior of the support.  Under tail regularity and a
stronger sparsity condition, we prove that uniformity extends to unbounded
intervals.  We also show that uniform limits fail at atoms of the degree
distribution, for instance in stochastic block models.  The analysis rests on
an exact multilinear representation of the conditional degree distribution,
which identifies its Hoeffding decomposition without Edgeworth expansions,
combined with a lattice version of Billingsley's tightness criterion and
moment inequalities of Efron--Stein type.
\end{abstract}

\section{Introduction}
The degree distribution is among the most fundamental descriptors of a large network. Empirically, it is the often first quantity examined in the study of real-world graphs, and the ability of a random graph model to reproduce observed degree behavior, such as the presence of high-degree hubs and heavy tails, is a primary criterion for assessing its realism. While degree distributions are of substantial applied interest, functional central limit theorems for the empirical degree distribution are rare, precluding natural inferential procedures such as Kolmogorov-Smirnov-type goodness-of-fit tests for degree functionals.

In this work, we establish conditions under which a central limit theorem for
the empirical degree distribution holds for sparse $L^p$ graphons in the sense
of \citet{borgs2019lp1}.  Informally, we show that the finite-dimensional
distributions of the empirical degree distribution, centered and scaled by
$\sqrt{n}$, converge to a Gaussian limit, and that under additional regularity
the convergence holds in $\ell^\infty$.  The most closely related work is that of
\citet{delmas2021asymptotic}, who establish a finite-dimensional central limit
theorem for dense, bounded graphons and note that their techniques are not
strong enough to yield tightness; we operate in the sparse regime, permit
unbounded kernels, and obtain the functional limit. \citet{Bickel-Chen-Levina-method-of-moments} show that the empirical distribution of degrees normalized by the average degree converges to the law of a latent quantity when $n\rho_n\to\infty$; we study its fluctuations.

Mathematically, this problem sits at the interface of lattice and continuum
behavior: the discreteness of the degrees does not disappear in the limit,
and it surfaces in two ways.  First, the aggregate effect of a single vertex
on the degrees of the others is governed by the binomial mass at one lattice
point, so the limiting covariance involves the density of the latent expected
degrees and is determined by a local limit theorem rather than a central
limit theorem.  Second, the process jumps at every lattice point, by amounts
that do not vanish relative to the length of the interval over which they
occur, so tightness cannot be established through a continuous modulus.
Increments must instead be controlled on the lattice itself, by mechanisms
that change with scale: a cancellation between adjacent binomial masses below
the bandwidth at which a degree fluctuates, and the regularity of the kernel
and of the degree density above it.  Controlling this interplay in the sparse,
unbounded regime is the central probabilistic challenge of the paper.  Our
main tool is an exact multilinear representation of the conditional degree
distribution, which identifies its Hoeffding decomposition without Edgeworth
expansions.


\section{Problem Setup}
We consider a sequence of adjacency matrices  $(A^{(n)})_{n \geq 1}$, where $A^{(n)} \in \{0,1\}^{n \times n}$ is generated by the sparse graphon model:
\begin{align}
\label{eq-sparse-graphon}
A_{ij}^{(n)} = A_{ji}^{(n)} = \mathbbm{1}(\eta_{ij} \leq \rho_n w(\xi_i, \xi_j) \wedge 1),
\end{align}
where $(\eta_{ij})_{1 \leq i<j \leq n}$ and $(\xi_i)_{i \in [n]}$ are independent collections of iid $\mathrm{Uniform}[0,1]$ random variables, $w:[0,1]^2 \mapsto \mathbb{R}^+$ is a symmetric function satisfying  $\|w\|_{L^1} = 1$, 
and  $\rho_n \rightarrow 0$ determines the sparsity level.  We assume no self-loops throughout, though it is not crucial.  While not needed to define our model, further moment conditions will be imposed on $w$ in the main theorems.   

Furthermore, let $D_i^{(n)} =\sum_{j \neq i} A_{ij}^{(n)}$ denote the degree of vertex $i$, $n' = n-1$, and $n'' = n-2$. For simplicity, in all the notation that follows, we will suppress dependence on $n$ when there is no ambiguity. Moreover, let $\mathfrak{D}_i = \mathfrak{D}(\xi_i) = E[w(\xi_i,\xi_j) \ | \ \xi_i]$ for any $j \neq i$.   The CDF $F(t)  =P(\mathfrak{D}(\xi) \leq t)$ is the degree distribution of primary interest and it is estimated by the empirical degree distribution $\hat{F}(t) = \frac{1}{n} \sum_{i=1}^n \mathbbm{1}(D_i/n'\rho_n \leq t)$, which has expectation $F^\star(t) := F_1(t')$, where $F_1(t') := P(D_1 \leq t')$, and $t' = \lfloor n' \rho_n t\rfloor$. 
We will focus on the limiting behavior of the random process $\mathbb{G}_n: \mathbb{R}^+ \mapsto \mathbb{R}$, where:
\begin{align}
\mathbb{G}_n(t) = \sqrt{n}[\hat{F}(t) - F(t)].
\end{align}
We will also study the closely related process $\mathring{\mathbb{G}}_n$, where: 
\begin{align}
\mathring{\mathbb{G}}_n(t) = \sqrt{n}[\hat{F}(t) - F^\star(t)].
\end{align}
Finite-dimensional convergence of $\mathring{\mathbb{G}}_n$ requires milder
sparsity conditions, since it is centered at its own expectation and is
therefore exactly mean-zero. Establishing convergence of $\mathbb{G}_n$
involves, in addition, showing that the bias $\sqrt{n}\,[F^\star(t) - F(t)]$ is
negligible, which requires a stronger sparsity condition. The centering by
$F(t)$ is nonetheless more natural: unlike $F^\star(t)$, this centering does not
depend on $n$ and reflects the degree structure of the underlying kernel,
making it the appropriate target for inference.


\section{Finite-Dimensional Central Limit Theorem}
\begin{theorem}[Finite-Dimensional Central Limit Theorem for $\mathring{\mathbb{G}}_n$]
\label{thm-fidi-clt}
Suppose that the following conditions hold:
\begin{enumerate}[leftmargin=1.5cm]
\item[M1.]\mylabel{M1}
\label{condition-m1}
 (Moment) 
 There exists some $\beta > 0$ such that $ w \in L^{2+\beta}([0,1]^2)$. 

\item[L1.]\mylabel{L1}
\label{condition-l1}
(Local Moment)
For any $0 < R < \infty$,
\begin{align}
\sup_{\mathfrak{D}(x) \leq R}E[w^{2+\beta}(x,\xi) \ | \ x] < \infty
\end{align}
\item[G1.]\mylabel{G1} (Graphon Regularity)
\label{condition-g1}
There exists $c > 0$ such that $w \geq c \ a.s.$  Moreover, $x \rightarrow w(x,\cdot)$ is $L^2$ continuous. 

\item[D1.]\mylabel{D1} (Density Regularity)
\label{condition-d1}
The law of $\mathfrak{D}$ is absolutely continuous with respect to the Lebesgue measure, with density $f_\mathfrak{D}$ that is continuous on the interior of the support. 
  
\item[S1.]\mylabel{S1}(Sparsity) 
\label{condition-s1}
The sparsity level satisfies:
\begin{align}
\label{eq-sparsity-new-2}
\frac{n\rho_n}{\log n}\to\infty,\quad
n\rho_n^{\,2+\beta}\to0 .
\end{align}
\end{enumerate}
Then,
\begin{align*}
\mathring{\mathbb{G}}_n \stackrel{fdd}{\rightsquigarrow}\mathbb{G},
\end{align*}
where $\mathbb{G}$ is a mean $0$ Gaussian process and covariance function $\Sigma(s,t) = \sum_{i=1}^3 \Sigma_i(s,t)$, where for $s,t \in \mathrm{supp} \  F$,
\begin{align}
\label{eq-cov-function}
\begin{split}
\Sigma_1(s,t) &= F(s) \wedge F(t) - F(s)F(t),
\\ \Sigma_2(s,t) &= f_\mathfrak{D}(s) f_\mathfrak{D}(t) \  \mathrm{Cov}(w(\mathfrak{D}^{-1}(s),\xi), w(\mathfrak{D}^{-1}(t),\xi) )
\\ \Sigma_3(s,t) &= f_{\mathfrak{D}}(t) \cdot \left\{tF(s) - \int_{\mathfrak{D}(\xi) \leq s}w(\mathfrak{D}^{-1}(t),\xi)  dP\right\} 
\\ & \ \ \ +  f_{\mathfrak{D}}(s) \cdot \left\{sF(t) - \int_{\mathfrak{D}(\xi) \leq t}w(\mathfrak{D}^{-1}(s),\xi)  dP\right\}. 
\end{split}
\end{align}    
\end{theorem}

\begin{remark}[Comparison with \citet{delmas2021asymptotic}]
Our covariance function is equivalent to the one established in the finite-dimensional central limit theorem of \citet{delmas2021asymptotic} for dense graphons.  These authors  assume the graphon $w$ to be three times continuously differentiable and bounded away from both $0$ and $1$, with a strictly positive degree derivative. These assumptions support the Taylor and bivariate Edgeworth expansions they employ, which break down in the sparse regime. Our strategy, like theirs, involves H\'ajek projections, but uses different tools to operate under sparsity and to relax the smoothness requirements. Importantly, we permit unbounded $w$, thereby accommodating heavy-tailed degrees.  Their setting of dense, bounded graphons also allows a proof architecture that avoids the difficulties created by unbounded kernels.  These remain unpleasant even within an Edgeworth-based approach, such as the one underlying an earlier, unpublished version of this work: clipping, vertices of atypically large expected degree, and local limit estimates that are not uniform in the latent position must all be handled, and in our setting they are handled within the multilinear framework.
\end{remark}

\begin{remark}[Existence of $\mathfrak{D}^{-1}$]
\label{remark-inverse}
Note that there exists a measure-preserving transformation $\psi:[0,1] \mapsto [0,1]$ such that $\mathfrak{D}\circ\psi$ is monotone non-decreasing.  Since the law of $\mathfrak{D}\circ\psi$ is absolutely continuous, it agrees almost everywhere with a strictly increasing function on the support of $F$; we identify $\mathfrak{D}$ with this representative, so $\mathfrak{D}^{-1}$ is well-defined on its support.  For $s$ in the support of $F$ we write $x_s := \mathfrak D^{-1}(s)$.  
\end{remark}

A central feature of our model is that it permits unbounded kernels $w$, and hence heavy-tailed degree distributions. The extent of heavy-tailedness in real-world networks is itself contested: while a substantial literature models empirical degree distributions as power laws
$f_{\mathfrak{D}}(x) \sim x^{-\gamma}$ (see, e.g. \citet{albert2002statistical,newman2003structure}) more recent work has questioned how prevalent genuine scale-free structure is (see, for instance, \citet{broido2019scale} and \citet{Voitalov_2019}). We note that the strict scale-free regime $\gamma \in (2,3)$, in which the degree distribution has
infinite variance, lies fundamentally outside the scope of our Gaussian central limit theorem, but our theorem does accommodate a broad range of heavy-tailed behavior.

In what follows, for notational simplicity, let $\check{w} = w \wedge \rho_n^{-1}$ denote the clipped graphon, which yields the same generative model as \eqref{eq-sparse-graphon}, and let  $\check{\mathfrak{D}}_i = E[\check{w}(\xi_i,\xi_j) \ | \ \xi_i]$ denote the clipped analogue of the latent degree parameter $\mathfrak D_i$. The clipped degree functional introduces some subtlety into the analysis, since a vertex with arbitrarily large $\mathfrak{D}_i$ may have a moderate clipped degree. It is also what permits the mild tail assumptions in our theorems. 

Furthermore, for any $i \in [n]$ and $J \subseteq [n]\setminus \{i\}$, let $D_{i\setminus J} := D_i-\sum_{j\in J}A_{ij}$ denote the degree of $i$ with the edges to $J$ removed, and write $F_i(t'\mid\mathcal G) := P(D_i\leq t'\mid\mathcal G)$ and $F_{i \setminus J}(t'\mid\mathcal G) := P(D_{i\setminus J} \leq t'\mid\mathcal G)$ for the conditional distribution functions given the latent positions $\mathcal G$ (defined below), and $F_{i\setminus J}(t')$ for the unconditional one. For simplicity, for any $j \neq i$ and  $J = \{j\}$, we write $D_{i\setminus j}$ and $F_{i \setminus j}$ for $D_{i\setminus\{j\}}$ and $F_{i \setminus \{j\}}$. 
 
\subsection{The Multilinear Representation and the Hoeffding Decomposition}
\label{sec-multilinear-hoeffding}

Our proof strategy involves H\'ajek projections, which we discuss briefly.  In what follows, for $S \subseteq [n]$ and $i \in [n]$, let:
\begin{align*}
\mathcal{G} = \sigma( \xi_1, \ldots, \xi_n), \quad  \mathcal{G}_S = \sigma( (\xi_i)_{i \in S}), \quad \mathcal{G}_i = \mathcal{G}_{\{i\}}. 
\end{align*}
Recall that a square-integrable $T_n$ that is $\mathcal{G}$-measurable
admits a unique Hoeffding decomposition
\citep{hoeffding1948class,efron1981jackknife,vanzwet1984berry}
\begin{align}
\label{eq-hoeffding-rep}
T_n =  \sum_{ S \subseteq [n]} \Pi(T_n;S),
\end{align}
where $\Pi(T_n;S)$ denotes the projection of $T_n$ onto
\begin{align*}
\mathcal{H}_S := \bigl\{ f \in L^2 \ : \ f \text{ is } \mathcal{G}_S 
\text{-measurable},\ E[f \mid \mathcal{G}_A] = 0 \text{ for all } A \subsetneq S \bigr\},
\end{align*}
and $\Pi(T_n \; ; \; \emptyset) = E(T_n)$. Taking $A = \emptyset$ shows each component is mean zero, and one can readily show that the spaces
$\mathcal{H}_S$ are mutually orthogonal.
Consequently, $\mathrm{Var}(T_n) = \sum_{\emptyset \neq S} \mathrm{Var}(\Pi(T_n;S))$,
and since $E[T_n \mid \mathcal{G}_S] = \sum_{U \subseteq S} \Pi(T_n;U)$, the H\'ajek projection $H(T_n)$, which is the sum of the first-order projections, satisfies: 
\begin{align}
\begin{split}
\label{eq-hajek-proj}
H(T_n) &:= \sum_{i=1}^{n} \Pi(T_n;\{i\})
 = \sum_{i=1}^{n} E[T_n \mid \xi_i] - (n-1)E[T_n],
 \\
\mathrm{Var}(H(T_n))  
 &= \mathrm{Var}(T_n) -\sum_{|S| \geq 2} \mathrm{Var}\bigl(\Pi(T_n;S)\bigr).
\end{split}
\end{align}

Thus, if the degenerate components are asymptotically negligible and $H(T_n)$
converges in distribution, $T_n$ has the same limit.  Here $H(T_n)$ is a
normalized sum of i.i.d.\ terms for which the Lindeberg condition can be verified;
the substance of the argument lies in the variance matching.

Our process depends on the edge variables $(\eta_{ij})$ as well as the latent
positions, so \eqref{eq-hoeffding-rep} does not apply directly.  Conditioning on $\mathcal{G}$ repairs this at the
cost of one term: each $\xi_i$ is $\mathcal{G}$-measurable, so the residual
$T_n - E[T_n \mid \mathcal{G}]$ is orthogonal to $L^2(\mathcal{G}_i)$ for any $i \in [n]$ and
$H(T_n) = H(E[T_n \mid \mathcal{G}] )$,
whence, in general, we have:
\begin{align}
\label{eq-gap-two-sources}
  \Var( T_n)-\Var(H(T_n))
 = \sum_{|S| \geq 2} \Var\Bigl(\Pi(E(T_n \ | \ \mathcal{G});S)\Bigr)
 \;+\; E\bigl[\Var(T_n \mid \mathcal{G})\bigr],
\end{align}
the two terms on the right of the display representing higher-order latent fluctuations and residual edge
randomness respectively. To establish the finite-dimensional CLT of interest, by Cram\'er-Wold device, it suffices to show that, for any $K \in \mathbb{N}$ and $\gamma \in \mathbb{R}^K$, $\sum_{k=1}^K \gamma_k \mathring{\mathbb{G}}_{n}(t_k) \rightsquigarrow \sum_{k=1}^K \gamma_k \mathbb{G}(t_k)$; thus, we consider Hoeffding decompositions for $T_n(\gamma, t) = \sum_{k=1}^K \gamma_k \mathring{\mathbb{G}}_{n}(t_k)$, where $K \in \mathbb{N}$ and $(\gamma_k ,t_k)_{k \in [K]}$ are arbitrary.  

For general symmetric statistics, computing the Hoeffding decomposition can be an arduous task and is often not needed for establishing a Normal limit via a H\'ajek projection argument.  However, a key feature of our problem is that $E(T_n(\gamma,t) \ | \ \mathcal{G})$ is a multilinear form, and its Hoeffding decomposition can be computed.   

Before stating this lemma, we introduce some additional notation.  Since $D_i \ | \ \xi_i \sim \mathrm{Binomial}(n', \rho_n \check{\mathfrak{D}}(\xi_i))$, notation specific to the binomial distribution will be useful. For $M \in \mathbb{N}$ and $p \in (0,1)$, let $b_M(m;p)$ and $\mathcal{B}_M(m;p)$ denote the $\mathrm{Binomial}(M,p)$ PMF and CDF evaluated at $m$, respectively.  Furthermore, for $f:\{0, \ldots, M\} \mapsto \mathbb{R}$, $d \geq 1$, and $m \geq 1$, let $\mathsf{D}f(m) := f(m)-f(m-1)$ denote the
backward difference, so that $\mathsf{D}\mathcal{B}_M = b_M$ and
$\mathsf{D}^{d}\mathcal{B}_M = \mathsf{D}^{d-1}b_M$.  For $d = 0$, $\mathsf{D}^0 f = f$. We have the following result:

\begin{lemma}[Multilinear Representation and Hoeffding Decomposition]
\label{lem-multilinear}
For every $s$ and any $i \in [n]$,
\begin{align}
\label{eq-multilinear-2}
F_i(s' \ | \ \mathcal{G}) = \sum_{S\subseteq[n]\setminus\{i\}}
 \underbrace{(-\rho_n)^{|S|}\prod_{j\in S} a(\xi_i,\xi_j) \cdot
 \mathsf D^{|S|}\mathcal{B}_{n'-|S|}\bigl(s'; \rho_n \check{\mathfrak{D}}(\xi_i)\bigr)}_{=:\,K_S(s)},
\end{align}
where $a(\xi_i,\xi_j) = \check{w}(\xi_i,\xi_j) - \check{\mathfrak{D}}(\xi_i)$ and the  $S=\emptyset$ term is   $\mathcal{B}_{n'}( s';\rho_n \check{\mathfrak{D}}(\xi_i))$.
Moreover, for any $S \subseteq [n]/\{i\}$,  
\begin{align}
\begin{split}
\label{eq-projection-form}
\Pi(F_i(s' \ | \ \mathcal{G})\; ; \; S) &=  E[K_S(s') \ | \ \mathcal{G}_S]
\\ \Pi(F_i(s' \ | \ \mathcal{G}) \; ; \;  S \cup \{i\}) &= K_S(s') - E[K_S(s') \ | \ \mathcal{G}_S].
\end{split}
\end{align}

\end{lemma}

\begin{proof}
We start by proving \eqref{eq-multilinear-2}. For notational simplicity, fix $i \in [n]$ and let $p_{j} = \rho_n \check{w}(\xi_i, \xi_j)$, and $\bar{p} = E[p_j \ | \ \xi_i] = \rho_n \check{\mathfrak{D}}(\xi_i)$. We will view $F_{i}(s' \ | \ \mathcal{G})$ as a function of parameters $p = (p_j)_{j \neq i}$.  
By conditioning on the  outcome of $A_{ij}$, it is clear that:
\begin{align*}
G(p) = p_j F_{i \setminus j}(s'-1 \ | \ \mathcal{G}) + (1-p_j) F_{i \setminus j}(s' \ | \ \mathcal{G}),
\end{align*}
where $F_{i \setminus j}(\cdot \ | \ \mathcal{G})$ does not depend on $p_j$.  Since for any $j \neq i$, 
\begin{align*}
\partial_{p_j}G  = -\mathsf DF_{i \setminus j}(s'  \ | \ \mathcal{G}),
\end{align*}
it is clear that $G$ is multilinear.  Iterating over distinct coordinates,
\begin{align*}
\partial_SG = (-1)^{|S|}\mathsf D^{|S|}F_{i \setminus S}(s' \ | \ \mathcal{G}).
\end{align*}
A multilinear function equals its Taylor expansion about any point, and the
expansion only contains mixed partial derivatives.  Thus, a Taylor expansion around $\bar{p}$ takes the form:
\begin{align*}
G(p) = \sum_{S\subseteq[n]\setminus\{i\}}\prod_{j\in S}(p_j-\bar p)\;
 \partial_SG(\bar p,\dots,\bar p).
\end{align*}
Since $G(\bar{p}, \ldots, \bar{p}) = \mathcal{B}_{n'}(s'; \rho_n \check{\mathfrak{D}}(\xi_i))$,  substituting $p_{j}-\bar p = \rho_n\bigl(\check w(\xi_i,\xi_j)-\check{\mathfrak
D}(\xi_i)\bigr) = \rho_na(\xi_i,\xi_j)$  and $\partial_SG|_{p = \bar{p}}$ gives \eqref{eq-multilinear-2}.  Note that $S = \emptyset$ corresponds to the value of the function at the expansion point, namely $G(\bar p \ldots, \bar p)$. 

We now prove \eqref{eq-projection-form}. For $S \subseteq [n]\setminus\{i\}$ define
\begin{align}
 \label{eq-h-def}
h_S := E\bigl[K_S(s) \mid \mathcal{G}_S\bigr],
\qquad
h_{S \cup \{i\}} := K_S(s) - E\bigl[K_S(s) \mid \mathcal{G}_S\bigr],
\end{align}
noting that every $S' \subseteq [n]$ is of exactly one of these two forms, so
that \eqref{eq-h-def} defines $h_{S'}$ for every index set.  We show that
$h_{S'} = \Pi(F_i(s' \mid \mathcal{G}); S')$ for each $S' \neq
\emptyset$, and that $h_\emptyset = E[F_i(s' \mid \mathcal{G})]$.
 
\emph{Centering.}  Fix $\xi_i$.  The factor
$\mathsf{D}^{|S|}\mathcal{B}_{n'-|S|}(s';\rho_n\check{\mathfrak{D}}(\xi_i))$ is
then deterministic, the variables $(\xi_j)_{j \neq i}$ are independent, and
$E[a(\xi_i,\xi_j) \mid \xi_i] = 0$ for every $j$.  Hence integrating out
any single $\xi_j$ with $j \in S$ annihilates the product, so
\begin{align}
\label{eq-K-centered}
E\bigl[K_S(s) \mid \mathcal{G}_{(S \setminus \{j\}) \cup \{i\}}\bigr] = 0
\qquad \text{for every } j \in S .
\end{align}
 
\emph{Membership.}  By \eqref{eq-K-centered} and the tower property,
$E[K_S(s) \mid \mathcal{G}_A] = 0$ for every $A \subsetneq S \cup \{i\}$ with
$S \not\subseteq A$; the same holds for $h_{S \cup \{i\}}$.  For $A = S$, we have
$E[h_{S \cup \{i\}} \mid \mathcal{G}_S] = 0$ by construction.  As
$h_{S\cup\{i\}}$ is $\mathcal{G}_{S \cup \{i\}}$-measurable, this gives
$h_{S \cup \{i\}} \in \mathcal{H}_{S \cup \{i\}}$.  Similarly $h_S$ is
$\mathcal{G}_S$-measurable, and for $A \subsetneq S$, \eqref{eq-K-centered} and
the tower property give $E[h_S \mid \mathcal{G}_A] = 0$, so
$h_S \in \mathcal{H}_S$.
 
\emph{Summation.}
\begin{align*}
\sum_{S' \subseteq [n]} h_{S'}
 &= \sum_{S \subseteq [n]\setminus\{i\}}
 \Bigl( E[K_S(s) \mid \mathcal{G}_S] + K_S(s) - E[K_S(s)\mid\mathcal{G}_S]\Bigr)
\\ &= \sum_{S \subseteq [n]\setminus\{i\}} K_S(s)
 = F_i(s' \mid \mathcal{G}).
\end{align*}
 
\emph{Identification.}  Fix $S_0 \neq \emptyset$ and apply $\Pi(\,\cdot\,;S_0)$,
the orthogonal projection onto $\mathcal{H}_{S_0}$, to the previous display.
The projection is linear and the sum is finite; $\Pi(h_{S'};S_0) = 0$ for
$S' \neq S_0$ because $h_{S'} \in \mathcal{H}_{S'}$ and the spaces
$\mathcal{H}_{S'}$ are mutually orthogonal, while $\Pi(h_{S_0};S_0) = h_{S_0}$
because a projection fixes its own subspace.  Hence
$\Pi(F_i(s'\mid\mathcal{G}); S_0) = h_{S_0}$, which is
\eqref{eq-projection-form}.  Taking $S_0 = \emptyset$, for which
$\mathcal{H}_\emptyset$ is the space of constants, gives
$h_\emptyset = E[K_\emptyset(s)] = E[F_i(s' \mid \mathcal{G})]$.

Thus, the claim follows. 
\end{proof}

For fixed edge probabilities, the identity of Lemma~\ref{lem-multilinear} is
closely related to the Krawtchouk expansion of the Poisson-binomial
distribution around a binomial distribution \citep{roos2000binomial}.  Our
observation is that, when the edge probabilities are random through
independent latent positions, the same identity is exactly the Hoeffding
decomposition of the conditional degree distribution, so that every
component is available in closed form.

We now provide a high-level roadmap of the proof based on the above Hoeffding
decomposition.  For simplicity, we examine the H\'ajek projection for a single
$s \in \{t_1,\ldots,t_K\}$; the general case follows by linearity of
$H(\cdot)$.  Recalling that $H(T_n)$ depends on $T_n$ only through
$E[T_n \mid \mathcal{G}]$, and that $\Pi(\,\cdot\,;S)$ is linear across the $n$
summands of $E[\mathring{\mathbb{G}}_{n}(s) \mid \mathcal{G}]$, we obtain
\begin{align*}
H\bigl(\mathring{\mathbb{G}}_{n}(s)\bigr)
&= \frac{1}{\sqrt n}\sum_{i=1}^{n}\Bigl\{
 \Pi\bigl(F_i(s'\mid\mathcal{G});\{i\}\bigr)
 + \sum_{j \neq i}\Pi\bigl(F_i(s'\mid\mathcal{G});\{j\}\bigr)\Bigr\}
\\ &= \frac{1}{\sqrt n}\sum_{i=1}^{n}\Bigl[
 \bigl\{F_i(s'\mid\mathcal{G}_i) - F_i(s')\bigr\}
 - \rho_n\sum_{j\neq i}
E\bigl[a(\xi_i,\xi_j)\,b_{n''}\bigl(s';\rho_n\check{\mathfrak{D}}(\xi_i)\bigr)
 \,\big|\,\mathcal{G}_j\bigr]\Bigr]
\\ &= \frac{1}{\sqrt n}\sum_{j=1}^{n}\bigl\{\psi_s(\xi_j)+\phi_s(\xi_j)\bigr\} .
\end{align*}
Here $F_i(s') = F^\star(s)$ for every $i$ by exchangeability, and
$F_i(s' \mid \mathcal{G}_i) = \mathcal{B}_{n'}(s';\rho_n
\check{\mathfrak{D}}(\xi_i))$, so that the two
first-order contributions are the functions
\begin{align}
\label{eq-psi-phi-def}
\psi_s(\xi) := \mathcal{B}_{n'}\bigl(s';\rho_n\check{\mathfrak{D}}(\xi)\bigr)
 - F^\star(s),
\quad
\phi_s(\xi) := -\,n'\rho_n\int_0^1 a(x,\xi)\,
 b_{n''}\bigl(s';\rho_n\check{\mathfrak{D}}(x)\bigr)\,dx,
\end{align}
evaluated at $\xi_j$.  The first measures the effect of a vertex's latent
position on its own degree; the second, the aggregate effect of that position on
the degrees of the remaining $n'$ vertices, each individually of order $\rho_n$.
Thus, it follows that:
\begin{align*}
& \ \mathrm{Cov}\bigl(H(\mathring{\mathbb{G}}_{n}(s)),\,
 H(\mathring{\mathbb{G}}_{n}(t))\bigr)
\\ = & \ \mathrm{Cov}\bigl(\psi_s(\xi),\psi_t(\xi)\bigr)
 + \mathrm{Cov}\bigl(\phi_s(\xi),\phi_t(\xi)\bigr)
 + \mathrm{Cov}\bigl(\psi_s(\xi),\phi_t(\xi)\bigr)
 + \mathrm{Cov}\bigl(\phi_s(\xi),\psi_t(\xi)\bigr)
\\ = & \ \Sigma_{n,1}(s,t) + \Sigma_{n,2}(s,t) + \Sigma_{n,3}(s,t),
 \qquad \text{say},
\end{align*}
with the last two covariances being grouped into $\Sigma_{n,3}$.

Since $H(\mathring{\mathbb{G}}_{n}(s))$ is a normalized sum of i.i.d.\
terms whose summands converge in $L^2$, the Lindeberg condition is immediate,
and a Gaussian limit follows from the Lindeberg--Feller theorem once the
variances converge.  The argument therefore reduces to showing that (1)
$\Sigma_{n,i} \to \Sigma_i$ for $i \in [3]$, (2) the residual network noise is
negligible, and (3) the higher-order Hoeffding components are negligible.

Step (1) is not a matter of routine approximation.  While the limit of
$\psi_s$ follows from concentration of $D_1$ about its conditional
mean, the limit of $\phi_s$ is governed by the binomial PMF at a parameter of
order $\rho_n \rightarrow 0$ and an argument of order $n\rho_n \rightarrow \infty$.  Rogozin's inequality
already shows that this PMF is $O((n\rho_n)^{-1/2})$ near the mean, so the
factor $n'\rho_n$ must be balanced against it, and delicate analysis is needed
to establish convergence to $f_{\mathfrak{D}}(s)$ once the latent variable is
integrated out.  This step may be viewed as an integrated local limit
statement, and it is in this way that the discreteness of the degrees survives
into the covariance of the limiting process.

Step (3) also requires substantial care.  By
Lemma~\ref{lem-multilinear} and an exchangeability argument discussed further in Section \ref{sec-degenerate} , the component indexed by a set of size $d$
carries a factor $\binom{n''}{d}\rho_n^{2d}$ against a $d$-fold product of
covariances of $a$ and two $d$-th order backward differences of
binomial distribution functions.  The differences are of order
$\sigma_*^{-d}$ with $\sigma_*^2 \asymp n\rho_n$, so the combinatorial and
lattice factors alone contribute $(d!)^{-1}(\rho_n/c)^d$. As result,  a spare power of
$\rho_n$ is available at each order; moreover the series is convergent so long as the covariance
factor remains bounded.  Whether it does is governed by
$\|a(x,\cdot)\|_{L^2}$, which is a truncated version of a product involving a kernel that is only assumed to have $2+\beta$ moments to begin with. 
The analysis involves partitioning the latent space into regions on which
different bounds for $\|a(x,\cdot)\|_{L^2}$ apply. The differences of
binomial distribution functions are controlled by Fourier inversion, uniformly
in the order $d$ and in the lattice point.

In the following subsection, we develop some tools that will be needed to carry out this strategy. 

\subsection{Binomial and Local Limit Estimates}
We begin by collecting some properties of the binomial distribution that will be useful in our analysis. The most non-standard item is \ref{binom-fact-c}, which provides a uniform bound on the $d$-fold backward difference via Fourier inversion. In what follows, let $\sigma^2(M,p) = Mp(1-p)$ and for $m \in \{0,1, \ldots, M\}$ let $x_m(M,p)$ denote the shifted and scaled lattice point $x_m(M,p) = (m-Mp)/\sigma(M,p)$. For ease of notation, we suppress the arguments $M$ and $p$ when there is no ambiguity.   
\begin{lemma}[Useful Properties of the Binomial Distribution]
\label{lem-binomial-facts}
The following statements hold for every $M\ge1$: 
\begin{enumerate}[(a)]
\item 
\label{binom-fact-a}
For $p\in(0,1/4]$ and $0\le m\le M$, there exists $C >0$ independent of $M$, $p$, and $m$ such that
$b_{M}(m;p) \le C\sigma^{-1} \ e^{-\zeta x_m^2}$ whenever $|m-Mp|\le Mp$.
\item \label{binom-fact-b} For $p\in(0,1/4]$, 
 $\mathcal{B}_{M}(m;p) \le e^{-Mp/8}$
whenever $m \leq Mp/2$  and $1-\mathcal{B}_{M}(m;p) \le e^{-Mp/8}$ whenever $m > 2Mp$.
\item \label{binom-fact-c} For $d \geq 1$, $\displaystyle\sup_m\bigl|\mathsf D^{d}\mathcal{B}_M(m;p)\bigr|
 \le \frac{\Gamma(d/2)}{2\pi\,(\zeta'\sigma^2)^{d/2}},$
  where $\zeta' := \frac{2}{\pi^2}$.  Consequently, \newline
$\sup_{m}|\mathsf D^{d}\mathcal{B}_M(m;p) | \leq  \Gamma_d \  \sigma^{-d}$, where $\Gamma_d  =(C_0^dd!)^{1/2}$ and $C_0 >0$ is independent of $M$, $d$, and $p$.  
\item \label{binom-fact-d}
$\sum_m|\mathsf D^{d}\mathcal{B}_M(m;p)|\le2^{d-1}$.
\item \label{binom-fact-e}
$\displaystyle\int_0^1b_M(m;p)\,dp = \frac1{M+1}$ for every
$0\le m\le M$.
\item \label{binom-fact-f}
$\displaystyle\partial_pb_M(m;p) = b_M(m;p)\,\frac{m-Mp}{p(1-p)}$.
\end{enumerate}
\end{lemma}

\begin{proof}
(a) Let $m_0 := \lfloor (M+1)p\rfloor$ denote the mode of $b_M(\cdot\,;p)$, and
note $|m_0 - Mp| \leq 1$.  Stirling's formula applied to the three factorials
in $\binom{M}{m_0}p^{m_0}(1-p)^{M-m_0}$ gives the standard bound at the mode,
\begin{align}
\label{eq-mode-bound}
b_M(m_0;p) \;\leq\; \frac{C_1}{\sigma},
\end{align}
for an absolute constant $C_1$; see, e.g., \citet[][Ch.~VII]{feller1971volume2}.
 
For the decay away from the mode, write the log-increment
\begin{align*}
r(m) := \log\frac{b_M(m;p)}{b_M(m-1;p)}
 = \log\frac{(M-m+1)\,p}{m\,(1-p)},
\qquad 1 \leq m \leq M,
\end{align*}
which is decreasing in $m$; this is precisely the log-concavity of the mass
function, and it is what allows the deviation bound to be obtained by
telescoping.  We claim that there is an absolute $\zeta_1 > 0$ with
\begin{align}
\label{eq-increment-bound}
r(m) \;\leq\; -\,\zeta_1\,\frac{m-m_0}{\sigma^2}
\qquad \text{for } m_0 < m \leq 2Mp,
\end{align}
and symmetrically $r(m) \geq \zeta_1(m_0-m)/\sigma^2$ for
$0 \leq m \leq m_0$.  Indeed, on this range
\begin{align*}
r(m) = \log\Bigl(1 - \frac{m - (M+1)p}{(M+1)(1-p)}\Bigr)
 - \log\Bigl(1 + \frac{m-(M+1)p}{(M+1)p}\Bigr),
\end{align*}
and since $m \leq 2Mp$ implies $M-m+1 \geq M(1-2p) + 1 \geq M/2 \cdot (1-p)$
and $m \asymp Mp$, both arguments are bounded away from $-1$ by an absolute
constant; using $\log(1+u) \asymp u$ on such a range, the two terms are each of
order $-(m-m_0)/\bigl(Mp(1-p)\bigr) = -(m-m_0)/\sigma^2$, which gives
\eqref{eq-increment-bound}.
 
Summing \eqref{eq-increment-bound} from $m_0$ to $m$, for
$m_0 < m \leq 2Mp$,
\begin{align*}
\log\frac{b_M(m;p)}{b_M(m_0;p)} = \sum_{k=m_0+1}^{m} r(k)
 \;\leq\; -\frac{\zeta_1}{\sigma^2}\sum_{k=m_0+1}^{m}(k-m_0)
 \;\leq\; -\frac{\zeta_1}{2}\,\frac{(m-m_0)^2}{\sigma^2},
\end{align*}
and since $|m_0-Mp| \leq 1$ we may replace $m_0$ by $Mp$ at the cost of an
absolute constant.  Combining with \eqref{eq-mode-bound} gives
\ref{binom-fact-a} on $Mp \leq m \leq 2Mp$; the argument for
$0 \leq m \leq Mp$ is identical with the inequalities reversed, and together
these cover the range $|m-Mp| \leq Mp$.
 
\ref{binom-fact-b} is the multiplicative Chernoff bound
$P(\mathrm{Bin}(M,p) \leq (1-\theta)Mp) \leq e^{-\theta^2Mp/2}$ evaluated at
$\theta = 1/2$.

\ref{binom-fact-c} 
For any $f$ supported on $\{0,1,\ldots\}$, reindexing gives
\begin{align*}
\sum_m \mathsf{D}f(m)z^m = (1-z)\sum_m f(m)z^m,
\end{align*}
so iterating and using that the generating function of $b_M(\cdot\,;p)$ is
$\sum_m b_M(m;p)z^m = (1-p+pz)^M$, we have
\begin{align*}
\sum_m \mathsf{D}^{d-1}b_M(m;p)\,z^m = (1-z)^{d-1}(1-p+pz)^M .
\end{align*}
Both sides are polynomials, so we may evaluate at $z = e^{iu}$ and extract
coefficients by Fourier inversion on the unit circle, yielding
\begin{align*}
\mathsf D^{d-1}b_M(m;p)
 = \frac1{2\pi}\int_{-\pi}^{\pi}(1-e^{iu})^{d-1}(1-p+pe^{iu})^{M}e^{-imu}\,du .
\end{align*}
Now $|1-e^{iu}|\le|u|$, and since $1-\cos u\ge2u^2/\pi^2$ on $[-\pi,\pi]$,
\begin{align*}
|1-p+pe^{iu}|^2 = 1-2p(1-p)(1-\cos u)
 \le\exp\bigl(-4p(1-p)u^2/\pi^2\bigr),
\end{align*}
whence $|1-p+pe^{iu}|^{M}\le\exp(-2\sigma^2u^2/\pi^2) =
\exp(-\zeta'\sigma^2u^2)$.  Therefore
\begin{align*}
\bigl|\mathsf D^{d-1}b_M(m;p)\bigr|
 \le\frac1{2\pi}\int_{\mathbb R}|u|^{d-1}e^{-\zeta'\sigma^2u^2}\,du
 = \frac{\Gamma(d/2)}{2\pi(\zeta'\sigma^2)^{d/2}},
\end{align*}
which is the stated bound.   The second claim follows from the Legendre duplication formula
$\Gamma(d/2)^2\le\Gamma(d/2)\Gamma(\tfrac{d+1}2) = \sqrt\pi\,2^{1-d}\Gamma(d)$ for $d \geq 3$ applied to the first claim; the case $d <3$ holds trivially.  

\ref{binom-fact-d}
$\sum_m|\mathsf Df(m)|\le2\sum_m|f(m)|$;
iterating $d-1$ times from $\sum_mb_M(m;p) = 1$ gives $2^{d-1}$.
\ref{binom-fact-e}- \ref{binom-fact-f} follow immediately by calculation.  
\end{proof}

We now state a lemma that establishes the integrated local limit theorem needed to carry out step (1) of our strategy.  Throughout, we call
\begin{align}
\label{eq-good-set}
\mathscr G_n := \bigl\{x\in[0,1]:\ \rho_n\mathfrak D(x)\leq\tfrac14\bigr\}
\end{align}
the \emph{good set} of $\mathfrak D$.  Since $\check{\mathfrak D}\leq\mathfrak D$
and $\check{\mathfrak D}\geq c$ by \ref{condition-g1}, on $\mathscr G_n$ both
$\rho_n\check{\mathfrak D}(x)$ and $\rho_n\mathfrak D(x)$ lie in
$[c\rho_n,\tfrac14]$, the range in which
Lemma~\ref{lem-binomial-facts}\ref{binom-fact-a} applies.  Its complement is
small: by Markov's inequality, Jensen's inequality and \ref{condition-m1},
\begin{align}
\label{eq-good-complement}
\mathrm{Leb}\bigl(\mathscr G_n^c\bigr) = P\bigl(\mathfrak D>(4\rho_n)^{-1}\bigr)
 \leq(4\rho_n)^{2+\beta}E\bigl[\mathfrak D^{2+\beta}\bigr]
 \leq C\rho_n^{2+\beta}.
\end{align}
We also record two facts used repeatedly.  First, with
$g(x) := \mathfrak D(x)-\check{\mathfrak D}(x) = E[(w(x,\xi)-\rho_n^{-1})_+]$,
Markov's inequality and \ref{condition-m1} give
\begin{align}
\label{eq-g-global}
\int_0^1 g(x)\,dx \;\leq\; \rho_n^{1+\beta}E\bigl[w^{2+\beta}\bigr]
 \;=\; C_w\,\rho_n^{1+\beta},
\qquad C_w := \|w\|_{L^{2+\beta}}^{2+\beta}.
\end{align}
Second, for every $p\in[c\rho_n,1]$ and every $m\leq2n''\rho_ns$,
\begin{align}
\label{eq-pmf-sup}
b_{n''}(m;p)\;\leq\;C\,(n\rho_n)^{-1/2},
\end{align}
by the mode bound of Lemma~\ref{lem-binomial-facts}\ref{binom-fact-a} when
$p\leq\frac14$, and by the lower-tail bound of
Lemma~\ref{lem-binomial-facts}\ref{binom-fact-b} when $p>\frac14$, in which case
$m\leq n''p/2$ for $n$ large and $b_{n''}(m;p)\leq e^{-n''/32}$.  The following lemma and arguments used in its proof to control truncation bias will play an important role throughout our argument.  

For $s$ in the interior of $\mathrm{supp}\,F$, define the measure $\mu_s$ on
$[0,1]$ by
\begin{align}
\label{eq-mu-s-def}
\mu_s(dx) := n'\rho_n\,b_{n''}\bigl(s';\rho_n\check{\mathfrak D}(x)\bigr)\,dx ,
\end{align}
so that $\phi_s(\xi) = -\int_0^1a(x,\xi)\,\mu_s(dx)$ by \eqref{eq-psi-phi-def}.
Since $D_{1\setminus2}$ given $\xi_1 = x$ is $\mathrm{Bin}(n'',\rho_n\check{\mathfrak D}(x))$,
the total mass is $\mu_s([0,1]) = n'\rho_n\,P(D_{1\setminus2}=s')$.

\begin{lemma}[Mass and localization of $\mu_s$]
\label{lem-mu-limit}
Suppose \ref{condition-m1}, \ref{condition-g1}, \ref{condition-d1} and
\ref{condition-s1} hold, and let $s$ lie in the interior of
$\mathrm{supp}\,F$.  Then
\begin{enumerate}[(a)]
\item
\label{localization-item-a}
$\mu_s([0,1]) = n'\rho_n\,P(D_{1\setminus2}=s')\to f_{\mathfrak D}(s)$;
\item
\label{localization-item-b}
for every $\eta>0$, $\;\mu_s\bigl(\{x:|\mathfrak D(x)-s|>\eta\}\bigr)\to0$.
\end{enumerate}
\end{lemma}

\begin{proof}
We split $\mu_s$ into its parts on the good set $\mathscr G_n$ and
on its complement.

\emph{Step 1: replacing $\check{\mathfrak D}$ by $\mathfrak D$ on
$\mathscr G_n$.}  By Lemma~\ref{lem-binomial-facts}\ref{binom-fact-f} and the
mean value theorem, for each $x\in\mathscr G_n$ there is $p_x$ between
$\rho_n\check{\mathfrak D}(x)$ and $\rho_n\mathfrak D(x)$, hence in
$[c\rho_n,\frac14]$, with
\begin{align}
\begin{split}
\label{eq-mvt-step}
\bigl|b_{n''}(s';\rho_n\check{\mathfrak D}(x))-b_{n''}(s';\rho_n\mathfrak D(x))\bigr|
 &\;=\; \rho_n\,g(x)\,\mathcal K(x),
\\ \text{where } \mathcal K(x) &:= \frac{b_{n''}(s';p_x)\,\bigl|s'-n''p_x\bigr|}{p_x(1-p_x)}.
\end{split}
\end{align}
We claim $\mathcal K\leq C/\rho_n$ on $\mathscr G_n$.  Since
$p_x(1-p_x)\geq3c\rho_n/4$, it suffices to bound the numerator by a constant,
and there are three cases.  If $|s'-n''p_x|\leq n''p_x$, then with
$x_{s'} := x_{s'}(n'',p_x)$, Lemma~\ref{lem-binomial-facts}\ref{binom-fact-a}
gives
$b_{n''}(s';p_x)|s'-n''p_x| = b_{n''}(s';p_x)\,\sigma(n'',p_x)\,|x_{s'}|
\leq C|x_{s'}|e^{-\zeta x_{s'}^2}\leq C'$.  If $s'\leq n''p_x/2$, then
$|s'-n''p_x|\leq n''p_x$ and the lower-tail bound gives
$b_{n''}(s';p_x)\leq e^{-n''p_x/8}$, so the numerator is at most
$\sup_{y>0}y\,e^{-y/8}\leq C$.  If $s'>2n''p_x$, then
$|s'-n''p_x|\leq s'\leq n'\rho_ns$, while the upper-tail bound gives
$b_{n''}(s';p_x)\leq1-\mathcal B_{n''}(s'-1;p_x)\leq e^{-n''p_x/8}
\leq e^{-n''c\rho_n/8}$, so the numerator is
$O\bigl(n\rho_ne^{-n''c\rho_n/8}\bigr) = o(1)$ by $n\rho_n/\log n\to\infty$.
Combining the claim with \eqref{eq-mvt-step} and \eqref{eq-g-global},
\begin{align}
\label{eq-step1}
n'\rho_n\!\int_{\mathscr G_n}\!
 \bigl|b_{n''}(s';\rho_n\check{\mathfrak D})-b_{n''}(s';\rho_n\mathfrak D)\bigr|\,dx
 \;\leq\; n'\rho_n^2\cdot\frac{C}{\rho_n}\int_0^1\! g
 \;\leq\; C\,n\rho_n^{\,2+\beta} \;\rightarrow\; 0 .
\end{align}

\emph{Step 2: the main term.}  Pushing forward the Lebesgue measure on
$[0,1]$ by $\mathfrak D$, which by Remark~\ref{remark-inverse} and
\ref{condition-d1} sends $dx$ to $f_{\mathfrak D}(\lambda)\,d\lambda$, and
substituting $p = \rho_n\lambda$, the good set corresponds exactly to
$p\in[0,\frac14]$, so that
\begin{align}
\label{eq-main-term}
n'\rho_n\!\int_{\mathscr G_n}\! b_{n''}(s';\rho_n\mathfrak D(x))\,dx
 = n'\!\int_0^{1/4}\! b_{n''}(s';p)\,f_{\mathfrak D}(p/\rho_n)\,dp
 = f_{\mathfrak D}(s)\frac{n'}{n''+1} + \mathrm{(I)} - \mathrm{(II)},
\end{align}
by Lemma~\ref{lem-binomial-facts}\ref{binom-fact-e}, where
\begin{align*}
\mathrm{(I)} := n'\!\int_0^{1/4}\! b_{n''}(s';p)
 \bigl[f_{\mathfrak D}(p/\rho_n)-f_{\mathfrak D}(s)\bigr]dp,
\qquad
\mathrm{(II)} := f_{\mathfrak D}(s)\,n'\!\int_{1/4}^{1}\! b_{n''}(s';p)\,dp .
\end{align*}
For $\mathrm{(II)}$, $b_{n''}(s';p)\leq e^{-n''/32}$
for $p>\frac14$ gives $\mathrm{(II)}\leq f_{\mathfrak D}(s)\,n'e^{-n''/32}\to0$.
For $\mathrm{(I)}$, fix $\eta>0$ with $[s-\eta,s+\eta]\subset\mathrm{supp}\,F$
and split at $|p/\rho_n-s|\leq\eta$.  On that set
$|f_{\mathfrak D}(p/\rho_n)-f_{\mathfrak D}(s)|\leq\omega(\eta)$, the modulus of
continuity of $f_{\mathfrak D}$ at $s$ supplied by \ref{condition-d1}, while
$n'\int_0^1b_{n''}(s';p)\,dp = n'/(n''+1)\leq2$, so that part is at most
$2\omega(\eta)$.  Off that set the Chernoff bounds of
Lemma~\ref{lem-binomial-facts}\ref{binom-fact-b} in both tails give
$b_{n''}(s';p)\leq e^{-\zeta_1n\rho_n\eta^2}$, while
$n'\int_0^{1/4}[f_{\mathfrak D}(p/\rho_n)+f_{\mathfrak D}(s)]\,dp
\leq n'\bigl(\rho_n+f_{\mathfrak D}(s)\bigr)$, so that part is
$O\bigl(n\,e^{-\zeta_1n\rho_n\eta^2}\bigr) = o(1)$ by
$n\rho_n/\log n\to\infty$.  Hence $\limsup_n|\mathrm{(I)}|\leq2\omega(\eta)$ for
every admissible $\eta$, and letting $\eta\downarrow0$ gives $\mathrm{(I)}\to0$.

\emph{Step 3: the complement of the good set.}  Since
$\rho_n\check{\mathfrak D}\geq c\rho_n$ and $s'\leq2n''\rho_ns$ for $n$ large,
\eqref{eq-pmf-sup} and \eqref{eq-good-complement} give
\begin{align}
\label{eq-far-mu}
\mu_s\bigl(\mathscr G_n^c\bigr)
 \;\leq\; n'\rho_n\cdot\mathrm{Leb}\bigl(\mathscr G_n^c\bigr)\cdot
 \frac{C}{\sqrt{n\rho_n}}
 \;\leq\; C\sqrt n\,\rho_n^{\,5/2+\beta}
 \;=\; C\bigl(n\rho_n^{\,5+2\beta}\bigr)^{1/2}\;\rightarrow\;0,
\end{align}
since $n\rho_n^{5+2\beta}\leq n\rho_n^{2+\beta}$.  Combining Steps~1--3 with
\eqref{eq-main-term} proves \ref{localization-item-a}.

\medskip
We now prove \ref{localization-item-b}.  Let $\eta>0$ and split
$\{|\mathfrak D-s|>\eta\}$ into three parts.  On $\mathscr G_n^c$ the
contribution is at most $\mu_s(\mathscr G_n^c)\to0$ by \eqref{eq-far-mu}.  On
$\mathscr G_n\cap\{|\check{\mathfrak D}-s|>\eta/2\}$, the Chernoff bounds in both
tails give $b_{n''}(s';\rho_n\check{\mathfrak D}(x))\leq
e^{-\zeta_1n\rho_n\eta^2/4}$, so this part contributes at most
$n'\rho_ne^{-\zeta_1n\rho_n\eta^2/4}\to0$.  On the remaining set
$\{|\mathfrak D-s|>\eta,\ |\check{\mathfrak D}-s|\leq\eta/2\}$ we have
$g\geq\eta/2$, so by Markov's inequality and \eqref{eq-g-global} its Lebesgue
measure is at most $2C_w\rho_n^{1+\beta}/\eta$; with \eqref{eq-pmf-sup}, its
contribution is at most
\begin{align*}
n'\rho_n\cdot\frac{2C_w\rho_n^{1+\beta}}{\eta}\cdot\frac{C}{\sqrt{n\rho_n}}
 \;=\; \frac{C}{\eta}\bigl(n\rho_n^{\,3+2\beta}\bigr)^{1/2}\;\rightarrow\;0,
\end{align*}
since $n\rho_n^{3+2\beta}\leq n\rho_n^{2+\beta}$.
\end{proof}

\subsection{Convergence of H\'ajek Projection Terms}
With the integrated local limit theorem in hand, we can now establish $L^2$ limits for H\'ajek projections, a key step for deriving the limiting distribution of the signal term.  We have the following result: 

\begin{lemma}[$L^2$ limits]
\label{lem-hajek-limits}
Under the assumptions \ref{condition-m1}, \ref{condition-g1}, \ref{condition-d1}, and   \ref{condition-s1}  for every $s$ in
the interior of $\mathrm{supp}\,F$,
\begin{align*}
\psi_s\to \psi_s^\infty :=\mathbbm 1(\mathfrak D\le s)-F(s)
\quad\text{and}\quad
\phi_s\to \phi_s^\infty:=  f_{\mathfrak D}(s)\bigl(s-w(x_s,\cdot)\bigr)
\quad\text{in }L^2 .
\end{align*}
\end{lemma}

\begin{proof}
\emph{The term $\psi_s$.}  Fix $x$ with $\mathfrak D(x)\ne s$.  
Given $\xi_1 = x$, the degree $D_1$ has the
$\mathrm{Bin}(n',\rho_n\check{\mathfrak D}(x))$ distribution, with mean
$n'\rho_n\check{\mathfrak D}(x)$ and variance at most its mean, while
$s'/(n'\rho_n)\to s$ and $\check{\mathfrak D}(x)\to\mathfrak D(x)$.
Chebyshev's inequality therefore gives
\begin{align*}
\mathcal B_{n'}\bigl(s';\rho_n\check{\mathfrak D}(x)\bigr)\to\mathbbm 1\bigl(\mathfrak D(x)<s\bigr).
\end{align*}
The set $\{\mathfrak D = s\}$ is Lebesgue-null by \ref{condition-d1}, so the
limit agrees almost everywhere with $\mathbbm 1(\mathfrak D(x)\leq s)$, and the
functions are bounded by $1$, so dominated convergence gives both the stated
$L^2$ limit and $F^\star(s)\to F(s)$.

\emph{The term $\phi_s$.}  Write:
\begin{align*}
\phi_s(\xi) = -\int_0^1a(x,\xi)\,\mu_s(dx)
 = -\int_0^1\bigl[a(x,\xi)-a(x_s,\xi)\bigr]\mu_s(dx)
 - a(x_s,\xi)\,\mu_s([0,1]).
\end{align*}
The second term converges to $-f_{\mathfrak D}(s)(w(x_s,\cdot)-s)$ in $L^2$: by
Lemma~\ref{lem-mu-limit}(a) the mass converges to $f_{\mathfrak D}(s)$, and
$a(x_s,\xi) = \check w(x_s,\xi)-\check{\mathfrak D}(x_s)\to w(x_s,\xi)-s$ in
$L^2$ by monotone convergence together with $E[w^2(x_s,\xi)]<\infty$ (implicit
in the $L^2$-continuity of \ref{condition-g1}) and $\mathfrak D(x_s) = s$.

For the first term, Minkowski's integral inequality in $L^2(d\xi)$ gives
\begin{align*}
\Bigl\|\int_0^1\bigl[a(x,\cdot)-a(x_s,\cdot)\bigr]\mu_s(dx)\Bigr\|_{L^2}
 \le \int_0^1\bigl\|a(x,\cdot)-a(x_s,\cdot)\bigr\|_{L^2}\,\mu_s(dx),
\end{align*}
which we split at $|\mathfrak D(x)-s|\le\eta$.  On that set, truncation being
$1$-Lipschitz,
\begin{align*}
\|a(x,\cdot)-a(x_s,\cdot)\|_{L^2}
 \le 2\|w(x,\cdot)-w(x_s,\cdot)\|_{L^2},
\end{align*}
which is small uniformly over the set when $\eta$ is small, by the
$L^2$-continuity in \ref{condition-g1} and continuity of $\mathfrak D^{-1}$ at
$s$; multiplying by the bounded mass $\mu_s([0,1])$ leaves a quantity that
vanishes as $\eta\downarrow0$.  On the complement, bound
$\|a(x,\cdot)-a(x_s,\cdot)\|\le\|a(x,\cdot)\|+\|a(x_s,\cdot)\|$ and use
$\|a(x,\cdot)\|\le\sqrt{\check{\mathfrak D}(x)/\rho_n}$ from
\eqref{eq-two-bounds}; the term involving the constant $\|a(x_s,\cdot)\|$
vanishes by Lemma~\ref{lem-mu-limit}(b).  It remains to show
\begin{align}
\label{eq-hajek-far}
\int_{\{|\mathfrak D-s|>\eta\}}\sqrt{\check{\mathfrak D}(x)/\rho_n}\;\mu_s(dx)\;\rightarrow\;0,
\qquad \mu_s(dx) = n'\rho_n\,b_{n''}\bigl(s';\rho_n\check{\mathfrak D}(x)\bigr)\,dx .
\end{align}
We split $\{|\mathfrak D-s|>\eta\}$ into two regions, according to whether the
clipped degree $\check{\mathfrak D}(x)$ is also far from $s$.

\emph{Region 1: $|\check{\mathfrak D}(x)-s|>\eta/2$.}  For $n$ large, $s'$ deviates
from the mean $n''\rho_n\check{\mathfrak D}(x)$ of $\mathrm{Bin}(n'',\rho_n\check{\mathfrak D}(x))$
by at least $\frac14n''\rho_n|\check{\mathfrak D}(x)-s|$, and the variance is at
most $n''\rho_n\check{\mathfrak D}(x)$, so Bernstein's inequality gives
\begin{align*}
b_{n''}\bigl(s';\rho_n\check{\mathfrak D}(x)\bigr)\leq
\begin{cases}
2e^{-cn\rho_n\eta^2}, & \check{\mathfrak D}(x)\leq2s,\\
2e^{-cn\rho_n\check{\mathfrak D}(x)}, & \check{\mathfrak D}(x)>2s,
\end{cases}
\end{align*}
for some $c>0$ depending on $s$.  Using $\check{\mathfrak D}\leq\rho_n^{-1}$ on
$\{\check{\mathfrak D}>2s\}$, the integrand in \eqref{eq-hajek-far} is therefore at
most
\begin{align*}
2n'\rho_n\sqrt{2s/\rho_n}\,e^{-cn\rho_n\eta^2}\leq Cn\rho_n^{1/2}e^{-cn\rho_n\eta^2}
\quad\text{or}\quad
2n'\rho_n\cdot\rho_n^{-1}e^{-2cn\rho_ns}\leq Cne^{-2cn\rho_ns}
\end{align*}
on the two parts.  Both bounds tend to zero because $n\rho_n/\log n\to\infty$.

\emph{Region 2: $|\mathfrak D(x)-s|>\eta$ and $|\check{\mathfrak D}(x)-s|\leq\eta/2$.}
Here the binomial mass is not small, and we use instead that the region is
small.  Since $\check{\mathfrak D}\leq\mathfrak D$, the two conditions force
$\mathfrak D(x)>s+\eta$, and hence $g(x) = \mathfrak D(x)-\check{\mathfrak D}(x)>\eta/2$.
By Markov's inequality and \eqref{eq-g-global}, the region has Lebesgue
measure at most $2C_w\rho_n^{1+\beta}/\eta$.  On it,
$\sqrt{\check{\mathfrak D}/\rho_n}\leq\sqrt{2s/\rho_n}$, and \eqref{eq-pmf-sup} gives
$b_{n''}(s';\rho_n\check{\mathfrak D}(x))\leq C(n\rho_n)^{-1/2}$.  Its contribution
to \eqref{eq-hajek-far} is therefore at most
\begin{align*}
n'\rho_n\cdot\sqrt{\frac{2s}{\rho_n}}\cdot\frac{C}{\sqrt{n\rho_n}}\cdot\frac{2C_w\rho_n^{1+\beta}}{\eta}
 \;\leq\;\frac C\eta\bigl(n\rho_n^{\,2+2\beta}\bigr)^{1/2}\;\rightarrow\;0,
\end{align*}
since $n\rho_n^{2+2\beta}\leq n\rho_n^{2+\beta}\to0$.  This proves
\eqref{eq-hajek-far}, and letting $\eta\downarrow0$ completes the proof.
\end{proof}

Now, the limiting variance associated with the H\'ajek projection is readily computed.  We have the following result:
\begin{corollary}[Convergence of Covariance Terms]
\label{cor-covariance}
Under conditions \ref{condition-m1}, \ref{condition-g1}, \ref{condition-d1}, and   \ref{condition-s1}, for any $s,t \in \text{support } F$ and $i \in [3]$,
\begin{align*}
\Sigma_{n,i}(s,t) \rightarrow \Sigma_i(s,t). 
\end{align*}
\end{corollary}

\begin{proof}
The convergences $\Sigma_{n,1}(s,t) \rightarrow \Sigma_1(s,t)$ and
$\Sigma_{n,2}(s,t) \rightarrow \Sigma_2(s,t)$ immediately follow from Cauchy-Schwarz and the $L^2$ convergence established in Lemma \ref{lem-hajek-limits}.  For the third claim, observe that:
\begin{align*}
E[\psi_s^\infty \phi_t^\infty] &= E[\mathbbm{1}(\mathfrak{D} \leq s) f_\mathfrak{D}(t) \cdot \{t- w(x_t, \cdot)\}]
\\&= f_\mathfrak{D}(t) \cdot \left\{tF(s) - \int_{\mathfrak{D} \leq s} w(\mathfrak{D}^{-1}(t), \xi) dP \right\}. 
\end{align*}
An analogous expression holds for $E[\psi_t^\infty \phi_s^\infty]$. Thus, the claim follows by Cauchy-Schwarz and $L^2$ convergence. 
\end{proof}

\subsection{Negligibility of the Higher-Order Components}
\label{sec-degenerate}

Throughout this section $s,t$ denote a generic pair of evaluation points and
$\bar t := \max_{k\leq K}t_k+1$; by bilinearity of the covariance it suffices
to argue for such a pair.  We write $R_n(s,t)$ for $n'$ times the sum of the
Hoeffding components of order $|S|\geq2$ appearing in
\eqref{eq-gap-two-sources}, and prove that $R_n(s,t)\to0$.

Our analysis has rather involved case analysis, which warrants additional discussion. The cases arise from two separate sources: the index set $S$ may meet the two factors being paired in four ways, and the domain of integration splits into three regions on which different bounds for $\|a(x,\cdot)\|_{L^2}$ are available.  We discuss each of these sources in more detail below.  

\paragraph{On the Four Set Types.}
The four set types arises from an asymmetry present in our functional that is not present with more classical functionals analyzed with Hoeffding decompositions, such as U-statistics.  
The
quantity being decomposed, $F_i(s'\mid\mathcal G)$, is a function of all
$n$ latent positions but is asymmetric in them: the index $i$ is
distinguished, since it is the vertex whose degree is counted. We refer to $i$ as the root of $F_i$. A given index set $S$ meets the root nodes in different ways, and the pairing
depends on whether $S$ contains the root of the first factor, of the second,
of both, or of neither.  Without loss of generality, one may use identical distributions of the nodes and bound the overall contribution of these higher order terms in terms of the root nodes $1$ and $2$. 

To see this, notice that, by orthogonality of the spaces $\mathcal{H}_S$, components indexed by distinct
sets do not interact, so
\begin{align*}
&\mathrm{Cov}\Bigl(\frac{1}{\sqrt n}\sum_{i=1}^{n}\sum_{|S| \geq 2}
 \Pi\bigl(F_i(s'\mid\mathcal G);S\bigr),\;
 \frac{1}{\sqrt n}\sum_{j=1}^{n}\sum_{|S| \geq 2}
 \Pi\bigl(F_j(t'\mid\mathcal G);S\bigr)\Bigr)
\\=& \frac{1}{n}\sum_{i,j \in [n]}\ \sum_{|S| \geq 2}
 \mathrm{Cov}\Bigl(\Pi\bigl(F_i(s');S\bigr),\,
 \Pi\bigl(F_j(t');S\bigr)\Bigr).
\end{align*}
For a fixed pair $i \neq j$, a permutation of $[n]$ carrying $(i,j)$ to $(1,2)$
also carries each $S$ to $\pi(S)$, and therefore permutes the inner sum among
its own terms; by exchangeability of the latent positions the inner sum is thus
the same for every ordered pair $i \neq j$, and likewise for every $i = j$.
Hence, the above display is equal to:
\begin{align*}
\sum_{|S| \geq 2}
 E\bigl[\Pi(F_1(s');S)\,\Pi(F_1(t');S)\bigr]
\;+\; \underbrace{n'\sum_{|S| \geq 2}
 E\bigl[\Pi(F_1(s');S)\,\Pi(F_2(t');S)\bigr]}_{=\,R_n(s,t)},
\end{align*}
where the first sum is handled in Lemma \ref{lem-diagonal}, and the four types of $S$ for the summands in $R_n(s,t)$ are handled in Lemma \ref{lem-degenerate}.

\paragraph{On the Three Regions.}
It is helpful to see informally what has to be estimated, since this makes the partition more clear.  Consider a set $S$ satisfying $\{1,2\} \subseteq S$, of size $d$, and
recall from Lemma~\ref{lem-multilinear} that each of the two paired components
carries $d$ centered factors $a$ and one $d$-th order backward
difference of a binomial distribution function.  Integrating out the $d$ latent
positions indexed by $S$ pairs the centered factors against one another through
\begin{align}
\label{eq-C-def}
C(x,y) := E_\xi\bigl[a(x,\xi)\,a(y,\xi)\bigr],
\qquad
|C(x,y)| \leq \bigl\|a(x,\cdot)\bigr\|_{L^2}\bigl\|a(y,\cdot)\bigr\|_{L^2},
\end{align}
so that the contribution of all sets of this type and size is
of the form
\begin{align*}
n\cdot\binom{n''}{d}\,\rho_n^{2d}\;\times\;
 \bigl(\text{a }d\text{-fold product of } C\bigr)\;\times\;
 \bigl(\text{two }d\text{-th order differences}\bigr),
\end{align*}
and the remaining types are similar.  The differences turn out to be of size
$\sigma_*^{-d}$, where $\sigma_*^2 \asymp n\rho_n$ is a lower bound for the
variance of the relevant binomial law; since
$\binom{n''}{d} \approx (n\rho_n)^{d}/(d!\,\rho_n^{d})$, the powers of $n$
cancel and the display behaves like a geometric series in $d$ of ratio
$\rho_n\|a(x,\cdot)\|_{L^2}^2$.  Everything therefore turns on the size
of $\|a(x,\cdot)\|_{L^2}$, for which two bounds are available:
\begin{align}
\label{eq-two-bounds}
\bigl\|a(x,\cdot)\bigr\|_{L^2}^2 \;\leq\; E\bigl[\check w^2(x,\xi)\bigr]
 \;\leq\; \min\Bigl(
 \underbrace{\textstyle\sup_{\{\mathfrak D\leq R\}}E\bigl[w^2(\cdot,\xi)\bigr]}
 _{\text{\ref{condition-l1}, if }\mathfrak D(x)\leq R},\ \
 \underbrace{\rho_n^{-1}\check{\mathfrak D}(x)}_{\check w\leq\rho_n^{-1}}\Bigr).
\end{align}
Where the expected degree is bounded the first applies, the ratio is
$O(\rho_n)$, and the series converges with a spare factor of $\rho_n$ at each
order.  Where it is not, only the second is available, and since it carries
the factor $\rho_n^{-1}$, the ratio is of order $\check{\mathfrak D}(x)$ and
the spare factor is lost; the series must then be summed by other means, and a
further distinction appears according to whether $\check{\mathfrak D}(x)$ is
itself large, in which case the binomial masses involved are exponentially
small, or not.
 
Accordingly, with $\bar t := \max_{k\leq K}t_k+1$, we partition $[0,1]$ into
\begin{align}
\label{eq-regions}
\mathscr A := \bigl\{\check{\mathfrak D}\leq3\bar t,\ \mathfrak D\leq6\bar t\bigr\},
\qquad
\mathscr B := \bigl\{\check{\mathfrak D}\leq3\bar t,\ \mathfrak D>6\bar t\bigr\},
\qquad
\mathscr C := \bigl\{\check{\mathfrak D}>3\bar t\bigr\} .
\end{align}
 The middle region is not vacuous: truncation can
reduce an arbitrarily large expected degree to a moderate one, so a bound on
$\check{\mathfrak D}$ carries no information about $\mathfrak D$, and on
$\mathscr B$ neither the first bound in \eqref{eq-two-bounds} nor the
exponential estimates available on $\mathscr C$ apply.  On $\mathscr A$ the
first bound and \ref{condition-l1} supply the two constants
\begin{align}
\label{eq-A-def}
C_a := \sup_{\{\mathfrak D\leq6\bar t\}}\bigl\|a(x,\cdot)\bigr\|_{L^2}<\infty,
\quad
g(x) := \mathfrak D(x)-\check{\mathfrak D}(x) \leq C_g\rho_n^{1+\beta}
 \quad\text{on }\mathscr A,
\end{align}
the second by Markov's inequality applied to \ref{condition-l1} with
$R = 6\bar t$, as in \eqref{eq-g-global}
.  In what follows, we collect various inequalities that will be useful in different regions.

\begin{lemma}[Region estimates]
\label{lem-regions}
For all $n$ large enough, depending only on $s,t$ and the constants in
\ref{condition-m1}--\ref{condition-d1}:
\begin{enumerate}[(a)]
\item\label{reg-univ} $\|a(x,\cdot)\|_{L^2}\leq\sqrt{\check{\mathfrak{D}}(x)/\rho_n}$
for every $x$, hence $\leq\sqrt{3\bar t/\rho_n}$ on
$\mathscr A\cup\mathscr B$;
\item\label{reg-I} $\mathrm{Leb}(\mathscr B)\leq C_w\rho_n^{1+\beta}/(3\bar t)$,
where
$C_w := \|w\|^{2+\beta}_{L^{2+\beta}}$ is the constant of
\eqref{eq-g-global};
\item\label{reg-F} $s'\leq n''\rho_n\check{\mathfrak{D}}(x)/2$ on $\mathscr C$;
\item\label{reg-p} $\rho_n\check{\mathfrak{D}}(x)\in[c\rho_n,\,1/4]$ on
$\mathscr A\cup\mathscr B$;
\item\label{reg-kernel} $\bigl|\mathsf D^{d}\mathcal B_{n'-d}
\bigl(s';\rho_n\check{\mathfrak{D}}(x)\bigr)\bigr|\leq\Gamma_d\,\sigma_*^{-d}$ on
$\mathscr A\cup\mathscr B$, for every $1\leq d\leq n''/2$, where
$\sigma_*^2 := n''c\rho_n/4$ and $\Gamma_d$ is as in
Lemma~\ref{lem-binomial-facts}\ref{binom-fact-c};
\item\label{reg-h} $\displaystyle\int_0^1 a(x,y)^2\,dx\;\leq\;
\frac{C}{\rho_n}\bigl(1+\check{\mathfrak{D}}(y)\bigr)$ for every $y$.
\end{enumerate}
\end{lemma}

\begin{proof}
For \ref{reg-univ}, $a(x,\cdot)$ is $\check w(x,\cdot)$ centered, so, using
$\check w\leq\rho_n^{-1}$,
\begin{align*}
\|a(x,\cdot)\|^2_{L^2}\leq E[\check w^2(x,\xi)]\leq
\rho_n^{-1}E[\check w(x,\xi)] = \check{\mathfrak{D}}(x)/\rho_n ;
\end{align*}
the second assertion follows from
$\check{\mathfrak{D}}\leq3\bar t$ on $\mathscr A\cup\mathscr B$.

For \ref{reg-I}, on $\mathscr B$ we have $g = \mathfrak D-\check{\mathfrak{D}}>3\bar t$,
so Markov's inequality and the global bias bound \eqref{eq-g-global} give
$\mathrm{Leb}(\mathscr B)\leq(3\bar t)^{-1}\int_0^1g\leq
C_w\rho_n^{1+\beta}/(3\bar t)$.

For \ref{reg-F}, on $\mathscr C$ we have
$n''\rho_n\check{\mathfrak{D}}(x)/2>n''\rho_n\cdot3\bar t/2>n'\rho_ns\geq s'$ for $n$
large, since $\bar t\geq s+1$.

For \ref{reg-p}, the lower bound is $\check{\mathfrak{D}}\geq c$, which holds because
$w\geq c$ by \ref{condition-g1} and $\rho_n^{-1}\geq c$ for $n$ large; the
upper bound is $\rho_n\check{\mathfrak{D}}\leq3\bar t\rho_n\leq1/4$ for $n$ large, since
$\bar t$ is a constant and $\rho_n\to0$.

Part \ref{reg-kernel} is Lemma~\ref{lem-binomial-facts}\ref{binom-fact-c},
whose hypothesis $p\leq1/4$ is supplied by \ref{reg-p}, together with the
variance bound
\begin{align*}
\sigma^2\bigl(n'-d,\,\rho_n\check{\mathfrak{D}}(x)\bigr)
 = (n'-d)\,\rho_n\check{\mathfrak{D}}(x)\bigl(1-\rho_n\check{\mathfrak{D}}(x)\bigr)
 \;\geq\; \frac{n''}{2}\cdot c\rho_n\cdot\frac34 \;\geq\; \sigma_*^2,
\end{align*}
valid for $d\leq n''/2$ by \ref{reg-p}.  Note that \ref{reg-kernel} is claimed
only on $\mathscr A\cup\mathscr B$: on $\mathscr C$ the parameter
$\rho_n\check{\mathfrak{D}}(x)$ need not be bounded away from $1$, and no Fourier bound
is used there.

For \ref{reg-h}, by $\check w\leq\rho_n^{-1}$ and the symmetry of $w$,
\begin{align*}
\int_0^1a(x,y)^2dx &\leq 2\int_0^1\check w(x,y)^2dx + 2E[\mathfrak D^2]
 \\&\leq \frac{2}{\rho_n}\int_0^1\check w(x,y)\,dx + 2E[\mathfrak D^2]
\leq \frac{2\,\check{\mathfrak{D}}(y)}{\rho_n} + C,
\end{align*}
with $E[\mathfrak D^2]\leq\|w\|^2_{L^2}<\infty$ by \ref{condition-m1}.
\end{proof}

\subsubsection*{Two summed bounds}
We now prove two bounds that are useful on $\mathscr B$ and $\mathscr C$.  Both bound the entire sum over $d$ at once, in terms of a quantity that is small on those regions: the Lebesgue measure of
$\mathscr B$, on $\mathscr C$, the conditional
probability that the degree falls at or below $s'$, which is exponentially
small since $s'$ lies deep in the lower tail.

\begin{lemma}[Summed bound, distribution function]
\label{lem-parseval-cdf}
For every $x\in[0,1]$ and every $s$,
\begin{align}
\label{eq-parseval}
V_s(x) := \sum_{d\geq1}\binom{n''}{d}\rho_n^{2d}
 \bigl\|a(x,\cdot)\bigr\|^{2d}
 \Bigl(\mathsf D^{d}\mathcal B_{n'-d}\bigl(s';\rho_n\check{\mathfrak D}(x)\bigr)\Bigr)^2
 \;\leq\; \min\Bigl(\tfrac14,\;
 \mathcal B_{n'}\bigl(s';\rho_n\check{\mathfrak D}(x)\bigr)\Bigr).
\end{align}
In particular $V_s(x)\leq e^{-n''\rho_n\check{\mathfrak D}(x)/8}$ on $\mathscr C$.
\end{lemma}

\begin{proof}
Fix $\xi_1=x$.  By Lemma~\ref{lem-multilinear} the summands $K_S(s)$ are
orthogonal given $\xi_1$, with $K_\emptyset(s)$ the conditional mean, and since
the $(\xi_j)_{j\in S}$ are i.i.d.\ and the difference factor is deterministic
given $\xi_1$,
\begin{align}
\label{eq-Ksq}
E\bigl[K_S(s)^2\mid\xi_1=x\bigr]
 = \rho_n^{2|S|}\bigl\|a(x,\cdot)\bigr\|^{2|S|}
 \Bigl(\mathsf D^{|S|}\mathcal B_{n'-|S|}\bigl(s';\rho_n\check{\mathfrak D}(x)\bigr)\Bigr)^2 .
\end{align}
Summing \eqref{eq-Ksq} over the $\binom{n''}{d}$ sets $S\subseteq\{3,\ldots,n\}$
with $|S|=d$, and over $d\geq1$, yields a subsum of nonnegative terms of
$\sum_{S\neq\emptyset}E[K_S(s)^2\mid\xi_1=x]
=\mathrm{Var}\bigl(F_1(s'\mid\mathcal G)\mid\xi_1=x\bigr)$.  Since
$F_1(s'\mid\mathcal G)\in[0,1]$, that variance is at most $1/4$, and
also at most $E[F_1(s'\mid\mathcal G)\mid\xi_1=x]
=\mathcal B_{n'}(s';\rho_n\check{\mathfrak D}(x))$, using $\mathrm{Var}(X)\leq E[X^2]\leq E[X]$ for
$X\in[0,1]$.  The final assertion follows from
Lemma~\ref{lem-regions}\ref{reg-F} and
Lemma~\ref{lem-binomial-facts}\ref{binom-fact-b}.
\end{proof}

\begin{lemma}[Summed bound, mass function]
\label{lem-parseval-pmf}
For every $x\in[0,1]$ and every $s$,
\begin{align}
\label{eq-parseval-pmf}
\begin{split}
W_s(x) &:= \sum_{d\geq1}\binom{n''}{d-1}\rho_n^{2(d-1)}
 \bigl\|a(x,\cdot)\bigr\|^{2(d-1)}
 \Bigl(\mathsf D^{d}\mathcal B_{n'-d}\bigl(s';\rho_n\check{\mathfrak D}(x)\bigr)\Bigr)^2
\\ &\;\leq\; b_{n''}\bigl(s';\rho_n\check{\mathfrak D}(x)\bigr).
\end{split}
\end{align}
\end{lemma}

\begin{proof}
Let $H(p):=P\bigl(\sum_{j\neq1,2}X_j=s'\bigr)$ as a function of
$p=(p_{1j})_{j\neq1,2}\in[0,1]^{n''}$, the $X_j$ being independent
$\mathrm{Bernoulli}(p_{1j})$ as in the proof of Lemma~\ref{lem-multilinear}.
Being a difference of two values of a distribution function, $H$ is multilinear,
with
\begin{align*}
\partial_{S'}H = (-1)^{|S'|}\mathsf D^{|S'|}
 \bigl(\mathsf D\mathcal B_{n''-|S'|}\bigr)(s')
 = (-1)^{|S'|}\mathsf D^{|S'|+1}\mathcal B_{n''-|S'|}(s')
\end{align*}
at the common parameter, where $n''-|S'|=n'-(|S'|+1)$.  Expanding about
$\rho_n\check{\mathfrak D}(x)$ and substituting
$p_{1j}-\rho_n\check{\mathfrak D}(x)=\rho_na(x,\xi_j)$,
\begin{align*}
P\bigl(D_{1\setminus2}=s'\mid\mathcal G\bigr)
 &= \sum_{S'\subseteq[n]\setminus\{1,2\}}L_{S'},
\\
L_{S'} &:= (-\rho_n)^{|S'|}\prod_{j\in S'}a(x,\xi_j)\,
 \mathsf D^{|S'|+1}\mathcal B_{n'-|S'|-1}\bigl(s';\rho_n\check{\mathfrak D}(x)\bigr),
\end{align*}
and, exactly as in Lemma~\ref{lem-multilinear}, the $L_{S'}$ are orthogonal
given $\xi_1=x$, with $E[L_{S'}^2\mid\xi_1=x]$ equal to the $d=|S'|+1$ summand
of \eqref{eq-parseval-pmf}.  Since
$P(D_{1\setminus2}=s'\mid\mathcal G)\in[0,1]$,
\begin{align*}
W_s(x) &= E\Bigl[P\bigl(D_{1\setminus2}=s'\mid\mathcal G\bigr)^2
 \Bigm|\xi_1=x\Bigr]
\\ &\leq E\Bigl[P\bigl(D_{1\setminus2}=s'\mid\mathcal G\bigr)\Bigm|\xi_1=x\Bigr]
 = b_{n''}\bigl(s';\rho_n\check{\mathfrak D}(x)\bigr).
\end{align*}
The result follows.
\end{proof}

Two consequences of \eqref{eq-parseval-pmf} are used below:
\begin{align}
\label{eq-W-bounds}
\begin{split}
&\sup_x W_s(x)^{1/2}\leq\frac{C}{(n\rho_n)^{1/4}},
\\
&\int_0^1W_s(x)\,dx \leq \int_0^1 b_{n''}\bigl(s';\rho_n\check{\mathfrak D}(x)\bigr)dx
 = \frac{\mu_s([0,1])}{n'\rho_n} \leq \frac{C}{n\rho_n},
\end{split}
\end{align}
the first by Lemma~\ref{lem-binomial-facts}\ref{binom-fact-a} at the mode, the
second by Lemma~\ref{lem-mu-limit}(a).

\subsubsection*{An integrated kernel bound on $\mathscr A$}
On $\mathscr A$ the two exact identities of the preceding subsection are
unavailable, since the quantities they bound are $O(1)$ there.  Applying
Lemma~\ref{lem-regions}\ref{reg-kernel} directly would bound each difference by
$\Gamma_d\sigma_*^{-d}$ uniformly and cost a factor $\mathrm{Leb}(\mathscr A)\asymp1$,
which discards the fact that the difference is not uniformly of that size: by
Lemma~\ref{lem-binomial-facts}\ref{binom-fact-a} it decays in the standardized
lattice distance $x_{s'}$, and is therefore appreciable only where
$\check{\mathfrak D}(x)$ lies within $O\bigl(\sigma_*/(n\rho_n)\bigr) =
O\bigl((n\rho_n)^{-1/2}\bigr)$ of $s$.  Integrating rather than taking a
supremum recovers that window, and so gains a factor $(n\rho_n)^{-1/2}$ over
the trivial bound.  We have the following bound:
 
\begin{lemma}
Under Condition \ref{condition-d1}, 
\label{lem-Jd}
for $1\leq d\leq\sigma_*$,
\begin{align}
\label{eq-Jd}
J_d(s) := \int_{\mathscr A}\bigl|\mathsf D^{d}\mathcal B_{n'-d}
 \bigl(s';\rho_n\check{\mathfrak D}(x)\bigr)\bigr|\,dx
 \;\leq\; C\Gamma_d\sigma_*^{-d}\Bigl[
 \frac{\sigma_*}{n\rho_n}\bigl(d+\sqrt{d\log(n\rho_n)}\bigr)
 + \rho_n^{1+\beta}\Bigr],
\end{align}
and $J_d(s)\leq\Gamma_d\sigma_*^{-d}$ for $\sigma_*<d\leq n''/2$.
\end{lemma}

\begin{proof}
Split $\mathscr A$ according to whether
\begin{align*}
\bigl|x_{s'}\bigl(n'-d,\rho_n\check{\mathfrak D}(x)\bigr)\bigr|\leq
R_d:=d/\sigma_*+K_d,\qquad K_d^2:=(2/\zeta)\log(2^d\sigma_*^d/\Gamma_d),
\end{align*}
so that $K_d\leq C\sqrt{d\log(n\rho_n)}$.  On that set bound the integrand by
$\Gamma_d\sigma_*^{-d}$ using Lemma~\ref{lem-regions}\ref{reg-kernel}; the set
is contained in $\{x\in\mathscr A:|\check{\mathfrak D}(x)-s|\leq h\}$ with
$h\asymp R_d\sigma_*/(n''\rho_n)$, and by \eqref{eq-A-def} together with the
pushforward of Lebesgue measure under $\mathfrak D$
(Remark~\ref{remark-inverse} and \ref{condition-d1}),
\begin{align}
\label{eq-measure-bound}
\mathrm{Leb}\bigl(\{x\in\mathscr A:|\check{\mathfrak D}(x)-s|\leq h\}\bigr)
 \leq \mathrm{Leb}\bigl(\{|\mathfrak D-s|\leq h+C_g\rho_n^{1+\beta}\}\bigr)
 \leq 2C_f\bigl(h+C_g\rho_n^{1+\beta}\bigr),
\end{align}
where $C_f := \sup\{f_{\mathfrak D}(\lambda):|\lambda-s|\leq\eta_0\}$ for a fixed
$\eta_0>0$ with $[s-\eta_0,s+\eta_0]$ in the interior of $\mathrm{supp}\,F$,
finite by \ref{condition-d1}.  This gives the two terms in \eqref{eq-Jd}.  Off that set, the $\ell^1$ bound
of Lemma~\ref{lem-binomial-facts}\ref{binom-fact-d} together with the local
bound \ref{binom-fact-a} gives
\begin{align*}
|\mathsf D^{d}\mathcal B_{n'-d}|\leq2^{d}C\sigma_*^{-1}e^{-\zeta(|x_{s'}|-d/\sigma_*)^2},
\end{align*}
whose integral over the complement is at most
\begin{align*}
2^dC'\sigma_*^{-1}e^{-\zeta K_d^2/2}\,\frac{\sigma_*}{n\rho_n}
\leq\Gamma_d\sigma_*^{-d}\,\frac{C'\sigma_*}{n\rho_n}
\end{align*}
by the choice of $K_d$.  Here the
measure of the set where $x_{s'}$ lies in a unit interval is bounded using the
density bound $C_f$, which is valid where $|\check{\mathfrak D}(x)-s|\leq\eta_0/2$;
on the rest of $\mathscr A$, the binomial Chernoff bounds make the integrand
at most $2^de^{-cn\rho_n\eta_0^2}$, which is negligible.  For
$\sigma_*<d\leq n''/2$ the bound is Lemma~\ref{lem-regions}\ref{reg-kernel}
together with $\mathrm{Leb}(\mathscr A)\leq1$.
\end{proof}

With appropriate estimates in hand for the three regions, we are now ready to show negligibility of the off-diagonal higher-order projection terms. 

\begin{lemma}[Negligibility of $R_n(s,t)$]
\label{lem-degenerate}
Under \ref{condition-m1}, \ref{condition-l1}, \ref{condition-g1}, \ref{condition-d1}, and \ref{condition-s1}, $R_n(s,t)\to0$ for every $s,t$ in
the interior of $\mathrm{supp}\,F$.
\end{lemma}

\begin{proof}
The components of the Hoeffding decomposition of order $|S|\geq2$ fall into
four types, according to how many of the indices $1$ and $2$ belong to $S$:
\begin{enumerate}[(i)]
\item $1\notin S$ and $2\notin S$;
\item $2\in S$, $1\notin S$;
\item $1\in S$, $2\notin S$;
\item $1\in S$ and $2\in S$.
\end{enumerate}
Types (ii) and (iii) are exchanged by swapping $(1,s)$ with $(2,t)$, so only
(i), (ii) and (iv) require argument.  In each case
$|C(x,y)|\leq\|a(x,\cdot)\|\,\|a(y,\cdot)\|$ and the double integral over
$(x,y)$ is split over the pairs of regions \eqref{eq-regions}.  Types (i) and
(iv) are symmetric under $(x,s)\leftrightarrow(y,t)$, so six pairs suffice for
them; type (ii) is not, since the two kernels have different orders and the
factor $a(x,y)$ is not symmetric, and there all nine pairs are treated, the
pairs $\mathscr A$--$\mathscr B$ and $\mathscr B$--$\mathscr A$ separately.
Write
\begin{align}
\begin{split}
\label{eq-ud-def}
u_d(x) &:= \binom{n''}{d}^{1/2}\rho_n^{d}\bigl\|a(x,\cdot)\bigr\|^{d}
 \bigl|\mathsf D^{d}\mathcal B_{n'-d}\bigl(s';\rho_n\check{\mathfrak{D}}(x)\bigr)\bigr|,
\\ \tilde u_d(x) &:= \binom{n''}{d-1}^{1/2}\rho_n^{d-1}
 \bigl\|a(x,\cdot)\bigr\|^{d-1}
 \bigl|\mathsf D^{d}\mathcal B_{n'-d}\bigl(s';\rho_n\check{\mathfrak{D}}(x)\bigr)\bigr|
 \end{split}
\end{align}
so that $\sum_{d\geq1}u_d(x)^2\leq V_s(x)$ by \eqref{eq-parseval} and
$\sum_{d\geq1}\tilde u_d(x)^2\leq W_s(x)$ by \eqref{eq-parseval-pmf}, and let
$v_d(y)$, $\tilde v_d(y)$ denote the same quantities at $t$ and $y$.

The termwise bounds below use Lemma~\ref{lem-regions}\ref{reg-kernel} and
Lemma~\ref{lem-Jd}, which require $d\leq n''/2$.  For $d>n''/2$ we use instead
the crude bound $|\mathsf D^d\mathcal B_{n'-d}|\leq2^{d-1}$ of
Lemma~\ref{lem-binomial-facts}\ref{binom-fact-d}.  Since
$\binom{n''}{d}\leq2^{n''}$, $\rho_n\|a(x,\cdot)\|^2\leq3\bar t$ on
$\mathscr A\cup\mathscr B$ by Lemma~\ref{lem-regions}\ref{reg-univ}, and
$\int\!\!\int a(x,y)^2\,dx\,dy\leq C\rho_n^{-1}$ by
Lemma~\ref{lem-regions}\ref{reg-h}, the terms with $d>n''/2$ of every type,
restricted to pairs of regions in $\mathscr A\cup\mathscr B$, are at most
$Cn\,2^{n''}(C\rho_n)^{n''/2-1}$, which is negligible.  Pairs meeting
$\mathscr C$ and the pair $\mathscr B$--$\mathscr B$ are handled by the summed
bounds \eqref{eq-parseval} and \eqref{eq-parseval-pmf}, which cover all $d$ at
once.  Below, $d\leq n''/2$ throughout the termwise bounds.

\medskip\noindent\emph{Type (i).}  With $d:=|S|\geq2$ the contribution is
\begin{align}
\label{eq-type-i}
\sum_{d\geq2}n'\binom{n''}{d}\rho_n^{2d}\int\!\!\int|C(x,y)|^{d}
 \bigl|\mathsf D^{d}\mathcal B_{n'-d}(s';\rho_n\check{\mathfrak{D}}(x))\bigr|
 \bigl|\mathsf D^{d}\mathcal B_{n'-d}(t';\rho_n\check{\mathfrak{D}}(y))\bigr|\,dy\,dx,
\end{align}
whose $d$-th integrand is at most $n'u_d(x)v_d(y)$.

For a pair of regions $\mathcal X\times\mathcal Y$ with
$\mathcal X,\mathcal Y\subseteq\mathscr B\cup\mathscr C$, the Cauchy--Schwarz
inequality in the summation index $d$ gives
$\sum_du_dv_d\leq(\sum_du_d^2)^{1/2}(\sum_dv_d^2)^{1/2}$, so that the
contribution is at most
\begin{align}
\label{eq-CS-in-d}
n'\int_{\mathcal X}V_s(x)^{1/2}dx\int_{\mathcal Y}V_t(y)^{1/2}dy .
\end{align}
By \eqref{eq-parseval} and Lemma~\ref{lem-regions}\ref{reg-I},
\begin{align*}
\mathscr B\text{--}\mathscr B:\quad
 &\tfrac14n'\mathrm{Leb}(\mathscr B)^2 \leq Cn\rho_n^{2+2\beta}
 \leq Cn\rho_n^{2+\beta}\to0,
\\
\mathcal X\text{--}\mathscr C\ \text{for any }\mathcal X:\quad
 &\tfrac12n'\!\int_{\mathscr C}e^{-n''\rho_n\check{\mathfrak{D}}/16}
 \leq n\,e^{-3n''\rho_n\bar t/16}\to0,
\end{align*}
the latter by $n\rho_n/\log n\to\infty$.

On $\mathscr A$--$\mathscr A$, \eqref{eq-CS-in-d} is not available, since
$V_s\asymp1$ and $\mathrm{Leb}(\mathscr A)\asymp1$.  By \eqref{eq-A-def} and
\eqref{eq-Jd}, for $2\leq d\leq\sigma_*$ the $d$-th term of
\eqref{eq-type-i} restricted to $\mathscr A\times\mathscr A$ is at most
\begin{align*}
\frac{n^{d+1}\rho_n^{2d}}{d!}C_a^{2d}J_d(s)J_d(t)
 \leq C\rho_n^{\,d-1}\Bigl(\frac{C_a^2C_0}{c}\Bigr)^{d}d\log(n\rho_n)
 + Cn\rho_n^{\,d+2+2\beta}\Bigl(\frac{C_a^2C_0}{c}\Bigr)^{d},
\end{align*}
using $(n'')^{d}\rho_n^{2d}\sigma_*^{-2d}=(4\rho_n/c)^{d}$ and
$\Gamma_d^2\leq C_0^dd!$, the absolute factor $4^d$ being absorbed into $C_0$.
Since $\rho_nC_a^2C_0/c\leq1/2$ for $n$ large, summing over
$2\leq d\leq\sigma_*$ gives $O(\rho_n\log(n\rho_n))+O(n\rho_n^{4+2\beta})\to0$.
For $d>\sigma_*$ the bound $J_d\leq\Gamma_d\sigma_*^{-d}$ gives
$n(\rho_nC_a^2C_0/c)^d$, summing to at most
$n(\rho_nC_a^2C_0/c)^{\sigma_*}\to0$ since
$\sigma_*\asymp\sqrt{n\rho_n}\gg\log n$.

On $\mathscr A$--$\mathscr B$, \eqref{eq-CS-in-d} would give only
$n\mathrm{Leb}(\mathscr B)\asymp n\rho_n^{1+\beta}$, and the termwise bound is used
instead.  By \eqref{eq-A-def} and Lemma~\ref{lem-regions}\ref{reg-univ},
$|C(x,y)|^d\leq C_a^d(3\bar t/\rho_n)^{d/2}$; integrating the $x$-kernel by
\eqref{eq-Jd} and the $y$-kernel by Lemma~\ref{lem-regions}\ref{reg-kernel}
times $\mathrm{Leb}(\mathscr B)$, the $d$-th term is at most
\begin{align*}
C\frac{n\rho_n^{1+\beta}\log(n\rho_n)}{\sqrt{n\rho_n}}\,d\,r^{\,d},
\qquad r := \frac{C_aC_0\sqrt{3\bar t}}{c}\,\rho_n^{1/2}\to0,
\end{align*}
the half power $\rho_n^{-d/2}$ leaving $\rho_n^{d/2}$ to spare; summing over
$d\geq2$ gives
\begin{align*}
O\bigl((n\rho_n^{2+2\beta})^{1/2}\log(n\rho_n)\bigr)\to0 .
\end{align*}

\medskip\noindent\emph{Type (ii).}  Here $S=\{2\}\cup S'$ with
$S'\subseteq\{3,\ldots,n\}$ and $d:=|S'|+1\geq2$; the count is
$\binom{n''}{d-1}$, the two kernels have orders $d$ and $d-1$, and there is an
extra factor $a(x,y)$.  Since
\begin{align*}
\binom{n''}{d-1}\rho_n^{2d-1}
=\rho_n\cdot\binom{n''}{d-1}^{1/2}\rho_n^{d-1}
\cdot\binom{n''}{d-1}^{1/2}\rho_n^{d-1},
\end{align*}
the $d$-th integrand is at most
\begin{align}
\label{eq-type-ii-factor}
n'\rho_n\,|a(x,y)|\;\tilde u_d(x)\,v_{d-1}(y),
\end{align}
pairing the order-$d$ kernel with the weight $\binom{n''}{d-1}$ against
\eqref{eq-parseval-pmf} and the order-$(d-1)$ kernel with the same weight
against \eqref{eq-parseval}.  Cauchy--Schwarz in $d$ therefore bounds the
contribution of $\mathcal X\times\mathcal Y$ by
\begin{align}
\label{eq-CS-type-ii}
n'\rho_n\int_{\mathcal X}\!\int_{\mathcal Y}|a(x,y)|\,
 W_s(x)^{1/2}V_t(y)^{1/2}\,dy\,dx .
\end{align}
This is effective whenever $\mathcal X,\mathcal Y\subseteq\mathscr B\cup\mathscr C$,
because $W_s$ is uniformly of order $(n\rho_n)^{-1/2}$ by \eqref{eq-W-bounds}.
On $\mathscr B$--$\mathscr B$, Cauchy--Schwarz and
Lemma~\ref{lem-regions}\ref{reg-univ}--\ref{reg-I} give
\begin{align*}
\int_{\mathscr B}\!\int_{\mathscr B}|a| \leq
\bigl((3\bar t/\rho_n)\mathrm{Leb}(\mathscr B)\bigr)^{1/2}\mathrm{Leb}(\mathscr B)
\leq C\rho_n^{1+3\beta/2},
\end{align*}
so that \eqref{eq-CS-type-ii} is at most
\begin{align*}
n\rho_n\cdot\frac{C}{(n\rho_n)^{1/4}}\cdot C\rho_n^{1+3\beta/2}
 \;=\; C\bigl(n\rho_n^{\,7/3+2\beta}\bigr)^{3/4} \;\to\;0,
\end{align*}
since $7/3+2\beta\geq2+\beta$.  For a pair meeting $\mathscr C$, either
$W_s$ or $V_t$ is at most $e^{-3n''\rho_n\bar t/16}$ there, while the
remaining factors are $O(\rho_n^{-1})$ by
Lemma~\ref{lem-regions}\ref{reg-h}, and the contribution is
$O\bigl(ne^{-3n''\rho_n\bar t/16}\bigr)\to0$.

It remains to treat $\mathscr A$--$\mathscr A$, $\mathscr A$--$\mathscr B$ and
$\mathscr B$--$\mathscr A$ (first coordinate $x$, second $y$), where
\eqref{eq-CS-type-ii} is too crude, since $V_t\asymp1$ and
$\mathrm{Leb}(\mathscr A)\asymp1$ on $\mathscr A$.  We use termwise bounds instead.
Write
$K^{(d)}_s(x) := \mathsf D^{d}\mathcal B_{n'-d}(s';\rho_n\check{\mathfrak{D}}(x))$,
$K^{(d-1)}_t(y) := \mathsf D^{d-1}\mathcal B_{n'-d+1}(t';\rho_n\check{\mathfrak{D}}(y))$,
$B_n := (3\bar t/\rho_n)^{1/2}$ and
$\ell_n := (\log(n\rho_n)/(n\rho_n))^{1/2}$.  For $2\le d\le\sigma_*$,
Lemma~\ref{lem-Jd} gives
\begin{align}
\label{eq-Jd-short}
J_d(s) \;\le\; C\,\Gamma_d\,\sigma_*^{-d}\,d\,\ell_n ,
\end{align}
since $\sigma_*/(n\rho_n)\asymp(n\rho_n)^{-1/2}\le\ell_n$ and
$\rho_n^{1+\beta}\le(n\rho_n)^{-1/2}$ for $n$ large, the latter because
$n\rho_n^{3+2\beta}\le n\rho_n^{2+\beta}\to0$.  In each case the $d$-th term is
\begin{align*}
T_d := n'\binom{n''}{d-1}\rho_n^{2d-1}&\int_{\mathcal X}\!\int_{\mathcal Y}
 |a(x,y)|\,\|a(x,\cdot)\|^{d-1}\|a(y,\cdot)\|^{d-1}
\\ &\qquad\times|K^{(d)}_s(x)|\,|K^{(d-1)}_t(y)|\,dy\,dx ,
\end{align*}
after bounding $|C(x,y)|^{d-1}$, and we use
$n'\binom{n''}{d-1}\rho_n^{2d-1}\sigma_*^{-(2d-1)}
\le C\,c^{-d}\,n^{1/2}\rho_n^{\,d-1/2}/(d-1)!$ together with
$\Gamma_d\Gamma_{d-1}\le C_0^{\,d}\,d!$.

\emph{$\mathscr A$--$\mathscr A$.}  Here
$\|a(x,\cdot)\|,\|a(y,\cdot)\|\le C_a$.  Integrating $y$ first by
Cauchy--Schwarz,
\begin{align*}
\int_{\mathscr A}|a(x,y)||K^{(d-1)}_t(y)|\,dy
\le\|a(x,\cdot)\|\,\|K^{(d-1)}_t\|_{L^2(\mathscr A)}
\le C_a\bigl(\Gamma_{d-1}\sigma_*^{-(d-1)}J_{d-1}(t)\bigr)^{1/2},
\end{align*}
and then $x$ by \eqref{eq-Jd-short},
\begin{align*}
T_d \;\le\; C\Bigl(\frac{C_a^2C_0}{c}\Bigr)^{d} d^{5/2}\,
 n^{1/2}\rho_n^{\,d-1/2}\,\ell_n^{3/2}
 \;=\; C\Bigl(\frac{C_a^2C_0}{c}\Bigr)^{d} d^{5/2}\,
 n^{-1/4}\rho_n^{\,d-5/4}\,\log^{3/4}(n\rho_n),
\end{align*}
so $\sum_{2\le d\le\sigma_*}T_d
= O\bigl((\rho_n^3\log^3(n\rho_n)/n)^{1/4}\bigr)\to0$, the series converging
because $\rho_nC_a^2C_0/c\le1/2$ for $n$ large.

\emph{$\mathscr A$--$\mathscr B$.}  Here $\|a(x,\cdot)\|\le C_a$ and
$\|a(y,\cdot)\|\le B_n$.  Bound
$|K^{(d-1)}_t(y)|\le\Gamma_{d-1}\sigma_*^{-(d-1)}$ pointwise,
$\int_{\mathscr B}|a(x,y)|dy\le C_a\mathrm{Leb}(\mathscr B)^{1/2}$ by
Cauchy--Schwarz, and integrate $x$ by \eqref{eq-Jd-short}:
\begin{align*}
T_d \;\le\; C\Bigl(\frac{C_aC_0}{c}\Bigr)^{d}B_n^{\,d-1}
 d^{2}\,n^{1/2}\rho_n^{\,d-1/2}\,\ell_n\,\mathrm{Leb}(\mathscr B)^{1/2}
 \;\le\; C\,r_n^{\,d-2}\,d^{2}\;\rho_n^{\,1+\beta/2}\log^{1/2}(n\rho_n),
\end{align*}
with $r_n := (C_aC_0/c)(3\bar t)^{1/2}\rho_n^{1/2}\to0$, using
$B_n=(3\bar t/\rho_n)^{1/2}$ and
$\mathrm{Leb}(\mathscr B)\le C\rho_n^{1+\beta}$; the sum over
$2\le d\le\sigma_*$ is $O\bigl(\rho_n^{1+\beta/2}\log^{1/2}(n\rho_n)\bigr)\to0$.

\emph{$\mathscr B$--$\mathscr A$.}  Here $\|a(x,\cdot)\|\le B_n$ and
$\|a(y,\cdot)\|\le C_a$.  Integrate $y$ first as in the
$\mathscr A$--$\mathscr A$ case,
\begin{align*}
\int_{\mathscr A}|a(x,y)||K^{(d-1)}_t(y)|\,dy
\le B_n\bigl(\Gamma_{d-1}\sigma_*^{-(d-1)}J_{d-1}(t)\bigr)^{1/2},
\end{align*}
then bound
$|K^{(d)}_s(x)|\le\Gamma_d\sigma_*^{-d}$ pointwise and integrate over
$\mathscr B$:
\begin{align*}
T_d \;\le\; C\Bigl(\frac{C_aC_0}{c}\Bigr)^{d}B_n^{\,d}
 d^{3/2}\,n^{1/2}\rho_n^{\,d-1/2}\,\ell_n^{1/2}\,\mathrm{Leb}(\mathscr B)
 \;\le\; C\,r_n^{\,d-2}\,d^{3/2}\;
 \bigl(n\rho_n^{\,5+4\beta}\log(n\rho_n)\bigr)^{1/4},
\end{align*}
and the sum over $2\le d\le\sigma_*$ tends to $0$ since
$n\rho_n^{5+4\beta}\le n\rho_n^{2+\beta}\to0$.

In all three cases the range $d>\sigma_*$ is handled as in type (i), with
$J_d\le\Gamma_d\sigma_*^{-d}$ and
$\|K^{(d-1)}_t\|_{L^2(\mathscr A)}\le\Gamma_{d-1}\sigma_*^{-(d-1)}$, giving a
contribution of at most $n^{1/2}(C\rho_n^{1/2})^{\sigma_*}\to0$.

\medskip\noindent\emph{Type (iv).}  Here $S=\{1,2\}\cup S'$ with
$d:=|S'|+1\geq1$, and the contribution is
\begin{align}
\label{eq-type-iv}
\begin{split}
\sum_{d\geq1}n'\binom{n''}{d-1}\rho_n^{2d}&\int\!\!\int
 \bigl|a(x,y)a(y,x)\bigr||C(x,y)|^{d-1}
\\ &\times\bigl|\mathsf D^{d}\mathcal B_{n'-d}(s';\rho_n\check{\mathfrak{D}}(x))\bigr|
 \bigl|\mathsf D^{d}\mathcal B_{n'-d}(t';\rho_n\check{\mathfrak{D}}(y))\bigr|\,dy\,dx .
\end{split}
\end{align}
Bounding $|C(x,y)|^{d-1}\leq\|a(x,\cdot)\|^{d-1}\|a(y,\cdot)\|^{d-1}$ and
applying Cauchy--Schwarz in $d$ against \eqref{eq-parseval-pmf} on both sides,
the inner sum is at most $\rho_n^2W_s(x)^{1/2}W_t(y)^{1/2}$, so the
contribution of $\mathcal X\times\mathcal Y$ is at most
\begin{align}
\label{eq-type-iv-CS}
n'\rho_n^2\int_{\mathcal X}\!\int_{\mathcal Y}
 \bigl|a(x,y)a(y,x)\bigr|W_s(x)^{1/2}W_t(y)^{1/2}dy\,dx ,
\end{align}
to which we apply $|a(x,y)a(y,x)|\leq\tfrac12(a(x,y)^2+a(y,x)^2)$ and
Lemma~\ref{lem-regions}\ref{reg-h}.

On $\mathscr B$--$\mathscr B$, by Lemma~\ref{lem-regions}\ref{reg-univ}--\ref{reg-I},
\begin{align*}
\int_{\mathscr B}\!\int_{\mathscr B}a(x,y)^2dy\,dx\leq
\int_{\mathscr B}\|a(x,\cdot)\|^2dx\leq(3\bar t/\rho_n)\mathrm{Leb}(\mathscr B)
\leq C\rho_n^{\beta},
\end{align*}
so \eqref{eq-type-iv-CS} and \eqref{eq-W-bounds} give
\begin{align*}
n'\rho_n^2\cdot\frac{C}{\sqrt{n\rho_n}}\cdot C\rho_n^{\beta}
 \leq C\bigl(n\rho_n^{\,3+2\beta}\bigr)^{1/2}\to0 .
\end{align*}
For a pair meeting $\mathscr C$, by Lemma~\ref{lem-regions}\ref{reg-F} and
Lemma~\ref{lem-binomial-facts}\ref{binom-fact-b},
\begin{align*}
W_t(y)\leq b_{n''}(t';\rho_n\check{\mathfrak{D}}(y))\leq
e^{-n''\rho_n\check{\mathfrak{D}}(y)/8}\leq e^{-3n''\rho_n\bar t/8},
\end{align*}
while the remaining factors
are $O(\rho_n^{-1})$ by Lemma~\ref{lem-regions}\ref{reg-h} and
\eqref{eq-W-bounds}; the contribution is
$O\bigl(ne^{-3n''\rho_n\bar t/16}\bigr)\to0$.

On $\mathscr A$--$\mathscr B$, Cauchy--Schwarz and \eqref{eq-W-bounds} give
\begin{align*}
\int_{\mathscr B}W_t^{1/2}\leq\Bigl(\int_0^1W_t\Bigr)^{1/2}
\mathrm{Leb}(\mathscr B)^{1/2}\leq C(n\rho_n)^{-1/2}\rho_n^{(1+\beta)/2}.
\end{align*}
Integrating $x$ first, Lemma~\ref{lem-regions}\ref{reg-h} gives
$\int_{\mathscr A}a(x,y)^2dx\leq C\rho_n^{-1}$, and
$\|a(y,\cdot)\|^2\leq3\bar t/\rho_n$ for $y\in\mathscr B$, so
\eqref{eq-type-iv-CS} is at most
\begin{align*}
n'\rho_n^2\cdot\sup_xW_s(x)^{1/2}\cdot\frac{C}{\rho_n}
 \int_{\mathscr B}W_t^{1/2}
 \leq Cn\rho_n\cdot\frac{1}{(n\rho_n)^{1/4}}\cdot
 \frac{\rho_n^{(1+\beta)/2}}{(n\rho_n)^{1/2}}
 = C\bigl(n\rho_n^{\,3+2\beta}\bigr)^{1/4}\to0 .
\end{align*}
On $\mathscr A$--$\mathscr A$ the termwise bound is used, as for type (i):
$\|a\|\leq C_a$ and $E[\check w^2(x,\xi)]\leq C$ there by
\ref{condition-l1}, while the weight $\binom{n''}{d-1}$ supplies one further
power of $\rho_n$ relative to type (i), giving
$C\rho_n^{\,d}(n\rho_n)^{-1/2}(C_a^2C_0/c)^dd\log(n\rho_n)$, which sums to
$O\bigl(\rho_n(n\rho_n)^{-1/2}\log n\bigr)\to0$.

\medskip
Collecting each of the bounds, it follows that $R_n(s,t)\to0$.
\end{proof}

The decomposition \eqref{eq-gap-two-sources} also produces a diagonal
contribution, namely the terms with $i = j$ in the double sum over roots.  It
carries no factor of $n$, and is correspondingly easier; we record it
separately.

\begin{lemma}[Diagonal contribution]
\label{lem-diagonal}
Under \ref{condition-m1}--\ref{condition-s1},
\begin{align*}
\sum_{|S|\geq2}E\bigl[\Pi\bigl(F_1(s'\mid\mathcal G);S\bigr)\,
 \Pi\bigl(F_1(t'\mid\mathcal G);S\bigr)\bigr] \;\longrightarrow\; 0
\end{align*}
for every $s,t$ in the interior of $\mathrm{supp}\,F$.
\end{lemma}

\begin{proof}
By the Cauchy--Schwarz inequality it suffices to show
\begin{align*}
\Delta_n(s) := \sum_{|S|\geq2}E\bigl[\Pi(F_1(s'\mid\mathcal G);S)^2\bigr]\to0 .
\end{align*}
By Lemma~\ref{lem-multilinear}, the components of $F_1(s'\mid\mathcal G)$
come in pairs indexed by $S\subseteq[n]\setminus\{1\}$:
$\Pi(F_1;S) = E[K_S(s)\mid\mathcal G_S]$ and
$\Pi(F_1;S\cup\{1\}) = K_S(s)-E[K_S(s)\mid\mathcal G_S]$, whose second moments
add up to $E[K_S(s)^2]$.  Every component of order at least two belongs to a
pair with $|S|\geq1$; note that this includes the components
$\Pi(F_1;\{1,j\})$, for which $|S| = 1$.  Hence, by \eqref{eq-Ksq},
\begin{align}
\label{eq-diagonal-reduction}
\Delta_n(s) \;\leq\; E\bigl[\widetilde V_s(\xi_1)\bigr],
\qquad
\widetilde V_s(x) := \sum_{d\geq1}\binom{n'}{d}\rho_n^{2d}
 \bigl\|a(x,\cdot)\bigr\|^{2d}
 \Bigl(\mathsf D^{d}\mathcal B_{n'-d}\bigl(s';\rho_n\check{\mathfrak{D}}(x)\bigr)\Bigr)^2 ,
\end{align}
the sum over $S$ of a given size having been collapsed as in
\eqref{eq-parseval}.  The argument of Lemma~\ref{lem-parseval-cdf}, applied to
all $S\subseteq[n]\setminus\{1\}$, shows that
$\widetilde V_s(x)\leq\mathrm{Var}(F_1(s'\mid\mathcal G)\mid\xi_1 = x)\leq
\min\{\tfrac14,\mathcal B_{n'}(s';\rho_n\check{\mathfrak D}(x))\}$.  We bound $\widetilde V_s$ on each of the three regions
\eqref{eq-regions}.

On $\mathscr A$ we use the termwise bounds.  By \eqref{eq-A-def} and
Lemma~\ref{lem-regions}\ref{reg-kernel}, together with the identity
$(n'')^{d}\rho_n^{2d}\sigma_*^{-2d} = (4\rho_n/c)^{d}$ and the bound
$\Gamma_d^2\leq C_0^{d}d!$,
\begin{align*}
\widetilde V_s(x) \;\leq\; \sum_{d\geq1}
 \Bigl(\frac{2\rho_nC_a^2C_0}{c}\Bigr)^{d} \;\leq\; C\rho_n ,
\end{align*}
using $(n')^d\leq2^d(n'')^d$, the series converging because
$2\rho_nC_a^2C_0/c\leq1/2$ for $n$ large; the terms with $d>n''/2$ are handled
by the crude bound of Lemma~\ref{lem-binomial-facts}\ref{binom-fact-d}, as in
the proof of Lemma~\ref{lem-degenerate}.

On $\mathscr B$ the termwise bounds are unavailable, even though
Lemma~\ref{lem-regions}\ref{reg-kernel} applies there: by
Lemma~\ref{lem-regions}\ref{reg-univ} the ratio of the resulting series is
$\rho_n\cdot(3\bar t/\rho_n)\,C_0/c = 3\bar tC_0/c$, which is not small.  We
use instead the exact bound $\widetilde V_s\leq1/4$, together with
Lemma~\ref{lem-regions}\ref{reg-I}:
\begin{align*}
E\bigl[\widetilde V_s(\xi_1)\mathbbm 1(\xi_1\in\mathscr B)\bigr]
 \;\leq\; \tfrac14\mathrm{Leb}(\mathscr B) \;\leq\; C\rho_n^{1+\beta} .
\end{align*}

On $\mathscr C$ we use the bound
$\widetilde V_s\leq\mathcal B_{n'}(s';\rho_n\check{\mathfrak{D}}(x))$, which by
Lemma~\ref{lem-regions}\ref{reg-F} and
Lemma~\ref{lem-binomial-facts}\ref{binom-fact-b} is at most
$e^{-n''\rho_n\check{\mathfrak{D}}(x)/8}\leq e^{-3n''\rho_n\bar t/8}$.

Collecting the three regions,
$\Delta_n(s) = O(\rho_n) + O(\rho_n^{1+\beta}) + O(e^{-3n''\rho_n\bar t/8})
\to0$.
\end{proof}

Now, it remains to establish step (2) of our strategy discussed in Section \ref{sec-multilinear-hoeffding}, pertaining to the claimed negligibility of network noise. 

\begin{lemma}[Negligibility of the network noise]
\label{lemma-sigma-n-0}
Assume \ref{condition-m1}, \ref{condition-g1} and \ref{condition-d1}, and that
$n\rho_n/\log n\to\infty$ and $n\rho_n^{2+\beta}\to0$.  Then, for all $s,t$ in
the interior of $\mathrm{supp}\,F$,
\begin{align*}
\Sigma_{n,0}(s,t) := n'\,E\bigl[\mathrm{Cov}\bigl(\mathbbm 1(D_1\leq s'),\,
 \mathbbm 1(D_2\leq t')\mid\mathcal G\bigr)\bigr]
 \rightarrow\; 0 .
\end{align*}
\end{lemma}
 
\begin{proof}
\emph{Step 1: the cancellation argument.}  For $i\neq j$ and an integer $k$,
let $f_{i\setminus j}(k\mid\mathcal G) := P(D_{i\setminus j} = k\mid\mathcal G)$.
Given $\mathcal G$, $D_1 = D_{1\setminus2}+A_{12}$ and $D_2 = D_{2\setminus1}+A_{12}$,
where the three summands are conditionally independent.  Conditioning further
on $A_{12}$ makes the two indicators independent, so
\begin{align*}
\mathrm{Cov}\bigl(\mathbbm 1(D_1\leq s'),\mathbbm 1(D_2\leq t')\bigm|\mathcal G\bigr)
 = \mathrm{Cov}\bigl(u(A_{12}),v(A_{12})\bigm|\mathcal G\bigr),
\end{align*}
where $u(\epsilon) := P(D_{1\setminus2}+\epsilon\leq s'\mid\mathcal G)$ and
$v(\epsilon) := P(D_{2\setminus1}+\epsilon\leq t'\mid\mathcal G)$.  For a
$\mathrm{Bernoulli}(p)$ variable $A$,
$\mathrm{Cov}(u(A),v(A)) = p(1-p)\bigl(u(1)-u(0)\bigr)\bigl(v(1)-v(0)\bigr)$, and
$u(1)-u(0) = -f_{1\setminus2}(s'\mid\mathcal G)$,
$v(1)-v(0) = -f_{2\setminus1}(t'\mid\mathcal G)$; the two signs cancel.  With
$p_{12} = \rho_n\check w_{12}$, $\check w_{12} := \check w(\xi_1,\xi_2)$,
\begin{align*}
\Sigma_{n,0}(s,t) = n'\rho_n\,E\Bigl[\check w_{12}\bigl(1-\rho_n\check w_{12}\bigr)
 f_{1\setminus2}(s'\mid\mathcal G)\,f_{2\setminus1}(t'\mid\mathcal G)\Bigr]\;\geq\;0,
\end{align*}
so that it suffices to bound the right-hand side above, and we may drop the
factor $1-\rho_n\check w_{12}\leq1$.  This is the same identity as
\eqref{eq-c-pairing} in the proof of tightness, there at the scale of a cell.
 
\emph{Step 2: anticoncentration on the $\xi_2$ side.}  Edges with
$\rho_nw(\xi_2,\xi_j)\geq1$ are deterministic given $\mathcal G$ and contribute
no variance, so the Kolmogorov--Rogozin inequality must be applied to a
sub-sum.  Let
\begin{align*}
\mathcal J := \bigl\{j\neq1,2 \;:\; \rho_nw(\xi_2,\xi_j)\leq\tfrac12\bigr\},
\end{align*}
so that $p_{2j}\in[c\rho_n,\tfrac12]$ for $j\in\mathcal J$ by
\ref{condition-g1}, whence $p_{2j}(1-p_{2j})\geq c\rho_n/2$.  Since the
concentration function of a sum of independent variables is at most that of any
sub-sum,
\begin{align}
\label{eq-rogozin}
\sup_k f_{2\setminus1}(k\mid\mathcal G)
 \;\leq\; \sup_k P\Bigl(\sum_{j\in\mathcal J}A_{2j} = k \Bigm| \mathcal G\Bigr)
 \;\leq\; \frac{C}{\sqrt{|\mathcal J|\,\rho_n}}
\qquad\text{on }\bigl\{|\mathcal J|\geq n''/2\bigr\},
\end{align}
by \citet[][Theorem 2]{rogozin-concentration} (see also \citealp{kesten1969sharper}), giving $C(n\rho_n)^{-1/2}$ there.
 
\emph{Step 3: the exceptional set.}  Given $\xi_2$, the indicators
$\mathbbm 1(j\notin\mathcal J)$, $j\neq1,2$, are i.i.d.\ Bernoulli with
parameter
$q(\xi_2) := P\bigl(w(\xi_2,\xi)>(2\rho_n)^{-1}\mid\xi_2\bigr)
\leq(2\rho_n)^{2+\beta}Z(\xi_2)$ by Markov's inequality, where
$Z(\xi_2) := E[w^{2+\beta}(\xi_2,\xi)\mid\xi_2]$ satisfies
$E[Z] = \|w\|^{2+\beta}_{L^{2+\beta}}<\infty$ by \ref{condition-m1}.  Put
\begin{align*}
\theta_n := \tfrac14(2\rho_n)^{-(2+\beta)},
\qquad
B := \bigl\{Z(\xi_2)\geq\theta_n\bigr\},
\end{align*}
so that $q(\xi_2)\leq1/4$ off $B$, and the multiplicative Chernoff bound gives
\begin{align}
\label{eq-J-tail}
P\bigl(|\mathcal J^c|>n''/2 \bigm| \xi_2\bigr) \;\leq\; e^{-n''/12}
\qquad\text{off } B .
\end{align}
 
\emph{Step 4: the $\xi_1$ side.}  Since $D_{1\setminus2}$ is a function of
$\xi_1$ and $(\xi_j,\eta_{1j})_{j\neq1,2}$, it is conditionally independent of
$\check w_{12}$ given $\xi_1$, so that
\begin{align}
\label{eq-xi1-side}
E\bigl[\check w_{12}\,f_{1\setminus2}(s'\mid\mathcal G)\bigm|\xi_1\bigr]
 = \check{\mathfrak{D}}(\xi_1)\,P\bigl(D_{1\setminus2}=s'\bigm|\xi_1\bigr)
 = \check{\mathfrak{D}}(\xi_1)\,b_{n''}\bigl(s';\rho_n\check{\mathfrak{D}}(\xi_1)\bigr),
\end{align}
the last equality since conditioning on the single
latent position $\xi_1$ makes the $n''$ indicators i.i.d.\ with common
parameter $\rho_n\check{\mathfrak{D}}(\xi_1)$, so the Poisson-binomial mass function is
averaged into a binomial one.  Consequently
\begin{align}
\label{eq-xi1-integral}
n'\rho_n\,E\bigl[\check w_{12}\,f_{1\setminus2}(s'\mid\mathcal G)\bigr]
 = \int_0^1\check{\mathfrak{D}}(x)\,\mu_s(dx) \;=\; O(1)
\end{align}
by Lemma~\ref{lem-mu-limit}, splitting the integral at $\check{\mathfrak{D}}(x)\leq2s$
and using, beyond it, the multiplicative Chernoff bound for the lower tail of
$\mathrm{Bin}(n'',\rho_n\check{\mathfrak D}(x))$, which holds for every value of
the parameter.
 
\emph{Step 5: combining.}  Let $G := B^c\cap\{|\mathcal J|\geq n''/2\}$.  On
$G$ we use \eqref{eq-rogozin} and then \eqref{eq-xi1-integral}:
\begin{align*}
n'\rho_n E\bigl[\check w_{12}f_{1\setminus2}(s'\mid\mathcal G)
 f_{2\setminus1}(t'\mid\mathcal G)\mathbbm 1_G\bigr]
 \;\leq\; \frac{C}{\sqrt{n\rho_n}}\,n'\rho_nE\bigl[\check w_{12}
 f_{1\setminus2}(s'\mid\mathcal G)\bigr]
 \;=\; O\bigl((n\rho_n)^{-1/2}\bigr).
\end{align*}
On $G^c$ we bound both mass functions by $1$ and use
$E[\check w_{12}\mid\xi_2] = \check{\mathfrak{D}}(\xi_2)\leq\mathfrak D(\xi_2)$.  By
\eqref{eq-J-tail} and conditional independence of $\check w_{12}$ from
$\mathcal J$ given $\xi_2$,
\begin{align*}
n'\rho_nE\bigl[\check w_{12}\mathbbm 1(B^c)\mathbbm 1(|\mathcal J^c|>n''/2)\bigr]
 \;\leq\; n'\rho_n\,e^{-n''/12}\,E[\mathfrak D] \;\longrightarrow\;0 .
\end{align*}
Finally, by conditional Jensen $\mathfrak D(\xi_2)\leq Z(\xi_2)^{1/(2+\beta)}$,
so that, $1/(2+\beta)-1<0$,
\begin{align*}
n'\rho_nE\bigl[\check w_{12}\mathbbm 1_B\bigr]
 &\;\leq\; n'\rho_n\,E\bigl[Z^{1/(2+\beta)}\mathbbm 1(Z\geq\theta_n)\bigr]
\\ &\;\leq\; n'\rho_n\,\theta_n^{\frac1{2+\beta}-1}E[Z]
 \;=\; O\bigl(n\rho_n^{\,2+\beta}\bigr),
\end{align*}
which tends to $0$ by hypothesis.  Collecting the three contributions gives
the statement.
\end{proof}

\subsection{Proof of the Fidi CLT}
We now assemble lemmas from the previous sections to prove Theorem \ref{thm-fidi-clt}. 
\begin{proof}[Proof of Theorem \ref{thm-fidi-clt}]
By Lemma~\ref{lemma-sigma-n-0} and dominated convergence applied to
$E[Q_s(1-Q_s)]$, where $Q_s := F_1(s'\mid\mathcal G) \to \mathbbm 1(\mathfrak D<s)$
almost surely, we have $E[\mathrm{Var}(T_n(\gamma,t)\mid\mathcal G)] \to 0$, so
that $T_n(\gamma,t) - E[T_n(\gamma,t)\mid\mathcal G] \to 0$ in $L^2$.  By
Lemmas~\ref{lem-degenerate} and \ref{lem-diagonal}, $E[T_n(\gamma,t)\mid\mathcal G] -
H(T_n(\gamma,t)) \to 0$ in $L^2$ as well.  Hence $T_n(\gamma,t)$ and
$H(T_n(\gamma,t))$ have the same limit in distribution, if either has one.  

Now, for any $j \in [K]$, let $Y_{n,i} = \sum_{k=1}^K \gamma_k \{\psi_{t_k}(\xi_i) + \phi_{t_k}(\xi_i) \}$.   By $L^2$ convergence implied by Lemma \ref{lem-hajek-limits}, it follows that $Y_{n,1}^2$ is uniformly integrable.  Since the summands are iid, this implies that the Lindeberg condition is satisfied.  Thus, $H(T_n(\gamma, t))$ tends toward a Normal limit, and Corollary \ref{cor-covariance}  provides the limiting variance, which matches with the variance implied by the posited limiting covariance function \eqref{eq-cov-function}.  The result follows.   

\end{proof}

\subsubsection{On the Choice of Centering}

Centering the empirical degree distribution at $F^\star(t) = P(D_i/(n'\rho_n) \leq t)$
yields a finite-dimensional central limit theorem under the mild sparsity condition
$n\rho_n/\log n\to\infty$. For statistical
applications the centering at $F(t) = P(\mathfrak D(\xi_i) \leq t)$ is more
natural, since it does not depend on $n$ and reflects the degree structure of
the underlying kernel.  The price is a stronger sparsity requirement and a smoothness condition on $f_\mathfrak{D}$, which the
following lemma quantifies.
\begin{lemma}[Centering]
\label{lem-centering}
Suppose \ref{condition-m1}, \ref{condition-g1} and \ref{condition-d1} hold, and
let $I\subset(0,\infty)$ be a compact interval.  Suppose that, for some
$\delta>0$ and $\theta\in(0,1]$, $f_{\mathfrak D}$ is $\theta$-H\"older
continuous on the $\delta$-neighborhood $I^\delta$ of $I$, and that
\begin{align}
\label{eq-centering-conditions}
\frac{n\rho_n}{\log n}\to\infty, \qquad
n\rho_n^{\,3/2+\beta}\to0, \qquad
n^{\theta}\rho_n^{\,1+\theta}\to\infty .
\end{align}
Then
\begin{align*}
\sup_{t\in I}\sqrt n\,\bigl|F^\star(t)-F(t)\bigr|\;\rightarrow\;0 .
\end{align*}
\end{lemma}

\begin{proof}
We begin by controlling the effect of truncation for $F^*(t)$. 
Write
\begin{align}
\label{eq-bias-decomp}
\sqrt n\,\bigl[F^{\star}(t)-F(t)\bigr]
 = \sqrt n\int_0^1\Bigl[\mathcal B_{n'}\bigl(t';\rho_n\check{\mathfrak D}(x)\bigr)
 - \mathbbm 1\bigl(\mathfrak D(x)\leq t\bigr)\Bigr]dx ,
\end{align}
Here, we use reasoning similar to Steps 1 and 3 of Lemma~\ref{lem-mu-limit} to consider the behavior on and off $\mathscr{G}_n$.
 
\emph{Clipping bias on $\mathscr{G}_n$.}  Since
$\partial_p\mathcal B_{n'}(m;p) = -n'\,b_{n''}(m;p)$, the mean value theorem
gives, for each such $x$, a $p_x$ between $\rho_n\check{\mathfrak D}(x)$ and
$\rho_n\mathfrak D(x)$, hence in $[c\rho_n,1/4]$ by \ref{condition-g1}, with
\begin{align*}
\Bigl|\mathcal B_{n'}\bigl(t';\rho_n\check{\mathfrak D}(x)\bigr)
 -\mathcal B_{n'}\bigl(t';\rho_n\mathfrak D(x)\bigr)\Bigr|
 = n'\rho_n\,g(x)\,b_{n''}(t';p_x)
 \;\leq\; \frac{C\,n\rho_n}{\sqrt{n\rho_n}}\,g(x),
\end{align*}
using $\sup_m b_{n''}(m;p)\leq C\sigma(n'',p)^{-1}\leq C'(n\rho_n)^{-1/2}$ from
Lemma~\ref{lem-binomial-facts}\ref{binom-fact-a}.  Hence, by the global bias
bound \eqref{eq-g-global},
\begin{align}
\label{eq-bias-clip}
\sqrt n\!\int_{\mathscr{G}_n}\!\Bigl|\mathcal B_{n'}\bigl(t';\rho_n\check{\mathfrak D}\bigr)
 -\mathcal B_{n'}\bigl(t';\rho_n\mathfrak D\bigr)\Bigr|dx
 \;\leq\; C\sqrt{n}\,\sqrt{n\rho_n}\int_0^1 g
 \;\leq\; C\,n\rho_n^{\,3/2+\beta} .
\end{align}
 
\emph{The region $\mathscr{G}_n^c$.}  Here the integrand in
\eqref{eq-bias-decomp} is a difference of two quantities in $[0,1]$, so by
\eqref{eq-good-complement} its contribution is at most
\begin{align*}
\sqrt n\,\mathrm{Leb}\bigl(\mathscr G_n^c\bigr) \;\leq\; C\sqrt n\,\rho_n^{\,2+\beta}
 \;=\; C\bigl(n\rho_n^{\,4+2\beta}\bigr)^{1/2} \;\rightarrow\; 0,
\end{align*}
since $n\rho_n^{4+2\beta}\leq n\rho_n^{3/2+\beta}$.  Both this bound and
\eqref{eq-bias-clip} are free of $t$, so they hold uniformly. 

Putting together the pieces, it follows that, uniformly in $t\in I$,
\begin{align}
\label{eq-bias-reduced}
\sqrt n\,\bigl[F^{\star}(t)-F(t)\bigr]
 = \sqrt n\int_{\mathscr G_n}\Bigl[\mathcal B_{n'}\bigl(t';\rho_n\mathfrak D(x)\bigr)
 - \mathbbm 1\bigl(\mathfrak D(x)\leq t\bigr)\Bigr]dx
 + O\bigl(n\rho_n^{\,3/2+\beta}\bigr) + o(1).
\end{align}
Pushing forward the Lebesgue measure on $[0,1]$ by $\mathfrak D$, which by
Remark~\ref{remark-inverse} and \ref{condition-d1} sends $dx$ to
$f(\lambda)\,d\lambda$, the good set corresponds to
$\lambda\in[0,L_n]$ with $L_n := (4\rho_n)^{-1}$, and the integral in
\eqref{eq-bias-reduced} equals
\begin{align}
\label{eq-bias-kernel}
\sqrt n\int_0^{L_n}K_t(\lambda)\,f(\lambda)\,d\lambda,
\qquad
K_t(\lambda) := \mathcal B_{n'}\bigl(t';\rho_n\lambda\bigr)-\mathbbm 1(\lambda\leq t).
\end{align}
We write $f(\lambda) = f(t)+[f(\lambda)-f(t)]$ and treat the two parts
separately.  Note that $L_n>\sup I$ for $n$ large.

\emph{The constant part: the lattice term.}  Summing
Lemma~\ref{lem-binomial-facts}\ref{binom-fact-e} over $k\leq m$ gives
\begin{align*}
\int_0^1\mathcal B_M(m;p)\,dp = \frac{m+1}{M+1} .
\end{align*}
For $p\geq\frac14$ we have
$t'\leq n'p/2$ for $n$ large, so
Lemma~\ref{lem-binomial-facts}\ref{binom-fact-b} gives
$\mathcal B_{n'}(t';p)\leq e^{-n'/32}$.  Hence, substituting $p = \rho_n\lambda$,
\begin{align*}
\int_0^{L_n}K_t(\lambda)\,d\lambda
 = \frac1{\rho_n}\int_0^{1/4}\mathcal B_{n'}(t';p)\,dp - t
 = \frac{t'+1-(n'+1)\rho_nt}{(n'+1)\rho_n} + O\bigl(\rho_n^{-1}e^{-n'/32}\bigr).
\end{align*}
Since $n'\rho_nt-1<t'\leq n'\rho_nt$, the numerator lies in
$(-\rho_nt,\,1]$, so that
\begin{align}
\label{eq-bias-lattice}
\sqrt n\,f(t)\Bigl|\int_0^{L_n}K_t(\lambda)\,d\lambda\Bigr|
 \;\leq\; \frac{C\,\sup_If}{\sqrt n\,\rho_n} + \text{negligible},
\end{align}
uniformly in $t\in I$.  This term is the discretization of the degree scale,
which oscillates with the fractional part of $n'\rho_nt$; integration over the
latent variable reduces it to the closed form above, so that no correction for
the lattice is needed.

\emph{The variation part: the smoothing bias.}  For $\lambda<t$,
$K_t(\lambda) = -P(X_\lambda>t')$, and for $\lambda>t$,
$K_t(\lambda) = P(X_\lambda\leq t')$, where
$X_\lambda\sim\mathrm{Bin}(n',\rho_n\lambda)$.  In both cases the event
requires $|X_\lambda-n'\rho_n\lambda|\geq n'\rho_n|\lambda-t|-1$, and
$\mathrm{Var}(X_\lambda)\leq n'\rho_n\lambda$.  By Bernstein's inequality,
\begin{align}
\label{eq-bernstein-K}
|K_t(\lambda)|\;\leq\;
\begin{cases}
C\exp\bigl(-\zeta(\lambda-t)^2/h_t^2\bigr), & |\lambda-t|\leq t/2,\\[2pt]
2\exp\bigl(-c\,n\rho_n|\lambda-t|\bigr), & |\lambda-t|>t/2,
\end{cases}
\end{align}
with $h_t^2 = t/(n''\rho_n)$ and absolute constants $C,\zeta,c>0$; in the first
case the denominator in Bernstein's bound is of order $n\rho_nt$, and the
discretization $|t'-n'\rho_nt|\leq1$ is absorbed into $C$.  We split the
integral at $|\lambda-t| = \delta$.  On $|\lambda-t|\leq\delta\leq t/2$, the
H\"older condition and \eqref{eq-bernstein-K} give
\begin{align*}
\int_{|\lambda-t|\leq\delta}|K_t(\lambda)|\,|f(\lambda)-f(t)|\,d\lambda
 \;\leq\; CK_f\int_{\mathbb R}|\lambda-t|^{\theta}
 e^{-\zeta(\lambda-t)^2/h_t^2}\,d\lambda
 \;=\; C'h_t^{\,1+\theta}.
\end{align*}
On $|\lambda-t|>\delta$, \eqref{eq-bernstein-K} gives
$|K_t(\lambda)|\leq\epsilon_n := C\exp(-c\,n\rho_n\delta^2/\sup I)$, while
$|f(\lambda)-f(t)|\leq f(\lambda)+f(t)$, so that part is at most
$\epsilon_n\bigl(1+L_n\sup_If\bigr)$, which is negligible since
$n\rho_n/\log n\to\infty$.  Hence, uniformly in $t\in I$,
\begin{align}
\begin{split}
\label{eq-bias-curvature}
\sqrt n\Bigl|\int_0^{L_n}K_t(\lambda)\bigl[f(\lambda)-f(t)\bigr]d\lambda\Bigr|
 &\;\leq\; C\sqrt n\,h_t^{\,1+\theta}+\text{negligible}
 \\&\;\leq\; C\sqrt n\,(n\rho_n)^{-(1+\theta)/2}+\text{negligible}.
\end{split}
\end{align}

\emph{Conclusion.}  Combining \eqref{eq-bias-reduced}--\eqref{eq-bias-curvature},
\begin{align}
\label{eq-bias-total}
\sup_{t\in I}\sqrt n\,\bigl|F^\star(t)-F(t)\bigr|
 \;\leq\; C\Bigl[n\rho_n^{\,3/2+\beta}
 + \frac1{\sqrt n\,\rho_n}
 + \sqrt n\,(n\rho_n)^{-(1+\theta)/2}\Bigr] + o(1).
\end{align}
The middle term is dominated by the last, since their ratio is
$(n\rho_n)^{-(1-\theta)/2}\leq1$.  The right side therefore tends to $0$
provided
\begin{align}
\label{eq-bias-conditions}
n\rho_n^{\,3/2+\beta}\to0
\qquad\text{and}\qquad
n^{\theta}\rho_n^{\,1+\theta}\to\infty,
\end{align}
the second being equivalent to $\sqrt n(n\rho_n)^{-(1+\theta)/2}\to0$. 
\end{proof}

\begin{remark}
\label{rem-bias-edge}
The H\"older condition on a neighborhood of $I$ cannot be dropped: if $f$ is
discontinuous at a point of $I$, the kernel $K_t$ straddles the discontinuity
for $t$ within $O(h_t)$ of it, and $\sqrt n[F^\star-F]$ is of order
$\sqrt n\,h_t\asymp\rho_n^{-1/2}$ there.  This is relevant in particular at
$\inf\mathrm{supp}\,F$, since $\mathfrak D\geq c>0$ by \ref{condition-g1}.  The
two conditions in \eqref{eq-bias-conditions} are compatible if and only if
$\theta(\tfrac12+\beta)>1$.
\end{remark}

\section{Tightness}
\label{sec-tight-compact}

The conditions needed for tightness are stronger than those for the
finite-dimensional theorem, and they are of two kinds.  On a compact interval
in the interior of the support, local regularity suffices, and for smooth
densities and kernels the sparsity condition is that of
Theorem~\ref{thm-fidi-clt}.  Uniformity over an unbounded range of degrees
requires, in addition, control of the tail of the degree distribution and a
stronger sparsity condition.  We state the compact result first, then its
extension, then discuss the assumptions, and close the subsection with the
strategy of the proof.

\subsection{Assumptions and results}

Let $d_0 := \inf\mathrm{supp}\,F$ and $d_1 := \sup\mathrm{supp}\,F\leq\infty$,
fix points $d_0<t_0^-<t_0<t_1<t_1^+<d_1$, and write $I := [t_0,t_1]$ for the
index set and $J := [t_0^-,t_1^+]$ for a neighborhood of it.  As in
Remark~\ref{remark-inverse}, $x_\lambda := \mathfrak D^{-1}(\lambda)$ and
$f := f_{\mathfrak D}$.

\paragraph{Compact intervals.}
\begin{enumerate}[leftmargin=1.5cm]
\item[D2.]\mylabel{D2}\label{condition-d2}
(Density) $f$ is bounded above and below by positive constants on $J$, and
$|f(\lambda)-f(\lambda')|\leq K_f|\lambda-\lambda'|^{\theta_f}$ for
$\lambda,\lambda'\in J$.
\item[G2.]\mylabel{G2}\label{condition-g2}
(Graphon regularity) $\|w(x_\lambda,\cdot)-w(x_{\lambda'},\cdot)\|_{L^2}\leq
K_w|\lambda-\lambda'|^{\theta_w}$ for $\lambda,\lambda'\in J$.
\item[L2.]\mylabel{L2}\label{condition-l2}
(Local moments) $\sup\{E[w^4(x,\xi)]:\mathfrak D(x)\in J\}<\infty$.  If
$\mathrm{supp}\,F$ is unbounded, in addition
$E[w^{2+\beta}(x,\xi)]\leq K\,\mathfrak D(x)^{2+\beta}$ for all $x$ with
$\mathfrak D(x)\geq t_0^-$.
\item[S2.]\mylabel{S2}\label{condition-s2}
(Sparsity) $\rho_n = n^{-b}$ for some $b\in(0,1)$, and with
$\vartheta := \min(\theta_w,\theta_f,1/2)>1/4$,
\begin{align*}
b\,(5-2\vartheta)>1-2\vartheta .
\end{align*}
\end{enumerate}

\begin{theorem}[Uniform CLT on compact intervals]
\label{thm-tightness-compact}
Under the conditions of Theorem~\ref{thm-fidi-clt} and
\ref{condition-d2}--\ref{condition-s2},
$\mathring{\mathbb G}_n\rightsquigarrow\mathbb G$ in $\ell^\infty(I)$, where
$\mathbb G$ is the Gaussian process of Theorem~\ref{thm-fidi-clt}, which has a
version with continuous sample paths on $I$.
\end{theorem}


\begin{remark}[Interpreting the sparsity condition]
\label{rem-s2}
Write $\rho_n = n^{-b}$.  Condition \ref{condition-s1} then reads
$b\in\bigl(\tfrac1{2+\beta},1\bigr)$, and \ref{condition-s2} adds the
requirement
\begin{align*}
b>\frac{1-2\vartheta}{5-2\vartheta},\qquad
\vartheta = \min\bigl(\theta_w,\theta_f,\tfrac12\bigr).
\end{align*}
This is a lower bound on $b$, so it asks the graph to be sparse enough.  It
is mild in three respects.  First, it is vacuous when $\vartheta = \frac12$,
that is, when $\lambda\mapsto w(x_\lambda,\cdot)$ and $f$ are H\"older of
order at least $\frac12$; Theorem~\ref{thm-tightness-compact} then holds
under exactly the sparsity condition of Theorem~\ref{thm-fidi-clt}.  Second,
the threshold never exceeds $\frac19$, its value as $\vartheta\downarrow\frac14$.
Third, it is implied by \ref{condition-s1} unless
$\beta>(3+2\vartheta)/(1-2\vartheta)\geq7$, so it restricts the sparsity
level beyond the finite-dimensional theorem only for kernels with many
moments and rough dependence on the latent degree.

Heuristically, the condition balances the roughness of the kernel against
the resolution of the lattice.  Denser graphs have smaller bandwidth
$h_t\asymp(n\rho_n)^{-1/2}$, so the increments resolve cells on which the
latent fluctuations are controlled only through the H\"older modulus, of
order $h_t^{\vartheta}$, while tightness requires a bound of order
$\Delta_A^{\vartheta'}$ with $\vartheta'>\frac14$.  When $\vartheta<\frac12$,
the gap between the two must be absorbed by the factor $\rho_n^2$ that
accompanies the interaction of two vertices through the edge
probabilities, and \ref{condition-s2} is exactly what makes this possible
(see \eqref{eq-c-PA}).  The condition arises from an operator-norm bound in
Step~Z4, and we do not expect it to be sharp.
\end{remark}

Combined with Lemma~\ref{lem-centering}, Theorem~\ref{thm-tightness-compact}
gives a uniform CLT for the process centered at $F$.

\begin{corollary}[Uniform CLT for $\mathbb G_n$]
\label{cor-Gn}
Suppose the conditions of Theorem~\ref{thm-tightness-compact} hold, and the
conditions of Lemma~\ref{lem-centering} hold on $I$.  Then
$\mathbb G_n\rightsquigarrow\mathbb G$ in $\ell^\infty(I)$.
\end{corollary}

\begin{proof}
Since $\mathbb G_n(t) = \mathring{\mathbb G}_n(t)+\sqrt n[F^\star(t)-F(t)]$ and the
second term is deterministic, Lemma~\ref{lem-centering} gives
$\sup_{t\in I}|\mathbb G_n(t)-\mathring{\mathbb G}_n(t)|\to0$.  The claim then
follows from Theorem~\ref{thm-tightness-compact} by Slutsky's lemma in
$\ell^\infty(I)$.
\end{proof}

As discussed after Theorem~\ref{thm-fidi-clt}, centering at $F$ is more natural
for inference, but it requires the stronger sparsity conditions of
Lemma~\ref{lem-centering}, which interact with the smoothness of the density.

\paragraph{Uniformity over the tail.}
Theorem~\ref{thm-tightness-compact} extends to the unbounded index set
$[t_0,\infty)$ if the regularity conditions hold uniformly in the tail, in a
scale-invariant form, and the sparsity condition is strengthened.  Let
$J_\infty := [t_0^-,\infty)$.  A condition below \emph{holds on $J_\infty$} if
its inequalities hold for all $\lambda,\lambda'\in J_\infty$ (and all $x$ with
$\mathfrak D(x)\in J_\infty$) satisfying the stated proximity; in particular a
jump of $f$ at $d_0$ is allowed.  We write $h_t := (t/(n''\rho_n))^{1/2}$ for
the bandwidth at $t$, $g := \mathfrak D-\check{\mathfrak D}$,
$g_t := \sup\{g(x):\mathfrak D(x)\leq2t\}$, and $\Upsilon_n(t) := g_t/h_t$ for the
\emph{clipping ratio}.

\begin{enumerate}[leftmargin=1.5cm]
\item[D3.]\mylabel{D3}\label{condition-d3}
(Density) $f$ is positive on $J_\infty$, and for $\lambda,\lambda'\in J_\infty$
with $|\lambda'-\lambda|\leq\lambda/2$,
\begin{align*}
&\text{(i)}\ c_df(\lambda)\leq f(\lambda')\leq C_df(\lambda); \qquad
\text{(ii)}\ |f(\lambda)-f(\lambda')|\leq K_f\,f(\lambda)
 \Bigl(\frac{|\lambda-\lambda'|}{\lambda}\Bigr)^{\theta_f};
\\
&\text{(iii)}\ c\,\lambda^{-\gamma_+}\leq f(\lambda)\leq C\lambda^{-\gamma_-},
 \qquad 1-F(\lambda)\leq C\lambda f(\lambda),
\end{align*}
with $\theta_f\in(1/4,1]$ and $\gamma_+\geq\gamma_->3+\beta$.

\item[G3.]\mylabel{G3}\label{condition-g3}
(Graphon regularity) For $\lambda,\lambda'\in J_\infty$ with
$|\lambda-\lambda'|\leq\lambda/2$,
\begin{align*}
\bigl\|w(x_\lambda,\cdot)-w(x_{\lambda'},\cdot)\bigr\|_{L^2}
 \leq K_w\,\lambda\Bigl(\frac{|\lambda-\lambda'|}{\lambda}\Bigr)^{\theta_w},
\qquad \theta_w\in(1/4,1].
\end{align*}

\item[L3.]\mylabel{L3}\label{condition-l3}
(Conditional moments) For $x$ with $\mathfrak D(x)\in J_\infty$,
\begin{align*}
E\bigl[w^{2+\beta}(x,\xi)\bigr]\leq K\,\mathfrak D(x)^{2+\beta}, \qquad
E\bigl[w^{4}(x,\xi)\bigr]\leq K\,\mathfrak D(x)^{4q_4},
\end{align*}
for some $q_4\geq1$ with $\gamma_-\geq2q_4$.

\item[S3.]\mylabel{S3}\label{condition-s3}
(Sparsity and tail cutoff) $\rho_n = n^{-b}$ with $b\in(0,1)$, and there is
$a\in(0,b/2)$ such that $\lambda_n := n^a$ satisfies, with
$\vartheta := \min(\theta_w,\theta_f,1/2)$,
\begin{align*}
&\text{(i)}\ \sqrt n\,P(\mathfrak D>\lambda_n/2)\to0; \qquad
\text{(ii)}\ \sup_{t_0\leq t\leq\lambda_n}\Upsilon_n(t)\,t\,h_t^{-\vartheta}\to0;
\qquad
\text{(iii)}\ n^{-1}\rho_n\,\lambda_n^{\,2\gamma_+}\to0;
\\
&\text{(iv)}\ \sup_{t_0\leq t\leq\lambda_n}
 n^{(1-2\vartheta)/4}\rho_n^{(9-2\vartheta)/4}\,t^{(11-2\vartheta)/4}f(t)^{1/2}\to0;
\qquad
\text{(v)}\ b\,(5-2\vartheta)>1-2\vartheta .
\end{align*}
\end{enumerate}

\begin{theorem}[Uniform CLT over the tail]
\label{thm-tightness-ring}
Assume the conditions of Theorem~\ref{thm-fidi-clt}, and
\ref{condition-d3}--\ref{condition-l3} and \ref{condition-s3} on $J_\infty$.
Then $\mathring{\mathbb G}_n\rightsquigarrow\mathbb G$ in
$\ell^\infty([t_0,\infty))$, and $\mathbb G$ has a version with continuous
sample paths on $[t_0,\infty)$.
\end{theorem}

As on compact intervals, the process centered at $F$ requires control of the
centering bias $\sqrt n(F^\star-F)$, now uniformly over the tail.

\begin{proposition}[Centering over the tail]
\label{prop-centering-tail}
Suppose \ref{condition-m1}, \ref{condition-g1} and \ref{condition-d1} hold,
that \ref{condition-d3} holds on $J_\infty$, that \ref{condition-s3}(i) holds for
a sequence $\lambda_n$ with $\rho_n\lambda_n\to0$, and that
\eqref{eq-centering-conditions} holds with $\theta = \theta_f$.  Then
\begin{align*}
\sup_{t\geq t_0}\sqrt n\,\bigl|F^\star(t)-F(t)\bigr|\;\rightarrow\;0 .
\end{align*}
\end{proposition}

The proof, in Appendix~\ref{app-centering-tail},
follows that of Lemma~\ref{lem-centering}.  Below the cutoff $\lambda_n$, the
clipping error is already uniform in $t$, and under \ref{condition-d3} the
lattice term and the smoothing bias $f(t)\,t^{-\theta_f}h_t^{1+\theta_f}$ are bounded
uniformly in $t$ by their values on compact sets.  Beyond the cutoff, both
$1-F$ and $1-F^\star$ are $o(n^{-1/2})$ by \ref{condition-s3}(i).  The
conditions of Proposition~\ref{prop-centering-tail} are those of
Lemma~\ref{lem-centering} with $\theta = \theta_f$; by
Remark~\ref{rem-bias-edge} they are compatible if and only if
$\theta_f(\tfrac12+\beta)>1$.

\begin{remark}[The sparsity window for centering over the tail]
\label{rem-centering-tail-window}
Centering at $F$ requires $\rho_n$ to be not too small,
$n^{\theta_f}\rho_n^{1+\theta_f}\to\infty$, that is $b<\theta_f/(1+\theta_f)\leq\frac12$
when $\rho_n = n^{-b}$, whereas uniformity over the tail requires $\rho_n$ to be
small, $b>b_{\min}$ in Lemma~\ref{lem-feasibility}.  The window is nonempty
only for kernels with enough moments.  For example, with
$\theta_f = 1$, $\vartheta = \frac12$ and $f(\lambda)\asymp\lambda^{-\gamma}$, the
conditions of Corollary~\ref{cor-Gn-tail} can be met for every
$\gamma>3+\beta$ when $\beta\geq0.84$, for some $\gamma$ when
$\frac12<\beta<0.84$, and for no $\gamma$ when $\beta\leq\frac12$; the last
restriction already applies on compact intervals.  For heavier-tailed kernels,
Theorem~\ref{thm-tightness-ring} still applies to the process
$\mathring{\mathbb G}_n$ centered at $F^\star$.
\end{remark}

\begin{corollary}[Uniform CLT for $\mathbb G_n$ over the tail]
\label{cor-Gn-tail}
Suppose the conditions of Theorem~\ref{thm-tightness-ring} and of
Proposition~\ref{prop-centering-tail} hold.  Then
$\mathbb G_n\rightsquigarrow\mathbb G$ in $\ell^\infty([t_0,\infty))$.
\end{corollary}

\begin{proof}
As in the proof of Corollary~\ref{cor-Gn},
$\sup_{t\geq t_0}|\mathbb G_n(t)-\mathring{\mathbb G}_n(t)| = \sup_{t\geq t_0}\sqrt n|F^\star(t)-F(t)|\to0$
by Proposition~\ref{prop-centering-tail}, and the claim follows from
Theorem~\ref{thm-tightness-ring} by Slutsky's lemma in
$\ell^\infty([t_0,\infty))$.
\end{proof}

\subsection{Discussion of the assumptions}

\paragraph{Consequences of the assumptions.}
Three consequences of \ref{condition-g3} and \ref{condition-l3} are used
throughout; the compact conditions give the same bounds for
$\mathfrak D(x)\in J$, with constants in place of the powers of
$\mathfrak D(x)$.  First, the proofs use the centered clipped kernel
$a(x,\xi) := \check w(x,\xi)-\check{\mathfrak D}(x)$ rather than $w$; since
$u\mapsto u\wedge\rho_n^{-1}$ is $1$-Lipschitz, \ref{condition-g3} continues
to hold for $\check w$, and for $a$ with constant $2K_w$.  Second, by Jensen's
inequality and \ref{condition-l3}, and by Markov's inequality applied to
$g(x) = E[(w(x,\xi)-\rho_n^{-1})_+]\leq\rho_n^{1+\beta}E[w^{2+\beta}(x,\xi)]$,
for every $x$ with $\mathfrak D(x)\in J_\infty$,
\begin{align}
\label{eq-a-growth}
\|a(x,\cdot)\|_{L^2}\leq K_2\,\mathfrak D(x), \qquad
\|a(x,\cdot)\|_{L^4}\leq K_4\,\mathfrak D(x)^{q_4}, \qquad
g(x)\leq K\rho_n^{1+\beta}\mathfrak D(x)^{2+\beta}.
\end{align}
Third, the clipping ratio, which measures the clipping bias at degree $t$ on
the scale of the bandwidth, satisfies
\begin{align}
\label{eq-chi-explicit}
\Upsilon_n(t)\leq K'\rho_n^{1+\beta}t^{2+\beta}\sqrt{n''\rho_n/t}
 = K'\,t^{\,3/2+\beta}\rho_n^{1+\beta}\sqrt{n''\rho_n},
\end{align}
and \ref{condition-s3}(ii) requires it to be negligible, uniformly below the
cutoff and relative to $h_t^{\vartheta}/t$.

\paragraph{The compact and the unbounded case.}
On a compact interval, \ref{condition-d3}(i)--(ii) reduce to positivity and
H\"older continuity of $f$, which is \ref{condition-d2}; \ref{condition-g3}
reduces to \ref{condition-g2}; and \ref{condition-l3} reduces to
\ref{condition-l1} together with a local fourth moment.
\ref{condition-s3}(ii) reduces to
$\rho_n^{1+\beta}(n\rho_n)^{(1+\vartheta)/2}\to0$, which is implied by
\ref{condition-s1}, since
$\rho_n^{1+\beta}(n\rho_n)^{(1+\vartheta)/2} =
(n\rho_n^{2+\beta})^{(1+\vartheta)/2}\rho_n^{(1+\beta)(1-\vartheta)/2}$.
\ref{condition-s3}(iv) reduces to
$n^{(1-2\vartheta)/4}\rho_n^{(9-2\vartheta)/4}\to0$, which is implied by
\ref{condition-s3}(v), that is \ref{condition-s2}.  The remaining parts of
\ref{condition-s3}, the cutoff (i) and the comparison (iii), concern only the
tail.  Thus Theorem~\ref{thm-tightness-compact} is
Theorem~\ref{thm-tightness-ring} restricted to a compact index set, with all
tail assumptions removed except the shape condition in \ref{condition-l2}.
With $\rho_n = n^{-b}$ and $\lambda_n = n^a$, every quantity in
\ref{condition-s3} is a constant multiple of a power of $n$, so each
convergence holds at a polynomial rate, and $a<b/2$ gives
$\rho_n\lambda_n^2\to0$.  For bounded $w$, $g\equiv0$ for $n$ large, and
\ref{condition-s3}(ii) is vacuous.  Remark~\ref{rem-s2} interprets
\ref{condition-s2} on compact intervals.


\begin{lemma}[Feasibility for power-law densities]
\label{lem-feasibility}
Suppose $c\lambda^{-\gamma}\leq f(\lambda)\leq C\lambda^{-\gamma}$ on $J_\infty$,
with $\gamma>3+\beta$, and let $\rho_n = n^{-b}$.  Put
$a_0 := 1/(2(\gamma-1))$,
\begin{align*}
b_* := \frac{\tfrac{1+\vartheta}2+\bigl(\tfrac52+\beta-\tfrac\vartheta2\bigr)a_0}
 {1+\beta+\tfrac{1+\vartheta}2},
\qquad
b_{\min} := \max\Bigl\{\frac1{2+\beta},\ \frac{1-2\vartheta}{5-2\vartheta},\ b_*\Bigr\}.
\end{align*}
If $b\in(b_{\min},1)$, then \ref{condition-s1} holds, and
\ref{condition-s3} holds with $\lambda_n = n^a$ for every $a>a_0$ close enough
to $a_0$.  Conversely, if $b\leq b_*$, no $\lambda_n = n^a$ satisfies both
\ref{condition-s3}(i) and the bound on \ref{condition-s3}(ii) supplied by
\eqref{eq-chi-explicit}.
\end{lemma}
The proof, in Appendix~\ref{app-feasibility}, shows that for power-law
densities \ref{condition-s3}(iii) and (iv) are implied by \ref{condition-s1},
and that the binding constraint is the tension between the cutoff (i) and the
clipping condition (ii).

\paragraph{Interpretation of \ref{condition-l3}.}
By the monotone reparametrization of Remark~\ref{remark-inverse}, we may take
$\mathfrak D$ nondecreasing, so that larger values of the latent coordinate
correspond to vertices of higher expected degree.  The first part of
\ref{condition-l3}, and the second part of \ref{condition-l2}, concern the
moment ratio
\begin{align*}
R(x) := \frac{E\bigl[w^{2+\beta}(x,\xi)\bigr]}{\mathfrak D(x)^{2+\beta}}\;\geq\;1,
\end{align*}
where the inequality is Jensen's, and require $R$ to be bounded on
$\{\mathfrak D\geq t_0^-\}$.  Writing $w(x,\cdot) = \mathfrak D(x)\,u(x,\cdot)$
with $E[u(x,\xi)] = 1$, we have $R(x) = E[u^{2+\beta}(x,\xi)]$.  So the
condition says that the \emph{normalized} connection profiles $u(x,\cdot)$
have uniformly bounded moments of order $2+\beta$.  In words, a vertex of high
expected degree is, up to scale, no more concentrated in its connections than
a typical vertex: its degree comes from many connections of comparable
strength rather than from a few exceptionally strong ones.  The degree itself
may be heavy-tailed, as in \ref{condition-d3}(iii); what is controlled is the
shape of each vertex's profile, not its scale.

This holds in degree-corrected models \citep{karrer2011stochastic}
\begin{align*}
w(x,y) = \theta(x)\,\theta(y)\,\varpi(x,y),
\end{align*}
where the degree parameter $\theta$ governs the propensity of a vertex to form
edges and $\varpi$ captures the association structure.  For the remaining
conditions: if $\theta(\xi)$ has a density then $\mathfrak D$ is absolutely
continuous, as \ref{condition-d1} requires; and \ref{condition-g1},
\ref{condition-g2} and \ref{condition-g3} hold when the association structure
does not vary across vertices of the same degree, for instance when $\varpi$
is constant.  For degree-corrected stochastic block models with several
blocks of overlapping degree ranges, \ref{condition-g1} generally fails.

\subsection{Strategy of the proof}
\label{sec-tight-strategy}
By Theorem~\ref{thm-fidi-clt} the finite-dimensional distributions converge,
so by the standard criterion \citep[Theorem~1.5.7]{van2016weak} it remains to
show asymptotic tightness: that the increments
$\mathring{\mathbb G}_n(t)-\mathring{\mathbb G}_n(r)$ are uniformly small over
short intervals.  We do this with Billingsley's chaining criterion, which asks
for a bound on the \emph{product} of two adjacent increments,
\begin{align}
\label{eq-grid-moment-informal}
E\Bigl[\bigl|\mathring{\mathbb G}_n(s)-\mathring{\mathbb G}_n(r)\bigr|^2
 \bigl|\mathring{\mathbb G}_n(t)-\mathring{\mathbb G}_n(s)\bigr|^2\Bigr]
 \;\leq\; C\,\bigl(F(t)-F(r)\bigr)^{2\nu},\qquad \nu>\tfrac12,
\end{align}
for $r<s<t$ (Proposition~\ref{prop-grid} below).  The product is what makes a
chaining argument work without any independence: the criterion needs only
moments, not a structure for the dependence.  By the Cauchy--Schwarz
inequality, \eqref{eq-grid-moment-informal} follows from a bound on each
increment separately, and this is the inequality to which the whole section
reduces: for every $r<t$ in $I$,
\begin{align}
\label{eq-target}
\bigl\|\mathring{\mathbb G}_n(t)-\mathring{\mathbb G}_n(r)\bigr\|_{L^4}
 \;\leq\; C\,\bigl(F(t)-F(r)\bigr)^{\vartheta'}
 \qquad\text{for some }\vartheta'>\tfrac14 .
\end{align}
Fourth moments, rather than second, are forced by the chaining: a
second-moment bound controls each increment but not the maximum over all of
them.  The threshold $\vartheta'>1/4$ is the point at which the chained bound converges.

Two features of the problem shape how \eqref{eq-target} is proved.  Since the
process is piecewise constant, with jumps of variance $\asymp f(t)/(n\rho_n)$,
\eqref{eq-target} can hold only for $r,t$ on the lattice of attainable degree
values, and the criterion must be posed there (Step~L).  And the increment has
two sources of randomness, the edges and the latent positions:
\begin{align}
\label{eq-split-informal}
\mathring{\mathbb G}_n(t)-\mathring{\mathbb G}_n(r) = B_A+Z_A,
\quad
E B_A^4\lesssim\bigl(F(t)-F(r)\bigr)^2,
\quad
E Z_A^4\lesssim\bigl(F(t)-F(r)\bigr)^{4\vartheta'},
\end{align}
where $B_A$ is the edge noise and $Z_A$ the latent part.  It is useful to
think of the latent part in turn as the sum of its H\'ajek projection and the
higher-order terms of its Hoeffding decomposition.  Our proof does not
separate these two explicitly, since $Z_A$ is bounded as a whole by moment
inequalities of Efron--Stein type (Steps~Z2 and~Z4), and $B_A$ by conditional
cumulants (Step~B); but the bounds organize along these lines, and each of the
three terms presents a different difficulty.

\emph{The H\'ajek part.}  The first-order Hoeffding component of an increment
consists of two terms: $\psi_A$, the effect of a vertex's own latent position
on its own degree, and $\phi_A$, the aggregate effect of that position on the
degrees of all other vertices.  The first behaves like an increment of a
classical empirical process, with variance at most $p_A\asymp F(t)-F(r)$.  The
second is where the lattice enters: $\phi_A$ is an integral of a difference of
binomial masses across the cell, and a bound in absolute value is too large
for cells below the bandwidth.  Its size, of order $\Xi_A^{1/2}$, follows only
from the cancellation of Step~W combined with the H\"older regularity of $f$
and of the kernel.  In terms of the target \eqref{eq-target}, this gives
$\|\phi_A\|_2\lesssim\Xi_A^{1/2}\lesssim(F(t)-F(r))^{\vartheta}$: the latent
fluctuation over a cell is controlled by the $\vartheta$-th power of its
$F$-measure, compared with the power $\frac12$ for the empirical-process part
$\psi_A$.  The H\"older regularity of the kernel and the density thus sets the
regularity of the limit, and this is why the assumptions require
$\vartheta>\frac14$, so that fourth moments of order
$(F(t)-F(r))^{4\vartheta'}$ with $4\vartheta'>1$ are available for chaining.

\emph{The higher-order Hoeffding terms.}  The degenerate components are
numerous and overlap, and the fourth-moment criterion requires more than
their variances.  Each additional order costs a factor $\rho_n$, so the
components form geometric series, but pairs of components interact in three
distinct ways (classes ($\alpha$)--($\gamma$) in Step~Z2), and one of them
requires the cancellation a second time.  At the level of fourth moments, the
Efron--Stein approach reduces everything to an exact operator identity
(Step~Z4), avoiding an explicit fourth-moment expansion of the Hoeffding
decomposition; this is the one place where the exponent drops from the
H\"older exponent $\vartheta$ of the assumptions to $\vartheta'$, a little above
$1/4$.  The sparsity condition \ref{condition-s2} enters in these two places.

\emph{The edge noise.}  Given the latent positions, $B_A$ is a normalized sum
of centered indicators $\mathbbm 1(D_i\in A)-q_i(A)$, dependent only through
shared edges.  Its size is governed by the probability of landing in the
cell: the conditional variance is essentially $n^{-1}\sum_iq_i(A)$, whose
mean is $p_A\asymp F(t)-F(r)$, exactly as for the binomial count of an
empirical process, and its fourth moment is of order
$(F(t)-F(r))^2+n^{-1}(F(t)-F(r))$, as for a binomial.  The dependence through
shared edges adds covariance and cumulant terms, but these are of the same
order: two vertices interact through a cell only if both respond to a shared
edge, and each response is damped by the factor $\mu_A$.  The term
$n^{-1}(F(t)-F(r))$, which forces the product form of the criterion for the
classical empirical process, is harmless here, because on the lattice
$F(t)-F(r)$ is bounded below by a negative power of $n$.

Inequality \eqref{eq-target} is proved as Lemma~\ref{lem-c-single}, and the
reduction from \eqref{eq-target} to the theorem is Steps~H and~L.

We prove Theorem~\ref{thm-tightness-compact} in the main text
(Section~\ref{sec-tight-compact-proof}).  The full proofs of the two latent
moment bounds, which are the technical core, are in
Appendix~\ref{app-latent}; routine estimates are deferred to
Appendix~\ref{app-compact-proofs}, and each deferred proof is replaced in the main text by a
one- or two-sentence summary.  The proof of
Theorem~\ref{thm-tightness-ring}, in Appendix~\ref{app-weighted}, follows the
same strategy, step by step, with the density, the degree and the bandwidth
carried as weights; it involves considerably more bookkeeping, but no new
idea beyond those described at the end of this subsection.  Table~\ref{tab-roadmap} summarizes the steps, the bound each supplies, and
where its proofs are given.  Several features make these bounds delicate.

\paragraph{The lattice.}  As explained in the introduction, increments are
controlled on lattice cells, by their $F$-measure, and the mechanism that
makes an increment small depends on whether the cell is longer than the
bandwidth $h_t\asymp(t/(n\rho_n))^{1/2}$ at which the degree of a vertex with
expected degree $t$ fluctuates: the regularity of $\lambda\mapsto w(x_\lambda,\cdot)$
and of $f$ if it is, a cancellation between adjacent values of a binomial
mass function if it lies within a single fluctuation window.  An exact identity, that $\int_0^1b_M(m;p)\,dp$ does
not depend on $m$, captures both mechanisms in a single bound. The identity has a simple probabilistic meaning: if $p$ is uniform on
$[0,1]$ and the count is $\mathrm{Bin}(M,p)$ given $p$, then the count is
uniform on $\{0,\ldots,M\}$, so mixing over the success probability erases the
lattice.  In our setting the density of latent expected degrees plays the role
of the uniform mixing distribution; it is nearly constant on the scale of the
bandwidth, so adjacent binomial masses nearly cancel, and what remains is
controlled by the variation of the density across the cell.  The cancellation
is used only for averages over the latent positions of other vertices: the
H\'ajek kernel describing the aggregate effect of one vertex on the degrees of
the others, and its higher-order analogues in the fourth-moment bounds.  The
fluctuations of an individual degree at a fixed latent position do not cancel;
they are controlled instead by localization, since only vertices whose
expected degree lies within a few bandwidths of the cell contribute.  On the
lattice, every cell below the cutoff has $F$-measure at least a negative power
of $n$, and this makes it possible to deduce Billingsley's product condition
from fourth-moment bounds on single cells.

\paragraph{Dependence.}  The modern tools for uniform central limit theorems,
such as symmetrization, bracketing and the associated maximal inequalities
(see, e.g., \citet{van2016weak}), rest on independence of the summands, or on
dependence that is weak in a quantifiable sense, such as mixing (see, e.g.,
\citet{rio2017asymptotic}).  Neither holds here: every pair of vertices shares
an edge and all vertices share the latent field, so the dependence among the
indicators does not decay.  Billingsley's criterion requires only moment
bounds on increments and is indifferent to the dependence structure that
produces them; the two sources of dependence are then handled by the
different tools described above.

\paragraph{Clipped hubs.}  A vertex whose expected degree exceeds the clipping
level $\rho_n^{-1}$ by a large factor can have a moderate clipped expected
degree while its degree, given the latent positions, is nearly
deterministic.  For such vertices the local estimates on which everything
rests fail.  They are rare, however, and we remove them as \emph{roots} of the
process: on an event of probability tending to one this changes the process
only by a deterministic shift that vanishes uniformly.  They continue to enter
as neighbors, where no local estimate is needed.

\paragraph{The tail.} Uniformity over $[t_0,\infty)$ is proved in
Appendix~\ref{app-weighted} by the same argument, with $f(t)$, $t$ and $h_t$
carried as weights in every bound and with the process shown to be flat
beyond the cutoff $\lambda_n$.  On a compact index set all weights are of
order one, and the argument is the one given in
Section~\ref{sec-tight-compact-proof}.

\begin{remark}[Where the additional care enters]
\label{rem-w-summary}
For reference, the places where the proof of Theorem~\ref{thm-tightness-ring}
(Appendix~\ref{app-weighted}) departs from the
compact one, and the assumption that pays for each:
\begin{enumerate}[(1)]
\item the polynomial lower bounds \eqref{eq-w-poly}, $\Delta_X\geq n^{-K_0}$ and
$nf(t)h_t\geq n^{c_1}$, which replace $\Delta_A\gtrsim(n\rho_n)^{-1}$ and are
used in every absorption of a negligible term and of a factor $n^{-1}$
(Lemmas~\ref{lem-w-sigma}, \ref{lem-w-edge}, \ref{lem-w-Z4}, \ref{lem-w-absorb});
they use \ref{condition-d3}(iii) and \ref{condition-s3}(iii);
\item one-sided cells in place of the splitting of long cells (Step~C,
Corollary~\ref{cor-w-long}), with the H\'ajek part bounded without
cancellation by $sf(s)$ (Proposition~\ref{prop-w-hajek}) and $p'_{O_s}\leq Csf(s)$
from \ref{condition-d3}(iii);
\item the scale-invariant H\"older moduli $(\bar\ell_X/t)^{\theta}$ of
\ref{condition-d3}(ii) and \ref{condition-g3}, which make the cancellation
lemma and the H\'ajek increment uniform in $t$;
\item the $t$-dependent clipping ratio $\Upsilon_n(t)$, controlled by
\ref{condition-s3}(ii) in Lemmas~\ref{lem-w-cancel}, \ref{lem-w-Z2} (class
($\beta$)) and \ref{lem-w-Z4};
\item the growth $t^{q_4}$ of fourth moments, \ref{condition-l3}, which enters
the remainder ratios of Step~Z4 and is matched against $f(t)$ through
constraint (iv) of \eqref{eq-w-theta-prime};
\item class ($\beta$) of Step~Z2, whose leading term is $v_n(t)$, hence
\ref{condition-s3}(iv); and the product $P_X$ of Step~Z4, hence
\ref{condition-s3}(v);
\item flatness beyond the cutoff, \ref{condition-s3}(i), in the proof of the
theorem.
\end{enumerate}
On a compact index set $f(t)\asymp t\asymp1$, items (1)--(5) are vacuous,
\ref{condition-s3}(iv) reduces to $v_n\to0$ and \ref{condition-s3}(v) is
\ref{condition-s2}, and the argument is that of
Section~\ref{sec-tight-compact-proof}.
\end{remark}
\begin{table}[t]
\caption{Structure of the proof of Theorem~\ref{thm-tightness-compact}.
Each step lists its main results, the bound it supplies, the earlier steps
and conditions it uses, and where its proofs are given (M: main text;
A: Appendix~\ref{app-latent}; S: Appendix~\ref{app-compact-proofs}).  Bounds are stated up to constants and negligible
terms.}
\label{tab-roadmap}
\small
\setlength{\tabcolsep}{4pt}
\begin{tabular}{@{}l >{\raggedright\arraybackslash}p{2.3cm} >{\raggedright\arraybackslash}p{5.6cm} >{\raggedright\arraybackslash}p{2.1cm} l@{}}
\hline
Step & Results & Main bound & Uses & Proofs\\
\hline
H & Lemmas~\ref{lem-c-hubs},~\ref{lem-c-hub-reduction}
  & hubs are rare; $\sup_t|\mathring{\mathbb G}_n-\mathbb G'_n|\to0$
  & \ref{condition-l1}, \ref{condition-l2}, \ref{condition-s1} & S, M\\
L & Prop.~\ref{prop-grid}
  & product moment bound on lattice cells $\Rightarrow$ weak convergence
  & --- & S\\
C & ---
  & long cells reduce to short ones & --- & M\\
W & Lemmas~\ref{lem-poisson-binomial}--\ref{lem-c-cancel}
  & $|g|\lesssim m_A\omega_A$;\ $n\rho_n|\int\mathsf K^{(d)}_A|\lesssim\sigma_t^{-d}\mu_A\bar\ell_A^{\vartheta}$
  & \ref{condition-d2}, \ref{condition-g2}, \ref{condition-l2} & S\\
K & Lem.~\ref{lem-c-kernels}, Prop.~\ref{prop-c-hajek}
  & $E\psi_A^2\lesssim\Delta_A$;\ $\|\phi_A\|_2\lesssim\Xi_A^{1/2}$
  & W & M\\
B & Lemmas~\ref{lem-c-pairing}--\ref{lem-c-edge}
  & $EB_A^4\lesssim\Delta_A^2+n^{-1}\|Z_A\|_2^2+n^{-1}\Delta_A$
  & W, K, Z2 & M, S\\
Z2 & Lemma~\ref{lem-c-Z2}
  & $\|Z_A\|_2^2\lesssim\Delta_A+\Xi_A$
  & W, K; \ref{condition-s2} & A\\
Z4 & Lemma~\ref{lem-c-Z4}
  & $\|Z_A\|_4\lesssim\Delta_A^{\vartheta'}$
  & W, K, Z2; \ref{condition-s2} & A\\
A & Lemma~\ref{lem-c-single}
  & $\|\Delta'_A\|_4\lesssim\Delta_A^{\vartheta'}$ $\Rightarrow$ Theorem~\ref{thm-tightness-compact}
  & H, L, C, B, Z2, Z4 & M\\
\hline
\end{tabular}
\end{table}

\subsection{Proof of Theorem~\ref{thm-tightness-compact}}
\label{sec-tight-compact-proof}

Table~\ref{tab-notation} collects the notation of the proof.

\begin{table}[t]
\caption{Notation for the proof of Theorem~\ref{thm-tightness-compact}, in order
of appearance.  The weighted proof of Theorem~\ref{thm-tightness-ring} uses the
same notation, with the additions listed at the start of
Section~\ref{app-weighted}.}
\label{tab-notation}
\small
\setlength{\tabcolsep}{4pt}
\begin{tabular}{@{}l >{\raggedright\arraybackslash}p{7.4cm} l@{}}
\hline
Symbol & Meaning & Defined in\\
\hline
\multicolumn{3}{@{}l}{\emph{Lattice and cells}}\\
$\lambda_m$, $t'$ & lattice point $m/(n'\rho_n)$; the integer $\lfloor n'\rho_nt\rfloor$ & Step~L\\
$A = A_{rt}$ & cell $\{r'+1,\ldots,t'\}$ between lattice points $r<t$ & \eqref{eq-c-cell}\\
$\ell_A$, $\Delta_A$ & length $t-r$; $F$-measure $F(t)-F(r)$ & \eqref{eq-c-cell}\\
$p_A$, $p'_A$ & $P(D_1\in A)$; the same on the hub-free event & \eqref{eq-c-cell}, \eqref{eq-c-BZ}\\
\multicolumn{3}{@{}l}{\emph{Scales}}\\
$h_t$, $\sigma_t$ & bandwidth $(t/(n''\rho_n))^{1/2}$; degree s.d.\ $n''\rho_nh_t$ & Conventions\\
$\mu_A$, $m_A$ & relative size $\min(1,\ell_A/h_t)$ of a mass difference; $\sigma_t^{-1}\mu_A$ & \eqref{eq-c-cell-quantities}\\
$\bar\ell_A$ & $\max(h_t,\ell_A)$, the scale of the H\"older error & \eqref{eq-c-cell-quantities}\\
$Q_A$, $\Xi_A$ & $h_t\mu_A^2$ (interactions); $\mu_A^2\bar\ell_A^{2\vartheta}$ (H\'ajek part) & \eqref{eq-c-cell-quantities}\\
$\vartheta$, $\vartheta'$ & $\min(\theta_w,\theta_f,\frac12)$; the exponent in the target, $>\frac14$ & \ref{condition-s2}, \eqref{eq-c-theta-prime}\\
$\mathfrak e_n$ & negligible quantity, $\leq n^{-K}$ for every $K$ & Conventions\\
\multicolumn{3}{@{}l}{\emph{Hubs and clipping}}\\
$g$, $\mathscr H$ & clipping bias $\mathfrak D-\check{\mathfrak D}$; hub set $\{g>\mathfrak D/4\}$ & Step~H\\
$\chi$, $\mathbb G'_n$ & indicator of a non-hub; the hub-free process & \eqref{eq-c-G-prime}\\
$\Upsilon_n$ & relative clipping error on the windows & \eqref{eq-c-clip}\\
\multicolumn{3}{@{}l}{\emph{The increment}}\\
$\Delta'_A$ & increment of $\mathbb G'_n$ over $A$ & \eqref{eq-c-BZ}\\
$q_i(A)$, $\bar p_A$ & $P(D_i\in A\mid\mathcal G)$; $n^{-1}\sum_i\chi_iq_i(A)$ & \eqref{eq-c-BZ}\\
$B_A$, $Z_A$ & edge noise; latent part, $\sqrt n(\bar p_A-p'_A)$ & \eqref{eq-c-BZ}\\
\multicolumn{3}{@{}l}{\emph{Localization and kernels}}\\
$z_s$, $\omega^{[j]}_A$ & standardized distance to an endpoint; Gaussian localization weight & \eqref{def-omega-e-a}\\
$e_A$ & distance in expected degree to the reference point $x_t$ & \eqref{def-omega-e-a}\\
$W_s$, $W_A$ & windows of width about $h_t$ around $s$ and around the cell & Step~W\\
$\mathsf K^{(d)}_A$ & $d$-th order mass difference across the cell & \eqref{eq-c-K}\\
$\psi_A$, $\phi_A$ & H\'ajek kernels: own position; effect on other degrees & \eqref{eq-c-hajek}\\
$g^{(j)}_i(A)$ & effect of the edge $\{i,j\}$ on $q_i(A)$ & \eqref{eq-c-pairing}\\
\multicolumn{3}{@{}l}{\emph{Moment bounds}}\\
$\bar\Sigma(A)$ & conditional pairing sum & \eqref{eq-c-Sigma-bar}\\
$\kappa$, $\kappa_4$ & conditional joint and fourth cumulants & Step~B\\
$\Phi_m$, $\Psi_m$ & resampling differences in the Efron--Stein bounds & \eqref{eq-c-Phi-Psi}\\
$\mathcal T_m$, $P_A$ & resampling operator; the product $\rho_n^2\mu_A^2\Xi_A$ & \eqref{eq-c-T-identity}, \eqref{eq-c-PA}\\
\hline
\end{tabular}
\end{table}

\paragraph{Lattice and cells.}
Since $\mathbbm 1(D_i\leq t')$ depends on $t$ only through
$t' = \lfloor n'\rho_nt\rfloor$, the process $\mathring{\mathbb G}_n$ is
constant between consecutive \emph{lattice points} $\lambda_m := m/(n'\rho_n)$,
$m\in\mathbb N$.  For lattice points $r = \lambda_{r'}<t = \lambda_{t'}$ the
corresponding \emph{cell} is $A = A_{rt} := \{r'+1,\ldots,t'\}$, and we write
\begin{align}
\label{eq-c-cell}
k_A := t'-r',\quad \ell_A := t-r = \frac{k_A}{n'\rho_n},\quad
\Delta_A := F(t)-F(r),\quad p_A := P(D_1\in A),
\end{align}
for its lattice length, its length, its $F$-measure and its probability.  The
increment over the cell is
$\mathring{\mathbb G}_n(t)-\mathring{\mathbb G}_n(r)
= n^{-1/2}\sum_{i=1}^n[\mathbbm 1(D_i\in A)-p_A]$.  For $i\neq j$ and
$J'\subseteq[n]\setminus\{i\}$ we use $D_{i\setminus J'}$ as in
Section~\ref{sec-multilinear-hoeffding}, i.e.\ the degree of $i$ with the edges
to $J'$ removed.

\paragraph{Conventions.}  Constants $C,c$ depend on $I$, $J$ and the constants
in the assumptions, but not on $n$ or on any cell.  Integrals over latent
positions are with respect to Lebesgue measure on $[0,1]$, and we omit the
differential when the variable of integration is clear.  A nonnegative
$\mathfrak e_n$ is \emph{negligible} if $\mathfrak e_n\leq n^{-K}$ for every
fixed $K$; Chernoff tails $e^{-cn\rho_n}$ are negligible since
$n\rho_n = n^{1-b}$.  For a lattice point $t\in I$, $h_t := (t/(n''\rho_n))^{1/2}$
and $\sigma_t := n''\rho_nh_t$, so that uniformly on $I$,
$h_t\asymp(n\rho_n)^{-1/2}$ and $\sigma_t\asymp(n\rho_n)^{1/2}$.  All cells are
$A = A_{rt}$ with $t_0\leq r<t\leq t_1$.  We classify a cell by comparing its
length with the two scales of the problem, the bandwidth $h_t$ and the
location $t$:
\begin{align}
\label{eq-c-cell-types}
\underbrace{\ell_A\leq h_t}_{\text{small regime}}
\qquad
\underbrace{h_t<\ell_A\leq t/2}_{\text{large regime}}
\qquad
\underbrace{\ell_A>t/2}_{\text{long cell}} .
\end{align}
Cells in the first two classes are called \emph{short}.  The comparison with
$h_t$ decides which mechanism makes the increment small: in the small regime
both endpoints lie in a single fluctuation window, and smallness comes from
cancellation between adjacent binomial masses, while in the large regime it
comes from the regularity of $f$ and of $\lambda\mapsto w(x_\lambda,\cdot)$.
The comparison with $t$ is a matter of bookkeeping: for a short cell
$r\geq t/2$, so the bandwidths, variances and densities at the two endpoints
are comparable and the regularity conditions apply across the cell, and every
estimate can be written at the single scale $t$.  Long cells need no
smallness, since their $F$-measure is bounded below, and they are reduced to
short ones in Step~C.  For a short cell let
\begin{align}
\begin{split}
\label{eq-c-cell-quantities}
m_A &:= \sigma_t^{-1}\mu_A,\quad \mu_A := \min(1,\ell_A/h_t),\quad
\bar\ell_A := \max(h_t,\ell_A),
\\ Q_A &:= h_t\mu_A^2,\qquad
\Xi_A := \mu_A^2\,\bar\ell_A^{\,2\vartheta}.
\end{split}
\end{align}
These quantities have simple interpretations.  The factor $\mu_A$ is the
relative size of a difference of mass functions across the cell: of order one
when the endpoints lie in different fluctuation windows, and of order
$\ell_A/h_t$ when they lie in the same window, where adjacent masses nearly
cancel.  Quantities involving two such differences, such as covariances
between vertices through a shared edge, are therefore of order
$Q_A = h_t\mu_A^2$, the factor $h_t$ being the measure of a window.  The H\'ajek
kernel of the cell involves one such difference, integrated against the
connection profiles; by the cancellation of Step~W, only the variation of the
density and of the profiles across the window survives, of order
$\bar\ell_A^{\,\vartheta}$, so that its squared norm is of order
$\Xi_A = \mu_A^2\bar\ell_A^{\,2\vartheta}$.

Since $f$ is bounded above and below on $J$ by \ref{condition-d2},
$c\ell_A\leq\Delta_A\leq C\ell_A$ for every cell.  This is the main point at
which the compact setting simplifies the argument: over the tail $f$ is not
bounded below, and the corresponding comparison \eqref{eq-w-Delta-ell} carries
the weight $f(t)$ and holds only for short cells.  In particular
\begin{align}
\label{eq-c-Q-Xi}
Q_A\leq C\Delta_A,\qquad \Xi_A\leq C\Delta_A^{2\vartheta},\qquad
\Delta_A\geq c\,(n'\rho_n)^{-1}\geq c\,n^{-(1-b)} ,
\end{align}
the second because $\mu_A^2h_t^{2\vartheta} = \ell_A^{2\vartheta}(\ell_A/h_t)^{2-2\vartheta}
\leq\ell_A^{2\vartheta}$ in the small regime, so that the H\'ajek part is
controlled by the $\vartheta$-th power of the $F$-measure of the cell.  We fix $\vartheta'\in(1/4,\vartheta]$
with
\begin{align}
\begin{split}
\label{eq-c-theta-prime}
\text{(i) } &(1-b)(4\vartheta'-1)<b/2,\quad
\text{(ii) }4\vartheta'-1<b(4\vartheta'+1),
\\ \text{(iii) }&\vartheta'\leq5/16,\quad
\text{(iv) }(2\vartheta'-\vartheta)(1-b)\leq2b ;
\end{split}
\end{align}
(i)--(iii) hold for $\vartheta'$ close to $1/4$, and (iv) can then be met if and
only if \ref{condition-s2} holds.  By (i) and \eqref{eq-c-Q-Xi}, for every cell,
\begin{align}
\label{eq-c-absorb}
\varepsilon\,\Delta_A\leq\Delta_A^{4\vartheta'}\ \text{ whenever }
0\leq\varepsilon\leq n^{-b/2},\qquad
\mathfrak e_n\leq\Delta_A^{4\vartheta'}\ \text{ for negligible }\mathfrak e_n ,
\end{align}
for $n\geq n_0$.  In particular this applies with $\varepsilon = \rho_n$ and
with $\varepsilon = n^{-1/2}$.
Moreover, since $\vartheta'>\frac14\geq\frac\vartheta2$, for every short cell
\begin{align}
\label{eq-c-Xi-Delta}
\Xi_A\leq C\,h_t^{\,2\vartheta-4\vartheta'}\,\Delta_A^{4\vartheta'} :
\end{align}
in the small regime $\Xi_A = \mu_A^2h_t^{2\vartheta}$, $\Delta_A\geq c\mu_Ah_t$ and
$\mu_A^{2-4\vartheta'}\leq1$; in the large regime $\Xi_A = \ell_A^{2\vartheta}$,
$\Delta_A\geq c\ell_A$ and $\ell_A^{2\vartheta-4\vartheta'}\leq h_t^{2\vartheta-4\vartheta'}$,
since $2\vartheta-4\vartheta'<0$.

\paragraph{Step H: hub reduction.}
Let $\mathscr H := \{x:\ g(x)>\mathfrak D(x)/4\}$ denote the values of $x$ corresponding to hub nodes, with
$g = \mathfrak D-\check{\mathfrak D}$.  The following lemma establishes that hubs are rare, which will allow us to argue that the empirical process is not perturbed much even if they are ignored. 

\begin{lemma}[Hubs]
\label{lem-c-hubs}
For $n\geq n_0$: $\mathrm{Leb}(\mathscr H)\leq C\rho_n^{2+\beta}$;
$\mathscr H\cap\{\mathfrak D\leq2t_1^+\} = \emptyset$; and
$g(x)\leq g_* := C\rho_n^{1+\beta}$ whenever $\mathfrak D(x)\leq2t_1^+$.  If
$\mathrm{supp}\,F$ is bounded, $\mathscr H = \emptyset$.
\end{lemma}

The proof is a Markov-inequality argument: on bounded sets of
$\mathfrak D$, \ref{condition-l1} makes the clipping bias $g$ of order
$\rho_n^{1+\beta}$, far below $\mathfrak D/4$, while \ref{condition-l2} forces
every hub to have $\mathfrak D\gtrsim\rho_n^{-1}$, a set of measure
$O(\rho_n^{2+\beta})$ by \ref{condition-m1}.  Details are in
Appendix~\ref{app-c-H}.

Let $\chi(x) := \mathbbm 1(x\notin\mathscr H)$ and $\chi_i := \chi(\xi_i)$.  We
remove hub vertices as \emph{roots} of the process: let
\begin{align}
\label{eq-c-G-prime}
F'_n(t) := E\bigl[\chi_1\mathbbm 1(D_1\leq t')\bigr],\qquad
\mathbb G'_n(t) := \frac1{\sqrt n}\sum_{i=1}^n\bigl[\chi_i\mathbbm 1(D_i\leq t')-F'_n(t)\bigr].
\end{align}
Hub vertices still enter $\mathbb G'_n$ as neighbors, through the edges of the
other vertices.

\begin{lemma}[Hub reduction]
\label{lem-c-hub-reduction}
$P(\xi_i\in\mathscr H\text{ for some }i)\leq Cn\rho_n^{2+\beta}\to0$, and
\begin{align*}
\sup_{t\geq0}\bigl|\mathring{\mathbb G}_n(t)-\mathbb G'_n(t)\bigr|\to0
\qquad\text{in probability.}
\end{align*}
Consequently $\mathbb G'_n\stackrel{fdd}{\rightsquigarrow}\mathbb G$, and
$\mathbb G'_n\rightsquigarrow\mathbb G$ in $\ell^\infty(I)$ if and only if
$\mathring{\mathbb G}_n\rightsquigarrow\mathbb G$ in $\ell^\infty(I)$.
\end{lemma}

\begin{proof}
The first claim is a union bound with Lemma~\ref{lem-c-hubs} and
\ref{condition-s1}.  Let $\mathscr E^{\mathscr H}_n$ be the event that no
$\xi_i$ lies in $\mathscr H$.  On $\mathscr E^{\mathscr H}_n$ all $\chi_i = 1$,
so for every $t$,
\begin{align*}
\mathring{\mathbb G}_n(t)-\mathbb G'_n(t) = \sqrt n\bigl(F'_n(t)-F^\star(t)\bigr)
 = -\sqrt n\,E\bigl[(1-\chi_1)\mathbbm 1(D_1\leq t')\bigr],
\end{align*}
and therefore
\begin{align*}
\sup_t\bigl|\mathring{\mathbb G}_n(t)-\mathbb G'_n(t)\bigr|
 \leq\sqrt n\,\mathrm{Leb}(\mathscr H)\leq C\bigl(n\rho_n^{4+2\beta}\bigr)^{1/2}\to0 .
\end{align*}
Hence
\begin{align*}
P\Bigl(\sup_t|\mathring{\mathbb G}_n-\mathbb G'_n|>\epsilon\Bigr)\leq
P\bigl((\mathscr E^{\mathscr H}_n)^c\bigr)+\mathbbm 1\Bigl(C\bigl(n\rho_n^{4+2\beta}\bigr)^{1/2}>\epsilon\Bigr)\to0 .
\end{align*}
The consequences follow from Theorem~\ref{thm-fidi-clt} and Slutsky's lemma in
$\ell^\infty(I)$.
\end{proof}

From here on we work with $\mathbb G'_n$.  For a cell $A$, let
$q_i(A) := P(D_i\in A\mid\mathcal G)$, $p'_A := E[\chi_1\mathbbm 1(D_1\in A)]$,
$\bar p_A := n^{-1}\sum_i\chi_iq_i(A)$, and decompose the hub-free increment as
\begin{align}
\begin{split}
\label{eq-c-BZ}
\Delta'_A &:= \mathbb G'_n(t)-\mathbb G'_n(r) = B_A+Z_A, \text{ where:}
\\
B_A &:= \frac1{\sqrt n}\sum_i\chi_i\bigl[\mathbbm 1(D_i\in A)-q_i(A)\bigr],\quad
Z_A := \sqrt n\,(\bar p_A-p'_A).
\end{split}
\end{align}
$B_A$ is the \emph{edge noise}: given $\mathcal G$, it is a centered function of
the independent edge variables.  $Z_A$ is the \emph{latent part}, a function of
$\xi_1,\ldots,\xi_n$.  Vertices in $\mathscr H$ never appear as roots in either
term.  We also write $\Upsilon_n := g_*/h_{t_0}\leq C\rho_n^{1+\beta}(n\rho_n)^{1/2}$
for the \emph{clipping ratio}.  Since $n\rho_n^{2+\beta}\to0$,
\begin{align}
\label{eq-c-clip}
\Upsilon_n\,(n\rho_n)^{\vartheta/2}\leq C\rho_n^{1+\beta}(n\rho_n)^{(1+\vartheta)/2}
 \leq C\bigl(n\rho_n^{2+\beta}\bigr)^{(1+\vartheta)/2}
 \rho_n^{(1+\beta)(1-\vartheta)/2}\to0 ,
\end{align}
so $\Upsilon_n\leq h_t^\vartheta\leq\bar\ell_A^\vartheta$ for $n\geq n_0$.

\paragraph{Step L: lattice criterion.}
The following is a lattice adaptation of a chaining argument of
\citet{billingsley1999convergence}; for extensions of such moment criteria
to general control measures, see \citet{bickel1971convergence}.  It is stated for $[t_0,\infty)$ with a
cutoff, so that it also serves the proof of Theorem~\ref{thm-tightness-ring};
here it is applied on $I$ with $T_n\equiv t_1$.

\begin{proposition}[Grid criterion with a right-endpoint cutoff]
\label{prop-grid}
Let $(\mathbb H_n)$ be processes on $[t_0,\infty)$ that are constant between
consecutive lattice points, and suppose there is a nondecreasing sequence
$T_n\in(t_0,\infty]$ with
$\sup_{\lambda_m\leq T_n}\bigl(F(\lambda_{m+1})-F(\lambda_m)\bigr)\to0$ such
that:
\begin{enumerate}[(i)]
\item $\mathbb H_n\stackrel{fdd}{\rightsquigarrow}\mathbb H$, where $\mathbb H$
has continuous sample paths;
\item there are $\nu>1/2$ and $C<\infty$ such that for all $n$ and all lattice
points $t_0\leq r\leq s\leq t\leq T_n$,
\begin{align}
\label{eq-grid-moment}
E\bigl[|\mathbb H_n(s)-\mathbb H_n(r)|^2\,|\mathbb H_n(t)-\mathbb H_n(s)|^2\bigr]
 \leq C\,\bigl(F(t)-F(r)\bigr)^{2\nu};
\end{align}
\item $\sup_{t\geq T_n}|\mathbb H_n(t)-\mathbb H_n(T_n)|\to0$ in probability.
\end{enumerate}
Then $\mathbb H_n\rightsquigarrow\mathbb H$ in $D[t_0,\infty)$, and hence in
$\ell^\infty[t_0,\infty)$.  If $T_n\equiv\infty$, condition (iii) is vacuous.
The same holds on $[t_0,t_1]$ with $T_n\equiv t_1$, in which case (iii) is void.
\end{proposition}

In the parametrization $u = F(t)$ the process is constant between the images
of lattice points, so Billingsley's fluctuation inequality applies to its
values on the grid and yields his modulus bound below $F(T_n)$; beyond $F(T_n)$
the modulus is controlled by (iii).  The proof is in
Appendix~\ref{app-c-L}.

\begin{remark}
\label{rem-why-lattice}
The criterion must be posed on the lattice.  For non-lattice $r<t$ straddling a
single lattice point, $t-r$ can be arbitrarily small while the increment has
variance $\asymp P(D_1 = t')\asymp f(t)/(n\rho_n)$, so no fixed continuous
modulus can satisfy \eqref{eq-grid-moment} for all $r<t$.  On the lattice,
$p_A\lesssim\Delta_A$ (Proposition~\ref{prop-c-hajek}), and increments are small
exactly when their $F$-length is.
\end{remark}

We apply Proposition~\ref{prop-grid} on $I$ with $T_n\equiv t_1$.  The grid
mesh is $\sup_m(F(\lambda_{m+1})-F(\lambda_m))\leq C/(n\rho_n)\to0$ by
\ref{condition-d2}.  Conditions (i) and (ii) are verified in Step~A.

\paragraph{Step C: reduction to short cells.}
Let $A = A_{rt}$ be a long cell.  Choose lattice points
$r = s_0<s_1<\dots<s_k = t$ with each $s_j$ the largest lattice point not
exceeding $\min(\tfrac32s_{j-1},t)$.  Then
\begin{align*}
\ell_{A_{s_{j-1}s_j}} = s_j-s_{j-1}\leq\tfrac12s_{j-1}\leq\tfrac12s_j,
\end{align*}
so each $A_{s_{j-1}s_j}$ is short, and since
$s_j\geq\tfrac32s_{j-1}-(n'\rho_n)^{-1}\geq\tfrac54s_{j-1}$ for $n\geq n_0$, we
have $k\leq\log_{5/4}(t_1/t_0)+2$.  By Minkowski's inequality,
\begin{align*}
\|\Delta'_A\|_{L^4}\leq\sum_{j=1}^k\|\Delta'_{A_{s_{j-1}s_j}}\|_{L^4},
\end{align*}
so a bound $\|\Delta'_B\|_{L^4}\leq C\Delta_B^{\vartheta'}$ for short cells $B$
gives $\|\Delta'_A\|_{L^4}\leq Ck\,\Delta_A^{\vartheta'}$.  \emph{From here on every
cell is short.}

\paragraph{Step W: windows, roots and cancellation.}

The edge noise is a normalized sum of centered indicators
$\mathbbm 1(D_v\in A)-q_v(A)$, so its moments involve products of such
indicators for several vertices.  Conditionally on the latent positions and on
the edges among those vertices, the remaining degree of each vertex is a
Poisson-binomial variable, independent across vertices, so the products
factorize into Poisson-binomial probabilities of landing in shifted copies of
the cell.  Averaging over the shared edges then produces differences of
Poisson-binomial mass functions, multiplied by the probabilities of those
edges.  These differences are not needed in the finite-dimensional theory,
where the edge noise only has to vanish at a fixed point, which follows from a
single concentration bound for the Poisson-binomial mass.  Tightness requires
a fourth-moment bound for the edge noise proportional to the $F$-measure of
each cell: the diagonal must be of the size of the cell, which requires
localization; the covariances involve first differences of mass functions
across the cell, which must be damped by $\mu_A$; and the fourth cumulant
involves differences of higher order, one for each shared edge in a
configuration.

Step~W therefore supplies two families of estimates.  Conditionally on all
latent positions, Lemmas~\ref{lem-poisson-binomial}--\ref{lem-c-loc} bound
differences of Poisson-binomial mass functions, pointwise and after
integration over a window; these bounds are used for the edge noise in
Step~B and for the conditional quantities of Step~K.  After averaging over the
latent positions of the other vertices, the multilinear representation
expresses everything in terms of binomial mass functions with the average
success probability; Lemma~\ref{lem-c-cancel} gives the binomial counterparts
of these bounds, together with the cancellation, and these are used for the
latent part in Steps~Z2 and~Z4.  The two families are used differently.  The
fourth cumulant involves at most four vertices, so the Poisson-binomial
estimates are needed only for differences of bounded order, and
Lemma~\ref{lem-c-roots} is stated for a fixed range of orders.  The Hoeffding
expansion in the latent positions has components of every order, one for each
set of vertices, so the binomial estimates of Lemma~\ref{lem-c-cancel} must
hold uniformly in the order, and the components are summed as geometric
series in $\rho_n$, with a truncation at high order.

We use two local estimates for sums of independent Bernoulli variables. The
first establishes that the $d$-th difference of a Poisson-binomial mass function
is uniformly of order $\Sigma^{-d/2}$, matching a binomial law with the same
variance.

\begin{lemma}[Poisson-binomial differences]
\label{lem-poisson-binomial}
Let $X = \sum_{j\in J'}X_j$ with independent $X_j\sim\mathrm{Bernoulli}(p_j)$, put
$\Sigma_{J'} := \sum_{j\in J'}p_j(1-p_j)$, and let $Y$ be independent of $X$.
Then for every $d\geq1$,
$\sup_m\bigl|\mathsf D^{d-1}P(X+Y=m)\bigr|\leq\Gamma_d\,\Sigma_{J'}^{-d/2}$, with
$\Gamma_d$ as in Lemma~\ref{lem-binomial-facts}\ref{binom-fact-c}.
\end{lemma}

The proof is the Fourier argument of
Lemma~\ref{lem-binomial-facts}\ref{binom-fact-c}, with the binomial
characteristic function replaced by a product of Bernoulli ones, in the
classical Fourier-analytic approach to local limit theorems
\citep{petrov1975sums}; see
Appendix~\ref{app-c-W}.

Near its mean, the $k$-th difference of a Poisson-binomial mass function is
of order $\Sigma^{-(k+1)/2}$ and decays like a Gaussian in the standardized
distance from the mean; the lemma below establishes a  local limit estimate that supports the rest of the proof.
\begin{lemma}[Localized Poisson-binomial differences]
\label{lem-pb-local}
Let $X = \sum_{j\in J'}X_j$ with independent $X_j\sim\mathrm{Bernoulli}(p_j)$,
$p_j\in[0,1]$, and put $\mu := EX$ and $\Sigma := \sum_{j\in J'}p_j(1-p_j)\geq1$.
For each $k\geq0$ there are constants $C_k,\zeta,c_0>0$, with $\zeta,c_0$
absolute and $C_k\leq C^{k+1}\Gamma_{k+1}$, such that for $|m-\mu|\leq c_0\Sigma$,
writing $x_m := (m-\mu)/\Sigma^{1/2}$,
\begin{align*}
\bigl|\mathsf D^{k}P(X = m)\bigr|
 \;\leq\; C_k\,\Sigma^{-(k+1)/2}\bigl(1+|x_m|\bigr)^{k}e^{-\zeta x_m^2}.
\end{align*}
\end{lemma}

The Fourier bound of Lemma~\ref{lem-poisson-binomial} is uniform in $m$ but
gives no decay away from the mean.  The proof therefore uses the classical
change of measure from large-deviation theory: the law is tilted
exponentially so that $m$ becomes its mean, where the uniform bound is sharp,
and the original mass equals the tilted one times $e^{-\Lambda^*(m)}$, where
$\Lambda^*$ is the Cram\'er rate function.  The Gaussian decay comes from
$\Lambda^*(m)\gtrsim x_m^2$, and the polynomial factor $(1+|x_m|)^k$ from the
small tilt parameter, which converts differences of the original law into
differences of the tilted one; see Appendix~\ref{app-c-W}.

Let $c_0,\zeta$ be the constants of Lemma~\ref{lem-pb-local},
$\zeta_* := \zeta/6$, and
\begin{align}
\label{eq-c-cw}
c_w := \min\Bigl\{\frac{c_0}{16},\ \frac12\Bigl(1-\frac{t_0^-}{t_0}\Bigr),\
 \frac14\Bigl(\frac{t_1^+}{t_1}-1\Bigr)\Bigr\}.
\end{align}
For a short cell $A = A_{rt}$ with endpoint set $S_A := \{r,t\}$, and
$z_s(x) := (\check{\mathfrak D}(x)-s)/h_s$, let
\begin{align}
\label{def-omega-e-a}
\omega^{[j]}_A(x) := \sum_{s\in S_A}\bigl(1+|z_s(x)|\bigr)^je^{-\zeta_*z_s(x)^2},
\quad
e_A(x) := \min_{s\in S_A}\bigl(|\mathfrak D(x)-s|+|s-t|\bigr),
\end{align}
and $\omega_A := \omega^{[1]}_A$. 
Both quantities are functions of a vertex's latent position $x$, through its
expected degree.  The weight $\omega^{[j]}_A(x)$ is a sum of Gaussian bumps,
one centred at each endpoint $s$ of the cell and of width $h_s$ in the
expected-degree scale.  It describes how strongly a vertex at position $x$
interacts with the cell: by the local limit estimates below, differences of
that vertex's conditional degree mass function at the endpoints are of order
$\sigma_t^{-(k+1)}\omega^{[k]}_A(x)$, so they are appreciable only when
$\check{\mathfrak D}(x)$ lies within a few bandwidths of an endpoint, and they
decay like a Gaussian in the standardized distance $z_s(x)$ beyond that.  The
polynomial factor $(1+|z_s|)^j$ accounts for $j$-th order differences, just as
the $j$-th derivative of a Gaussian density carries a polynomial of degree $j$.
Integrating such a weight over $x$ yields the measure of a window of width
$h_t$ around the cell (Lemma~\ref{lem-c-loc}).

The quantity $e_A(x)$ measures how far a vertex at $x$ is, in expected degree,
from the reference point $x_t$ of the cell, with the distance taken through the
nearest endpoint.  It controls how much the connection profile of that vertex
can differ from the profile at $x_t$: by \ref{condition-g2},
$\|a(x,\cdot)-a(x_t,\cdot)\|_{L^2}\leq Ce_A(x)^{\theta_w}$
(see \eqref{eq-c-a-window} below).  On the support of the weights,
$e_A\lesssim\bar\ell_A(1+|z_s|)$, so $e_A^\theta$ is of order
$\bar\ell_A^{\,\theta}$, and products $e_A^\theta\,\omega^{[j]}_A$ record the
H\"older error at the scale of the cell.

The \emph{windows} are
$W_s := \{x\notin\mathscr H:\ |\check{\mathfrak D}(x)-s|\leq c_ws\}$ and
$W_A := W_r\cup W_t$. Every kernel below depends on a vertex's latent position $x$ through the gap
between $\check{\mathfrak D}(x)$ and an endpoint $s$ of the cell.  The windows
$W_s$ mark the gaps of relative size at most $c_w$: inside them the local
bounds of Lemma~\ref{lem-pb-local} apply and $\mathfrak D(x)\in J$, so the
assumptions can be used, while outside them every kernel is negligible by
Chernoff bounds.  The finer localization, on the scale $h_s$ of the bandwidth,
is carried by the weights $\omega^{[j]}_A$.  

For $x\in W_A$ and $n\geq n_0$, Lemma~\ref{lem-c-hubs}
gives $\mathfrak D(x)\leq\check{\mathfrak D}(x)+g_*\leq(1+c_w)t_1+g_*$, so
$\mathfrak D(x)\in J$.  Off $\mathscr H$ we also have
$\mathfrak D\leq\frac43\check{\mathfrak D}$, and $0\leq\mathfrak D-\check{\mathfrak D}\leq g_*$
on $\{\mathfrak D\leq2t_1^+\}$.  Hence, by Remark~\ref{remark-inverse},
\ref{condition-d1} and \ref{condition-d2}, for $s\in I$, $0\leq u\leq c_ws$ and
$n\geq n_0$,
\begin{align}
\label{eq-c-window-measure}
\mathrm{Leb}\bigl(\{x\notin\mathscr H:\ |\check{\mathfrak D}(x)-s|\leq u\}\bigr)
 \leq\mathrm{Leb}\bigl(\{|\mathfrak D-s|\leq u+g_*\}\bigr)
 \leq2\bigl(\textstyle\sup_Jf\bigr)(u+g_*) .
\end{align}
Consequently, by \ref{condition-l1}, \ref{condition-g2}
(applied twice, at $\mathfrak D(x),s$ and at $s,t$) and \ref{condition-l2}, for
$x\in W_A$,
\begin{align}
\label{eq-c-a-window}
\|a(x,\cdot)\|_{L^2}+\|a(x,\cdot)\|_{L^4}\leq C,\qquad
\|a(x,\cdot)-a(x_t,\cdot)\|_{L^2}\leq C\,e_A(x)^{\theta_w},
\end{align}
using that $a(x,\cdot)$ is the centering of $\check w(x,\cdot)$ and that
$u\mapsto u\wedge\rho_n^{-1}$ is $1$-Lipschitz.  The point $x_t$ is not a hub,
by Lemma~\ref{lem-c-hubs}.

Most estimates below bound a quantity attached to a vertex by a power of
$\sigma_t^{-1}$ times a localization weight $\omega^{[j]}_A$ evaluated at that
vertex's latent position.  The next lemma converts such weights into the
measure of a window of width $h_t$ around the cell, which is how the scale of
the cell enters every bound.

\begin{lemma}[Localization integrals]
\label{lem-c-loc}
For fixed $J_0$ there are $C,n_0$ such that for $n\geq n_0$, every short cell
$A$, $0\leq j\leq J_0$ and $\theta\in[0,1]$:
\begin{align*}
\int_{[0,1]\setminus\mathscr H}e_A(x)^\theta\,\omega^{[j]}_A(x)\,dx
 \leq C\,h_t\,\bar\ell_A^{\,\theta},
\end{align*}
where $\omega^{[j]}_A(x)$ and $e_A(x)^\theta$ are defined in  \eqref{def-omega-e-a} and $\bar\ell_A^{\,\theta}$ is defined in \eqref{eq-c-cell-quantities}. The same holds with $\zeta_*$ replaced by any fixed $\zeta'\in(0,\zeta_*]$.
\end{lemma}

The integral is computed over dyadic shells around each endpoint: by
\eqref{eq-c-window-measure} the $k$-th shell has measure $O(2^kh_t)$, on which the
Gaussian factor in $\omega^{[j]}_A$ is $e^{-\zeta_*4^{k-1}}$, and outside the
windows the weights are negligible.  See Appendix~\ref{app-c-W}.

The basic local estimate concerns a vertex whose degree is counted.  Given the
latent positions, its degree is Poisson-binomial, and the next lemma shows
that, outside a negligible event, the differences of its mass function over a
cell are as small and as localized as those of a binomial law with the same
mean.

\begin{lemma}[Roots]
\label{lem-c-roots}
For $i\in[n]$ and $I'\ni i$ with $|I'|\leq4$, let
$\mu_{i,I'} := \sum_{l\notin I'}p_{il}$ and
$\Sigma_{i,I'} := \sum_{l\notin I'}p_{il}(1-p_{il})$.  Let $\mathscr E_i$ be the
event that for all such $I'$:
\begin{enumerate}[(E1)]
\setcounter{enumi}{-1}
\item $|\mu_{i,I'}-n''\rho_n\check{\mathfrak D}(\xi_i)|\leq\frac{c_w}4n''\rho_n
\max\{\check{\mathfrak D}(\xi_i),t_0\}$;
\item for every lattice point $s\in I$ with
$|\check{\mathfrak D}(\xi_i)-s|\leq c_ws$:
$|\mu_{i,I'}-n''\rho_n\check{\mathfrak D}(\xi_i)|\leq\sigma_s$ and
$\Sigma_{i,I'}\geq\mu_{i,I'}/4$.
\end{enumerate}
Then (a) $\sup_{x\notin\mathscr H}P(\mathscr E_i^c\mid\xi_i = x)$ is negligible;
and (b) on $\mathscr E_i\cap\{\xi_i\notin\mathscr H\}$, for every short cell $A$,
every $I'$ as above, $1\leq k\leq4$ and $|j|\leq3$, with
$P_{I'}(m) := P(D_{i\setminus I'} = m\mid\mathcal G)$,
\begin{align*}
\Bigl|\sum_{m\in A}\mathsf D^kP_{I'}(m+j)\Bigr|
 \leq C\,m_A\,\sigma_t^{-(k-1)}\,\omega^{[k]}_A(\xi_i)+\mathfrak e_n .
\end{align*}
\end{lemma}

Part (a) follows from Bernstein's inequality for the conditional mean and
variance of the degree, a Chernoff bound for the number of nearly
deterministic edges, and a union bound over the polynomially many $I'$ and
lattice points.  For (b), on the windows Lemma~\ref{lem-pb-local} applies with
variance $\asymp\sigma_s^2$, and off them every mass is negligible; the bound is
then summed over the $k_A$ points of a cell in the small regime, and
telescoped to the two endpoints in the large regime.  See
Appendix~\ref{app-c-W}.

Every quantity whose moments are bounded below is a polynomial with at most
$n^4$ terms in the $p_{ij}$, conditional mass functions and values
$a(\xi_i,\xi_j)$, all bounded by $n^C$.  The event
$\bigcup_i(\{\xi_i\notin\mathscr H\}\cap\mathscr E_i^c)$ therefore contributes
negligibly to every moment, and we suppress it.

When the positions of the neighbors are averaged out, the conditional mass
function of a degree becomes binomial, and the resulting kernels integrate
almost to zero over the latent space, because $\int_0^1b_M(m;p)\,dp$ does not
depend on $m$ (equivalently, a binomial count whose success probability is
drawn uniformly from $[0,1]$ is uniformly distributed, so mixing removes the
lattice).  This cancellation, up to errors controlled by the regularity
of $f$ and by clipping, is what makes increments over short cells small.

For $d\geq0$ and $|M-n''|\leq8$ define:
\begin{align}
\begin{split}
\label{eq-c-K}
 \mathsf K^{(d)}_{A,M}(p) &:= (-1)^d[\mathsf D^db_M(r';p)-\mathsf D^db_M(t';p)],
\\ \mathsf K^{(d)}_A(x) &:= \mathsf K^{(d)}_{A,M}(\rho_n\check{\mathfrak D}(x)).
\end{split}
\end{align}
With $d = 0$ and $M = n''$,
$n'\rho_n\mathsf K^{(0)}_A(x) = n'\rho_n[b_{n''}(r';\rho_n\check{\mathfrak D}(x))
-b_{n''}(t';\rho_n\check{\mathfrak D}(x))]$ is the density of the signed measure
that defines the H\'ajek kernel $\phi_A$ in Step~K.
The cases $d\geq1$ control the higher-order Hoeffding terms in Steps~Z2
and~Z4: by Lemma~\ref{lem-multilinear}, the Hoeffding component of a
conditional mass difference indexed by a set $S$ of latent positions carries
the difference $\mathsf K^{(|S|)}$ of order $|S|$.  Each additional order
costs a factor $\sigma_t^{-1}$ in Lemma~\ref{lem-c-cancel}, and together with
the factor $\rho_n$ per additional latent position this makes the sums over
the order geometric in $\rho_n$.

The binomial estimates of Section~\ref{sec-degenerate} do not suffice here,
for three reasons.  They bound differences at a single lattice point, whereas
the increments involve differences between the two endpoints of a cell, which
are smaller by the factor $\mu_A$ when the cell lies within one fluctuation
window.  They control the higher-order components only after integration over
fixed regions, whereas the fourth-moment bounds require pointwise bounds that
retain the Gaussian localization, so that they can be multiplied and then
integrated.  And they establish convergence of the H\'ajek kernels, whereas
tightness requires the cancellation to hold at a rate, uniformly over cells.
Lemma~\ref{lem-c-cancel} provides these refinements.

\begin{lemma}[Cancellation]
\label{lem-c-cancel}
Fix $d_1$.  For $n\geq n_0$, $0\leq d\leq d_1$, $|M-n''|\leq8$, every short cell
$A$ and $\theta\in[0,1]$:
\begin{enumerate}[(a)]
\item $|\mathsf K^{(d)}_A(x)|\leq C\sigma_t^{-(d+1)}\mu_A\,\omega^{[d+1]}_A(x)+\mathfrak e_n$
for all $x$, and the left side is negligible for $x\notin W_A\cup\mathscr H$;
\item $\int_{[0,1]\setminus\mathscr H}e_A^\theta|\mathsf K^{(d)}_A|^j\leq
C\sigma_t^{-j(d+1)}\mu_A^jh_t\bar\ell_A^\theta$ for $j\in\{1,2\}$;
\item $n'\rho_n\bigl|\int_{[0,1]\setminus\mathscr H}\mathsf K^{(d)}_A(x)\,dx\bigr|
\leq C\sigma_t^{-d}\mu_A\bigl(\bar\ell_A^{\,\theta_f}+\Upsilon_n\bigr)+\mathfrak e_n
\leq C\sigma_t^{-d}\mu_A\bar\ell_A^{\,\vartheta}+\mathfrak e_n$.
\end{enumerate}
\end{lemma}

Parts (a) and (b) are the binomial analogues of Lemma~\ref{lem-c-roots}(b)
and Lemma~\ref{lem-c-loc}.  Part (c) is the heart of the matter.  After
replacing $\check{\mathfrak D}$ by $\mathfrak D$ at a cost $\Upsilon_n$, and pushing
Lebesgue measure forward by $\mathfrak D$, the integral becomes
\begin{align}
\label{eq-c-cancel-idea}
\int_0^{2t}\mathsf K^{(d)}_{A,M}(\rho_n\lambda)\,f(\lambda)\,d\lambda
 = f(t)\int_0^{2t}\mathsf K^{(d)}_{A,M}(\rho_n\lambda)\,d\lambda
 +\int_0^{2t}\mathsf K^{(d)}_{A,M}(\rho_n\lambda)\bigl[f(\lambda)-f(t)\bigr]d\lambda .
\end{align}
The first integral vanishes up to a negligible term, by the identity
$\int_0^1b_M(m;p)\,dp = 1/(M+1)$ of
Lemma~\ref{lem-binomial-facts}\ref{binom-fact-e}; Remark~\ref{rem-cancel-signs}
explains why, for every order $d$.  The second is bounded by the H\"older
modulus of $f$ across the window, which gives $\bar\ell_A^{\,\theta_f}$.  The complete proof is in Appendix~\ref{app-c-W}.

\begin{remark}[The cancellation]
\label{rem-cancel-signs}
For $d=0$ the cancellation has a direct interpretation.  Fix the latent
positions and add the edge $\{i,j\}$.  Since $D_i = D_{i\setminus j}+A_{ij}$, and
shifting by one moves the cell $\{r'+1,\ldots,t'\}$ to $\{r',\ldots,t'-1\}$, the
conditional probability that $D_i$ lies in $A$ changes by
\begin{align*}
P(D_{i\setminus j}=r'\mid\mathcal G)-P(D_{i\setminus j}=t'\mid\mathcal G).
\end{align*}
This is positive when the degree of $i$ tends to lie below the cell, so that
the edge pushes it into the cell, and negative when it tends to lie above the
cell, so that the edge pushes it out.  Averaging over the positions of the
other vertices replaces it, for a vertex at position $x$, by
\begin{align*}
b_{n''}\bigl(r';\rho_n\check{\mathfrak D}(x)\bigr)-b_{n''}\bigl(t';\rho_n\check{\mathfrak D}(x)\bigr)
 = \mathsf K^{(0)}_A(x),
\end{align*}
and the H\'ajek kernel $\phi_A$ of Lemma~\ref{lem-c-kernels}, the aggregate
effect of one vertex on the degrees of all others, integrates these signed
contributions over $x$ against the connection profile.  After the profile is
evaluated at the reference point $x_t$ (Proposition~\ref{prop-c-hajek}), what
remains is, up to clipping,
\begin{align*}
&n'\rho_n\int_0^1\Bigl[b_{n''}\bigl(r';\rho_n\check{\mathfrak D}(x)\bigr)-b_{n''}\bigl(t';\rho_n\check{\mathfrak D}(x)\bigr)\Bigr]dx
\\&\qquad = n'\rho_n\int\Bigl[b_{n''}(r';\rho_n\lambda)-b_{n''}(t';\rho_n\lambda)\Bigr]f(\lambda)\,d\lambda .
\end{align*}
If $f$ were equal to a constant $c$ near $t$, the substitution $p = \rho_n\lambda$
would give
\begin{align*}
n'c\int_0^1\Bigl[b_{n''}(r';p)-b_{n''}(t';p)\Bigr]dp
 = n'c\Bigl[\frac1{n''+1}-\frac1{n''+1}\Bigr] = 0 :
\end{align*}
the vertices pushed into the cell from below exactly balance those pushed
out of it above.  In general the integral equals
$n'\rho_n\int[b_{n''}(r';\rho_n\lambda)-b_{n''}(t';\rho_n\lambda)][f(\lambda)-f(t)]\,d\lambda$,
which is small by the H\"older continuity of $f$ across the window.  Without
the cancellation, the integral could only be bounded by
$n'\rho_n\int|b_{n''}(r';\rho_n\lambda)-b_{n''}(t';\rho_n\lambda)|f(\lambda)\,d\lambda$,
which is of order $\mu_A$; the cancellation gains the factor
$\bar\ell_A^{\,\theta_f}$.

For $d\geq1$ the same identity applies term by term.  Since
\begin{align*}
\mathsf D^db_M(m;p) = \sum_{j=0}^d(-1)^j\binom dj b_M(m-j;p),
\end{align*}
integrating in $p$ gives
\begin{align*}
\int_0^1\mathsf D^db_M(m;p)\,dp = \frac1{M+1}\sum_{j=0}^d(-1)^j\binom dj = 0 .
\end{align*}
In words, mixing over a flat density makes the masses constant in $m$, and
every finite difference of a constant vanishes.
\end{remark}

\paragraph{Step K: the H\'ajek increment.}
We now introduce analogs of the H\'ajek projections defined for the half-intervals in Section \ref{sec-multilinear-hoeffding} to short cells below.  
\begin{lemma}[H\'ajek kernels]
\label{lem-c-kernels}
For a short cell $A$, the H\'ajek projection of $\Delta'_A$ is
$n^{-1/2}\sum_j[\psi_A(\xi_j)+\phi_A(\xi_j)]$, with
\begin{align}
\label{eq-c-hajek}
\psi_A(\xi) := \chi(\xi)\,\mathcal P_A(\xi)-p'_A,\qquad
\phi_A(\xi) := n'\rho_n\int_{[0,1]\setminus\mathscr H}a(x,\xi)\,
 \mathsf K^{(0)}_A(x)\,dx ,
\end{align}
where $\mathcal P_A(\xi) := P(D_1\in A\mid\xi_1 = \xi)$.  Moreover
$E\psi_A = E\phi_A = 0$.
\end{lemma}

\begin{proof}
By \eqref{eq-hajek-proj} it suffices to compute $E[\Delta'_A\mid\xi_j]$.
Since $E[\chi_i\mathbbm 1(D_i\in A)\mid\xi_j]$ equals
$E[\chi_iq_i(A)\mid\xi_j]$,
\begin{align*}
\sqrt n\,E[\Delta'_A\mid\xi_j] = \bigl(E[\chi_jq_j(A)\mid\xi_j]-p'_A\bigr)
 +\sum_{i\neq j}\bigl(E[\chi_iq_i(A)\mid\xi_j]-p'_A\bigr).
\end{align*}
The first term is $\psi_A(\xi_j)$.  For $i\neq j$, since
$E[\chi_i\mathcal P_A(\xi_i)] = p'_A$,
\begin{align*}
E[\chi_iq_i(A)\mid\xi_j]-p'_A = E\bigl[\chi(\xi_i)\bigl(q_i(A)-\mathcal P_A(\xi_i)\bigr)\bigm|\xi_j\bigr].
\end{align*}
By Lemma~\ref{lem-multilinear}, applied to $F_i(t'\mid\mathcal G)$ and
$F_i(r'\mid\mathcal G)$, $q_i(A)-\mathcal P_A(\xi_i)$ is the sum over
$S\neq\emptyset$ of the differences of the components $K_S$ at $t'$ and $r'$.
Each is centered in every $\xi_l$, $l\in S$, given the other positions.  Hence
only $S = \{j\}$ survives conditioning on $\xi_j$, and since
$\mathsf D\mathcal B_{n''} = b_{n''}$, its contribution is
\begin{align*}
-\rho_na(\xi_i,\xi_j)\bigl[b_{n''}(t';\rho_n\check{\mathfrak D}(\xi_i))
-b_{n''}(r';\rho_n\check{\mathfrak D}(\xi_i))\bigr] = \rho_na(\xi_i,\xi_j)\mathsf K^{(0)}_A(\xi_i).
\end{align*}
Integrating $\xi_i$ against $\chi$ and summing over the $n'$ indices $i\neq j$
gives $\phi_A(\xi_j)$.  Finally $E\psi_A = 0$ by definition of $p'_A$, and
$E\phi_A = 0$ because $E_\xi[a(x,\xi)] = 0$.
\end{proof}

\begin{proposition}[H\'ajek increment, compact case]
\label{prop-c-hajek}
For $n\geq n_0$ and every short cell $A$:
$E\psi_A^2\leq p'_A\leq C\Delta_A$ and
$\|\phi_A\|_{L^2}\leq C\,\Xi_A^{1/2}+\mathfrak e_n\leq C\Delta_A^{\vartheta}$.
\end{proposition}

\begin{proof}
First, $E\psi_A^2\leq E[\chi P(D_1\in A\mid\xi_1)]=p'_A$.  By
Lemma~\ref{lem-binomial-facts}\ref{binom-fact-a} on the windows and Chernoff
bounds off them, for each $m\in A$,
\begin{align*}
E\bigl[\chi\, b_{n'}(m;\rho_n\check{\mathfrak D}(\xi))\bigr]\leq
C\sigma_t^{-1}h_t\leq C(n\rho_n)^{-1},
\qquad\text{so}\qquad
p'_A\leq\frac{Ck_A}{n\rho_n}\leq C\ell_A\leq C\Delta_A .
\end{align*}
For $\phi_A$, write
\begin{align*}
\phi_A(\xi) &= n'\rho_n\int_{[0,1]\setminus\mathscr H}\bigl[a(x,\xi)-a(x_t,\xi)\bigr]
 \mathsf K^{(0)}_A(x)\,dx+a(x_t,\xi)\,M^\circ_A,
\\ M^\circ_A &:= n'\rho_n\int_{[0,1]\setminus\mathscr H}\mathsf K^{(0)}_A(x)\,dx .
\end{align*}
By Minkowski's inequality, \eqref{eq-c-a-window} and
Lemma~\ref{lem-c-cancel}(a),(b) with $\theta = \theta_w$ (the part off $W_A$ is
negligible), the first term has $L^2$ norm at most
\begin{align*}
Cn\rho_n\cdot\sigma_t^{-1}\mu_Ah_t\bar\ell_A^{\theta_w}\leq C\mu_A\bar\ell_A^\vartheta .
\end{align*}
By Lemma~\ref{lem-c-cancel}(c) with $d = 0$, the second is at most
$C\mu_A\bar\ell_A^\vartheta+\mathfrak e_n$.  Thus
$\|\phi_A\|\leq C\Xi_A^{1/2}+\mathfrak e_n$, and \eqref{eq-c-Q-Xi} completes the
proof.
\end{proof}

The same argument covers both regimes: in the small regime $\mu_A = \ell_A/h_t$
records the first-difference cancellation between the two windows, and in the
large regime $\bar\ell_A = \ell_A$ records the H\"older moduli of
$\lambda\mapsto w(x_\lambda,\cdot)$ and of $f$ across the cell.

\paragraph{Step B: edge noise.}
Recall from \eqref{eq-c-BZ} that $B_A$, the edge noise, is a centered
function of the edge variables given the latent positions.
Step~B is the counterpart, at the scale of a cell, of
Lemma~\ref{lemma-sigma-n-0}, with two differences.  The edge noise is no
longer negligible: for a cell no wider than the bandwidth, whether a degree
falls in the cell remains random given the latent positions, so the
conditional variance of $B_A$ is of the same order as $\Delta_A$, and what is
needed is a bound of the correct order rather than negligibility.  And the
criterion of Step~L requires fourth moments, so in addition to the signed
pairing, which bounds the conditional covariances, we need the fourth
conditional cumulant of $B_A$.
For $i\neq j$ let $g^{(j)}_i(A) := P(D_{i\setminus j} = r'\mid\mathcal G)-
P(D_{i\setminus j} = t'\mid\mathcal G)$.

\begin{lemma}[Signed pairing]
\label{lem-c-pairing}
For $i\neq j$,
\begin{align}
\label{eq-c-pairing}
\mathrm{Cov}\bigl(\mathbbm 1(D_i\in A),\mathbbm 1(D_j\in A)\bigm|\mathcal G\bigr)
 = p_{ij}(1-p_{ij})\,g^{(j)}_i(A)\,g^{(i)}_j(A),
\end{align}
and, for $\xi_i\notin\mathscr H$, $|g^{(j)}_i(A)|\leq Cm_A\omega_A(\xi_i)+\mathfrak e_n$.
\end{lemma}

\begin{proof}
Given $\mathcal G$, $D_i = D_{i\setminus j}+A_{ij}$ and
$D_j = D_{j\setminus i}+A_{ij}$, with the three summands conditionally
independent.  Conditioning further on $A_{ij}$ makes the two indicators
independent, so
\begin{align*}
\mathrm{Cov}\bigl(\mathbbm 1(D_i\in A),\mathbbm 1(D_j\in A)\bigm|\mathcal G\bigr)
 = \mathrm{Cov}\bigl(u(A_{ij}),v(A_{ij})\bigm|\mathcal G\bigr),
\end{align*}
where $u(\epsilon) := P(D_{i\setminus j}+\epsilon\in A\mid\mathcal G)$ and $v$ is
defined likewise.  For a Bernoulli$(p)$ variable,
$\mathrm{Cov}(u(A),v(A)) = p(1-p)(u(1)-u(0))(v(1)-v(0))$, and
\begin{align*}
u(1)-u(0) = P(D_{i\setminus j} = r'\mid\mathcal G)-P(D_{i\setminus j} = t'\mid\mathcal G) = g^{(j)}_i(A).
\end{align*}
The bound is Lemma~\ref{lem-c-roots}(b) with $k = 1$ and $I' = \{i,j\}$, since
$g^{(j)}_i(A) = -\sum_{m\in A}\mathsf DP(D_{i\setminus j} = m\mid\mathcal G)$.
\end{proof}

Let
\begin{align}
\label{eq-c-Sigma-bar}
\bar\Sigma(A) := \frac1n\sum_{i\neq j}\chi_i\chi_j\,p_{ij}\,
 \bigl|g^{(j)}_i(A)\bigr|\bigl|g^{(i)}_j(A)\bigr| .
\end{align}

\begin{lemma}[Edge covariance]
\label{lem-c-sigma}
$E\bar\Sigma(A)\leq CQ_A\leq C\Delta_A$ and $\|\bar\Sigma(A)\|_{L^2}\leq C\Delta_A$.
\end{lemma}

The first moment follows by bounding one factor $|g|$ by $Cm_A$ and
integrating the other against the edge weight with Lemma~\ref{lem-c-loc}.
The second moment follows from the Efron--Stein inequality
\citep{efron1981jackknife}, separating the
direct effect of resampling $\xi_m$ on the terms rooted at $m$ from its
indirect effect through the edge probabilities $p_{im}$.  See
Appendix~\ref{app-c-B}.

Fix a set $V$ of $k\leq4$ distinct vertices.  For $v\in V$ write
$D_v = D^{\rm ext}_v+s_v$ with $s_v := \sum_{u\in V\setminus\{v\}}A_{uv}\in\{0,\ldots,3\}$,
so that $D^{\rm ext}_v = D_{v\setminus(V\setminus\{v\})}$.  Let $E_v$ be the set
of edges of $V$ at $v$, $\mathcal I := \sigma(A_e:e\subseteq V)$,
$\pi_v(s) := P(D^{\rm ext}_v+s\in A\mid\mathcal G)$ and $Y_e := A_e-p_e$.  Given
$\mathcal G$, the $D^{\rm ext}_v$, $v\in V$, are independent of each other and of
$\mathcal I$.

We write $\kappa(X_1,\ldots,X_k\mid\mathcal G)$ for the joint cumulant of
$X_1,\ldots,X_k$ conditional on $\mathcal G$, and $\kappa_4(X\mid\mathcal G)$
for the fourth conditional cumulant of $X$; for joint cumulants of products,
see \citet{leonov1959method}.  We work with
cumulants rather than moments because the fourth cumulant isolates the part
of the fourth moment not determined by the conditional variance, which the
signed pairing already controls, and because a joint cumulant vanishes
whenever its arguments split into independent groups, so that only connected
configurations of vertices contribute.  The resulting bound is of Rosenthal
type: the fourth moment splits into a Gaussian part controlled by the
conditional variance and a non-Gaussian part, the fourth cumulant, which
under independence would reduce to a sum over single vertices and here is a
sum over connected configurations.  Bounding cumulants by sums over connected
configurations, controlled by spanning trees, is in the spirit of cumulant
bounds for sums with dependency graphs
\citep{janson1988normal,feray2018weighted} and of moment computations for
random graphs \citep{janson2000random}; the new ingredient here is that each
summand is a threshold function of many edges, which we handle by
conditioning on the internal edges of a configuration and bounding the
resulting block coefficients by local limit estimates.

To bound the fourth cumulant of $B_A$, we expand it over quadruples of
vertices.  The indicators attached to distinct vertices are dependent only
through the edges among those vertices, and conditionally on $\mathcal G$ and
$\mathcal I$ the indicator of $v$ depends on its internal edges only through
their number:
\begin{align*}
P\bigl(D_v\in A\bigm|\mathcal G,\mathcal I\bigr) = \pi_v(s_v),
\qquad s_v = \sum_{e\in E_v}A_e .
\end{align*}
Since each $A_e$ takes the values $0$ and $1$, this is multilinear in the
internal edges, and expanding it in the centered edges gives
\begin{align*}
\pi_v(s_v) = \sum_{T\subseteq E_v}c_T\prod_{e\in T}Y_e,
\qquad
c_T = E\Bigl[\nabla^{|T|}\pi_v\Bigl(\sum_{e\in E_v\setminus T}A_e\Bigr)\Bigm|\mathcal G\Bigr],
\end{align*}
where $\nabla\pi(s) := \pi(s+1)-\pi(s)$ is the forward difference.  Since
$\pi_v(s) = \sum_{m\in A}P(D^{\rm ext}_v = m-s\mid\mathcal G)$, each
difference of $\pi_v$ is a difference of the Poisson-binomial mass function of
$D^{\rm ext}_v$ across the cell, so $c_T$ is a difference of order $|T|$, of
the kind controlled by Lemma~\ref{lem-c-roots}.  Each additional edge thus
costs a factor $\sigma_t^{-1}$ and keeps the coefficient localized at the
cell, which is what the next lemma records.

\begin{lemma}[Block coefficients]
\label{lem-c-block}
Let $v\in V$ with $\xi_v\notin\mathscr H$, work on $\mathscr E_v$, and let $P$ be
a polynomial of degree at most $4$.  Then
$P(\pi_v(s_v)) = \sum_{T\subseteq E_v}c_T\prod_{e\in T}Y_e$, where
\begin{align*}
|c_T|\leq C\,m_A\,\sigma_t^{-(|T|-1)}\,\omega^{[|T|]}_A(\xi_v)+\mathfrak e_n
\qquad(T\neq\emptyset),
\end{align*}
with $C$ depending only on the coefficients of $P$; and $c_\emptyset = 0$ if
$P(\pi) = \pi-q_v(A)$.
\end{lemma}

The coefficient $c_T$ is a $|T|$-fold forward difference of $P\circ\pi_v$ in
the number $s_v$ of internal edges, and the differences of $\pi_v$ are
exactly the cell differences bounded in Lemma~\ref{lem-c-roots}(b); a
general polynomial $P$ is handled by the discrete product rule.  See
Appendix~\ref{app-c-B}.

In the expansion of the cumulant, a configuration of $k$ distinct vertices
contributes only if the corresponding edges connect all $k$ of them, since otherwise
the vertices split into conditionally independent groups, and only if every
internal edge that appears does so in the expansions of both of its
endpoints, since a centered edge appearing once averages to zero.  With blocks
$T_v\subseteq E_v$ at the vertices $v$ of the configuration, this forces
\begin{align*}
\sum_{v}|T_v|\;\geq\;2(k-1),
\end{align*}
so that configurations with more vertices carry differences of higher total
order.  Such a configuration contains a spanning tree.  The $n^k$ choices of vertices are offset by a
factor $\rho_n$ for each of the $k-1$ tree edges, and the coefficients of
Lemma~\ref{lem-c-block} get smaller with each additional edge, so that every
additional vertex costs a factor $\mu_A\leq1$.  Making this count rigorous
requires summing over the latent positions of vertices that are coupled
through the edge weights, while retaining the localization at every vertex,
and this is what the next lemma does.

\begin{lemma}[Tree sums]
\label{lem-c-tree}
For a tree $F$ on $\{1,\ldots,k\}$, $k\leq4$, and $d_a\in\{1,\ldots,9\}$,
\begin{align*}
\sum_{v_1,\ldots,v_k\text{ distinct}}E\Bigl[\prod_{\{a,b\}\in F}\rho_n\check w(\xi_{v_a},\xi_{v_b})
 \prod_{a=1}^k\chi(\xi_{v_a})\,\omega^{[d_a]}_A(\xi_{v_a})\Bigr]\leq Cn^k\rho_n^{k-1}h_t .
\end{align*}
\end{lemma}

The positions are integrated from the leaves of the tree inwards, using only
$\int_0^1\check w(x,y)\,dy = \check{\mathfrak D}(x)$ and absorbing each factor
$\check{\mathfrak D}$ into the localization weight at the parent; the root is
integrated last with Lemma~\ref{lem-c-loc}.  See Appendix~\ref{app-c-B}.

\begin{lemma}[Edge-noise moments]
\label{lem-c-edge}
For every short cell $A$ and $n\geq n_0$,
\begin{align*}
E[B_A^2\mid\mathcal G]&\leq\bar p_A+\bar\Sigma(A),\qquad
E[\kappa_4(B_A\mid\mathcal G)]\leq Cn^{-1}\Delta_A+\mathfrak e_n,\qquad
\\E[B_A^4] &\leq C\bigl(\Delta_A^2+n^{-1}\|Z_A\|_{L^2}^2+n^{-1}\Delta_A\bigr)+\mathfrak e_n .
\end{align*}
\end{lemma}

\begin{proof}
\emph{Second moment.}  By \eqref{eq-c-pairing},
\begin{align*}
E[B_A^2\mid\mathcal G] = \frac1n\sum_i\chi_iq_i(1-q_i)+\frac1n\sum_{i\neq j}
\chi_i\chi_jp_{ij}(1-p_{ij})g^{(j)}_ig^{(i)}_j\leq\bar p_A+\bar\Sigma(A).
\end{align*}

\emph{Reduction of the cumulant.}  With
$a_v := \chi_v[\mathbbm 1(D_v\in A)-q_v(A)]$, multilinearity of joint
cumulants gives
\begin{align*}
\kappa_4(B_A\mid\mathcal G) = \frac1{n^2}\sum_{v_1,\ldots,v_4}
\kappa(a_{v_1},\ldots,a_{v_4}\mid\mathcal G).
\end{align*}
Fix a tuple, let $V$ be its set of distinct vertices and $k := |V|$.  The tuple
contributes only if $\xi_v\notin\mathscr H$ for all $v\in V$.  Given
$(\mathcal G,\mathcal I)$, the $a_v$ at distinct vertices are functions of
disjoint sets of external edges, hence independent.  By the law of total
cumulance \citep{brillinger1969calculation}, the cumulant is the sum over
partitions $\pi$ of the four slots of
\begin{align*}
\kappa\bigl(\kappa(a_\beta\mid\mathcal G,\mathcal I):\beta\in\pi\bigm|\mathcal G\bigr).
\end{align*}
The inner cumulant vanishes unless all slots of $\beta$ belong to one vertex
$v$, in which case it equals $P_{|\beta|}(\pi_v(s_v))$.  Here
$P_1(\pi) = \pi-q_v(A)$, and $P_2,P_3,P_4$ are the cumulants of a
Bernoulli$(\pi)$ variable, polynomials with bounded coefficients.  Let
$b := |\pi|\geq k$.

\emph{One vertex.}  The only contribution is $P_4(q_v)$, and $|P_4(q)|\leq q$;
the $n$ tuples contribute at most
\begin{align*}
n^{-1}E[\chi_1q_1(A)] = n^{-1}p'_A\leq Cn^{-1}\Delta_A,
\end{align*}
by Proposition~\ref{prop-c-hajek}.

\emph{Several vertices.}  Let $k\geq2$.  Expanding each $P_{|\beta|}(\pi_v(s_v))$
by Lemma~\ref{lem-c-block}, the outer cumulant becomes
\begin{align*}
\sum_{(T_\beta)}\prod_\beta c_{T_\beta}\,\kappa\Bigl(\prod_{e\in T_\beta}Y_e:\beta\in\pi\Bigm|\mathcal G\Bigr).
\end{align*}
A term vanishes in each of three cases:
\begin{enumerate}[(V1)]
\item some $T_\beta = \emptyset$, since a joint cumulant of order at least two
with a constant argument is zero;
\item some edge lies in exactly one $T_\beta$, since every moment in the
moment--cumulant formula containing that argument factorizes through
$EY_e = 0$;
\item the blocks split into two groups with disjoint edge sets, since the groups
are then independent.
\end{enumerate}
For a surviving term, $|\kappa(\cdots)|\leq C\prod_{e\in U}p_e$ with
$U := \bigcup_\beta T_\beta$, because $|EY_e^j|\leq p_e$ for $j\geq2$.  By (V1)
and (V3), $U$ is connected and spans $V$, so it contains a spanning tree $F$;
by (V2), $\sum_\beta|T_\beta|\geq2|U|\geq2(k-1)$.  By Lemma~\ref{lem-c-block},
\begin{align*}
\prod_\beta|c_{T_\beta}|\leq Cm_A^b\sigma_t^{\,b-\sum|T_\beta|}
 \prod_{v\in V}\omega^{[d_v]}_A(\xi_v)\leq Cm_A^k\sigma_t^{\,2-k}
 \prod_{v\in V}\omega^{[d_v]}_A(\xi_v),
\end{align*}
where $d_v := \sum_{\beta\text{ at }v}|T_\beta|\in[1,9]$, using
$(m_A\sigma_t)^b\leq(m_A\sigma_t)^k$.  Also
$\prod_{e\in U}p_e\leq\prod_{e\in F}\rho_n\check w_e$.  The numbers of tuple
patterns, partitions, families $(T_\beta)$ and trees are bounded, so by
Lemma~\ref{lem-c-tree} the tuples with $|V| = k$ contribute at most
\begin{align*}
Cn^{-2}m_A^k\sigma_t^{\,2-k}\,n^k\rho_n^{k-1}h_t
 = C\Bigl(\frac n{n''}\Bigr)^{k-1}\frac{\mu_A^kh_t}n\leq C\frac{Q_A}n\leq Cn^{-1}\Delta_A ,
\end{align*}
using $\sigma_t^2 = n''\rho_nt$, $t\leq t_1$ and $\mu_A\leq1$.

\emph{Fourth moment.}  Since $B_A$ is centered given $\mathcal G$,
\begin{align*}
E[B_A^4] = E\bigl[3E(B_A^2\mid\mathcal G)^2+\kappa_4(B_A\mid\mathcal G)\bigr].
\end{align*}
Also $E(B_A^2\mid\mathcal G)^2\leq2\bar p_A^2+2\bar\Sigma(A)^2$,
$E\bar p_A^2 = (p'_A)^2+n^{-1}\|Z_A\|_2^2$, $p'_A\leq C\Delta_A$, and
$E\bar\Sigma(A)^2\leq C\Delta_A^2$ by Lemma~\ref{lem-c-sigma}.
\end{proof}

\paragraph{Step Z2: latent part, second moment.}
Recall from \eqref{eq-c-BZ} that $Z_A = \sqrt n(\bar p_A-p'_A)$ is the latent
part of the increment, a function of the latent positions alone.
Let $\xi'_m$ be an independent copy of $\xi_m$, and $q^{(m)}_i$, $Z^{(m)}_A$ the
values with $\xi_m$ replaced.  Since $q_i(A)$ is affine in $p_{im}$, and
$g^{(m)}_i(A)$ does not involve $\xi_m$, for $i\neq m$
\begin{align}
\label{eq-c-exact}
q_i(A)-q^{(m)}_i(A) = \rho_n\bigl[a(\xi_i,\xi_m)-a(\xi_i,\xi'_m)\bigr]g^{(m)}_i(A),
\end{align}
so that
\begin{align}
\label{eq-c-Z-diff}
\sqrt n\bigl(Z_A-Z^{(m)}_A\bigr) = \bigl[\chi_mq_m-\chi'_mq^{(m)}_m\bigr]
+\rho_n\sum_{i\neq m}\chi_i\bigl[a(\xi_i,\xi_m)-a(\xi_i,\xi'_m)\bigr]g^{(m)}_i .
\end{align}
Let
\begin{align}
\label{eq-c-Phi-Psi}
\Phi_m := \rho_n\sum_{i\neq m}\chi_ia(\xi_i,\xi_m)g^{(m)}_i(A),\qquad
\Psi_m := E_\xi\Bigl[\Bigl(\rho_n\sum_{i\neq m}\chi_ia(\xi_i,\xi)g^{(m)}_i(A)\Bigr)^2\Bigr],
\end{align}
and $\bar g(x) := E[g^{(1)}_2\mid\xi_2 = x] = \mathsf K^{(0)}_A(x)$, with
$M = n''$.

\begin{lemma}[Latent second moment]
\label{lem-c-Z2}
$\|Z_A\|_{L^2}^2\leq2p'_A+2E\Phi_1^2$ and
\begin{align}
\label{eq-c-Z2}
E\Phi_1^2\leq C\bigl[\Xi_A+\rho_nQ_A+\rho_n\Xi_A+v_n\Delta_A\bigr]+\mathfrak e_n,
\qquad v_n := n^{(1-2\vartheta)/4}\rho_n^{(9-2\vartheta)/4}\to0 .
\end{align}
In particular $\|Z_A\|_{L^2}^2\leq C(\Delta_A+\Xi_A)$.  The same bound holds for
$X^{(-m)}$ of Step~Z4, with three vertices removed; for its centered part,
\begin{align*}
E\|X^{(-m)}-EX^{(-m)}\|^2\leq C[\rho_nQ_A+\rho_n\Xi_A+v_n\Delta_A]+\mathfrak e_n .
\end{align*}
\end{lemma}

\begin{proof}[Outline of the proof]
By the Efron--Stein inequality, exchangeability and \eqref{eq-c-Z-diff},
\begin{align}
\label{eq-c-Z2-reduction}
\|Z_A\|_{L^2}^2 = \mathrm{Var}(Z_A)\leq E\bigl[(\chi_1q_1-\chi'_1q^{(1)}_1)^2\bigr]+2E\Phi_1^2
\leq2p'_A+2E\Phi_1^2 .
\end{align}
To bound $E\Phi_1^2$, replace each $g^{(1)}_i$ by its conditional mean
$\bar g(\xi_i)$ given $\xi_i$:
\begin{align*}
\Phi_1 = \underbrace{\rho_n\sum_{i\neq1}\chi_ia(\xi_i,\xi_1)\bar g(\xi_i)}_{\Phi^{\rm m}}
 +\underbrace{\rho_n\sum_{i\neq1}\chi_ia(\xi_i,\xi_1)\bigl[g^{(1)}_i-\bar g(\xi_i)\bigr]}_{\Phi^{\rm f}} .
\end{align*}
Given $\xi_1$, the \emph{main part} $\Phi^{\rm m}$ is a sum of i.i.d.\ terms
whose mean is essentially the H\'ajek kernel $\phi_A(\xi_1)$, so by
Proposition~\ref{prop-c-hajek},
\begin{align*}
E(\Phi^{\rm m})^2\leq\|\phi_A\|^2+n''\rho_n^2E\bigl[\chi\|a(\xi,\cdot)\|^2\bar g(\xi)^2\bigr]
 \leq C\Xi_A+C\rho_nQ_A .
\end{align*}
For the \emph{fluctuation part} $\Phi^{\rm f}$, expand each
$g^{(1)}_i-\bar g(\xi_i)$ into its Hoeffding components by
Lemma~\ref{lem-multilinear}.  After truncation of the high-order components,
which costs $o(\Delta_A)$, a pair of components contributes to $E(\Phi^{\rm f})^2$
only in one of three classes ($\alpha$), ($\beta$), ($\gamma$), according to
how the two index sets overlap, and
\begin{align}
\label{eq-c-Z2-classes}
(\alpha)\ \leq C\rho_n\,\Xi_A,\qquad
(\beta)\ \leq C v_n\,\Delta_A,\qquad
(\gamma)\ \leq C\rho_n^2\,Q_A ,
\end{align}
with a diagonal term at most $C\rho_nQ_A$.  Classes ($\alpha$) and ($\gamma$)
are geometric series in $\rho_n$.  Class ($\beta$) is the delicate one: a
single cancellation, Lemma~\ref{lem-c-cancel}(c), is needed for the lower-order
kernel, and its leading term relative to $\Delta_A$ is $v_n$, which tends to
$0$ by \ref{condition-s2}.  The complete proof is in Appendix~\ref{app-c-Z2}.
\end{proof}

\paragraph{Step Z4: latent part, fourth moment.}

\begin{lemma}[Latent fourth moment]
\label{lem-c-Z4}
$\|Z_A\|_{L^4}\leq C\Delta_A^{\vartheta'}$ for every short cell $A$ and $n\geq n_0$.
\end{lemma}

\begin{proof}[Outline of the proof]
By \citet[Theorem~2]{boucheron2005moment} with $q = 4$,
\begin{align}
\label{eq-c-BBLM}
\|Z_A\|_{L^4}^2\leq C\|V\|_{L^2},\qquad
V := \sum_mE'_m\bigl[(Z_A-Z^{(m)}_A)^2\bigr]\leq4(V_1+V_2)+2(V_3+V_4),
\end{align}
with $V_1 := n^{-1}\sum_m\chi_mq_m^2$, $V_2$ its resampled analogue,
$V_3 := n^{-1}\sum_m\Phi_m^2$ and $V_4 := n^{-1}\sum_m\Psi_m$.  By
Lemma~\ref{lem-c-Z2}, $EV\leq C(\Delta_A+\Xi_A)\leq C\Delta_A^{2\vartheta}$, so it
suffices to show $\mathrm{Var}(V_k)\leq C\Delta_A^{4\vartheta'}$ for each $k$.  For
$V_1,V_2$ this is a direct Efron--Stein computation.  For $V_3$ and $V_4$,
after removing vertices $1$ and $2$ at a cost $o(\Delta_A^{4\vartheta'})$, both
reduce to $\mathrm{Var}(\Psi^*)$, where
\begin{align*}
\Psi^* := \|X\|^2_{L^2},\qquad X := \rho_n\sum_{i\geq3}\chi_ia(\xi_i,\cdot)\,g^{(1,2)}_i .
\end{align*}
We bound $\mathrm{Var}(\Psi^*)$ by the Efron--Stein inequality once more.  By
multilinearity, resampling $\xi_m$ changes $X$ exactly by
$Y(\xi_m)+V(\xi_m)$, a \emph{direct} term from the summand $i = m$ and an
\emph{indirect} term through the edge probabilities $p_{im}$, and the latter
satisfies the exact identity
\begin{align}
\label{eq-c-T-identity}
\mathrm{Var}_y\bigl\langle X^{(-m)},V(y)\bigr\rangle = \rho_n^4\bigl\|\mathcal T_mX^{(-m)}\bigr\|^2,
\qquad
\mathcal T_m := \sum_{i\neq m}\chi_i\delta_i\,a_i\otimes a_i ,
\end{align}
with $\delta_i$ a second difference of the conditional mass function and
$a_i := a(\xi_i,\cdot)$.  The operator formulation is natural because the
indirect effect of resampling $\xi_m$ is linear in the connection profile of
the new position: its inner product with $X^{(-m)}$ is a linear functional of
$a(\cdot,y)$, so its variance over $y$ is a quadratic form, namely the squared
norm of $X^{(-m)}$ under an operator built from the profiles $a_i$ of the
other vertices.  This separates the randomness of the operator from that of
the vector it acts on, and replaces a four-index expansion of the variance by
bounds of the form $\|\mathcal T_m\|\,\|X^{(-m)}\|$, in which the operator's
mean, its fluctuation and the vector can each be controlled by the estimates
already available.  The operator $\mathcal T_m$ is split into its mean,
bounded with the cancellation of Lemma~\ref{lem-c-cancel}(c) at $d = 1$, and
its fluctuation, bounded in Hilbert--Schmidt norm without cancellation and in
operator norm on an event of probability close to one.  All terms are
$O(\Delta_A^{4\vartheta'})$ except possibly the product
$P_A := \rho_n^2\mu_A^2\Xi_A$, which by \eqref{eq-c-Xi-Delta} and $\mu_A\leq1$ satisfies
\begin{align}
\label{eq-c-PA}
\frac{P_A}{\Delta_A^{4\vartheta'}}\leq C\rho_n^2h_t^{2\vartheta-4\vartheta'}
 \leq Cn^{2\vartheta'-\vartheta}\rho_n^{2+2\vartheta'-\vartheta}\leq C
\end{align}
by \eqref{eq-c-theta-prime}(iv).  This is the only place where
\ref{condition-s2} is used beyond $v_n\to0$.  The complete proof is in
Appendix~\ref{app-c-Z4}.
\end{proof}

\paragraph{Step A: assembly.}

\begin{lemma}[Single-cell bound]
\label{lem-c-single}
$\|\Delta'_A\|_{L^4}\leq C\Delta_A^{\vartheta'}$ for every cell $A\subseteq I$ and
$n\geq n_0$.
\end{lemma}

\begin{proof}
By Step~C it suffices to treat short cells.  Then
$\|\Delta'_A\|_4\leq\|B_A\|_4+\|Z_A\|_4$.  By Lemma~\ref{lem-c-edge},
Lemma~\ref{lem-c-Z2} and \eqref{eq-c-absorb},
\begin{align*}
EB_A^4\leq C\bigl(\Delta_A^2+n^{-1}\Delta_A+n^{-1}\Xi_A\bigr)+\mathfrak e_n\leq C\Delta_A^{4\vartheta'} .
\end{align*}
This uses \eqref{eq-c-Xi-Delta} and
$h_t^{-(4\vartheta'-2\vartheta)}\asymp(n\rho_n)^{2\vartheta'-\vartheta}\leq n^{1/2}$, which give
\begin{align*}
n^{-1}\Xi_A\leq Cn^{-1}h_t^{2\vartheta-4\vartheta'}\Delta_A^{4\vartheta'}
\leq Cn^{-1/2}\Delta_A^{4\vartheta'} .
\end{align*}
Finally apply Lemma~\ref{lem-c-Z4}.
\end{proof}

\begin{proof}[Proof of Theorem~\ref{thm-tightness-compact}]
By Step~H it suffices to prove $\mathbb G'_n\rightsquigarrow\mathbb G$ in
$\ell^\infty(I)$.  We verify Proposition~\ref{prop-grid} on $I$ with
$T_n\equiv t_1$ and $\nu := 2\vartheta'>1/2$ (Step~L).

\emph{Condition (i).}  The fidi convergence of $\mathbb G'_n$ follows from
Theorem~\ref{thm-fidi-clt} and Step~H.  For continuity of the limit, let
$s<t$ in $I$, and let $A_n$ be the cell between $\lambda_{s'}$ and
$\lambda_{t'}$.  By Lemma~\ref{lem-c-single} and Fatou's lemma,
\begin{align*}
E|\mathbb G(t)-\mathbb G(s)|^4\leq C\bigl(F(t)-F(s)\bigr)^{4\vartheta'},
\end{align*}
and the Kolmogorov--Chentsov theorem in the parametrization $u = F(t)$ gives a
continuous version.

\emph{Condition (ii).}  For lattice points $r<s<t$ in $I$, by Cauchy--Schwarz
and Lemma~\ref{lem-c-single},
\begin{align*}
E\bigl[(\Delta'_{A_{rs}})^2(\Delta'_{A_{st}})^2\bigr]\leq C(\Delta_{A_{rs}}\Delta_{A_{st}})^{2\vartheta'}
\leq C\bigl(F(t)-F(r)\bigr)^{4\vartheta'} .
\end{align*}
\end{proof}

\subsection{Empirical Processes With Estimated Sparsity}
In statistical applications, $\rho_n$ is typically an unknown nuisance parameter that must be estimated. Let $\hat\rho_n := \frac{1}{{n \choose 2}}\sum_{1 \leq i < j \leq n }A_{ij}$; \citet{Bickel-Chen-Levina-method-of-moments} show that $\sqrt{n}(\hat{\rho}_n/\rho_n - 1)$ satisfies a central limit theorem if $n \rho_n \rightarrow \infty$. Define the
\emph{feasible process}
\begin{align}
\label{eq-feasible-process}
\hat{\mathbb G}_n(t) := \sqrt n\Bigl[\frac1n\sum_{i=1}^n
 \mathbbm 1\bigl(D_i\leq\lfloor n'\hat\rho_nt\rfloor\bigr)-F(t)\Bigr],
\end{align}
in which degrees are normalized by $n'\hat\rho_n$ in place of $n'\rho_n$.  For
$x\in[0,1]$ let $(T_w\mathfrak D)(x) := \int_0^1w(x,y)\mathfrak D(y)\,dy$, the
expected number of paths of length two from a vertex with latent position $x$
(in units of $\rho_n^2n^2$), and recall $x_s := \mathfrak D^{-1}(s)$.  The following result establishes uniform CLTs for $\hat{\mathbb{G}}_n$, on
compact intervals and over the tail, using an argument that exploits the
stochastic equicontinuity of $\mathbb{G}_n$.

\begin{corollary}[Estimated sparsity]
\label{cor-feasible}
Let $\hat{\mathbb G} := \mathbb G + 2\,\iota f\,Z$, where $\iota(t) := t$,
$f := f_{\mathfrak D}$, and $(\mathbb G,Z)$ is jointly Gaussian with
$\mathrm{Var}(Z) = \mathrm{Var}(\mathfrak D)$ and
$\mathrm{Cov}(\mathbb G(s),Z) = C(s)$, where:
\begin{align}
\label{eq-feasible-cross}
C(s) := E\bigl[(\mathfrak D(\xi)-1)\,\mathbbm 1(\mathfrak D(\xi)\leq s)\bigr]
 - f(s)\bigl[(T_w\mathfrak D)(x_s)-s\bigr].
\end{align}
Then $\hat{\mathbb G}$ is a centered Gaussian process with covariance
\begin{align*}
\hat\Sigma(s,t) = \Sigma(s,t)+2sf(s)\,C(t)+2tf(t)\,C(s)
+4\,st\,f(s)f(t)\,\mathrm{Var}(\mathfrak D),
\end{align*}
and the following hold.
\begin{enumerate}[(a)]
\item \emph{Compact intervals.}  If the conditions of Corollary~\ref{cor-Gn}
hold on a compact interval $I$, and $I_0$ is a compact interval contained in
the interior of $I$, then $\hat{\mathbb G}_n\rightsquigarrow\hat{\mathbb G}$ in
$\ell^\infty(I_0)$, and $\hat{\mathbb G}$ has continuous sample paths on $I_0$.
\item \emph{Unbounded intervals.}  If the conditions of
Corollary~\ref{cor-Gn-tail} hold, then for every $t_0'>t_0$,
$\hat{\mathbb G}_n\rightsquigarrow\hat{\mathbb G}$ in
$\ell^\infty([t_0',\infty))$, and $\hat{\mathbb G}$ has continuous sample paths on
$[t_0',\infty)$.
\end{enumerate}
\end{corollary}

\begin{proof}
\emph{Part (a).}  Since $\lfloor n'\hat\rho_nt\rfloor = \lfloor n'\rho_n(t\hat\rho_n/\rho_n)\rfloor$,
the feasible process is the original one evaluated at a random point, plus a
deterministic shift:
\begin{align}
\label{eq-feasible-split}
\hat{\mathbb G}_n(t)
 = \mathbb G_n\bigl(t\hat\rho_n/\rho_n\bigr)
 + \sqrt n\bigl[F(t\hat\rho_n/\rho_n)-F(t)\bigr].
\end{align}

\emph{The estimator $\hat\rho_n$.}  Write
$\hat\rho_n/\rho_n-1 = U_n+V_n+(E\check w-1)$, where
\begin{align*}
U_n := \frac1{nn'}\sum_{i\neq j}\bigl\{\check w(\xi_i,\xi_j)-E\check w\bigr\},
\qquad
V_n := \frac1{nn'\rho_n}\sum_{i\neq j}\bigl\{A_{ij}-\rho_n\check w(\xi_i,\xi_j)\bigr\}
\end{align*}
are a centered $U$-statistic and the edge noise.  Given the latent positions,
$V_n$ is a sum of independent centered terms, so
\begin{align*}
E[V_n^2]\leq\frac{2E\check w}{nn'\rho_n} = O\bigl((n^2\rho_n)^{-1}\bigr),
\end{align*}
and $\sqrt nV_n\to0$ in $L^2$.  The H\'ajek projection of $U_n$ is
$(2/n)\sum_j\{\check{\mathfrak D}(\xi_j)-E\check w\}$, and its remainder, a
degenerate $U$-statistic, has variance $O(n^{-2}E\check w^2)$.  Finally
$0\leq1-E\check w = \int g\leq C_w\rho_n^{1+\beta}$ by \eqref{eq-g-global} and
$E(\mathfrak D-\check{\mathfrak D})^2\to0$ by dominated convergence, so that
\begin{align}
\label{eq-rho-hat}
\sqrt n\Bigl(\frac{\hat\rho_n}{\rho_n}-1\Bigr)
 = \frac2{\sqrt n}\sum_{j=1}^n\bigl\{\mathfrak D(\xi_j)-1\bigr\}+o_P(1),
\end{align}
using $\sqrt n\rho_n^{1+\beta}\to0$, which is implied by
$n\rho_n^{3/2+\beta}\to0$.

\emph{Joint finite-dimensional convergence.}  Apply the Cram\'er--Wold argument
of the proof of Theorem~\ref{thm-fidi-clt} to
$\sum_ka_k\mathbb G_n(t_k)+b\sqrt n(\hat\rho_n/\rho_n-1)$.  By
\eqref{eq-rho-hat}, the additional term contributes
$2b\{\mathfrak D(\xi_j)-1\}$ to the i.i.d.\ summands of the H\'ajek projection
and an $o_P(1)$ remainder.  The Lindeberg condition for the enlarged summands
holds since $\mathfrak D(\xi)$ has a finite moment of order $2+\beta$ by
\ref{condition-m1}.  Hence
$(\mathbb G_n(t_1),\ldots,\mathbb G_n(t_K),\sqrt n(\hat\rho_n/\rho_n-1))$
converges jointly to a Gaussian vector $(\mathbb G(t_1),\ldots,\mathbb G(t_K),2Z)$,
where $\mathrm{Var}(Z) = \mathrm{Var}(\mathfrak D)$ and
$\mathrm{Cov}(\mathbb G(s),Z) = E[(\psi^\infty_s+\phi^\infty_s)(\xi)\,
\{\mathfrak D(\xi)-1\}]$.  With
$\psi^\infty_s = \mathbbm 1(\mathfrak D\leq s)-F(s)$,
$\phi^\infty_s = f(s)(s-w(x_s,\cdot))$ and $E\mathfrak D = 1$,
\begin{align*}
E\bigl[\psi^\infty_s(\mathfrak D-1)\bigr]
 &= E\bigl[(\mathfrak D-1)\mathbbm 1(\mathfrak D\leq s)\bigr],
\\ E\bigl[\phi^\infty_s(\mathfrak D-1)\bigr]
 &= -f(s)\Bigl(E\bigl[w(x_s,\xi)\mathfrak D(\xi)\bigr]-s\Bigr),
\end{align*}
and the sum is $C(s)$ in \eqref{eq-feasible-cross}.

\emph{The random evaluation point.}  By \eqref{eq-rho-hat},
$t_k\hat\rho_n/\rho_n\to t_k$ in probability, so it lies in $I$ with
probability tending to one.  By Corollary \ref{cor-Gn}, $\mathbb G_n$ is
asymptotically uniformly equicontinuous on $I$, and therefore
$\mathbb G_n(t_k\hat\rho_n/\rho_n)-\mathbb G_n(t_k)\to0$ in probability.

\emph{The shift.}  Since $F$ is differentiable at $t_k$ with derivative
$f(t_k)$ by \ref{condition-d1}, and $\sqrt n(\hat\rho_n/\rho_n-1) = O_P(1)$,
the delta method gives
$\sqrt n[F(t_k\hat\rho_n/\rho_n)-F(t_k)]
= t_kf(t_k)\sqrt n(\hat\rho_n/\rho_n-1)+o_P(1)$.

Combining the three displays with \eqref{eq-feasible-split},
$\hat{\mathbb G}_n(t_k) = \mathbb G_n(t_k)+t_kf(t_k)\sqrt n(\hat\rho_n/\rho_n-1)
+o_P(1)$, which converges jointly in $k$ to $\mathbb G(t_k)+2t_kf(t_k)Z$.

\emph{Uniformity.}  By  Corollary \ref{cor-Gn}, $\mathbb G_n$ is
asymptotically tight in $\ell^\infty(I)$, and
$\sqrt n(\hat\rho_n/\rho_n-1)$ is asymptotically tight in $\mathbb R$ by
\eqref{eq-rho-hat}; hence the pair is asymptotically tight in
$\ell^\infty(I)\times\mathbb R$, and by the joint finite-dimensional convergence
above, $(\mathbb G_n,\sqrt n(\hat\rho_n/\rho_n-1))\rightsquigarrow(\mathbb G,2Z)$
there.  Next, $\sup_{t\in I_0}|t\hat\rho_n/\rho_n-t|\leq(\sup I_0)
|\hat\rho_n/\rho_n-1|\to0$ in probability, so $t\hat\rho_n/\rho_n\in I$ for all
$t\in I_0$ with probability tending to one, and asymptotic uniform
equicontinuity of $\mathbb G_n$ on $I$ gives
$\sup_{t\in I_0}|\mathbb G_n(t\hat\rho_n/\rho_n)-\mathbb G_n(t)|\to0$ in
probability.  Finally, by the mean value theorem and continuity of $f$ on $I$,
with modulus of continuity $\omega_f$,
\begin{align*}
& \sup_{t\in I_0}\Bigl|\sqrt n\bigl[F(t\hat\rho_n/\rho_n)-F(t)\bigr]
 - tf(t)\sqrt n\Bigl(\frac{\hat\rho_n}{\rho_n}-1\Bigr)\Bigr|
 \\ \leq \ & \ (\sup I_0)\,\sqrt n\Bigl|\frac{\hat\rho_n}{\rho_n}-1\Bigr|\,
 \omega_f\Bigl((\sup I_0)\Bigl|\frac{\hat\rho_n}{\rho_n}-1\Bigr|\Bigr)
 \rightarrow0
\end{align*}
in probability.  By \eqref{eq-feasible-split},
$\hat{\mathbb G}_n = \mathbb G_n+\iota f\,\sqrt n(\hat\rho_n/\rho_n-1)+o_P(1)$
uniformly on $I_0$, and the claim follows from the continuous mapping theorem
applied to $(h,z)\mapsto h+\iota fz$, which is continuous from
$\ell^\infty(I)\times\mathbb R$ to $\ell^\infty(I_0)$ since $\iota f$ is bounded
on $I_0$.

\emph{Part (b).}  The expansion \eqref{eq-rho-hat}, the joint
finite-dimensional convergence and the identification of $C(s)$ do not involve
the index set, and \eqref{eq-rho-hat} uses only $\sqrt n\rho_n^{1+\beta}\to0$,
which is implied by \eqref{eq-centering-conditions}.  Write
$u_n := \hat\rho_n/\rho_n$.  Three steps of part (a) change.

\emph{Tightness.}  By Corollary~\ref{cor-Gn-tail}, $\mathbb G_n$ is
asymptotically tight in $\ell^\infty([t_0,\infty))$, so, as in part (a),
\begin{align*}
\bigl(\mathbb G_n,\sqrt n(u_n-1)\bigr)\rightsquigarrow(\mathbb G,2Z)
\qquad\text{in }\ell^\infty([t_0,\infty))\times\mathbb R .
\end{align*}

\emph{The random evaluation point.}  Let $d_F(s,t) := |F(s)-F(t)|$.  The set
$[t_0,\infty)$ is totally bounded under $d_F$, and by the proof of
Theorem~\ref{thm-tightness-ring} the limit $\mathbb G$ has a version whose
paths are uniformly $d_F$-continuous.  Hence $\mathbb G_n$ is asymptotically
uniformly $d_F$-equicontinuous in probability
\citep[Theorem~1.5.7 and Addendum~1.5.8]{van2016weak}.  Fix
$\varepsilon\in(0,\tfrac12)$ with $(1-\varepsilon)t_0'\geq t_0$.  On the event
$\{|u_n-1|\leq\varepsilon\}$, whose probability tends to one, $tu_n\geq t_0$ for
every $t\geq t_0'$, and by \ref{condition-d3}(i),
\begin{align*}
\sup_{t\geq t_0'}d_F(tu_n,t) = \sup_{t\geq t_0'}\Bigl|\int_t^{tu_n}f(\lambda)\,d\lambda\Bigr|
 \leq C_d\sup_{t\geq t_0}tf(t)\,|u_n-1|\;\rightarrow\;0
\end{align*}
in probability, using $\sup_{t\geq t_0}tf(t)<\infty$, which follows from
\ref{condition-d3}(i) and (iii).  Therefore
$\sup_{t\geq t_0'}|\mathbb G_n(tu_n)-\mathbb G_n(t)|\to0$ in probability.

\emph{The shift.}  On the same event, \ref{condition-d3}(ii) at the base point
$t$ gives $|f(\lambda)-f(t)|\leq K_ff(t)|u_n-1|^{\theta_f}$ for $\lambda$ between
$t$ and $tu_n$, so
\begin{align*}
\sup_{t\geq t_0'}\sqrt n\Bigl|F(tu_n)-F(t)-tf(t)(u_n-1)\Bigr|
 &\leq\sup_{t\geq t_0'}\sqrt n\Bigl|\int_t^{tu_n}\bigl(f(\lambda)-f(t)\bigr)d\lambda\Bigr|
\\ &\leq K_f\sup_{t\geq t_0}tf(t)\;\sqrt n\,|u_n-1|^{1+\theta_f}
 = O_P\bigl(n^{-\theta_f/2}\bigr).
\end{align*}

\emph{Conclusion.}  By \eqref{eq-feasible-split},
$\hat{\mathbb G}_n = \mathbb G_n+\iota f\,\sqrt n(u_n-1)+o_P(1)$ uniformly on
$[t_0',\infty)$.  Since $\iota f$ is bounded on $[t_0,\infty)$, the map
$(h,z)\mapsto(h+\iota fz)|_{[t_0',\infty)}$ is continuous from
$\ell^\infty([t_0,\infty))\times\mathbb R$ to $\ell^\infty([t_0',\infty))$, and the
continuous mapping theorem gives the claim.  The sample paths of
$\hat{\mathbb G}$ are continuous because those of $\mathbb G$ are and $f$ is
continuous on $J_\infty$ by \ref{condition-d3}(ii).
\end{proof}

\section{On the Failure of Uniform CLTs for SBMs}
\label{sec-sbm}

The uniform results of the previous sections assume that the law of
$\mathfrak D(\xi)$ is absolutely continuous.  This condition cannot simply be
dropped.  If $\mathfrak D(\xi)$ has an atom, the variance of $\mathbb G_n(t)$
at the atom diverges, at the rate $n\rho_n$, while away from the atoms it
remains bounded.  The mechanism is visible in the H\'ajek kernel itself.  The
measure $\mu_t$ of Lemma~\ref{lem-mu-limit} has total mass
$n'\rho_nP(D_{1\setminus2} = t')$; when $\mathfrak D$ has a density, this
mass converges to $f_{\mathfrak D}(t)$, but at an atom of mass $\pi$ it is of
order $\pi\sqrt{n\rho_n}$, since a positive fraction of vertices then have
expected degree exactly at the lattice point.  We demonstrate this for
stochastic block models.

For a $K$-block SBM, let $Z_1,\ldots,Z_n\stackrel{iid}{\sim}
\mathrm{Categorical}(\pi)$ denote the block memberships and, for
$1\leq i<j\leq n$,
\begin{align*}
A_{ij} = \mathbbm 1\bigl(\eta_{ij}\leq\rho_nB_{Z_iZ_j}\wedge1\bigr),
\end{align*}
with $(\eta_{ij})$ i.i.d.\ $\mathrm{Uniform}[0,1]$ and $B\in[0,\infty)^{K\times K}$
symmetric with $\sum_{v,v'}\pi_vB_{vv'}\pi_{v'} = 1$.  This is the graphon model
\eqref{eq-sparse-graphon} with $w(x,y) = B_{z(x)z(y)}$, where $z(x)$ is the
block of $x$ under a partition of $[0,1]$ into intervals of lengths $\pi_v$, so
that $Z_i = z(\xi_i)$ and $\mathfrak D(\xi) = \mathfrak D(Z) :=
\sum_{v}\pi_vB_{Zv}$.  The law of $\mathfrak D$ is supported on the finite set
$\{\mathfrak D(v):v\in[K]\}$.  Since $B$ is bounded, $\rho_nB_{vv'}\leq\frac14$
for $n$ large, so no clipping occurs and $\check{\mathfrak D} = \mathfrak D$.

\begin{proposition}[Variance explosion at atoms]
\label{prop-sbm-variance}
Let $\rho_n\to0$ with $n\rho_n\to\infty$.  Let $u\in[K]$ satisfy
$\pi_u>0$, $\mathfrak D(u)>0$ and $\mathfrak D(u)\neq\mathfrak D(v)$ for all
$v\neq u$ with $\pi_v>0$, and suppose
\begin{align}
\label{eq-sbm-row-variance}
\mathrm{Var}\bigl(B_{uZ}\bigr) = \sum_{v}\pi_v\bigl(B_{uv}-\mathfrak D(u)\bigr)^2>0 .
\end{align}
Then at $t = \mathfrak D(u)$,
\begin{align}
\label{eq-sbm-variance}
\mathrm{Var}\bigl(\mathbb G_n(t)\bigr) = \mathrm{Var}\bigl(\mathring{\mathbb G}_n(t)\bigr)
 \;\geq\; \frac{\pi_u^2\,\mathrm{Var}(B_{uZ})}{2\pi\,\mathfrak D(u)}\;n\rho_n\,\bigl(1+o(1)\bigr).
\end{align}
\end{proposition}

\begin{proof}
The two processes differ by a deterministic shift, so their variances agree.
The H\'ajek projection $H(\mathring{\mathbb G}_n(t))$ is the $L^2$ projection of
$\mathring{\mathbb G}_n(t)$ onto sums of functions of single latent positions,
so $\mathrm{Var}(\mathring{\mathbb G}_n(t))\geq\mathrm{Var}(H(\mathring{\mathbb G}_n(t)))$.
The representation \eqref{eq-psi-phi-def} of the projection is algebraic and
holds for any $w$ with $\rho_nw\leq1$; it gives
$H(\mathring{\mathbb G}_n(t)) = n^{-1/2}\sum_j(\psi_t+\phi_t)(\xi_j)$, a
normalized sum of i.i.d.\ centered terms, so that
\begin{align}
\label{eq-sbm-hajek-lower}
\mathrm{Var}\bigl(\mathbb G_n(t)\bigr)\;\geq\;E\bigl[(\psi_t+\phi_t)^2\bigr]
 \;\geq\;\bigl(\|\phi_t\|_{L^2}-\|\psi_t\|_{L^2}\bigr)_+^2
 \;\geq\;\bigl(\|\phi_t\|_{L^2}-1\bigr)_+^2,
\end{align}
the last step since $|\psi_t|\leq1$.  It remains to bound $\|\phi_t\|_{L^2}$
from below.

For the SBM, $a(x,\xi) = B_{z(x)Z}-\mathfrak D(z(x))$ with $Z = z(\xi)$, and the
binomial parameter $\rho_n\mathfrak D(x)$ depends on $x$ only through its block.
Integrating block by block, $\phi_t$ is the finite sum
\begin{align}
\label{eq-sbm-phi}
\phi_t(\xi) = -\sum_{v\in[K]}c_{n,v}\bigl(B_{vZ}-\mathfrak D(v)\bigr),
\qquad
c_{n,v} := n'\rho_n\,\pi_v\,b_{n''}\bigl(t';\rho_n\mathfrak D(v)\bigr).
\end{align}

\emph{The blocks $v\neq u$.}  Since $\mathfrak D(v)\neq t$, there is
$\epsilon_0>0$ with $|\mathfrak D(v)-t|\geq\epsilon_0$ for all such $v$ with
$\pi_v>0$, and $t' = n'\rho_nt+O(1)$ is separated from the mean
$n''\rho_n\mathfrak D(v)$ by at least $\epsilon_0n\rho_n/2$ for $n$ large.  By
Lemma~\ref{lem-binomial-facts}\ref{binom-fact-b} in both tails,
$b_{n''}(t';\rho_n\mathfrak D(v))\leq2e^{-c\,n\rho_n}$ for some $c>0$ (and it
vanishes if $\mathfrak D(v) = 0$), so $c_{n,v}\leq Cn\rho_ne^{-c\,n\rho_n}\to0$
because $n\rho_n\to\infty$.  As $B$ is bounded, these terms contribute $o(1)$ to
$\phi_t$, uniformly.

\emph{The block $u$.}  Here $p := \rho_nt\to0$, the binomial
$\mathrm{Bin}(n'',p)$ has variance $\sigma^2 = n''p(1-p)\to\infty$, and the
lattice point $t'$ satisfies $|t'-n''p|\leq|t'-n'\rho_nt|+\rho_nt\leq2$, so its
standardized value tends to $0$.  By the de Moivre--Laplace local limit theorem
\citep[Ch.~VII]{feller1971volume2},
$b_{n''}(t';p) = (2\pi\sigma^2)^{-1/2}(1+o(1))$, hence
\begin{align*}
c_{n,u} = \frac{n'\rho_n\,\pi_u}{\sqrt{2\pi\,n''\rho_nt}}\bigl(1+o(1)\bigr)
 = \pi_u\sqrt{\frac{n\rho_n}{2\pi\,\mathfrak D(u)}}\,\bigl(1+o(1)\bigr).
\end{align*}

Combining, $\phi_t = -c_{n,u}(B_{uZ}-\mathfrak D(u))+r_n$ with
$\|r_n\|_\infty\to0$, and since $E[B_{uZ}] = \mathfrak D(u)$,
\begin{align*}
\|\phi_t\|_{L^2}^2 = c_{n,u}^2\,\mathrm{Var}(B_{uZ})+O(c_{n,u}\|r_n\|_\infty)+o(1)
 = \frac{\pi_u^2\,\mathrm{Var}(B_{uZ})}{2\pi\,\mathfrak D(u)}\,n\rho_n\,\bigl(1+o(1)\bigr),
\end{align*}
which diverges by \eqref{eq-sbm-row-variance}.  Substituting into
\eqref{eq-sbm-hajek-lower} gives \eqref{eq-sbm-variance}.
\end{proof}

Several comments are in order.

\begin{remark}[What the proposition says about continuity]
\label{rem-sbm-continuity}
At an atom the H\'ajek kernel itself diverges: the coefficient $c_{n,u}$ is the
contribution of block $u$ to the mass of $\mu_t$, which in the absolutely
continuous case converges to the density $f_{\mathfrak D}(t)$ by
Lemma~\ref{lem-mu-limit}.  An atom corresponds to an infinite density, and the
divergence of $c_{n,u}$ at rate $\sqrt{n\rho_n}$ is the finite-$n$ form of that.
Since $\mathrm{Var}(\mathbb G_n(t))$ remains bounded at points away from the
atoms (Theorem~\ref{thm-sbm-donsker} below), the variance is unbounded on any
index set containing an atom, and no uniform CLT can hold there.
\end{remark}


Away from the atoms, by contrast, the empirical degree distribution is
asymptotically indistinguishable from the empirical distribution of the
expected degrees, and a Donsker theorem holds with no condition beyond
$n\rho_n\gg\log n$.  Let $d_1<\cdots<d_L$ denote the distinct values of
$\mathfrak D(v)$ over blocks with $\pi_v>0$, fix $\epsilon>0$, and let
\begin{align*}
\mathcal T_\epsilon := [0,\infty)\setminus\bigcup_{l=1}^L(d_l-\epsilon,d_l+\epsilon).
\end{align*}

\begin{theorem}[A Donsker theorem for sparse SBMs]
\label{thm-sbm-donsker}
Suppose that $\rho_n\to0$ and $n\rho_n/\log n\to\infty$.  Then, for
every $\epsilon>0$,
\begin{align*}
\mathbb G_n\rightsquigarrow\mathbb B\circ F\qquad\text{in }\ell^\infty(\mathcal T_\epsilon),
\end{align*}
where $\mathbb B$ is a standard Brownian bridge, so that the limit has
covariance $F(s)\wedge F(t)-F(s)F(t)$.
\end{theorem}

\begin{proof}
Let $\mathbb G^Z_n(t) := n^{-1/2}\sum_i[\mathbbm 1(\mathfrak D(Z_i)\leq t)-F(t)]$,
the empirical process of the i.i.d.\ variables $\mathfrak D(Z_i)$.  By the
classical Donsker theorem for indicators of half-lines,
$\mathbb G^Z_n\rightsquigarrow\mathbb B\circ F$ in $\ell^\infty(\mathbb R)$, hence
in $\ell^\infty(\mathcal T_\epsilon)$.  It therefore suffices to show
$\sup_{t\in\mathcal T_\epsilon}|\mathbb G_n(t)-\mathbb G^Z_n(t)|\to0$ in
probability.

For $t\in\mathcal T_\epsilon$ and each $i$, $\mathfrak D(Z_i)$ is a point
$d_l$ at distance at least $\epsilon$ from $t$.  Hence, the indicators
$\mathbbm 1(D_i\leq t')$ and $\mathbbm 1(\mathfrak D(Z_i)\leq t)$ can differ for
\emph{some} $t\in\mathcal T_\epsilon$ only on the single event
$E_i := \{|D_i/(n'\rho_n)-\mathfrak D(Z_i)|\geq\epsilon/2\}$, for $n$ large.  Consequently,
\begin{align*}
\sup_{t\in\mathcal T_\epsilon}\bigl|\mathbb G_n(t)-\mathbb G^Z_n(t)\bigr|
 \;\leq\;\frac1{\sqrt n}\sum_{i=1}^n\mathbbm 1(E_i).
\end{align*}
To bound $P(E_i)$, write $D_i/(n'\rho_n)-\mathfrak D(Z_i)$ as the sum of
$n'^{-1}\sum_{j\neq i}(B_{Z_iZ_j}-\mathfrak D(Z_i))$, an average of bounded
i.i.d.\ centered variables given $Z_i$, and
$(n'\rho_n)^{-1}\sum_{j\neq i}(A_{ij}-\rho_nB_{Z_iZ_j})$, a centered sum of
independent Bernoulli variables given $\mathbf Z$ with mean at most
$Cn\rho_n$.  By Hoeffding's inequality for the first and Bernstein's
inequality for the second,
$P(E_i)\leq2e^{-c\,n\epsilon^2}+2e^{-c\,n\rho_n\epsilon^2}$ for some $c>0$
depending only on $B$.  Hence
\begin{align*}
E\Bigl[\frac1{\sqrt n}\sum_{i=1}^n\mathbbm 1(E_i)\Bigr]
 \;\leq\;4\sqrt n\,e^{-c\,n\rho_n\epsilon^2}\;\longrightarrow\;0,
\end{align*}
since $n\rho_n/\log n\to\infty$, and Markov's inequality completes the proof.
\end{proof}

\section{Conclusion}
In this work, we have established central limit theorems for the empirical degree
distribution of sparse graphons with unbounded kernels. The finite-dimensional
result holds throughout the regime $n\rho_n\gg\log n$, subject only to an upper
bound on $\rho_n$ governed by the global moment of the kernel, and its limiting
covariance coincides with that of the dense case; under additional regularity
the convergence is functional.  The analysis rests on an exact multilinear
representation of the conditional degree distribution, which identifies the
Hoeffding components explicitly and so avoids the Edgeworth expansions on which
earlier work relied. 

Since the validity of subsampling for sparse graphons requires only the
existence of a limiting distribution, exactly as in the i.i.d.\ setting
\citep{lunde2023subsampling, politis1999subsampling}, our results already permit inference on the degree distribution,
including Kolmogorov-Smirnov-type tests based on the uniform result.  It would
be of interest to study alternative resampling schemes, such as bootstrap
procedures, which may improve on the rate of convergence of subsampling and
avoid the choice of a subsample size.  Two regimes lie outside the scope of
this work and appear to call for different limits. In the ultra-sparse regime $n\rho_n\to c\in(0,\infty)$, degrees no longer
concentrate about their expected values, so that the natural target is a mixed
Poisson law \citep{bollobas2007phase}, and the randomness of the edges
contributes to the fluctuations at the same order as that of the latent
positions.  The form of the limit in this regime appears to be open beyond
homogeneous models.  For heavy-tailed degree distributions with infinite variance,
such as power laws with exponent in $(2,3)$, the Gaussian approximation likely fails
outright.  In both cases non-Gaussian limits seem likely, although whether a
limiting distribution can be identified in a tractable form remains to be seen.

\subsection*{Acknowledgements}
Robert Lunde acknowledges support from NSF grant DMS-2515923. The arguments here refine and extend an earlier, unpublished version of this
work, which was based on smoothed lattice Edgeworth expansions
\citep{DEHEUVELS1989282} and yielded a finite-dimensional central limit
theorem under $\rho_n\gg n^{-2/3}$.  The authors discussed the proof strategy with Claude Fable~5.1
and Claude Opus~5.5.  The models suggested the
multilinear decomposition and drafted portions of the tightness arguments.
All arguments have been verified by the authors, who take full responsibility
for the results.

\bibliographystyle{plainnat}
\bibliography{mybib}

\appendix
\counterwithin{equation}{section}
\counterwithin{definition}{section}
\counterwithin{theorem}{section}
\counterwithin{proposition}{section}
\counterwithin{remark}{section}

\section{Proofs of the latent moment bounds}
\label{app-latent}

This appendix contains the complete proofs of Lemmas~\ref{lem-c-Z2} and
\ref{lem-c-Z4}, whose outlines are given in
Section~\ref{sec-tight-compact-proof}.  Notation and conventions are those of
Section~\ref{sec-tight-compact-proof}.

\subsection{Step Z2}
\label{app-c-Z2}

\begin{proof}[Proof of Lemma~\ref{lem-c-Z2}]
\emph{Reduction.}  The Efron--Stein inequality, exchangeability and
\eqref{eq-c-Z-diff} give
\begin{align*}
\mathrm{Var}(Z_A)\leq E\bigl[(q_1-q^{(1)}_1)^2\bigr]+2E\Phi_1^2,
\end{align*}
since $\frac12(T_1+\Phi_1-\Phi_1')^2\leq T_1^2+(\Phi_1-\Phi_1')^2$.  Write
$e_i := g^{(1)}_i-\bar g(\xi_i)$ and $\Phi_1 = \Phi^{\rm m}+\Phi^{\rm f}$ with
\begin{align*}
\Phi^{\rm m} := \rho_n\sum_i\chi_ia(\xi_i,\xi_1)\bar g(\xi_i),\qquad
\Phi^{\rm f} := \rho_n\sum_i\chi_ia(\xi_i,\xi_1)e_i .
\end{align*}

\emph{Main part.}  Given $\xi_1$, $\Phi^{\rm m}$ is a sum of i.i.d.\ terms with
mean $(n''/n')\phi_A(\xi_1)$, so by Proposition~\ref{prop-c-hajek},
\eqref{eq-c-a-window}, Lemma~\ref{lem-c-cancel}(b), and
$n\rho_n^2m_A^2h_t\asymp\rho_nQ_A$,
\begin{align*}
E(\Phi^{\rm m})^2\leq\|\phi_A\|^2+n''\rho_n^2E\bigl[\chi\|a(\xi,\cdot)\|^2\bar g(\xi)^2\bigr]
\leq C\Xi_A+C\rho_nQ_A .
\end{align*}

\emph{Fluctuation part: structure.}  With $C(x,y) := \langle a(x,\cdot),a(y,\cdot)\rangle$,
\begin{align*}
E(\Phi^{\rm f})^2 = n''\rho_n^2E\bigl[\chi_2\|a(\xi_2,\cdot)\|^2e_2^2\bigr]
+n''(n''-1)\rho_n^2E\bigl[\chi_2\chi_3C(\xi_2,\xi_3)e_2e_3\bigr] .
\end{align*}
The diagonal is at most $C\rho_nQ_A$, since $|e_2|\leq Cm_A\omega_A(\xi_2)$.
For the off-diagonal term, Lemma~\ref{lem-multilinear}, applied to $g^{(1)}_2$
given $\xi_2$, gives the components
\begin{align*}
\Pi(g^{(1)}_2;S) = \rho_n^{|S|}\prod_{j\in S}a(\xi_2,\xi_j)K_{|S|}(\xi_2),
\qquad S\subseteq[n]\setminus\{1,2\},
\end{align*}
with $K_d := \mathsf K^{(d)}_A$ and $M = n''-d$.  A pair of components
$(S,S')$ of $(e_2,e_3)$ contributes only if
$S\setminus\{3\} = S'\setminus\{2\} =: T$.  This leaves three classes, with
$d := |T|$ and $\Gamma_k := C^k$:
\begin{align*}
(\alpha)\quad&\rho_n^{2d}\iint\Gamma_{d+1}K_d(x)K_d(y),\qquad d\geq1;\\
(\beta)\quad&\rho_n^{2d+1}\iint a(x,y)\Gamma_{d+1}K_{d+1}(x)K_d(y),\qquad d\geq1,\ \text{and its mirror};\\
(\gamma)\quad&\rho_n^{2d+2}\iint a(x,y)a(y,x)\Gamma_{d+1}K_{d+1}(x)K_{d+1}(y),\qquad d\geq0.
\end{align*}
Here $x,y$ range over $[0,1]\setminus\mathscr H$, and contributions off $W_A$ are
negligible.  There are at most $n^d/d!$ sets $T$, and the prefactor is
$n''(n''-1)\rho_n^2$.

\emph{Bounds used.}  For $x,y\in W_A$, by \eqref{eq-c-a-window}, with
$\varepsilon(x) := e_A(x)^{\theta_w}$:
\begin{enumerate}[($\Gamma$1)]
\item $|\Gamma_k|\leq C^k$;
\item $|\Gamma_k(x,y)-\Gamma_k(x_t,y)|\leq kC^k\varepsilon(x)$;
\item the mixed difference of $\Gamma_k$ in $(x,y)$ is at most
$k^2C^k\varepsilon(x)\varepsilon(y)$, since $C$ is bilinear in
$(a(x,\cdot),a(y,\cdot))$;
\item $\int a(x',y)^2dx'\leq C$ and
$\int(a(x',y)-a(x',y_t))^2dx'\leq C\varepsilon(y)^2$, because
$a(x',y) = \check w(y,x')-\check{\mathfrak D}(x')$.
\end{enumerate}
The latent part involves only binomial mass functions, so the bounds below are
the binomial counterparts, from Lemma~\ref{lem-c-cancel}, of the roots and
localization estimates used for the edge noise.  By Lemma~\ref{lem-c-cancel}:
\begin{enumerate}[(K1)]
\item $|K_d|\leq C\sigma_t^{-(d+1)}\mu_A$ and
$\int K_d^2\leq C\sigma_t^{-2(d+1)}\mu_A^2h_t$;
\item $\int\varepsilon|K_d|\leq C\sigma_t^{-(d+1)}\mu_Ah_t\bar\ell_A^\vartheta$;
\item $|\int K_d|\leq C\sigma_t^{-(d+1)}\mu_Ah_t\bar\ell_A^\vartheta+\mathfrak e_n$ for $d\geq1$.
\end{enumerate}

\emph{Truncation.}  Choose $d_0$ with $n^{3+K}\rho_n^{d_0/4}\to0$, where
$\Delta_A\geq n^{-K}$.  The components with $|S|>d_0+1$ have conditional
second moment at most $C\sigma_t^{-2}(C\rho_n)^{d_0/2}$: use
Lemma~\ref{lem-binomial-facts}\ref{binom-fact-c} for $|S|\leq n''/2$,
\ref{binom-fact-d} beyond, $\Gamma_{d+1}^2\leq C_0^{d+1}(d+1)!$ and
$\binom{n''}d\leq n^d/d!$.  By Cauchy--Schwarz their total contribution is
$o(\Delta_A)$.

\emph{Class ($\alpha$).}  Split $\Gamma_{d+1}$ into its mixed difference and
three terms involving $x_t$ or $y_t$, and use ($\Gamma$1)--($\Gamma$3), (K2),
(K3).  Since $n\rho_n^2/\sigma_t^2\asymp\rho_n$ and $n\rho_nh_t\asymp\sigma_t$,
the contribution is at most
\begin{align*}
Cn^{d+2}\rho_n^{2d+2}\sigma_t^{-2d-2}\mu_A^2h_t^2\bar\ell_A^{2\vartheta}
\leq C(C\rho_n)^d\,\Xi_A .
\end{align*}
Summed over $d\geq1$ this is at most $C\rho_n\Xi_A$.

\emph{Class ($\beta$).}  Write
\begin{align*}
a(x,y)\Gamma_{d+1}(x,y) &= \bigl[a(x,y)-a(x,y_t)\bigr]\Gamma_{d+1}(x,y)
+a(x,y_t)\bigl[\Gamma_{d+1}(x,y)-\Gamma_{d+1}(x,y_t)\bigr]
\\&\quad+a(x,y_t)\Gamma_{d+1}(x,y_t).
\end{align*}
For the first two terms, bound the $x$-integral by Cauchy--Schwarz with
($\Gamma$4), ($\Gamma$2) and (K1), and the $y$-integral by (K2).  The third
factorizes, and is bounded with (K3) for $K_d$ ($d\geq1$) and Cauchy--Schwarz in
$x$.  The $d$-th term is at most
\begin{align*}
B_d := Cn^{d+2}\rho_n^{2d+3}\sigma_t^{-(2d+3)}\mu_A^2h_t^{3/2}\bar\ell_A^\vartheta,
\qquad B_{d+1}\leq C\rho_nB_d .
\end{align*}
In the small regime, $\mu_A^2h_t\leq\ell_A\leq C\Delta_A$ and $\bar\ell_A = h_t$.
In the large regime, $\mu_A = 1$, $\bar\ell_A = \ell_A\geq h_t$ and
$\Delta_A\geq c\ell_A$.  In both cases
\begin{align*}
\frac{B_1}{\Delta_A}\leq Cn^{1/2}\rho_n^{5/2}h_t^{1/2+\vartheta} = Cv_n,
\end{align*}
which tends to $0$ by \ref{condition-s2}, since
$(1-2\vartheta)/(9-2\vartheta)<(1-2\vartheta)/(5-2\vartheta)<b$.

\emph{Class ($\gamma$).}  By $|a(x,y)a(y,x)|\leq\frac12(a(x,y)^2+a(y,x)^2)$,
($\Gamma$1), ($\Gamma$4) and (K1), the $d$-th term is at most
\begin{align*}
Cn^{d+2}\rho_n^{2d+4}\sigma_t^{-2d-4}\mu_A^2h_t\leq C(C\rho_n)^{d+2}Q_A .
\end{align*}

\emph{Conclusion.}  Collecting,
\begin{align*}
E\Phi_1^2\leq C\bigl[\Xi_A+\rho_nQ_A+\rho_n\Xi_A+v_n\Delta_A\bigr]+\mathfrak e_n,
\end{align*}
and $p'_A\leq C\Delta_A$ gives the bound on $\|Z_A\|_2^2$.  For $X^{(-m)}$,
write its centered part as $(X^{\rm m}-EX^{\rm m})+X^{\rm f}$.  The first piece
is bounded by the variance term of the main part ($C\rho_nQ_A$), and the second
by the fluctuation part.
\end{proof}

\subsection{Step Z4}
\label{app-c-Z4}

\begin{proof}[Proof of Lemma~\ref{lem-c-Z4}]
By \eqref{eq-c-BBLM} and Lemma~\ref{lem-c-Z2}, it suffices to show
$\mathrm{Var}(V_k)\leq C\Delta_A^{4\vartheta'}$ for $k = 1,\ldots,4$; here we used
$E'_m[(X_m-X'_m)^2] = \Phi_m^2+\Psi_m$.  Write $M_4 := \|\Phi_1\|_{L^4}$.  By
Minkowski's inequality in $L^4(d\xi_1)$, \eqref{eq-c-a-window} and
$|g|\leq Cm_A\omega_A$,
\begin{align*}
M_4\leq C\rho_n m_A\,nh_t\leq C\mu_A,\qquad\text{and}\qquad \Delta_A\geq ch_t\mu_A .
\end{align*}

\emph{$V_1,V_2$.}  By Efron--Stein and \eqref{eq-c-exact},
$\mathrm{Var}(V_1)\leq Cn^{-1}(p'_A+\mu_A^2)$.  Here
$n^{-1}p'_A\leq\Delta_A^{4\vartheta'}$, and
\begin{align*}
n^{-1}\mu_A^2 = \frac{Q_A}{nh_t}\leq C\Delta_An^{-1/2}\leq\Delta_A^{4\vartheta'}
\end{align*}
by \eqref{eq-c-absorb}.  The same holds for $V_2$.

\emph{$V_3$: reduction.}  First,
\begin{align*}
\mathrm{Var}(V_3)&\leq n^{-1}E\Phi_1^4+|\mathrm{Cov}(\Phi_1^2,\Phi_2^2)|,
\\
n^{-1}M_4^4&\leq Cn^{-1}\mu_A^4\leq Cn^{-1}h_t^{-4\vartheta'}\Delta_A^{4\vartheta'}
\leq C\rho_n^{2\vartheta'}\Delta_A^{4\vartheta'},
\end{align*}
since $2\vartheta'\leq1$.  Let $g^{(1,2)}_i$ be computed with vertices $1,2$
removed, and
\begin{align*}
\Phi^*_l := \rho_n\sum_{i\geq3}\chi_ia(\xi_i,\xi_l)g^{(1,2)}_i,\qquad l = 1,2 .
\end{align*}
Given $\xi_3,\ldots,\xi_n$, $\Phi^*_1$ and $\Phi^*_2$ are i.i.d., so
\begin{align*}
\mathrm{Cov}(\Phi_1^{*2},\Phi_2^{*2}) = \mathrm{Var}(\Psi^*),\qquad
\Psi^* := \|X\|^2,\quad X := \rho_n\sum_{i\geq3}\chi_ia(\xi_i,\cdot)g^{(1,2)}_i .
\end{align*}
The remainder is
\begin{align*}
r_1 := \Phi_1-\Phi^*_1 = \underbrace{\rho_n\chi_2a(\xi_2,\xi_1)g^{(1)}_2}_{r^{(1)}}
+\underbrace{\rho_n\sum_{i\geq3}\chi_ia(\xi_i,\xi_1)p_{i2}\delta_i}_{r^{(2)}},
\end{align*}
where $\delta_i$ is the slope of $g^{(1)}_i$ in $p_{i2}$, with
$|\delta_i|\leq Cm_A\sigma_t^{-1}\omega^{[2]}_A(\xi_i)$.  Using $|a|\leq\rho_n^{-1}$
for $r^{(1)}$, and Minkowski's inequality in $L^4(d\xi_1d\xi_2)$ with
$\|p_{i2}\|_{L^4(d\xi_2)}\leq\min(C\rho_n,1)$ for $r^{(2)}$,
\begin{align*}
\|r^{(1)}\|_4\leq C\rho_n^{1/2}m_Ah_t^{1/4},\qquad
\|r^{(2)}\|_4\leq C \rho_n \,\mu_Ah_t .
\end{align*}
Since
$|\mathrm{Cov}(\Phi_1^2,\Phi_2^2)-\mathrm{Var}(\Psi^*)|\leq C(\|r_1\|_4M_4^3+\|r_1\|_4^4)$,
and $\Delta_A^{4\vartheta'}\geq c(h_t\mu_A)^{4\vartheta'}$, the four remainder
terms satisfy
\begin{align*}
\|r^{(1)}\|_4M_4^3&\leq C\rho_n^{1/2}h_t^{5/4-4\vartheta'}\Delta_A^{4\vartheta'},&
\|r^{(2)}\|_4M_4^3&\leq Cn^{(4\vartheta'-1)/2}\rho_n^{(4\vartheta'+1)/2}\Delta_A^{4\vartheta'},\\
\|r^{(1)}\|_4^4&\leq C\rho_n^2m_A^4h_t\leq Cn^{-2}\Delta_A,&
\|r^{(2)}\|_4^4&\leq Ch_t^{4-4\vartheta'}\Delta_A^{4\vartheta'} .
\end{align*}
All are $o(\Delta_A^{4\vartheta'})$, by \eqref{eq-c-theta-prime}(ii), (iii) and
\eqref{eq-c-absorb}.  It remains to bound $\mathrm{Var}(\Psi^*)$.  The same
reduction gives
\begin{align*}
\mathrm{Var}(V_4)\leq n^{-1}E\Phi_1^4+\mathrm{Var}(\Psi^*)+o(\Delta_A^{4\vartheta'}).
\end{align*}

\emph{$\mathrm{Var}(\Psi^*)$: exact Efron--Stein decomposition.}  For $m\geq3$,
let $X^{(-m)}$ be $X$ with vertex $m$ removed (from the sum and from every
$g$).  By multilinearity,
\begin{align*}
X &= X^{(-m)}+Y(\xi_m)+V(\xi_m),\qquad
Y(y) := \rho_n\chi(y)a(y,\cdot)g_m(y),
\\
V(y) &:= \rho_n^2\sum_{i\neq m}\chi_ia(\xi_i,\cdot)\check w(\xi_i,y)\delta_i,
\end{align*}
where $\delta_i$ is now the slope in $p_{im}$, and $X^{(-m)}$, $\delta_i$ do
not involve $\xi_m$.  The Efron--Stein inequality gives
\begin{align*}
\mathrm{Var}(\Psi^*)\leq\sum_m\Bigl(&16E\bigl[E_y\langle X^{(-m)},Y(y)\rangle^2\bigr]
 +16E\bigl[\mathrm{Var}_y\langle X^{(-m)},V(y)\rangle\bigr]
\\ &+4E\|Y(\xi_m)+V(\xi_m)\|^4\Bigr) .
\end{align*}

\emph{Direct part.}  $E_y\langle X^{(-m)},Y(y)\rangle^2\leq C\rho_n^2m_A^2h_t\|X^{(-m)}\|^2$,
so by Lemma~\ref{lem-c-Z2} it contributes at most
\begin{align*}
Cn\rho_n^2m_A^2h_tE\|X^{(-m)}\|^2\leq C\rho_nQ_A(\Delta_A+\Xi_A).
\end{align*}

\emph{Indirect part: an exact identity.}  Since
$\check w(x,y) = \check{\mathfrak D}(x)+a(x,y)$, \eqref{eq-c-T-identity} holds.
With $\bar{\mathcal T} := E\mathcal T_m$, $\tilde\varphi := EX^{(-m)}$ and
$\tilde X := X^{(-m)}-\tilde\varphi$,
\begin{align}
\label{eq-app-T-split}
\|\mathcal T_mX^{(-m)}\|^2\leq3\|\bar{\mathcal T}X^{(-m)}\|^2
+3\|(\mathcal T_m-\bar{\mathcal T})\tilde\varphi\|^2+3\|(\mathcal T_m-\bar{\mathcal T})\tilde X\|^2 .
\end{align}

\emph{Mean operator.}  We have
$\bar{\mathcal T} = (n-3)\int_{[0,1]\setminus\mathscr H}\bar\delta\,a(x)\otimes a(x)\,dx$
with $\bar\delta = \mathsf K^{(1)}_A$ and $M = n-4$.  Split
$a(x)\otimes a(x) = a_t\otimes a_t+[a(x)\otimes a(x)-a_t\otimes a_t]$.  The
first part is bounded by Lemma~\ref{lem-c-cancel}(c) with $d = 1$, and the
second by \eqref{eq-c-a-window} and Lemma~\ref{lem-c-cancel}(b).  This gives
\begin{align*}
\|\bar{\mathcal T}\|_{\rm op}\leq\frac{C\mu_A\bar\ell_A^\vartheta}{\rho_n\sigma_t},
\qquad\text{hence}\qquad
n\rho_n^4\|\bar{\mathcal T}\|_{\rm op}^2\leq C\rho_n\Xi_A,
\end{align*}
and the first term of \eqref{eq-app-T-split} contributes at most
$C\rho_n\Xi_A(\Delta_A+\Xi_A)$.

\emph{Hilbert--Schmidt bound, no cancellation.}  Split
$\mathcal T_m-\bar{\mathcal T}$ into the centered i.i.d.\ sum
$\sum_i\chi_i\bar\delta(\xi_i)a_i\otimes a_i-\bar{\mathcal T}$ and
$\mathcal T^{\rm f} := \sum_i\chi_i(\delta_i-\bar\delta(\xi_i))a_i\otimes a_i$.
The off-diagonal part of
$E\|\mathcal T^{\rm f}\|_{\rm HS}^2 = \sum_{i,k}E[\cdots C(\xi_i,\xi_k)^2]$ has
the three classes of Lemma~\ref{lem-c-Z2}, with $C^2$ in place of $C$, kernels
of one higher order, and prefactor $n\rho_n^4$ in place of $n^2\rho_n^2$.
Bounding them as there, without cancellation,
\begin{align*}
n\rho_n^4E\|\mathcal T_m-\bar{\mathcal T}\|_{\rm HS}^2
 \leq C\rho_n^2\bigl(Q_A+\mu_A^2\bigr)+\mathfrak e_n .
\end{align*}
Class ($\alpha$) gives $C\rho_n^{d+1}\mu_A^2$, class ($\beta$) gives
$C\rho_n^2\mu_A^2(\rho_n^3/n)^{1/4}$, class ($\gamma$) gives
$C\rho_n^{d+3}Q_A$, and the i.i.d.\ part and the diagonal give
$C\rho_n^2Q_A$.  Since $\|\tilde\varphi\|^2\leq C\Xi_A$, by
Proposition~\ref{prop-c-hajek} with $n-3$ trials, the second term of
\eqref{eq-app-T-split} contributes at most $C\Xi_A\rho_n^2(Q_A+\mu_A^2)$.

\emph{Operator bound against $\tilde X$.}  Let
$N := \sum_{i\neq m}\chi_i\omega^{[4]}_A(\xi_i)$.  Then $EN\leq C_2nh_t$ and
$nh_t\geq n^{1/2}$, so by Bernstein's inequality
$\mathscr G := \{N\leq2C_2nh_t\}$ has negligible complement.  On $\mathscr G$,
\begin{align*}
\|\mathcal T_m\|_{\rm op}+\|\bar{\mathcal T}\|_{\rm op}\leq Cm_A\sigma_t^{-1}nh_t,
\qquad\text{so}\qquad
n\rho_n^4\|\mathcal T_m-\bar{\mathcal T}\|_{\rm op}^2\leq C\rho_n\mu_A^2 .
\end{align*}
With Lemma~\ref{lem-c-Z2}, the third term of \eqref{eq-app-T-split} contributes
at most $C\rho_n\mu_A^2[\rho_nQ_A+\rho_n\Xi_A+v_n\Delta_A]+\mathfrak e_n$.

\emph{The product $P_A := \rho_n^2\mu_A^2\Xi_A$.}  This appears in both
fluctuation terms.  By \eqref{eq-c-Xi-Delta} and $\mu_A\leq1$,
\begin{align*}
\frac{P_A}{\Delta_A^{4\vartheta'}}\leq C\rho_n^2h_t^{2\vartheta-4\vartheta'}
\leq Cn^{2\vartheta'-\vartheta}\rho_n^{2+2\vartheta'-\vartheta}\leq C
\end{align*}
by \eqref{eq-c-theta-prime}(iv), which is \eqref{eq-c-PA}.

\emph{Quartic part.}  $E\|Y(\xi_m)\|^4\leq C\rho_n^4m_A^4h_t$, and
\begin{align*}
\|V(y)\|\leq C\rho_n^2m_A\sigma_t^{-1}S(y),\qquad
S(y) := \sum_i\chi_i\check w(\xi_i,y)\omega^{[3]}_A(\xi_i).
\end{align*}
By Rosenthal's inequality \citep{rosenthal1970subspaces}, Minkowski's inequality with \ref{condition-l2}, and
$\check w\leq\rho_n^{-1}$, $ES(\xi_m)^4\leq C(n^4h_t^4+n\rho_n^{-2}h_t)$.  Hence
\begin{align*}
\sum_m4E\|Y+V\|^4\leq C(\rho_n^2/n+\rho_n^3)\Delta_A .
\end{align*}

\emph{Conclusion.}  Collecting,
\begin{align*}
\mathrm{Var}(\Psi^*)\leq C\rho_n(\Delta_A+\Xi_A)^2+C\rho_n\Delta_A+CP_A+\mathfrak e_n
\leq C\Delta_A^{4\vartheta'} .
\end{align*}
This uses $(\Delta_A+\Xi_A)^2\leq C\Delta_A^{4\vartheta}$, \eqref{eq-c-absorb}
with $\varepsilon = \rho_n$, and \eqref{eq-c-PA}.  Therefore
$\|V\|_2\leq EV+\sum_k\mathrm{Var}(V_k)^{1/2}\leq C\Delta_A^{2\vartheta'}$.
\end{proof}

\section{Deferred proofs for Section~\ref{sec-tight-compact}}
\label{app-compact-proofs}

This appendix collects the proofs of Section~\ref{sec-tight-compact-proof} that consist of routine, though careful, estimates.  The
proofs of the latent moment bounds, Lemmas~\ref{lem-c-Z2} and
\ref{lem-c-Z4}, are in Appendix~\ref{app-latent}.  The subsections follow the steps of the main text; notation and
conventions are those of Section~\ref{sec-tight-compact-proof}.

\subsection{Step H}
\label{app-c-H}

\begin{proof}[Proof of Lemma~\ref{lem-c-hubs}]
For fixed $R$, Markov's inequality and \ref{condition-l1} give
\begin{align*}
g(x)\leq\rho_n^{1+\beta}E\bigl[w^{2+\beta}(x,\xi)\bigr]\leq C_R\rho_n^{1+\beta}
\qquad\text{on }\{\mathfrak D\leq R\},
\end{align*}
which is below $c/4\leq\mathfrak D/4$ by \ref{condition-g1} for $n$ large.
With $R = 2t_1^+$ this gives the last two claims, and with $R = d_1$ the claim
for bounded support.  If the support is unbounded and
$\mathfrak D(x)\geq t_0^-$, \ref{condition-l2} gives
\begin{align*}
g(x)\leq K\rho_n^{1+\beta}\mathfrak D(x)^{2+\beta},
\end{align*}
so $x\in\mathscr H$ forces $\mathfrak D(x)\geq c_H\rho_n^{-1}$, and Markov's
inequality with \ref{condition-m1} gives $\mathrm{Leb}(\mathscr H)\leq C\rho_n^{2+\beta}$.
\end{proof}

\subsection{Step L}
\label{app-c-L}

\begin{proof}[Proof of Proposition~\ref{prop-grid}]
Reparametrize by $u = F(t)\in[F(t_0),1]$.  Since $f>0$ on $J$
(respectively $J_\infty$), $F$ is continuous and strictly increasing there, so
this is a homeomorphism onto its image, a subinterval of $[0,1]$, with
$t = \infty$ corresponding to $u = 1$; the results of
\citet[Section~10 and Theorem~13.5]{billingsley1999convergence}, stated for
$[0,1]$, apply on any compact subinterval after an affine change of
variables, which changes the constant in the moment condition but not its
form.  The process is then constant
between the images $u_m := F(\lambda_m)$ of lattice points, and (ii) reads
\begin{align*}
E\bigl[|\mathbb H_n(u_j)-\mathbb H_n(u_i)|^2|\mathbb H_n(u_k)-\mathbb H_n(u_j)|^2\bigr]
\leq C(u_k-u_i)^{2\nu}
\end{align*}
for grid indices $i\leq j\leq k$ with $u_k\leq F(T_n)$.

\emph{Billingsley's fluctuation inequality.}  Let $S_1,\ldots,S_M$ be random
variables and $w_1,\ldots,w_M$ nonnegative weights such that
\begin{align*}
E\bigl[|S_j-S_i|^2|S_k-S_j|^2\bigr]\leq\Bigl(\sum_{i<l\leq k}w_l\Bigr)^{2\nu}
\qquad\text{for all }i\leq j\leq k .
\end{align*}
Then, by \citet[Section~10]{billingsley1999convergence}, for every $\lambda>0$,
\begin{align*}
P\Bigl(\max_{i\leq k\leq j}\min\bigl(|S_k-S_i|,|S_j-S_k|\bigr)\geq\lambda\Bigr)
\leq K_\nu\lambda^{-4}\Bigl(\sum_{i<l\leq j}w_l\Bigr)^{2\nu},
\end{align*}
with $K_\nu$ depending only on $\nu$.

\emph{Below the cutoff.}  Apply this to $S_m := \mathbb H_n(u_m)$ with
$w_m := C^{1/(2\nu)}(u_m-u_{m-1})$.  Following the proof of
\citet[Theorem~13.5]{billingsley1999convergence}, for every $\delta$ larger
than the grid mesh,
\begin{align*}
P\bigl(w''_{[F(t_0),F(T_n)]}(\mathbb H_n,\delta)\geq\epsilon\bigr)
\leq K\epsilon^{-4}\delta^{2\nu-1},
\end{align*}
where $w''$ is Billingsley's modulus (13.6).  Since $\mathbb H_n$ is constant
between grid points, no additional increments occur, and the mesh tends to
$0$ by assumption.

\emph{Beyond the cutoff.}  For a triple $u_1\leq u_2\leq u_3$ with
$u_2\leq F(T_n)<u_3$ and $u_3-u_1\leq\delta$,
$\min(|\mathbb H_n(u_2)-\mathbb H_n(u_1)|,|\mathbb H_n(u_3)-\mathbb H_n(u_2)|)$
is at most the same quantity for the triple $(u_1,u_2,F(T_n))$ plus
$\sup_{u\geq F(T_n)}|\mathbb H_n(u)-\mathbb H_n(F(T_n))|$, and triples with
$u_1>F(T_n)$ are bounded by twice that supremum.  Hence
\begin{align*}
w''_{[F(t_0),1]}(\mathbb H_n,\delta)\leq w''_{[F(t_0),F(T_n)]}(\mathbb H_n,\delta)
+2\sup_{u\geq F(T_n)}\bigl|\mathbb H_n(u)-\mathbb H_n(F(T_n))\bigr|,
\end{align*}
and the second term is $o_P(1)$ by (iii).

\emph{Conclusion.}  Together with (i) and continuity of $\mathbb H$, the
conditions of \citet[Theorem~13.3]{billingsley1999convergence} hold, so
$\mathbb H_n\rightsquigarrow\mathbb H$ in $D$; as the limit is continuous,
convergence holds in the uniform topology.
\end{proof}

\subsection{Step W}
\label{app-c-W}

\begin{proof}[Proof of Lemma~\ref{lem-poisson-binomial}]
The generating function of $\mathsf D^{d-1}P(X=\cdot)$ is the polynomial
$(1-z)^{d-1}\prod_{j\in J'}(1-p_j+p_jz)$, and
\begin{align*}
|1-p_j+p_je^{iu}|^2\leq\exp\bigl(-4p_j(1-p_j)u^2/\pi^2\bigr)\qquad\text{on }[-\pi,\pi],
\end{align*}
so the proof of Lemma~\ref{lem-binomial-facts}\ref{binom-fact-c} applies with
$\sigma^2$ replaced by $\Sigma_{J'}$.  Adding $Y$ convolves the mass function
with a probability measure, which commutes with $\mathsf D$ and does not
increase the supremum norm.
\end{proof}

\begin{proof}[Proof of Lemma~\ref{lem-pb-local}]
\emph{Tilting.}  Let $\Lambda(\lambda) := \log Ee^{\lambda X}$ and let $\lambda$
solve $\Lambda'(\lambda) = m$.  Since $\Lambda''\asymp\Sigma$ on $|\lambda|\leq1$,
\begin{align*}
|\lambda|\leq C|m-\mu|/\Sigma\leq C|x_m|\Sigma^{-1/2}.
\end{align*}
Under the tilted law $P_\lambda(X = j) := e^{\lambda j-\Lambda(\lambda)}P(X = j)$,
$X$ is again Poisson-binomial, with parameters
$p_je^\lambda/(1-p_j+p_je^\lambda)$, whose variances lie within a factor
$e^{\pm2|\lambda|}$ of $p_j(1-p_j)$.  Hence its variance $\Sigma_\lambda$
satisfies $\Sigma/2\leq\Sigma_\lambda\leq2\Sigma$ for $c_0$ small.

\emph{Differences under the tilt.}  With $S$ the shift $(Su)(j) := u(j-1)$,
\begin{align*}
\bigl((I-S)P\bigr)(m) = e^{\Lambda(\lambda)-\lambda m}
\bigl((I-e^{\lambda}S)P_\lambda\bigr)(m),
\end{align*}
and since $I-e^\lambda S = (I-S)+(1-e^\lambda)S$ and the two operators commute,
\begin{align*}
\mathsf D^kP(X = m) = e^{\Lambda(\lambda)-\lambda m}
 \sum_{l=0}^{k}\binom kl(1-e^{\lambda})^{l}
 \bigl((I-S)^{k-l}S^lP_\lambda\bigr)(m).
\end{align*}
By Lemma~\ref{lem-poisson-binomial} applied to $P_\lambda$, the $l$-th term is
at most
\begin{align*}
\binom kl(C|\lambda|)^l\Gamma_{k-l+1}\Sigma_\lambda^{-(k-l+1)/2}
\leq\binom kl C^l\Gamma_{k-l+1}|x_m|^l\Sigma^{-(k+1)/2}2^{(k+1)/2}.
\end{align*}

\emph{The Gaussian factor.}  Finally
$\Lambda(\lambda)-\lambda m = -\Lambda^*(m)\leq-\zeta x_m^2$ for
$|m-\mu|\leq c_0\Sigma$, by the quadratic lower bound on the Cram\'er rate
function of a sum of independent Bernoulli variables in the moderate-deviation
range.
\end{proof}

\begin{proof}[Proof of Lemma~\ref{lem-c-loc}]
Fix $s\in S_A$.  On $W_s$, by \eqref{eq-c-clip},
\begin{align*}
|\mathfrak D-s|\leq h_s|x_s|+g_*\leq h_t(1+|x_s|),
\qquad\text{so}\qquad
e_A\leq3\bar\ell_A(1+|x_s|).
\end{align*}
Split $W_s$ into $\{|\check{\mathfrak D}-s|\leq h_s\}$ and the shells
\begin{align*}
\bigl\{2^{k-1}h_s<|\check{\mathfrak D}-s|\leq2^kh_s\bigr\},\qquad k\geq1 .
\end{align*}
By \eqref{eq-c-window-measure}, the $k$-th shell has measure at most $C2^kh_t$,
and on it the integrand is at most
$(3\bar\ell_A)^\theta(1+2^k)^{j+1}e^{-\zeta_*4^{k-1}}$; summing over $k$ gives
$Ch_t\bar\ell_A^\theta$.

A term of $\omega^{[j]}_A$ indexed by $s_1$ at a point outside
$W_{s_1}\cup\mathscr H$ has $|x_{s_1}|\geq c_w\sigma_{s_1}$, so it is at most
$Ce^{-cn\rho_n}$.  With $e_A\leq1+\mathfrak D+t_1$ its integral is negligible,
hence at most $h_t\bar\ell_A^\theta$.
\end{proof}

\begin{proof}[Proof of Lemma~\ref{lem-c-roots}]
(a) Fix $x\notin\mathscr H$ and $I'$, and condition on $\xi_i = x$.  The
$p_{il} = \rho_n\check w(x,\xi_l)$, $l\notin I'$, are i.i.d.\ in $[0,1]$ with
mean $\rho_n\check{\mathfrak D}(x)$, so
\begin{align*}
\bigl|E\mu_{i,I'}-n''\rho_n\check{\mathfrak D}(x)\bigr|\leq2\rho_n\check{\mathfrak D}(x),
\qquad
\mathrm{Var}(\mu_{i,I'})\leq n\rho_n\check{\mathfrak D}(x),
\end{align*}
as $\check w\leq\rho_n^{-1}$.  Bernstein's inequality gives (E0) up to
probability $2e^{-cn\rho_nt_0}$.

For (E1), $|\check{\mathfrak D}(x)-s|\leq c_ws$ gives $x\in W_s$, so
$E[\check w(x,\xi)^2]\leq C$ by \eqref{eq-c-a-window}, and
$\mathrm{Var}(\mu_{i,I'})\leq Cn\rho_n^2$.  Bernstein's inequality at deviation
$\sigma_s/2$ gives probability at most
\begin{align*}
2\exp\bigl(-c\min\{\rho_n^{-1},\sigma_s\}\bigr),
\end{align*}
which is negligible.  For the second part of (E1), write
\begin{align*}
\mu_{i,I'}-\Sigma_{i,I'} = \sum_{l\notin I'}p_{il}^2\leq S+\mu_{i,I'}/2,
\qquad S := \#\{l\notin I':p_{il}>1/2\},
\end{align*}
where $ES\leq n(2\rho_n)^{2+\beta}E[w^{2+\beta}(x,\xi)]\leq Cn\rho_n^{2+\beta}$.
With $a := n''\rho_ns/16$ the Chernoff bound gives $P(S\geq a)\leq e^{-a}$.  On
the complementary event and the deviation event,
\begin{align*}
\mu_{i,I'}\geq(1-c_w)n''\rho_ns-\sigma_s\geq\tfrac12n''\rho_ns,
\end{align*}
so $S<\mu_{i,I'}/8$ and $\Sigma_{i,I'}\geq\mu_{i,I'}/4$.  A union bound over
at most $4n^3$ sets $I'$ and $n$ lattice points proves (a).

(b) Fix $I'$ and write $\mu,\Sigma$; given $\mathcal G$, $D_{i\setminus I'}$ is
Poisson-binomial with these parameters.  Let $s\in S_A$ and
$|m-s'|\leq\frac{c_w}4n''\rho_ns$; this covers every point used below.

\emph{Case 1: $|\check{\mathfrak D}(x)-s|\leq c_ws$.}  By (E1),
$\mu\in[\frac12n''\rho_ns,2n''\rho_ns]$, $\Sigma\in[\sigma_s^2/8,2\sigma_s^2]$, and
\begin{align*}
|m-\mu|&\leq|m-s'|+|s'-n''\rho_ns|+n''\rho_n|s-\check{\mathfrak D}(x)|+\sigma_s
\\&\leq\tfrac{c_w}4n''\rho_ns+2+c_wn''\rho_ns+\sigma_s
 \leq2c_wn''\rho_ns\leq\tfrac{c_0}8n''\rho_ns\leq c_0\Sigma .
\end{align*}
So Lemma~\ref{lem-pb-local} gives, with $y_m := (m-\mu)/\Sigma^{1/2}$,
\begin{align*}
|\mathsf D^kP_{I'}(m)|\leq C\Sigma^{-(k+1)/2}(1+|y_m|)^ke^{-\zeta y_m^2}.
\end{align*}
Since
\begin{align*}
m-\mu = (m-s')+(s'-n''\rho_ns)-\sigma_sz_s(x)+\bigl(n''\rho_n\check{\mathfrak D}(x)-\mu\bigr),
\end{align*}
for $|m-s'|\leq L\sigma_s$ we get $|y_m+\tau z_s(x)|\leq C_L$ with
$\tau := \sigma_s/\Sigma^{1/2}\in[2^{-1/2},\sqrt8]$, hence
$y_m^2\geq x_s^2/4-C_L^2$ and
\begin{align*}
|\mathsf D^kP_{I'}(m)|\leq C\sigma_s^{-(k+1)}(1+|x_s|)^ke^{-\zeta_*x_s^2}.
\end{align*}

\emph{Case 2: $|\check{\mathfrak D}(x)-s|>c_ws$.}  By (E0),
$|m-\mu|\geq c\max\{\mu,n''\rho_ns\}$ for $n$ large, and the Chernoff bound for
the Poisson-binomial law makes every $P_{I'}(m-j')$, $j'\leq k$, negligible.

\emph{Assembly.}  In the small regime all $m+j$ with $m\in A$ lie within
$2\sigma_t$ of $t'$; summing Case~1 over the $k_A$ points gives
\begin{align*}
Ck_A\sigma_t^{-(k+1)}\omega^{[k]}_A = Cm_A\sigma_t^{-(k-1)}\omega^{[k]}_A .
\end{align*}
In the large regime, telescope to
$\mathsf D^{k-1}P_{I'}(t'+j)-\mathsf D^{k-1}P_{I'}(r'+j)$ and bound each end at
its own endpoint, using $\sigma_r\geq\sigma_t/\sqrt2$ and $m_A = \sigma_t^{-1}$.
\end{proof}

\begin{proof}[Proof of Lemma~\ref{lem-c-cancel}]
(a) Let $\lambda>0$ and $p := \rho_n\lambda$.  If $|\lambda-s|\leq c_ws$ for
some $s\in S_A$, the argument of Lemma~\ref{lem-c-roots}(b), Case~1, applied
to $b_M(\cdot\,;p)$ (mean $Mp$, variance $\asymp\sigma_s^2$, no exceptional
event) and the same assembly give the bound with $z_s(\lambda) := (\lambda-s)/h_s$.
Otherwise the binomial Chernoff bounds make every relevant mass negligible.
Take $\lambda = \check{\mathfrak D}(x)$.

(b) By (a),
$|\mathsf K^{(d)}_A|^j\leq C\sigma_t^{-j(d+1)}\mu_A^j\omega^{[j(d+1)]}_A+\mathfrak e_n$;
apply Lemma~\ref{lem-c-loc}.

(c) Write $\mathsf K := \mathsf K^{(d)}_{A,M}$.

\emph{Step 1: localization.}  Off $\{\mathfrak D\leq2t\}$ and off $\mathscr H$,
$\check{\mathfrak D}\geq\frac34\mathfrak D>\frac32t$, so the integrand is
negligible there by (a).

\emph{Step 2: clipping.}  On $\{\mathfrak D\leq2t\}$,
$0\leq\mathfrak D-\check{\mathfrak D}\leq g_*$.  Since
$\partial_pb_M(m;p) = -M\,\mathsf Db_{M-1}(m;p)$ and $\mathsf D$ commutes with
$\partial_p$,
\begin{align*}
\bigl|\partial_p\mathsf K^{(d)}_{A,M}\bigr| = M\bigl|\mathsf K^{(d+1)}_{A,M-1}\bigr| .
\end{align*}
Apply the mean value theorem and (a) at a point $\lambda^*(x)$ between
$\check{\mathfrak D}(x)$ and $\mathfrak D(x)$, with
$|z_s(\lambda^*)-z_s(x)|\leq g_*/h_s\leq1$.  By Lemma~\ref{lem-c-loc} with
decay $\zeta_*/2$, $n\rho_nh_t\asymp\sigma_t$, $\rho_nM\asymp\sigma_t^2$ and
$g_*\leq C\Upsilon_n h_t$,
\begin{align*}
n'\rho_n\int_{\{\mathfrak D\leq2t\}}
 \bigl|\mathsf K(\rho_n\check{\mathfrak D})-\mathsf K(\rho_n\mathfrak D)\bigr|
 &\leq Cn\rho_n\cdot\rho_ng_*M\cdot\sigma_t^{-(d+2)}\mu_Ah_t+\mathfrak e_n
 \\&\leq C\sigma_t^{-d}\mu_A\Upsilon_n+\mathfrak e_n .
\end{align*}

\emph{Step 3: main term.}  By Remark~\ref{remark-inverse} and
\ref{condition-d1},
\begin{align*}
\int_{\{\mathfrak D\leq2t\}}\mathsf K(\rho_n\mathfrak D(x))\,dx =
\int_0^{2t}\mathsf K(\rho_n\lambda)f(\lambda)\,d\lambda .
\end{align*}
Let $U := \bigcup_{s\in S_A}[(1-c_w)s,(1+c_w)s]\subseteq J$; off $U$ the
integrand is negligible.  On $U$, \ref{condition-d2} gives
$|f(\lambda)-f(t)|\leq Ce(\lambda)^{\theta_f}$ with
$e(\lambda) := \min_s(|\lambda-s|+|s-t|)$, and the shell argument of
Lemma~\ref{lem-c-loc} in the variable $\lambda$ gives
\begin{align*}
n'\rho_n\Bigl|\int_U\mathsf K(\rho_n\lambda)\bigl[f(\lambda)-f(t)\bigr]d\lambda\Bigr|\leq
C\sigma_t^{-d}\mu_A\bar\ell_A^{\theta_f}.
\end{align*}
Finally, up to a negligible term,
\begin{align*}
\int_U\mathsf K(\rho_n\lambda)\,d\lambda = \frac1{\rho_n}\int_0^1\mathsf K(p)\,dp = 0 :
\end{align*}
for $d\geq1$ because $\int_0^1b_M(m;p)\,dp = 1/(M+1)$ does not depend on $m$
(Lemma~\ref{lem-binomial-facts}\ref{binom-fact-e}), and for $d = 0$ because
the two endpoint terms cancel.  The second inequality in (c) uses
$\theta_f\geq\vartheta$, $\bar\ell_A\leq C$ and \eqref{eq-c-clip}.
\end{proof}

\subsection{Step B}
\label{app-c-B}

\begin{proof}[Proof of Lemma~\ref{lem-c-sigma}]
\emph{First moment.}  The factor $g^{(2)}_1$ involves neither $\xi_2$ nor
$\eta_{12}$.  Bounding $|g^{(1)}_2|\leq Cm_A$ and integrating $\xi_2$ first,
by Lemma~\ref{lem-c-loc},
\begin{align*}
E\bar\Sigma(A)\leq Cn\rho_nm_A\,E\bigl[\chi_1\check{\mathfrak D}(\xi_1)|g^{(2)}_1|\bigr]
\leq Cn\rho_nm_A^2h_t = C\frac{n\rho_n}{\sigma_t^2}\mu_A^2h_t\leq CQ_A .
\end{align*}

\emph{Variance.}  By the Efron--Stein inequality,
\begin{align*}
\mathrm{Var}(\bar\Sigma)\leq\frac12\sum_mE\bigl[(\bar\Sigma-\bar\Sigma^{(m)})^2\bigr],
\end{align*}
where $\bar\Sigma^{(m)}$ has $\xi_m$ replaced by an independent copy.  Write
$c_{ij} := \chi_i\chi_jp_{ij}|g^{(j)}_i||g^{(i)}_j|$ for the summands, which
are symmetric in $(i,j)$.

\emph{Direct effect} (terms with $m\in\{i,j\}$).  Each contains the factor
$|g^{(j)}_m|\leq Cm_A\omega_A(\xi_m)$, and $\sum_jp_{mj}\leq Cn\rho_n$ on
$\mathscr E_m$, so the change is at most $C\rho_nm_A^2\omega_A(\xi_m)$, plus the
same for $\xi'_m$.  Its contribution to the variance is at most
\begin{align*}
Cn\rho_n^2m_A^4h_t = C\mu_A^4h_t/n\leq CQ_A/n\leq C\Delta_A^2,
\end{align*}
using $\Delta_A\geq n^{-1}$ from \eqref{eq-c-Q-Xi}.

\emph{Indirect effect} (through $p_{im}$ in $g^{(j)}_i$, $i,j\neq m$).  The
function $g^{(j)}_i$ is affine in $p_{im}$ with slope a second difference
$\delta$, and $|\delta|\leq Cm_A\sigma_t^{-1}\omega^{[2]}_A(\xi_i)$ by
Lemma~\ref{lem-c-roots}(b).  So the change is at most
\begin{align*}
C\rho_n^2m_A^2\sigma_t^{-1}S_m,\qquad
S_m := \sum_i(\check w_{im}+\check w'_{im})\omega^{[2]}_A(\xi_i).
\end{align*}
Given $\xi_m$, $S_m$ has mean at most $2n\int\check w(x,\xi_m)\omega^{[2]}_A(x)\,dx$
and variance at most $nE[\check w^2\omega^{[2]2}_A]$.  By Minkowski's inequality
and \eqref{eq-c-a-window}, $ES_m^2\leq Cn^2h_t^2+Cnh_t\leq Cn^2h_t^2$.  The
contribution is at most
\begin{align*}
Cn\cdot\rho_n^4m_A^4\sigma_t^{-2}n^2h_t^2\leq C\rho_n\mu_A^4h_t^2 = C\rho_nQ_A^2
\leq C\rho_n\Delta_A^2 .
\end{align*}
\end{proof}

\begin{proof}[Proof of Lemma~\ref{lem-c-block}]
\emph{Multilinear expansion.}  Given $\mathcal G$, $P(\pi_v(s_v))$ is a
function of the independent Bernoulli variables $(A_e)_{e\in E_v}$, hence
multilinear in them.  Expanding about $(p_e)$ as in the proof of
Lemma~\ref{lem-multilinear} gives the representation with
$c_T = E[\partial_TP(\pi_v(s_v))\mid\mathcal G]$, where $\partial_T$ is the
$|T|$-fold difference in the coordinates in $T$, with the other coordinates
Bernoulli$(p_e)$.  Since the function depends on $(A_e)$ only through $s_v$,
$\partial_T$ is the forward difference $\Delta^{|T|}$ in $s_v$, so
\begin{align*}
|c_T|\leq\max_{0\leq s\leq3}\bigl|\Delta^{|T|}[P\circ\pi_v](s)\bigr| .
\end{align*}

\emph{Differences of $\pi_v$.}  As
$\Delta_sP(D^{\rm ext}_v = m-s) = -\mathsf DP(D^{\rm ext}_v = \cdot)(m-s)$,
\begin{align*}
\Delta^k\pi_v(s) = (-1)^k\sum_{m\in A}\mathsf D^kP\bigl(D^{\rm ext}_v = m-s\bigm|\mathcal G\bigr),
\end{align*}
which by Lemma~\ref{lem-c-roots}(b), with $I' = V$ and $j = -s$, is at most
$Cm_A\sigma_t^{-(k-1)}\omega^{[k]}_A(\xi_v)+\mathfrak e_n$.

\emph{General $P$.}  Iterating
$\Delta(uw)(s) = \Delta u(s)\,w(s+1)+u(s)\,\Delta w(s)$ writes
$\Delta^k[P\circ\pi_v]$ as a bounded combination of products
$\prod_{l\leq\ell}\Delta^{k_l}\pi_v$ (at shifted arguments), with $\ell\geq1$,
$k_l\geq1$ and $\sum k_l = k$, times factors in $[0,1]$.  Each product is at
most
\begin{align*}
Cm_A^\ell\sigma_t^{-(k-\ell)}\omega^{[k]}_A\leq Cm_A\sigma_t^{-(k-1)}\omega^{[k]}_A .
\end{align*}
This uses $m_A\sigma_t = \mu_A\leq1$ and $\prod_l\omega^{[k_l]}_A\leq C\omega^{[k]}_A$:
a cross term between the two endpoints is bounded by the term at either one,
since $(1+|y|)^je^{-\zeta_*y^2}\leq C$.  Finally, if $P(\pi) = \pi-q_v(A)$,
then $c_\emptyset = E[\pi_v(s_v)\mid\mathcal G]-q_v(A) = 0$.
\end{proof}

\begin{proof}[Proof of Lemma~\ref{lem-c-tree}]
For distinct indices the positions are i.i.d., so the left side is at most
$n^k$ times the integral over $[0,1]^k$.  Root $F$ at $1$ and define
recursively, for a vertex $u$ with children $c$,
\begin{align*}
\Phi_u(x) := \chi(x)\omega^{[d_u]}_A(x)\prod_c\int_0^1\rho_n\check w(x,y)\Phi_c(y)\,dy .
\end{align*}
We claim that, for the subtree $S_u$ of $u$ and $J_u := \sum_{v\in S_u}(d_v+1)$,
\begin{align*}
\Phi_u\leq C\rho_n^{|S_u|-1}\chi\,\omega^{[J_u]}_A .
\end{align*}
This is trivial for a leaf.  For an internal $u$, the induction hypothesis,
$\omega\leq C$ and $\int\check w(x,y)\,dy = \check{\mathfrak D}(x)$ give
\begin{align*}
\int\rho_n\check w(x,y)\Phi_c(y)\,dy\leq C\rho_n^{|S_c|}\check{\mathfrak D}(x),
\end{align*}
and $\check{\mathfrak D}(x)\omega^{[d]}_A(x)\leq C\omega^{[d+1]}_A(x)$, because
$\check{\mathfrak D}\leq s+h_s|x_s|\leq C(1+|x_s|)$.  Integrating $\Phi_1$ with
Lemma~\ref{lem-c-loc} proves the lemma.  Only $\int\check w(x,y)\,dy$ enters,
so neighbor positions in $\mathscr H$ need no treatment.
\end{proof}

\section{Proof of Theorem~\ref{thm-tightness-ring}}
\label{app-weighted}

The proof follows Section~\ref{sec-tight-compact-proof} step by step and
keeps its step names.  Three ingredients are added.  First, every bound must
be uniform in the location of the cell, with $f(t)$, $t$ and $h_t$ entering
as weights; in the tail the $F$-measure of a cell can be far smaller than
$n^{-1}$, so terms of that order, harmless on compact sets, must be absorbed
through the weighted form of the bounds.  Second, long cells can no longer be
split into boundedly many short ones (Step~C of the compact proof fails), and
are written instead as differences of one-sided cells $\{m>s'\}$, which behave
like cells below the bandwidth with a single window.  Third, beyond the cutoff
$\lambda_n$ the process is shown to be flat, using the tail of the degree
distribution; this is condition~(iii) of Proposition~\ref{prop-grid}.  No new
idea is involved beyond these, but the bookkeeping is heavier.  At the start
of each step we say which parts are the compact argument carried through with
weights and which parts need additional care.

\subsection*{Standing conventions}

Notation is as in Table~\ref{tab-notation}, with four additions: the
one-sided cells $O_s$ and the class $\mathscr C_n$ of
Definition~\ref{def-w-onesided}; the domain $J_\infty = [t_0^-,\infty)$; the
weighted cell quantities $Q_X$, $\Lambda_X$, $\Xi_X$ of
\eqref{eq-w-cell-quantities}; and the relative H\"older error $\varepsilon_X$ of
\eqref{eq-w-a-window}.

The conditions of Theorem~\ref{thm-fidi-clt} are in force, together with
\ref{condition-d3}--\ref{condition-l3} on $J_\infty = [t_0^-,\infty)$ and
\ref{condition-s3}, so that $\rho_n = n^{-b}$ and $\lambda_n = n^a$ with
$0<a<b/2$.  Constants $C,c$ depend on $t_0,t_0^-$ and the constants in the
assumptions, never on $n$ or on a cell.  A nonnegative $\mathfrak e_n$ is
\emph{negligible} if $\mathfrak e_n\leq n^{-K}$ for every fixed $K$; Chernoff
tails $e^{-cn\rho_nt_0}$ are negligible since $n\rho_n = n^{1-b}$.  As in the
compact proof, $a(x,\xi) := \check w(x,\xi)-\check{\mathfrak D}(x)$,
$g := \mathfrak D-\check{\mathfrak D}$, $g_t := \sup\{g(x):\mathfrak D(x)\leq2t\}$
and $\Upsilon_n(t) := g_t/h_t$.

\paragraph{The polynomial regime.}
By \ref{condition-d3}(iii), $f(t)\geq c\lambda_n^{-\gamma_+}$ for
$t\in[t_0,\lambda_n]$, so every quantity below that is a power of $n$, $\rho_n$,
$\lambda_n$ and $f(t)$ is bounded below by a fixed negative power of $n$.  We
record two instances.  With $K_0 := 1-b+a\gamma_+$,
\begin{align}
\begin{split}
\label{eq-w-poly}
\Delta_X& \geq n^{-K_0}\quad\text{for every cell }X\in\mathscr C_n\text{ (defined below)},
\\ nf(t)h_t&\geq n^{c_1}\quad\text{for some }c_1>0\text{ and all }t\in[t_0,\lambda_n].
\end{split}
\end{align}
The second holds because
\begin{align*}
nf(t)h_t\geq c\,n\lambda_n^{-\gamma_+}h_{t_0}\asymp n^{(1+b)/2-a\gamma_+},
\end{align*}
and \ref{condition-s3}(iii) says $2a\gamma_+<1+b$.  Also, by
\ref{condition-s3}(ii) and $t\,h_t^{-\vartheta}\geq t_0^{1-\vartheta}h_t^{-\vartheta}$,
for $n\geq n_0$ and uniformly in $t\in[t_0,\lambda_n]$,
\begin{align}
\label{eq-w-chi}
\Upsilon_n(t)\leq t_0^{\vartheta-1}\Bigl(\frac{h_t}{t}\Bigr)^{\vartheta}\leq1,
\qquad g_t = \Upsilon_n(t)h_t\leq h_t .
\end{align}

\paragraph{The exponent $\vartheta'$.}
We fix $\vartheta'\in(1/4,\vartheta]$ satisfying
\begin{align}
\label{eq-w-theta-prime}
\begin{split}
&\text{(i) } K_0(4\vartheta'-1)\leq b/2;\qquad
\text{(ii) } \gamma_-(2-4\vartheta')\geq2-\vartheta-2\vartheta';\qquad
\text{(iii) } 4\vartheta'-1<b(4\vartheta'+1);\\
&\text{(iv) } \vartheta'\leq5/16,\quad
\gamma_-(4-4\vartheta')\geq\max\{6-\vartheta-2\vartheta',\,5q_4-1\},\quad
\gamma_-(13/4-4\vartheta')\geq3q_4;\\
&\text{(v) } (2\vartheta'-\vartheta)(1-b)\leq2b .
\end{split}
\end{align}
Constraints (i)--(iv) hold for all $\vartheta'$ close enough to $1/4$: at
$\vartheta' = 1/4$ each is strict, because $\gamma_->3+\beta$ and
$\gamma_-\geq2q_4$.  Constraint (v) can be met by some $\vartheta'>1/4$ if and
only if $(\tfrac12-\vartheta)(1-b)<2b$, which is \ref{condition-s3}(v); it is
the analogue of \eqref{eq-c-theta-prime}(iv) and is used in exactly one place
(Step~Z4).  Constraints (ii) and (iv) have no counterpart in the compact
proof: they bound the powers of $t$ and $f(t)$ that the weights produce, and
they are the price of uniformity in the tail.  By (i) and
\eqref{eq-w-poly}, for every $X\in\mathscr C_n$ and $n\geq n_0$,
\begin{align}
\label{eq-w-absorb}
\varepsilon\,\Delta_X\leq\Delta_X^{4\vartheta'}\ \text{ whenever }
0\leq\varepsilon\leq n^{-b/2},\qquad
\mathfrak e_n\leq\Delta_X^{4\vartheta'}\ \text{ for negligible }\mathfrak e_n ;
\end{align}
in particular this applies with $\varepsilon = n^{-1}$, with
$\varepsilon = \rho_n\lambda_n\leq\rho_n^{1/2} = n^{-b/2}$, and hence with
$\varepsilon = \rho_nt$ for $t\leq\lambda_n$.

\subsection*{Lattice, cells, and weighted cell quantities}

Lattice points $\lambda_m = m/(n'\rho_n)$, cells $A = A_{rt}$, and
$k_A,\ell_A,\Delta_A$ are as in \eqref{eq-c-cell}; $h_t := (t/(n''\rho_n))^{1/2}$
and $\sigma_t := n''\rho_nh_t = (n''\rho_nt)^{1/2}$, so that
\begin{align*}
\frac{h_t}{t} = \sigma_t^{-1},\qquad \frac{k_A}{\sigma_t} = \frac{\ell_A}{h_t}.
\end{align*}
Unlike the compact case, $h_t$ and $\sigma_t$ now vary by unbounded factors
across the index set, and every bound must track them.  Cells are classified
as in \eqref{eq-c-cell-types}: a \emph{short} cell, with $\ell_A\leq t/2$, is in
the \emph{small regime} if $\ell_A\leq h_t$ and in the \emph{large regime}
otherwise, and a cell with $\ell_A>t/2$ is \emph{long}.  As on compact
intervals, the short/long distinction keeps all estimates for a cell at the
single scale $t$.  Over the tail it matters more, since a long cell may span
unboundedly many scales; long cells are handled through one-sided cells in
Step~C below.  We write $\bar\ell_A := \max(h_t,\ell_A)$, $\mu_A := \min(1,\ell_A/h_t)$ and
$m_A := \sigma_t^{-1}\mu_A$.  For a short cell, \ref{condition-d3}(i) on
$[r,t]\subset[t/2,t]$ gives
\begin{align}
\label{eq-w-Delta-ell}
c\,f(t)\ell_A\leq\Delta_A\leq C\,f(t)\ell_A .
\end{align}

\begin{definition}[One-sided cells]
\label{def-w-onesided}
For a lattice point $s\in[t_0,\lambda_n]$ let $O_s := \{s'+1,s'+2,\ldots\}$, so
that $\mathbbm 1(D_i\in O_s) = \mathbbm 1(D_i>s')$.  The \emph{scale} of $O_s$
is $s$, its endpoint set is $S_{O_s} := \{s\}$, and
\begin{align*}
m_{O_s} := \sigma_s^{-1},\qquad \mu_{O_s} := 1,\qquad \bar\ell_{O_s} := h_s,\qquad
\Delta_{O_s} := s\,f(s).
\end{align*}
A short cell $A_{rt}$ has scale $t$ and endpoint set $S_A := \{r,t\}$.  Let
$\mathscr C_n$ be the set of short cells with right endpoint in
$[t_0,\lambda_n]$ together with the $O_s$, $s\in[t_0,\lambda_n]$.  For
$X\in\mathscr C_n$ with scale $t$ we set
\begin{align}
\label{eq-w-cell-quantities}
Q_X := f(t)\,h_t\,\mu_X^2,\quad
\Lambda_X := \Bigl(\frac{\bar\ell_X}{t}\Bigr)^{\vartheta}+\Upsilon_n(t),\quad
\Xi_X := \begin{cases} f(t)^2t^2\mu_X^2\Lambda_X^2, & X\text{ a short cell},\\
 (sf(s))^2, & X = O_s .\end{cases}
\end{align}
\end{definition}

The quantities $Q_X$ and $\Xi_X$ are the weighted forms of
$Q_A = h_t\mu_A^2$ and $\Xi_A = \mu_A^2\bar\ell_A^{2\vartheta}$ of
\eqref{eq-c-cell-quantities}: on a compact set $f(t)\asymp t\asymp1$ and
$\Upsilon_n(t)\leq h_t^\vartheta$, so they reduce to the compact ones.  By
\eqref{eq-w-chi}, for every $X\in\mathscr C_n$,
\begin{align*}
\Lambda_X\leq C\Bigl(\frac{\bar\ell_X}{t}\Bigr)^{\vartheta},\qquad
\Xi_X\geq\tfrac14f(t)^2t^2\mu_X^2\Lambda_X^2
\end{align*}
(for $O_s$ because $\Lambda_{O_s}\leq2$).

\begin{lemma}[Weighted cell comparisons]
\label{lem-w-cell}
For $n\geq n_0$ and every $X\in\mathscr C_n$ with scale $t$:
\begin{enumerate}[(a)]
\item $\Delta_X\geq c\,f(t)h_t\mu_X$, $\ Q_X\leq C\Delta_X$, and
$\Xi_X\leq C\Delta_X^{2\vartheta}$;
\item for a long cell $A = A_{rt}$ with $t\leq\lambda_n$,
$\Delta_A\geq c\max\{\Delta_{O_r},\Delta_{O_t}\}$;
\item for every lattice point $s\in[t_0,\lambda_n]$,
$1-F(s)\leq C\,sf(s)\leq C$, and the grid mesh satisfies
$F(\lambda_{m+1})-F(\lambda_m)\leq C_df(\lambda_m)/(n'\rho_n)\to0$
uniformly over $\lambda_m\leq\lambda_n$.
\end{enumerate}
\end{lemma}

\begin{proof}
(a) For a short cell, \eqref{eq-w-Delta-ell} gives
\begin{align*}
\Delta_A\geq cf(t)\ell_A\geq cf(t)h_t\mu_A,\qquad
Q_A = f(t)h_t\mu_A^2\leq f(t)\min(\ell_A,h_t)\leq C\Delta_A .
\end{align*}
For $\Xi_A$, by \eqref{eq-w-chi},
$\Xi_A\leq Cf(t)^2t^{2-2\vartheta}\mu_A^2\bar\ell_A^{2\vartheta}$.  In the small
regime $\mu_A^2\bar\ell_A^{2\vartheta} = \ell_A^{2\vartheta}(\ell_A/h_t)^{2-2\vartheta}\leq\ell_A^{2\vartheta}$,
and in the large regime it equals $\ell_A^{2\vartheta}$, so
\begin{align*}
\frac{\Xi_A}{\Delta_A^{2\vartheta}}\leq C\bigl(tf(t)\bigr)^{2-2\vartheta}\leq C
\end{align*}
by \ref{condition-d3}(iii).  For $O_s$, $\Delta_{O_s} = sf(s)\geq f(s)h_s$ since
$h_s\leq s$, $Q_{O_s} = f(s)h_s\leq\Delta_{O_s}$, and
$\Xi_{O_s} = (sf(s))^2\leq C(sf(s))^{2\vartheta}$ since $sf(s)\leq C$.

(b) If $\ell_A>t/2$ then $r<t/2$, so by \ref{condition-d3}(i),
\begin{align*}
\Delta_A\geq F(t)-F(t/2)\geq c\,tf(t),\qquad \Delta_A\geq F(2r)-F(r)\geq c\,rf(r).
\end{align*}

(c) The first claim is \ref{condition-d3}(iii) with $sf(s)\leq Cs^{1-\gamma_-}\leq C$;
the mesh bound is \ref{condition-d3}(i) with $f\leq C$.
\end{proof}

\paragraph{Step H: hub reduction.}
\emph{Analogous.}  The hub set $\mathscr H := \{g>\mathfrak D/4\}$ and the
hub-free process $\mathbb G'_n$ of \eqref{eq-c-G-prime} are unchanged, and
Lemma~\ref{lem-c-hub-reduction} applies verbatim, since its proof involves no
compactness.  \emph{Additional care.}  The reference positions
$x_\lambda = \mathfrak D^{-1}(\lambda)$ are now needed for all $\lambda\leq2\lambda_n$,
and the clipping bias $g$ must be controlled on $\{\mathfrak D\leq2t\}$ for
every $t\leq\lambda_n$, with a $t$-dependent bound; this is what the clipping
ratio $\Upsilon_n(t)$ and \ref{condition-s3}(ii) provide.

\begin{lemma}[Hubs, weighted]
\label{lem-w-hubs}
For $n\geq n_0$: $\mathrm{Leb}(\mathscr H)\leq C_H\rho_n^{2+\beta}$;
$\mathscr H\cap\{\mathfrak D\leq2\lambda_n\} = \emptyset$, so that no
$x_\lambda$ with $\lambda\leq2\lambda_n$ is a hub; and $g(x)\leq g_t = \Upsilon_n(t)h_t$
whenever $\mathfrak D(x)\leq2t$, $t\leq\lambda_n$.  Off $\mathscr H$,
$\tfrac34\mathfrak D\leq\check{\mathfrak D}\leq\mathfrak D$.
\end{lemma}

\begin{proof}
If $\mathfrak D(x)\in J_\infty$, \eqref{eq-a-growth} shows that $x\in\mathscr H$
forces
\begin{align*}
\mathfrak D(x)^{1+\beta}>(4K\rho_n^{1+\beta})^{-1},\qquad\text{that is,}\qquad
\mathfrak D(x)\geq c_H\rho_n^{-1}.
\end{align*}
If $\mathfrak D(x)<t_0^-$, Markov's inequality and \ref{condition-l1} with
$R = t_0^-$ give $g(x)\leq C\rho_n^{1+\beta}<c/4\leq\mathfrak D(x)/4$ for $n$
large, by \ref{condition-g1}.  Hence $\mathscr H\subseteq\{\mathfrak D\geq c_H\rho_n^{-1}\}$,
whose measure is at most $C_H\rho_n^{2+\beta}$ by Markov's inequality and
\ref{condition-m1}, and which is disjoint from $\{\mathfrak D\leq2\lambda_n\}$
because $\rho_n\lambda_n\to0$.  The bound on $g$ is \eqref{eq-a-growth}.
\end{proof}

\paragraph{Hub-free increments.}
For $X\in\mathscr C_n$ or a long cell, let $q_i(X) := P(D_i\in X\mid\mathcal G)$,
$p'_X := E[\chi_1\mathbbm 1(D_1\in X)]$, $\bar p_X := n^{-1}\sum_i\chi_iq_i(X)$, and
\begin{align}
\begin{split}
\label{eq-w-BZ}
\Delta'_X &:= \frac1{\sqrt n}\sum_i\chi_i\bigl[\mathbbm 1(D_i\in X)-p'_X\bigr] = B_X+Z_X,
\\
B_X &:= \frac1{\sqrt n}\sum_i\chi_i\bigl[\mathbbm 1(D_i\in X)-q_i(X)\bigr],
\qquad
Z_X := \sqrt n\,(\bar p_X-p'_X),
\end{split}
\end{align}
exactly as in \eqref{eq-c-BZ}; $\mathbb G'_n(t)-\mathbb G'_n(r) = \Delta'_{A_{rt}}$.
For $i\neq j$ let
\begin{align*}
g^{(j)}_i(A) := P(D_{i\setminus j} = r'\mid\mathcal G)-P(D_{i\setminus j} = t'\mid\mathcal G)
\end{align*}
for a short cell and $g^{(j)}_i(O_s) := P(D_{i\setminus j} = s'\mid\mathcal G)$,
and define, with $\chi_i$ inserted at every root,
\begin{align}
\begin{split}
\label{eq-w-Sigma-Phi}
\bar\Sigma(X) &:= \frac1n\sum_{i\neq j}\chi_i\chi_j\,p_{ij}\bigl|g^{(j)}_i(X)\bigr|\bigl|g^{(i)}_j(X)\bigr|,
\qquad
\Phi_m := \rho_n\sum_{i\neq m}\chi_i\,a(\xi_i,\xi_m)\,g^{(m)}_i(X),
\\ \Psi_m &:= E_\xi\Bigl[\Bigl(\rho_n\sum_{i\neq m}\chi_ia(\xi_i,\xi)g^{(m)}_i(X)\Bigr)^2\Bigr].
\end{split}
\end{align}

\paragraph{Step L: lattice criterion.}
\emph{Analogous.}  Proposition~\ref{prop-grid} is applied to $\mathbb G'_n$ on
$[t_0,\infty)$ with $T_n := \lambda_n$ and $\nu := 2\vartheta'$.  The grid mesh
below the cutoff tends to zero by Lemma~\ref{lem-w-cell}(c).
\emph{Additional care.}  Condition~(iii) of the proposition, flatness beyond
$\lambda_n$, is no longer void; it is verified in Step~A from
\ref{condition-s3}(i).  Condition~(ii) is required only for lattice points
below $\lambda_n$, which is why every lemma below is stated for cells with
scale in $[t_0,\lambda_n]$.

\paragraph{Step C: reduction to short and one-sided cells.}
\emph{Additional care.}  In the compact proof a long cell is cut into a
bounded number of short cells, a bound that is unavailable on
$[t_0,\lambda_n]$.  Instead, for a long cell $A = A_{rt}$ we use
$\mathbbm 1(D_i\in A) = \mathbbm 1(D_i>r')-\mathbbm 1(D_i>t')$, whence
\begin{align}
\label{eq-w-long-split}
\Delta'_A = \Delta'_{O_r}-\Delta'_{O_t},\qquad
B_A = B_{O_r}-B_{O_t},\qquad Z_A = Z_{O_r}-Z_{O_t},
\end{align}
and $\Delta_A\geq c\max\{\Delta_{O_r},\Delta_{O_t}\}$ by
Lemma~\ref{lem-w-cell}(b).  A one-sided cell has a single window, at $s$, and
no cancellation between two endpoints: it behaves like a cell in the small
regime with $\mu = 1$, and its H\'ajek part is bounded without cancellation
by $sf(s)$ (Step~K).  All lemmas in Steps~W--Z4 are therefore stated for
$X\in\mathscr C_n$, that is for short cells and one-sided cells alike, and the
proofs indicate the (few) places where the one-sided case differs.  Long cells
are then handled at the end (Step~A, Corollary~\ref{cor-w-long}).  The
covariance identity of Lemma~\ref{lem-c-pairing} holds verbatim for $O_s$,
since $g^{(j)}_i(O_s) = u(1)-u(0)$ for $u(\epsilon) := P(D_{i\setminus j}+\epsilon>s'\mid\mathcal G)$.

\paragraph{Step W: windows, roots and cancellation.}
\emph{Analogous.}  Lemmas~\ref{lem-poisson-binomial} and \ref{lem-pb-local}
are used unchanged.  The windows, the localization weights
$\omega^{[j]}_X$ and the three lemmas of this step (localization integrals,
roots, cancellation) are those of the compact proof with $f(t)$, $t$ and
$h_t$ carried as weights.  \emph{Additional care.}
\begin{enumerate}[(1)]
\item The measure of a window of width $h$ at $s$ is now $\asymp f(s)h$ rather
than $\asymp h$, and $\check{\mathfrak D}$ on a window is of order $t$, not
order one.  So every localization integral carries $f(t)h_t$, every power of
$\check{\mathfrak D}$ costs a power of $t$ and is absorbed into the weight by
Lemma~\ref{lem-w-loc}(a), and \eqref{eq-c-a-window} is replaced by the
scale-invariant \eqref{eq-w-a-window}.
\item In the roots lemma, (E1) must hold for every lattice point up to
$\lambda_n$, with the second moment of $\check w$ of order $s^2$; the
Bernstein exponent is $\min\{(\rho_ns)^{-1},\sigma_s\}$, and
$(\rho_ns)^{-1}\geq n^{b/2}$ uses $\rho_n\lambda_n^2\to0$.
\item In the cancellation lemma the H\"older error of $f$ is the relative one
of \ref{condition-d3}(ii), $(\bar\ell_X/t)^{\theta_f}$, and the clipping error
is $\Upsilon_n(t)$ rather than the fixed $\Upsilon_n$ of \eqref{eq-c-clip}; for
one-sided cells the $d = 0$ cancellation is unavailable, since
$\int_0^1b_M(s';p)\,dp\neq0$.
\end{enumerate}

Let $c_0,\zeta$ be the constants of Lemma~\ref{lem-pb-local},
$\zeta_* := \zeta/6$, and
\begin{align}
\label{eq-w-cw}
c_w := \min\Bigl\{\frac{c_0}{32},\ \frac12\Bigl(1-\frac{t_0^-}{t_0}\Bigr)\Bigr\},
\end{align}
so that $(1-c_w)s\geq t_0^-$ for every $s\geq t_0$; this is \eqref{eq-c-cw}
without the upper constraint, which is not needed on an unbounded index set.
For $X\in\mathscr C_n$ with scale $t$, $s\in S_X$ and $x\in[0,1]$, write
$z_s(x) := (\check{\mathfrak D}(x)-s)/h_s$,
\begin{align}
\label{eq-w-omega}
\begin{split}
\omega^{[j]}_X(x) &:= \sum_{s\in S_X}\bigl(1+|z_s(x)|\bigr)^je^{-\zeta_*z_s(x)^2},
\qquad \omega_X := \omega^{[1]}_X,
\\
e_X(x) &:= \min_{s\in S_X}\bigl(|\mathfrak D(x)-s|+|s-t|\bigr),
\end{split}
\end{align}
These are the quantities of \eqref{def-omega-e-a}, with the same
interpretation (see the discussion following \eqref{def-omega-e-a}); only the
scale $t$, the bandwidth $h_t$ and, through the windows, the density $f(t)$
now vary with the cell.  Let $W_s := \{x\notin\mathscr H:\ |\check{\mathfrak D}(x)-s|\leq c_ws\}$,
$W_X := \bigcup_{s\in S_X}W_s$.  For $x\in W_s$ with $s\in[t_0,\lambda_n]$,
Lemma~\ref{lem-w-hubs} gives
\begin{align*}
t_0^-\leq(1-c_w)s\leq\check{\mathfrak D}(x)\leq\mathfrak D(x)\leq\tfrac43\check{\mathfrak D}(x)\leq2t,
\end{align*}
so $\mathfrak D(x)\in J_\infty$ and $g(x)\leq g_t\leq h_t$; and by
\ref{condition-d3}(i), for $0\leq u\leq c_ws$,
\begin{align}
\begin{split}
\label{eq-w-window-measure}
\mathrm{Leb}\bigl(\{x\notin\mathscr H:\ |\check{\mathfrak D}(x)-s|\leq u\}\bigr)
 &\leq\mathrm{Leb}\bigl(\{|\mathfrak D-s|\leq u+g_t\}\bigr)
 \\&\leq2C_d\,f(s)\,(u+g_t)\leq C\,f(t)\,(u+g_t).
 \end{split}
\end{align}
Consequently, with $\varepsilon_X(x) := (e_X(x)/t)^{\theta_w}$, for $x\in W_X$,
\begin{align}
\label{eq-w-a-window}
\|a(x,\cdot)\|_{L^2}\leq K_2\mathfrak D(x)\leq Ct,\quad
\|a(x,\cdot)\|_{L^4}\leq Ct^{q_4},\quad
\|a(x,\cdot)-a(x_t,\cdot)\|_{L^2}\leq C\,t\,\varepsilon_X(x) .
\end{align}
The first two are \eqref{eq-a-growth}.  For the third, suppose first that the
minimum in $e_X(x)$ is attained at an $s$ with $x\in W_s$.  Then
$|\mathfrak D(x)-s|\leq c_ws+h_t\leq s/2$ and $|s-t|\leq t/2$, so
\ref{condition-g3} applied twice (it transfers to $\check w$ and to $a$ with
constant $2K_w$) gives
\begin{align*}
\|a(x,\cdot)-a(x_t,\cdot)\|\leq2K_w\bigl[s^{1-\theta_w}|\mathfrak D(x)-s|^{\theta_w}+t^{1-\theta_w}|s-t|^{\theta_w}\bigr]
\leq4K_wt\,\varepsilon_X(x).
\end{align*}
Otherwise $e_X(x)\geq c_ws-h_t\geq c_wt/3$, so $\varepsilon_X(x)\geq c$ and the
bound follows from $\|a(x,\cdot)\|+\|a(x_t,\cdot)\|\leq Ct$.  (The point $x_t$
is not a hub by Lemma~\ref{lem-w-hubs}.)

\begin{lemma}[Localization integrals, weighted]
\label{lem-w-loc}
Fix $J_0$.  There are $C,n_0$ such that for $n\geq n_0$, every
$X\in\mathscr C_n$ with scale $t$, $0\leq j\leq J_0$, $p\in[0,J_0]$ and
$\theta\in[0,1]$:
\begin{enumerate}[(a)]
\item $\check{\mathfrak D}(x)^p\,\omega^{[j]}_X(x)\leq t^p\,\omega^{[j+\lceil p\rceil]}_X(x)$
for every $x$;
\item $\displaystyle\int_{[0,1]\setminus\mathscr H}e_X(x)^\theta\,\omega^{[j]}_X(x)\,dx
\leq C\,f(t)\,h_t\,\bar\ell_X^{\,\theta}$.
\end{enumerate}
The same holds with $\zeta_*$ replaced by any fixed $\zeta'\in(0,\zeta_*]$,
with $C$ depending on $\zeta'$.
\end{lemma}

\begin{proof}
For $s\in S_X$ we have $s\in[t/2,t]$, hence $h_t/\sqrt2\leq h_s\leq h_t$ and
$f(s)\leq C_df(t)$.

(a) Since $h_s\leq s\leq t$, termwise
$\check{\mathfrak D}\leq s+h_s|x_s|\leq t(1+|x_s|)$.

(b) Fix $s\in S_X$.  On $W_s$, $g\leq h_t$, so
\begin{align*}
e_X\leq|\mathfrak D-s|+|s-t|\leq h_s|x_s|+h_t+\ell_A\mathbbm 1(s = r)\leq3\bar\ell_X(1+|x_s|).
\end{align*}
Split $W_s$ into $\{|\check{\mathfrak D}-s|\leq h_s\}$ and the shells
$\{2^{k-1}h_s<|\check{\mathfrak D}-s|\leq2^kh_s\}$, $2^kh_s\leq c_ws$.  By
\eqref{eq-w-window-measure} the $k$-th shell has measure at most
$Cf(t)2^kh_t$, and on it the integrand is at most
$(3\bar\ell_X)^\theta(1+2^k)^{j+1}e^{-\zeta_*4^{k-1}}$; summing gives
$Cf(t)h_t\bar\ell_X^\theta$.

A term of $\omega^{[j]}_X$ indexed by $s_1$ at a point
$x\notin W_{s_1}\cup\mathscr H$ has
\begin{align*}
|x_{s_1}|\geq c_ws_1/h_{s_1} = c_w\sigma_{s_1}\geq c_w\sigma_t/\sqrt2,
\end{align*}
so it is at most $Ce^{-cn\rho_nt_0}$, while $e_X^\theta\leq1+\mathfrak D+t$;
its integral is at most $Ce^{-cn\rho_nt_0}(2+E\mathfrak D+\lambda_n)$,
negligible, hence at most $f(t)h_t\bar\ell_X^\theta$ by \eqref{eq-w-poly}.
\end{proof}

\begin{lemma}[Roots, weighted]
\label{lem-w-roots}
For $i\in[n]$ and $I'\ni i$ with $|I'|\leq4$, let
$\mu_{i,I'} := \sum_{l\notin I'}p_{il}$ and
$\Sigma_{i,I'} := \sum_{l\notin I'}p_{il}(1-p_{il})$.  Let $\mathscr E_i$ be the
event that for all such $I'$:
\begin{enumerate}[(E1)]
\setcounter{enumi}{-1}
\item $|\mu_{i,I'}-n''\rho_n\check{\mathfrak D}(\xi_i)|\leq\frac{c_w}4n''\rho_n
\max\{\check{\mathfrak D}(\xi_i),t_0\}$;
\item for every lattice point $s\in[t_0,\lambda_n]$ with
$|\check{\mathfrak D}(\xi_i)-s|\leq c_ws$:
$|\mu_{i,I'}-n''\rho_n\check{\mathfrak D}(\xi_i)|\leq\sigma_s$ and
$\Sigma_{i,I'}\geq\mu_{i,I'}/4$.
\end{enumerate}
Then (a) $\sup_{x\notin\mathscr H}P(\mathscr E_i^c\mid\xi_i = x)$ is negligible;
and (b) on $\mathscr E_i\cap\{\xi_i\notin\mathscr H\}$, for every
$X\in\mathscr C_n$ with scale $t$, every $I'$ as above, $1\leq k\leq4$ and
$|j|\leq3$, with $P_{I'}(m) := P(D_{i\setminus I'} = m\mid\mathcal G)$,
\begin{align*}
\Bigl|\sum_{m\in A}\mathsf D^kP_{I'}(m+j)\Bigr|
 &\leq C\,m_A\,\sigma_t^{-(k-1)}\,\omega^{[k]}_A(\xi_i)+\mathfrak e_n
 \qquad(X = A\text{ a short cell}),\\
\bigl|\mathsf D^{k-1}P_{I'}(s'+j)\bigr|
 &\leq C\,m_{O_s}\,\sigma_s^{-(k-1)}\,\omega^{[k]}_{O_s}(\xi_i)+\mathfrak e_n
 \qquad(X = O_s),
\end{align*}
with $C$ and $\mathfrak e_n$ not depending on $i$, $X$, $I'$, $k$, $j$.
\end{lemma}

\begin{proof}
(a) Fix $x\notin\mathscr H$ and $I'$, and condition on $\xi_i = x$.  The
$p_{il} = \rho_n\check w(x,\xi_l)$, $l\notin I'$, are i.i.d.\ in $[0,1]$ with
mean $\rho_n\check{\mathfrak D}(x)$, so
\begin{align*}
\bigl|E\mu_{i,I'}-n''\rho_n\check{\mathfrak D}(x)\bigr|\leq2\rho_n\check{\mathfrak D}(x),
\qquad
\mathrm{Var}(\mu_{i,I'})\leq n\rho_n\check{\mathfrak D}(x).
\end{align*}
(E0) follows from Bernstein's inequality exactly as in
Lemma~\ref{lem-c-roots}(a), up to probability $2e^{-cn\rho_nt_0}$.

(E1): suppose $|\check{\mathfrak D}(x)-s|\leq c_ws$ for a lattice point
$s\in[t_0,\lambda_n]$.  Then $x\in W_s$, so $\mathfrak D(x)\in J_\infty$, and
\ref{condition-l3} with Jensen's inequality give
\begin{align*}
E\bigl[\check w(x,\xi)^2\bigr]\leq E\bigl[w^{2+\beta}(x,\xi)\bigr]^{2/(2+\beta)}\leq Cs^2,
\qquad\text{whence}\qquad
\mathrm{Var}(\mu_{i,I'})\leq Cn\rho_n^2s^2 .
\end{align*}
Bernstein's inequality with deviation $\sigma_s/2$ (which exceeds
$2\rho_n\check{\mathfrak D}(x)$ for $n$ large) gives
\begin{align*}
P\bigl(|\mu_{i,I'}-E\mu_{i,I'}|>\sigma_s/2\bigr)
 &\leq2\exp\Bigl(-\frac{\sigma_s^2}{8(Cn\rho_n^2s^2+\sigma_s/6)}\Bigr)
\\ &\leq2\exp\bigl(-c\min\{(\rho_ns)^{-1},\sigma_s\}\bigr),
\end{align*}
using $\sigma_s^2 = n''\rho_ns$.  Here $(\rho_ns)^{-1}\geq(\rho_n\lambda_n)^{-1}\geq n^{b/2}$
(as $\rho_n\lambda_n^2\leq1$) and $\sigma_s\geq cn^{(1-b)/2}$, so this is
negligible.

For $\Sigma_{i,I'}$, write
\begin{align*}
\mu_{i,I'}-\Sigma_{i,I'} = \sum_{l\notin I'}p_{il}^2\leq S+\mu_{i,I'}/2,
\qquad S := \#\{l\notin I':p_{il}>1/2\},
\end{align*}
which given $\xi_i = x$ is binomial with mean
\begin{align*}
ES\leq n(2\rho_n)^{2+\beta}E[w^{2+\beta}(x,\xi)]\leq Cn\rho_n^{2+\beta}s^{2+\beta}.
\end{align*}
With $a := n''\rho_ns/16$, the Chernoff bound $P(S\geq a)\leq(e\,ES/a)^a$ and
the bound $e\,ES/a\leq C(\rho_ns)^{1+\beta}\leq e^{-1}$ give
$P(S\geq a)\leq e^{-a}$.  On
$\{S<a\}$ and the deviation event,
\begin{align*}
\mu_{i,I'}\geq(1-c_w)n''\rho_ns-\sigma_s\geq\tfrac12n''\rho_ns,
\end{align*}
so $S<\mu_{i,I'}/8$ and $\Sigma_{i,I'}\geq\mu_{i,I'}/4$.  A union bound over at
most $4n^3$ sets $I'$ and at most $n'\rho_n\lambda_n+1\leq n$ lattice points
proves (a).

(b) Work on $\mathscr E_i$ with $\xi_i = x\notin\mathscr H$, fix $I'$, and
write $\mu,\Sigma$.  Let $s\in S_X$; we bound $\mathsf D^kP_{I'}(m)$ for
$|m-s'|\leq\frac{c_w}4n''\rho_ns$, which covers all points used below.

\emph{Case 1: $|\check{\mathfrak D}(x)-s|\leq c_ws$.}  By (E1),
$\mu\in[\frac12n''\rho_ns,2n''\rho_ns]$, $\Sigma\in[\sigma_s^2/8,2\sigma_s^2]$, and
\begin{align*}
|m-\mu|\leq\tfrac{c_w}4n''\rho_ns+2+c_wn''\rho_ns+\sigma_s\leq2c_wn''\rho_ns
\leq\tfrac{c_0}{16}n''\rho_ns\leq c_0\Sigma
\end{align*}
for $n$ large, so Lemma~\ref{lem-pb-local} gives, with
$y_m := (m-\mu)/\Sigma^{1/2}$,
\begin{align*}
|\mathsf D^kP_{I'}(m)|\leq C\Sigma^{-(k+1)/2}(1+|y_m|)^ke^{-\zeta y_m^2}.
\end{align*}
Since $n''\rho_nh_s = \sigma_s$,
\begin{align*}
m-\mu = (m-s')+(s'-n''\rho_ns)-\sigma_sz_s(x)+\bigl(n''\rho_n\check{\mathfrak D}(x)-\mu\bigr),
\end{align*}
so for $|m-s'|\leq L\sigma_s$, $|y_m+\tau z_s(x)|\leq C_L$ with
$\tau := \sigma_s/\Sigma^{1/2}\in[2^{-1/2},\sqrt8]$, hence
$y_m^2\geq x_s^2/4-C_L^2$ and
\begin{align*}
|\mathsf D^kP_{I'}(m)|\leq C\sigma_s^{-(k+1)}(1+|z_s(x)|)^ke^{-\zeta_*z_s(x)^2}.
\end{align*}

\emph{Case 2: $|\check{\mathfrak D}(x)-s|>c_ws$.}  By (E0), as in
Lemma~\ref{lem-c-roots}(b), $|m-\mu|\geq c\max\{\mu,n''\rho_ns\}$ for $n$ large,
and Bernstein's inequality for the Poisson-binomial law gives
$P_{I'}(m)\leq2e^{-cn\rho_nt_0}$; with
$|\mathsf D^kP_{I'}(m)|\leq2^k\max_{j'\leq k}P_{I'}(m-j')$ this is negligible.

\emph{Assembly.}  Identical to Lemma~\ref{lem-c-roots}(b): in the small regime
sum Case~1 at $s = t$ over the $k_A\leq\sigma_t$ points; in the large regime
telescope to the two ends and bound each at its own endpoint, using
$\sigma_r\geq\sigma_t/\sqrt2$; for $O_s$ the bound is the single-point bound
at $s$.
\end{proof}

Every quantity whose moments are bounded in this appendix is a polynomial
with at most $n^4$ terms in the $p_{ij}$, conditional mass functions, and
values $a(\xi_i,\xi_j)$, all bounded by $n^C$ in absolute value.  The event
$\bigcup_i(\{\xi_i\notin\mathscr H\}\cap\mathscr E_i^c)$ therefore contributes
negligibly to every moment, and we suppress it; by \eqref{eq-w-BZ}, roots in
$\mathscr H$ do not occur.

For $X\in\mathscr C_n$, an integer $M$ with $|M-n''|\leq d_1$ ($d_1$ a fixed
constant to be chosen in Step~Z2) and $d\geq0$, let
\begin{align}
\label{eq-w-K}
\begin{split}
\mathsf K^{(d)}_{X,M}(p) &:=
\begin{cases}
(-1)^d\bigl[\mathsf D^db_M(r';p)-\mathsf D^db_M(t';p)\bigr], & X = A_{rt},\\
\mathsf D^db_M(s';p), & X = O_s,
\end{cases}
\\
\mathsf K^{(d)}_X(x) &:= \mathsf K^{(d)}_{X,M}\bigl(\rho_n\check{\mathfrak D}(x)\bigr).
\end{split}
\end{align}
With $d = 0$ and $M = n''$, $n'\rho_n\mathsf K^{(0)}_X$ is the density of the
signed measure defining the H\'ajek kernel $\phi_X$ in Step~K, and the
Hoeffding components of $g^{(1)}_2(X)$ in Step~Z2 carry $\mathsf K^{(d)}_X$ with
$M = n''-d$.

\begin{lemma}[Cancellation, weighted]
\label{lem-w-cancel}
Fix $d_1$.  There are $C,n_0$ such that for $n\geq n_0$, $0\leq d\leq d_1$,
$|M-n''|\leq d_1$, and $X\in\mathscr C_n$ with scale $t$:
\begin{enumerate}[(a)]
\item for every $\lambda\in(0,\rho_n^{-1}]$, with $z_s(\lambda) := (\lambda-s)/h_s$,
\begin{align*}
|\mathsf K^{(d)}_{X,M}(\rho_n\lambda)|\leq C\sigma_t^{-(d+1)}\mu_X
\sum_{s\in S_X}(1+|z_s(\lambda)|)^{d+1}e^{-\zeta_*z_s(\lambda)^2}+\mathfrak e_n,
\end{align*}
and the left side is negligible if $|\lambda-s|>c_ws$ for every $s\in S_X$;
in particular $|\mathsf K^{(d)}_X(x)|\leq C\sigma_t^{-(d+1)}\mu_X\omega^{[d+1]}_X(x)+\mathfrak e_n$
for all $x$;
\item for $j\in\{1,2\}$ and $\theta\in[0,1]$,
$\displaystyle\int_{[0,1]\setminus\mathscr H}e_X^\theta|\mathsf K^{(d)}_X|^j
\leq C\sigma_t^{-j(d+1)}\mu_X^j\,f(t)\,h_t\,\bar\ell_X^{\,\theta}$;
\item if $d\geq1$, or $d = 0$ and $X$ is a short cell,
\begin{align}
\label{eq-w-cancel}
n'\rho_n\Bigl|\int_{[0,1]\setminus\mathscr H}\mathsf K^{(d)}_X(x)\,dx\Bigr|
 \leq C\,f(t)\,\sigma_t^{-d}\,\mu_X\,\Lambda_X+\mathfrak e_n .
\end{align}
\end{enumerate}
\end{lemma}

\begin{proof}
(a) Put $p := \rho_n\lambda$ and let $s\in S_X$.  If $|\lambda-s|\leq c_ws$,
then $p\leq2\rho_n\lambda_n\leq1/4$ and $\Sigma := Mp(1-p)\in[\sigma_s^2/2,2\sigma_s^2]$,
and the proof of Lemma~\ref{lem-w-roots}(b), Case~1, applies verbatim to the
binomial law $b_M(\cdot\,;p)$ (mean $Mp$, within $d_1\rho_n\lambda\leq1$ of
$n''\rho_n\lambda$, and no exceptional event).  The assembly over the cell is
as there.  If $|\lambda-s|>c_ws$ for every $s\in S_X$, the binomial Chernoff
bounds in both tails, valid for every $p\in[0,1]$, make every $b_M(m-j;p)$
with $|m-s'|\leq k_A+d_1$ at most $2e^{-cn\rho_nt_0}$.  The statement for $x$
is the case $\lambda = \check{\mathfrak D}(x)$.

(b) By (a), using $(\sum_su_s)^2\leq2\sum_su_s^2$ and
$((1+|y|)^{d+1}e^{-\zeta_*y^2})^2\leq(1+|y|)^{2d+2}e^{-\zeta_*y^2}$,
\begin{align*}
|\mathsf K^{(d)}_X|^j\leq C\sigma_t^{-j(d+1)}\mu_X^j\omega^{[j(d+1)]}_X+\mathfrak e_n .
\end{align*}
Lemma~\ref{lem-w-loc}(b) bounds the main term, and the negligible term
integrates against $e_X^\theta\leq1+\mathfrak D+t$ to a negligible quantity,
absorbed by \eqref{eq-w-poly}.

(c) Write $\mathsf K := \mathsf K^{(d)}_{X,M}$.

\emph{Step 1: localization.}  On $\{\mathfrak D>2t\}\setminus\mathscr H$,
$\check{\mathfrak D}\geq\frac34\mathfrak D>\frac32t\geq(1+c_w)s$ for every
$s\in S_X$, so the integrand is negligible there by (a).

\emph{Step 2: clipping.}  On $\{\mathfrak D\leq2t\}$,
$0\leq\mathfrak D-\check{\mathfrak D}\leq g_t = \Upsilon_n(t)h_t$.  Since
$\partial_pb_M(m;p) = -M\,\mathsf Db_{M-1}(m;p)$ and $\mathsf D$ commutes with
$\partial_p$, $|\partial_p\mathsf K^{(d)}_{X,M}| = M|\mathsf K^{(d+1)}_{X,M-1}|$.
By the mean value theorem there is
$\lambda^*(x)\in[\check{\mathfrak D}(x),\mathfrak D(x)]$ with
\begin{align*}
\bigl|\mathsf K(\rho_n\check{\mathfrak D}(x))-\mathsf K(\rho_n\mathfrak D(x))\bigr|
\leq\rho_ng_tM\bigl|\mathsf K^{(d+1)}_{X,M-1}(\rho_n\lambda^*(x))\bigr| .
\end{align*}
We apply (a) at $\lambda^*(x)$.  Since $|z_s(\lambda^*)-z_s(x)|\leq h_t/h_s\leq\sqrt2$,
\begin{align*}
e^{-\zeta_*z_s(\lambda^*)^2}\leq e^{2\zeta_*}e^{-\zeta_*z_s(x)^2/2},
\end{align*}
and the right side is at most
\begin{align*}
C\rho_ng_tM\sigma_t^{-(d+2)}\mu_X\omega^{[d+2],\flat}_X(x)+\mathfrak e_n,
\end{align*}
with $\flat$ denoting decay $\zeta_*/2$.  By Lemma~\ref{lem-w-loc}(b),
$n'\rho_nh_t = (n'/n'')\sigma_t$, $M\rho_n\leq2\sigma_t^2/t$ and
$h_t/t = \sigma_t^{-1}$,
\begin{align*}
n'\rho_n\int_{\{\mathfrak D\leq2t\}}\bigl|\mathsf K(\rho_n\check{\mathfrak D})-\mathsf K(\rho_n\mathfrak D)\bigr|
 &\leq C\,\sigma_t\cdot\frac{\sigma_t^2}t\cdot\Upsilon_n(t)h_t\cdot\sigma_t^{-(d+2)}\mu_Xf(t)+\mathfrak e_n
\\ &= Cf(t)\sigma_t^{-d}\mu_X\Upsilon_n(t)+\mathfrak e_n .
\end{align*}

\emph{Step 3: main term.}  By Remark~\ref{remark-inverse} and
\ref{condition-d1},
\begin{align*}
\int_{\{\mathfrak D\leq2t\}}\mathsf K(\rho_n\mathfrak D(x))\,dx
= \int_0^{2t}\mathsf K(\rho_n\lambda)f(\lambda)\,d\lambda .
\end{align*}
Let $U := \bigcup_{s\in S_X}[(1-c_w)s,(1+c_w)s]\subseteq J_\infty\cap[0,2t]$; off
$U$ the integrand is negligible by (a), and $\int f\leq1$.  For $\lambda\in U$
with $|\lambda-s|\leq c_ws$, \ref{condition-d3}(ii) twice (at $(s,\lambda)$ and
$(t,s)$) and \ref{condition-d3}(i) give
\begin{align*}
|f(\lambda)-f(t)|\leq Cf(t)\Bigl(\frac{e(\lambda)}t\Bigr)^{\theta_f},
\qquad e(\lambda) := |\lambda-s|+|s-t| .
\end{align*}
By (a) and the shell argument of Lemma~\ref{lem-w-loc}(b) with Lebesgue
measure $d\lambda$ in place of $dx$ (the shells now have length $2^kh_s$ and no
factor $f$),
\begin{align*}
n'\rho_n\Bigl|\int_U\mathsf K(\rho_n\lambda)[f(\lambda)-f(t)]\,d\lambda\Bigr|
 &\leq Cn'\rho_nf(t)t^{-\theta_f}\sigma_t^{-(d+1)}\mu_Xh_t\bar\ell_X^{\theta_f}
\\ &\leq Cf(t)\sigma_t^{-d}\mu_X\Bigl(\frac{\bar\ell_X}t\Bigr)^{\theta_f}
 \leq Cf(t)\sigma_t^{-d}\mu_X\Bigl(\frac{\bar\ell_X}t\Bigr)^{\vartheta},
\end{align*}
using $\bar\ell_X\leq t$ and $\theta_f\geq\vartheta$.  Finally, up to a
negligible term by (a),
\begin{align*}
\int_U\mathsf K(\rho_n\lambda)\,d\lambda = \frac1{\rho_n}\int_0^1\mathsf K(p)\,dp = 0 :
\end{align*}
for $d\geq1$ because $\int_0^1b_M(m;p)dp = 1/(M+1)$ does not depend on $m$
(Lemma~\ref{lem-binomial-facts}\ref{binom-fact-e}), and for $d = 0$ and
$X = A_{rt}$ because the two endpoint terms cancel.  This is where $d = 0$
requires a short cell: for $O_s$ there is a single endpoint and no
cancellation.  Steps~1--3 give \eqref{eq-w-cancel}.
\end{proof}

\paragraph{Step K: the H\'ajek increment.}
\emph{Analogous.}  The H\'ajek kernels and the proof of
Proposition~\ref{prop-c-hajek} go through with weights: the $\psi$-part is
bounded by $p'_X$, and the $\phi$-part splits into a H\"older term, now of
size $tf(t)\mu_X(\bar\ell_X/t)^{\vartheta}$, and a mass term controlled by
Lemma~\ref{lem-w-cancel}(c).  \emph{Additional care.}  $p'_X\leq C\Delta_X$
needs \ref{condition-d3}(iii) for one-sided cells, and the H\'ajek part of a
one-sided cell is bounded without cancellation, which is why
$\Xi_{O_s} = (sf(s))^2$.

\begin{lemma}[H\'ajek kernels, weighted]
\label{lem-w-kernels}
For $X\in\mathscr C_n$ let $\mathcal P_X(\xi) := P(D_1\in X\mid\xi_1 = \xi)$.
The H\'ajek projection of $\Delta'_X$ is
$n^{-1/2}\sum_j[\psi_X(\xi_j)+\phi_X(\xi_j)]$, where
\begin{align}
\label{eq-w-hajek}
\psi_X(\xi) := \chi(\xi)\,\mathcal P_X(\xi)-p'_X,\quad
\phi_X(\xi) := n'\rho_n\int_{[0,1]\setminus\mathscr H}a(x,\xi)\,\mathsf K^{(0)}_X(x)\,dx
\quad(M = n''),
\end{align}
and $E\psi_X = E\phi_X = 0$, $E\psi_X^2\leq p'_X$.
\end{lemma}

\begin{proof}
For a short cell this is Lemma~\ref{lem-c-kernels}.  For $O_s$,
$q_i(O_s) = 1-F_i(s'\mid\mathcal G)$, whose $S = \{j\}$ component is, by
Lemma~\ref{lem-multilinear},
\begin{align*}
+\rho_na(\xi_i,\xi_j)b_{n''}\bigl(s';\rho_n\check{\mathfrak D}(\xi_i)\bigr) = \rho_na(\xi_i,\xi_j)\mathsf K^{(0)}_{O_s}(\xi_i);
\end{align*}
the rest of the argument is unchanged.
\end{proof}

\begin{proposition}[H\'ajek increment, weighted]
\label{prop-w-hajek}
For $n\geq n_0$ and every $X\in\mathscr C_n$ with scale $t$:
$E\psi_X^2\leq p'_X\leq C\Delta_X$ and
$\|\phi_X\|_{L^2}\leq C\,\Xi_X^{1/2}+\mathfrak e_n\leq C\Delta_X^{\vartheta}$.
\end{proposition}

\begin{proof}
\emph{The bound $p'_X\leq C\Delta_X$.}  For a short cell,
$p'_A = \sum_{m\in A}E[\chi b_{n'}(m;\rho_n\check{\mathfrak D}(\xi))]$.  By
Lemma~\ref{lem-w-cancel}(a) with $d = 0$ applied to the single-point kernel
$b_{n'}(m;\cdot)$ (each $m\in A$ lies within $k_A+1$ of an endpoint) together
with Lemma~\ref{lem-w-loc}(b), each term is at most $C\sigma_t^{-1}f(t)h_t$;
hence, by \eqref{eq-w-Delta-ell},
\begin{align*}
p'_A\leq Ck_A\sigma_t^{-1}f(t)h_t = Cf(t)\ell_A\cdot\frac{n'\rho_nh_t}{\sigma_t}\leq C\Delta_A .
\end{align*}
For $O_s$, by the binomial Chernoff bound for $D_1$ given $\xi_1$ with
$\check{\mathfrak D}(\xi_1)\leq(1-c_w)s$,
\begin{align*}
p'_{O_s}\leq P(D_1>s')&\leq P\bigl(\check{\mathfrak D}(\xi_1)>(1-c_w)s\bigr)+\mathfrak e_n
\\ &\leq1-F\bigl((1-c_w)s\bigr)+\mathfrak e_n\leq C\,sf(s) = C\Delta_{O_s},
\end{align*}
by \ref{condition-d3}(iii) at $(1-c_w)s\in J_\infty$ and \ref{condition-d3}(i).

\emph{Short cells.}  With $M^\circ_X := n'\rho_n\int_{[0,1]\setminus\mathscr H}\mathsf K^{(0)}_X(x)\,dx$, write
\begin{align*}
\phi_X(\xi) = n'\rho_n\int_{[0,1]\setminus\mathscr H}\bigl[a(x,\xi)-a(x_t,\xi)\bigr]\mathsf K^{(0)}_X(x)\,dx
 + a(x_t,\xi)\,M^\circ_X .
\end{align*}
By Minkowski's inequality, \eqref{eq-w-a-window} and
Lemma~\ref{lem-w-cancel}(a),(b) with $\theta = \theta_w$ (the part off $W_X$
is negligible, since there $\mathsf K^{(0)}_X$ is negligible and
$\|a(x,\cdot)\|\leq\rho_n^{-1}$), the first term has $L^2$ norm at most
\begin{align*}
Cn\rho_n\cdot t^{1-\theta_w}\sigma_t^{-1}\mu_Xf(t)h_t\bar\ell_X^{\theta_w}
\leq C\,tf(t)\mu_X\Bigl(\frac{\bar\ell_X}t\Bigr)^{\vartheta},
\end{align*}
using $n\rho_nh_t\asymp\sigma_t$ and $\theta_w\geq\vartheta$.  By
Lemma~\ref{lem-w-cancel}(c) with $d = 0$ and $\|a(x_t,\cdot)\|\leq Ct$, the
second is at most $Ctf(t)\mu_X\Lambda_X+\mathfrak e_n$.  Hence
\begin{align*}
\|\phi_X\|\leq Ctf(t)\mu_X\Lambda_X+\mathfrak e_n = C\Xi_X^{1/2}+\mathfrak e_n,
\end{align*}
and Lemma~\ref{lem-w-cell}(a) completes the proof.

\emph{One-sided cells.}  Here no cancellation is available.  By Minkowski's
inequality,
\begin{align*}
\|\phi_{O_s}\|\leq n'\rho_n\int_{[0,1]\setminus\mathscr H}\|a(x,\cdot)\|\,\bigl|\mathsf K^{(0)}_{O_s}(x)\bigr|\,dx .
\end{align*}
On $W_s$, $\|a(x,\cdot)\|\leq Cs$ by \eqref{eq-w-a-window}, and
Lemma~\ref{lem-w-cancel}(b) with $d = 0$, $j = 1$, $\theta = 0$ gives
$n'\rho_n\int|\mathsf K^{(0)}_{O_s}|\leq Cn'\rho_n\sigma_s^{-1}f(s)h_s\leq Cf(s)$;
off $W_s$ the kernel is negligible and $\|a\|\leq\rho_n^{-1}$.  Hence
$\|\phi_{O_s}\|\leq Csf(s)+\mathfrak e_n = C\Xi_{O_s}^{1/2}+\mathfrak e_n$.
\end{proof}

\paragraph{Step B: edge noise.}
\emph{Analogous.}  The signed pairing (Lemma~\ref{lem-c-pairing}), the block
coefficients (Lemma~\ref{lem-c-block}) and the cumulant expansion of
Lemma~\ref{lem-c-edge} are unchanged in structure; the tree sums
(Lemma~\ref{lem-c-tree}) carry $(\rho_nt)^{k-1}f(t)h_t$ in place of
$\rho_n^{k-1}h_t$.  \emph{Additional care.}  In the compact proof the
variance of $\bar\Sigma(A)$ was closed with $Q_A/n\leq\Delta_A/n\leq\Delta_A^2$,
which uses $\Delta_A\gtrsim n^{-1}$ and fails in the tail.  Here the direct
Efron--Stein term is instead compared with $\Delta_X^2$ through the factor
$(nf(t)h_t)^{-1}$, which tends to zero by \eqref{eq-w-poly}; this is the
first place where \ref{condition-s3}(iii) is used.  The same device closes
the several-vertex part of the fourth cumulant.

Let $\delta$ denote, as in the compact proof, the slope of $g^{(j)}_i(X)$ in
$p_{im}$ (a second difference; Lemma~\ref{lem-w-roots}(b) with $k = 2$ and
$I' = \{i,j,m\}$), so that on $\mathscr E_i\cap\{\xi_i\notin\mathscr H\}$,
\begin{align}
\label{eq-w-g-delta}
|g^{(j)}_i(X)|\leq Cm_X\,\omega_X(\xi_i)+\mathfrak e_n,\qquad
|\delta|\leq Cm_X\sigma_t^{-1}\omega^{[2]}_X(\xi_i)+\mathfrak e_n .
\end{align}

\begin{lemma}[Edge covariance, weighted]
\label{lem-w-sigma}
For $n\geq n_0$ and $X\in\mathscr C_n$ with scale $t$,
$E\bar\Sigma(X)\leq CQ_X\leq C\Delta_X$ and $\|\bar\Sigma(X)\|_{L^2}\leq C\Delta_X$.
\end{lemma}

\begin{proof}
\emph{First moment.}  Exactly as in Lemma~\ref{lem-c-sigma}, with
$\check{\mathfrak D}(\xi_1)\omega_X(\xi_1)\leq Ct\,\omega^{[2]}_X(\xi_1)$
(Lemma~\ref{lem-w-loc}(a)) and Lemma~\ref{lem-w-loc}(b):
\begin{align*}
E\bar\Sigma(X)\leq Cn\rho_nt\,m_X^2f(t)h_t = C\frac{n\rho_nt}{\sigma_t^2}\mu_X^2f(t)h_t\leq CQ_X .
\end{align*}

\emph{Variance.}  By the Efron--Stein inequality in the latent positions,
with the direct and indirect effects as in Lemma~\ref{lem-c-sigma}.

\emph{Direct effect.}  Changing $\xi_m$ alters the terms with $m\in\{i,j\}$
by at most $C\rho_nt\,m_X^2\omega_X(\xi_m)$ (using $\sum_jp_{mj}\leq Cn\rho_nt$
on the window of $m$), so the direct part contributes at most
$Cn\rho_n^2t^2m_X^4f(t)h_t$.  In the small regime,
\begin{align*}
n\rho_n^2t^2m_A^4f(t)h_t = n\rho_n^2t^2\sigma_t^{-4}\Bigl(\frac{\ell_A}{h_t}\Bigr)^4f(t)h_t
= \frac{C\ell_A^4f(t)}{nh_t^3}\leq\frac{C\Delta_A^2}{nf(t)h_t};
\end{align*}
in the large regime and for $O_t$ it is $Cf(t)h_t/n\leq C\Delta_X^2/(nf(t)h_t)$
as well, by Lemma~\ref{lem-w-cell}(a).  By \eqref{eq-w-poly},
$(nf(t)h_t)^{-1}\to0$ uniformly, so the direct part is $o(\Delta_X^2)$.

\emph{Indirect effect.}  By \eqref{eq-w-g-delta} the change through $p_{im}$
in $g^{(j)}_i$, $i,j\neq m$, is at most
\begin{align*}
C\rho_n^2t\,m_X^2\sigma_t^{-1}\sum_i\chi_i(\check w_{im}+\check w'_{im})\omega^{[2]}_X(\xi_i).
\end{align*}
Given $\xi_m$, the sum has mean at most $2n\int\check w(x,\xi_m)\omega^{[2]}_X(x)dx$
and variance at most $nE[\check w^2\omega^{[2]2}_X]$; by Minkowski's
inequality, \eqref{eq-w-a-window} and Lemma~\ref{lem-w-loc}, its second
moment is at most $Cn^2t^2(f(t)h_t)^2+Cnt^2f(t)h_t\leq Cn^2t^2(f(t)h_t)^2$,
using $nf(t)h_t\geq1$.  Summing over $m$, the indirect part contributes at most
\begin{align*}
Cn^3\rho_n^4t^4m_X^4\sigma_t^{-2}(f(t)h_t)^2
 = C\bigl(n\rho_ntm_X^2f(t)h_t\bigr)^2\frac{n\rho_n^2t^2}{\sigma_t^2}
 \leq C\rho_nt\,Q_X^2\leq C\rho_nt\,\Delta_X^2 .
\end{align*}
Hence $E\bar\Sigma(X)^2\leq(E\bar\Sigma(X))^2+\mathrm{Var}(\bar\Sigma(X))\leq C\Delta_X^2$.
\end{proof}

Fix $X\in\mathscr C_n$ with scale $t$ and a set $V$ of $k\leq4$ distinct
vertices; $D^{\rm ext}_v$, $s_v$, $E_v$, $\mathcal I$, $\pi_v(s) := P(D^{\rm ext}_v+s\in X\mid\mathcal G)$
and $Y_e := A_e-p_e$ are as in Step~B of the compact proof.

\begin{lemma}[Block coefficients, weighted]
\label{lem-w-block}
Let $v\in V$ with $\xi_v\notin\mathscr H$, work on $\mathscr E_v$, and let $P$ be
a polynomial of degree at most $4$.  Then
$P(\pi_v(s_v)) = \sum_{T\subseteq E_v}c_T\prod_{e\in T}Y_e$ with
$|c_T|\leq C\,m_X\,\sigma_t^{-(|T|-1)}\,\omega^{[|T|]}_X(\xi_v)+\mathfrak e_n$ for
$T\neq\emptyset$, and $c_\emptyset = 0$ if $P(\pi) = \pi-q_v(X)$.
\end{lemma}

\begin{proof}
As Lemma~\ref{lem-c-block}, with Lemma~\ref{lem-w-roots}(b) in place of
Lemma~\ref{lem-c-roots}(b).  For $O_s$,
$\pi_v(s) = P(D^{\rm ext}_v>s'-s\mid\mathcal G)$, so by telescoping to
$+\infty$,
$\Delta^k\pi_v = \pm\mathsf D^{k-1}P(D^{\rm ext}_v = s'-s\mid\mathcal G)$, and the
second bound of Lemma~\ref{lem-w-roots}(b) applies.
\end{proof}

\begin{lemma}[Tree sums, weighted]
\label{lem-w-tree}
Let $F$ be a tree on $\{1,\ldots,k\}$, $k\leq4$, and $d_a\in\{1,\ldots,9\}$.
For $X\in\mathscr C_n$ with scale $t$ and $n\geq n_0$,
\begin{align*}
\sum_{v_1,\ldots,v_k\text{ distinct}}
E\Bigl[\prod_{\{a,b\}\in F}\rho_n\check w(\xi_{v_a},\xi_{v_b})
 \prod_{a=1}^k\chi(\xi_{v_a})\,\omega^{[d_a]}_X(\xi_{v_a})\Bigr]
 \leq C\,n^k(\rho_nt)^{k-1}f(t)\,h_t .
\end{align*}
\end{lemma}

\begin{proof}
As Lemma~\ref{lem-c-tree}: rooting $F$ at $1$ and integrating leaf-first,
by induction
\begin{align*}
\int\rho_n\check w(x,y)\Phi_c(y)\,dy\leq C(\rho_nt)^{|S_c|-1}\rho_n\check{\mathfrak D}(x),
\end{align*}
and $\check{\mathfrak D}(x)^r\omega^{[d]}_X(x)\leq t^r\omega^{[d+r]}_X(x)$
by Lemma~\ref{lem-w-loc}(a); the factor $t$ per edge is the only change.
Integrating the root with Lemma~\ref{lem-w-loc}(b) gives $f(t)h_t$.  Only
$\int\check w(x,y)dy$ enters, so neighbor positions in $\mathscr H$ need no
treatment.
\end{proof}

\begin{lemma}[Edge-noise moments, weighted]
\label{lem-w-edge}
For $X\in\mathscr C_n$ and $n\geq n_0$,
\begin{align}
\begin{split}
\label{eq-w-edge}
E[B_X^2\mid\mathcal G]&\leq\bar p_X+\bar\Sigma(X),\qquad
E\bigl[\kappa_4(B_X\mid\mathcal G)\bigr]\leq Cn^{-1}\Delta_X+\mathfrak e_n,
\\ E\bigl[B_X^4\bigr]&\leq C\bigl(\Delta_X^2+n^{-1}\|Z_X\|_{L^2}^2+n^{-1}\Delta_X\bigr)+\mathfrak e_n .
\end{split}
\end{align}
\end{lemma}

\begin{proof}
The second moment and the reduction of the cumulant by the law of total
cumulance \citep{brillinger1969calculation} are as in Lemma~\ref{lem-c-edge};
we only recompute the count.

\emph{One vertex.}  $|P_4(q)|\leq q$ gives at most
$n^{-1}p'_X\leq Cn^{-1}\Delta_X$ by Proposition~\ref{prop-w-hajek}.

\emph{Several vertices.}  For $k\geq2$, a surviving term has, by
Lemma~\ref{lem-w-block} and $\sum_\beta|T_\beta|\geq2(k-1)$,
\begin{align*}
\prod_\beta|c_{T_\beta}|\leq Cm_X^k\sigma_t^{2-k}\prod_{v\in V}\omega^{[d_v]}_X(\xi_v),
\qquad d_v := \sum_{\beta\text{ at }v}|T_\beta|\leq9,
\end{align*}
and $\prod_{e\in U}p_e\leq\prod_{e\in F}\rho_n\check w_e$ for a spanning tree
$F$.  By Lemma~\ref{lem-w-tree} the tuples with $|V| = k$ contribute at most
\begin{align*}
Cn^{-2}m_X^k\sigma_t^{2-k}\,n^k(\rho_nt)^{k-1}f(t)h_t
 = C\Bigl(\frac n{n''}\Bigr)^{k-1}\frac{\mu_X^kf(t)h_t}n\leq C\frac{Q_X}n\leq Cn^{-1}\Delta_X,
\end{align*}
using $n^{k-2}(\rho_nt)^{k-1}\sigma_t^{2-k} = (n/n'')^{k-1}\sigma_t^k/n$, $\mu_X\leq1$
and Lemma~\ref{lem-w-cell}(a).

\emph{Fourth moment.}  As in Lemma~\ref{lem-c-edge}, with
$E\bar p_X^2 = (p'_X)^2+n^{-1}\|Z_X\|_2^2$, $p'_X\leq C\Delta_X$ and
$E\bar\Sigma(X)^2\leq C\Delta_X^2$ by Lemma~\ref{lem-w-sigma}.
\end{proof}

\paragraph{Step Z2: latent part, second moment.}
\emph{Analogous.}  The Efron--Stein reduction, the exact resampling identity
\eqref{eq-c-exact} (which holds after multiplication by $\chi_i$, since
$\chi_i$ is a function of $\xi_i$ alone), the split of $\Phi_1$ into a main
part and a fluctuation part, and the three classes ($\alpha$), ($\beta$),
($\gamma$) of paired Hoeffding components are those of Lemma~\ref{lem-c-Z2}.
\emph{Additional care.}
\begin{enumerate}[(1)]
\item The profile norms $\|a(x,\cdot)\|$ on a window are of order $t$, not
one, so the geometric ratio of the component series is $C\rho_nt$ rather than
$C\rho_n$, and the truncation level $d_0$ is chosen against $\rho_n\lambda_n$.
\item The H\"older moduli of $C(x,y)$ are the relative ones,
$\varepsilon_X(x) = (e_X(x)/t)^{\theta_w}$, carrying powers of $t$ that must be
matched against $f(t)$ through \ref{condition-d3}(iii).
\item Class ($\beta$) is where \ref{condition-s3}(iv) is born: its leading
term relative to $\Delta_X$ is exactly $v_n(t)$, and the clipping part of
$\Lambda_X$ contributes $v_n(t)\Upsilon_n(t)th_t^{-\vartheta}$, controlled by
\ref{condition-s3}(ii).
\end{enumerate}

Let $\xi'_m$ be an independent copy of $\xi_m$ and $q^{(m)}_i$, $Z^{(m)}_X$ the
values with $\xi_m$ replaced.  As in \eqref{eq-c-exact}, for $i\neq m$,
$\chi_i[q_i(X)-q^{(m)}_i(X)] = \rho_n\chi_i[a(\xi_i,\xi_m)-a(\xi_i,\xi'_m)]g^{(m)}_i(X)$,
so that
\begin{align}
\label{eq-w-Z-diff}
\begin{split}
\sqrt n\,\bigl(Z_X-Z^{(m)}_X\bigr) &= \bigl[\chi_mq_m(X)-\chi'_mq^{(m)}_m(X)\bigr]
\\ &\quad+\rho_n\sum_{i\neq m}\chi_i\bigl[a(\xi_i,\xi_m)-a(\xi_i,\xi'_m)\bigr]g^{(m)}_i(X),
\end{split}
\end{align}
with $\chi'_m := \chi(\xi'_m)$.  Let
$\bar g(x) := E[g^{(1)}_2(X)\mid\xi_2 = x] = \mathsf K^{(0)}_X(x)$ with $M = n''$,
and $e_i := g^{(1)}_i(X)-\bar g(\xi_i)$.

\begin{lemma}[Latent second moment, weighted]
\label{lem-w-Z2}
For $n\geq n_0$ and $X\in\mathscr C_n$ with scale $t$,
$\|Z_X\|_{L^2}^2\leq2p'_X+2E\Phi_1^2$ and
\begin{align}
\label{eq-w-Z2}
\begin{split}
E\Phi_1^2&\leq C\bigl[\Xi_X+\rho_nt\,Q_X+\rho_nt\,\Xi_X+\bar s_n\Delta_X\bigr]+\mathfrak e_n,
\\
\bar s_n &:= \sup_{t_0\leq t\leq\lambda_n}v_n(t)\bigl[1+\Upsilon_n(t)\,t\,h_t^{-\vartheta}\bigr]\to0,
\end{split}
\end{align}
where $v_n(t) := n^{(1-2\vartheta)/4}\rho_n^{(9-2\vartheta)/4}t^{(11-2\vartheta)/4}f(t)^{1/2}$
is the quantity of \ref{condition-s3}(iv).  In particular
$\|Z_X\|_{L^2}^2\leq C(\Delta_X+\Xi_X)$.  The same bound holds for the
$L^2(d\xi)$-valued sum $X^{(-m)} := \rho_n\sum_{i\notin\{1,2,m\}}\chi_ia(\xi_i,\cdot)g^{(1,2,m)}_i$
of Step~Z4 (three vertices removed), and for its centered part,
\begin{align*}
E\|X^{(-m)}-EX^{(-m)}\|^2\leq C[\rho_ntQ_X+\rho_nt\Xi_X+\bar s_n\Delta_X]+\mathfrak e_n .
\end{align*}
\end{lemma}

\begin{proof}
\emph{Reduction.}  By the Efron--Stein inequality, exchangeability and
\eqref{eq-w-Z-diff},
\begin{align*}
\mathrm{Var}(Z_X)\leq E\bigl[(\chi_1q_1-\chi'_1q_1^{(1)})^2\bigr]+E\bigl[(\Phi_1-\Phi'_1)^2\bigr]
\leq2p'_X+2E\Phi_1^2,
\end{align*}
where $\Phi'_1$ is $\Phi_1$ with $\xi'_1$; the cross term $E[\Phi_1\Phi'_1]$
vanishes because $E_\xi a(x,\xi) = 0$ and the $g^{(1)}_i$ do not involve
$\xi_1$.  Split $\Phi_1 = \Phi^{\rm m}+\Phi^{\rm f}$ with
\begin{align*}
\Phi^{\rm m} := \rho_n\sum_{i\neq1}\chi_ia(\xi_i,\xi_1)\bar g(\xi_i),\qquad
\Phi^{\rm f} := \rho_n\sum_{i\neq1}\chi_ia(\xi_i,\xi_1)e_i .
\end{align*}

\emph{Main part.}  Given $\xi_1$, $\Phi^{\rm m}$ is a sum of i.i.d.\ terms with
mean $\phi_X(\xi_1)$, so by Proposition~\ref{prop-w-hajek},
\eqref{eq-w-a-window} and Lemma~\ref{lem-w-cancel}(b) with $d = 0$, $j = 2$,
\begin{align*}
E(\Phi^{\rm m})^2&\leq\|\phi_X\|^2+n\rho_n^2E\bigl[\chi\|a(\xi,\cdot)\|^2\bar g(\xi)^2\bigr]
\\ &\leq C\Xi_X+Cn\rho_n^2t^2\sigma_t^{-2}\mu_X^2f(t)h_t\leq C\Xi_X+C\rho_ntQ_X .
\end{align*}

\emph{Fluctuation part: structure.}  Integrating out $\xi_1$,
\begin{align*}
E(\Phi^{\rm f})^2 &= n'\rho_n^2E\bigl[\chi_2\|a(\xi_2,\cdot)\|^2e_2^2\bigr]
+n'(n'-1)\rho_n^2E\bigl[\chi_2\chi_3\,C(\xi_2,\xi_3)\,e_2e_3\bigr],
\\ C(x,y) &:= \langle a(x,\cdot),a(y,\cdot)\rangle .
\end{align*}
By \eqref{eq-w-g-delta} and Lemma~\ref{lem-w-cancel}(a),
$|e_2|\leq Cm_X\omega_X(\xi_2)$, so the diagonal term is at most
$C\rho_ntQ_X$.  For the off-diagonal term, Lemma~\ref{lem-multilinear}
applied to $g^{(1)}_2$ given $\xi_2$ gives the components
\begin{align*}
\Pi(g^{(1)}_2;S) = \rho_n^{|S|}\prod_{j\in S}a(\xi_2,\xi_j)K_{|S|}(\xi_2),
\qquad S\subseteq[n]\setminus\{1,2\},
\end{align*}
with $K_d := \mathsf K^{(d)}_X$ and $M = n''-d$, each centered in every
$\xi_j$, $j\in S$.  A pair of components $(S,S')$ of $(e_2,e_3)$ contributes
only if $S\setminus\{3\} = S'\setminus\{2\} =: T$, leaving the three classes of
Lemma~\ref{lem-c-Z2}, with $d := |T|$ and $\Gamma_k(x,y) := C(x,y)^k$:
\begin{align*}
(\alpha)\quad&\rho_n^{2d}\iint\Gamma_{d+1}K_d(x)K_d(y),\qquad d\geq1;\\
(\beta)\quad&\rho_n^{2d+1}\iint a(x,y)\Gamma_{d+1}K_{d+1}(x)K_d(y),\qquad d\geq1,\ \text{and its mirror};\\
(\gamma)\quad&\rho_n^{2d+2}\iint a(x,y)a(y,x)\Gamma_{d+1}K_{d+1}(x)K_{d+1}(y),\qquad d\geq0.
\end{align*}
Here $x,y$ range over $[0,1]\setminus\mathscr H$, there are at most $n^d/d!$
sets $T$, and the prefactor is $n'(n'-1)\rho_n^2$.

\emph{Preliminary bounds.}  For $x,y\in W_X$ and $1\leq k\leq d_0+2$, by
\eqref{eq-w-a-window} and bilinearity of $C$:
\begin{enumerate}[($\Gamma$1)]
\item $|\Gamma_k(x,y)|\leq(Ct)^{2k}$;
\item $|\Gamma_k(x,y)-\Gamma_k(x_t,y)|\leq k(Ct)^{2k}\varepsilon_X(x)$, and
symmetrically in $y$;
\item $|\Gamma_k(x,y)-\Gamma_k(x_t,y)-\Gamma_k(x,y_t)+\Gamma_k(x_t,y_t)|
\leq k^2(Ct)^{2k}\varepsilon_X(x)\varepsilon_X(y)$, where $y_t := x_t$;
\item $\int_0^1a(x',y)^2dx'\leq Ct^2$ and
$\int_0^1(a(x',y)-a(x',y_t))^2dx'\leq(Ct\,\varepsilon_X(y))^2$.
\end{enumerate}
For ($\Gamma$3), with $u := C(x,y)$, $u_1 := C(x_t,y)$, $u_2 := C(x,y_t)$,
$u_{12} := C(x_t,y_t)$ and $P(\alpha,\beta) := \sum_{j<k}\alpha^j\beta^{k-1-j}$,
the mixed difference equals
\begin{align*}
(u-u_1-u_2+u_{12})P(u,u_1)+(u_2-u_{12})\bigl[P(u,u_1)-P(u_2,u_{12})\bigr],
\end{align*}
with $u-u_1-u_2+u_{12} = \langle a(x)-a(x_t),a(y)-a(y_t)\rangle$.  For
($\Gamma$4), $a(x',y) = \check w(y,x')-\check{\mathfrak D}(x')$ by symmetry, so
$\int a(x',y)^2dx'\leq2\|\check w(y,\cdot)\|^2+2E\mathfrak D^2\leq Ct^2$, and
$a(x',y)-a(x',y_t) = \check w(y,x')-\check w(y_t,x')$ is bounded as in
\eqref{eq-w-a-window}.  By Lemma~\ref{lem-w-cancel}, using
$(n'\rho_n)^{-1}\sigma_t^{-d} = (n''/n')h_t\sigma_t^{-(d+1)}$:
\begin{enumerate}[(K1)]
\item $\sup_{x\notin\mathscr H}|K_d(x)|\leq C\sigma_t^{-(d+1)}\mu_X+\mathfrak e_n$ and
$\int_{[0,1]\setminus\mathscr H}K_d^2\leq C\sigma_t^{-2(d+1)}\mu_X^2f(t)h_t$;
\item the integrals satisfy
\begin{align*}
\int_{[0,1]\setminus\mathscr H}\varepsilon_X|K_d|&\leq C\sigma_t^{-(d+1)}\mu_Xf(t)h_t\Bigl(\frac{\bar\ell_X}t\Bigr)^\vartheta,
\\
\int_{[0,1]\setminus\mathscr H}|K_d|&\leq C\sigma_t^{-(d+1)}\mu_Xf(t)h_t ;
\end{align*}
\item for $d\geq1$,
$\displaystyle\Bigl|\int_{[0,1]\setminus\mathscr H}K_d\Bigr|\leq C\sigma_t^{-(d+1)}\mu_Xf(t)h_t\Lambda_X+\mathfrak e_n$.
\end{enumerate}
$K_d$ is negligible on $[0,1]\setminus(\mathscr H\cup W_X)$ and
$|\Gamma_k|\leq\rho_n^{-2k}$ everywhere, so contributions with $x$ or
$y\notin W_X$ are negligible and are omitted below.

\emph{Truncation.}  Choose $d_0$ so large that
$n^{3+K_0}(C_*\rho_n\lambda_n)^{d_0/4}\to0$, with $C_*$ below, and set
$d_1 := d_0+2$ in \eqref{eq-w-K}; this is possible since
$\rho_n\lambda_n\leq n^{-b/2}$, and it is the one place where the truncation
depends on the cutoff.  Split $e_i = e^\leq_i+e^>_i$ with
$e^\leq_i := \sum_{1\leq|S|\leq d_0+1}\Pi(g^{(1)}_i;S)$.  Given $\xi_i$ the
components are orthogonal, so
\begin{align*}
E\bigl[(e^\leq_i)^2\mid\xi_i\bigr]\leq E\bigl[e_i^2\mid\xi_i\bigr]\leq E\bigl[|g^{(1)}_i|\bigm|\xi_i\bigr],
\end{align*}
which is $b_{n''}(r';\rho_n\check{\mathfrak D}(\xi_i))+b_{n''}(t';\rho_n\check{\mathfrak D}(\xi_i))$
(one point for $O_s$), negligible for $\xi_i\notin W_X$; this handles the
off-window positions, where $\|a(\xi_i,\cdot)\|$ is not bounded.  For
$x\in W_X$,
\begin{align*}
E\bigl[(e^>_i)^2\bigm|\xi_i = x\bigr]
 = \sum_{d>d_0+1}\binom{n''}{d}\rho_n^{2d}\|a(x,\cdot)\|^{2d}K_d(x)^2 .
\end{align*}
For $d\leq n''/2$, Lemma~\ref{lem-binomial-facts}\ref{binom-fact-c} gives
$|K_d(x)|\leq2\Gamma_{d+1}(\sigma_t/3)^{-(d+1)}$ (as $\sigma(n''-d,p)\geq\sigma_t/3$
on the window), and with $\Gamma_{d+1}^2\leq C_0^{d+1}(d+1)!$ and
$\binom{n''}d\leq n^d/d!$ the $d$-th term is at most
$C\sigma_t^{-2}(d+1)(C_*\rho_nt)^d$; for $d>n''/2$,
$|K_d|\leq2^{d+1}$ (Lemma~\ref{lem-binomial-facts}\ref{binom-fact-d}) and
$\binom{n''}d\rho_n^{2d}(Ct)^{2d}\leq(C\rho_n^2t^2)^d$, which is smaller.  Hence
$\|\chi_ie^>_i\|_{L^2}^2\leq C\sigma_t^{-2}(C_*\rho_nt)^{d_0/2}+\mathfrak e_n$,
and by Cauchy--Schwarz and ($\Gamma$1),
\begin{align*}
n'(n'-1)\rho_n^2\bigl|E[\chi_2\chi_3C(\xi_2,\xi_3)(e_2e_3-e^\leq_2e^\leq_3)]\bigr|
 &\leq Cn^2\rho_n^2t^2\|\chi_2e_2\|_2\|\chi_2e^>_2\|_2+\mathfrak e_n
 \\&\leq Cn^3(C_*\rho_n\lambda_n)^{d_0/4}+\mathfrak e_n ,
\end{align*}
which is $o(\Delta_X)$ uniformly by \eqref{eq-w-poly}.  In
$E[\chi_2\chi_3Ce^\leq_2e^\leq_3]$ only the classes with $d\leq d_0+1$ occur,
and \eqref{eq-w-K} covers $M = n''-d$ for these.

\emph{Class ($\alpha$).}  Write
\begin{align*}
\Gamma_{d+1}(x,y) &= \mathsf M(x,y)+\bigl[\Gamma_{d+1}(x_t,y)-\Gamma_{d+1}(x_t,y_t)\bigr]
\\ &\quad+\bigl[\Gamma_{d+1}(x,y_t)-\Gamma_{d+1}(x_t,y_t)\bigr]+\Gamma_{d+1}(x_t,y_t),
\end{align*}
with $\mathsf M$ the mixed difference of ($\Gamma$3).  The four terms integrate
against $K_d(x)K_d(y)$ to at most
\begin{align*}
Ct^{2d+2}\Bigl(\int\varepsilon_X|K_d|\Bigr)^2,\qquad
Ct^{2d+2}\Bigl|\int K_d\Bigr|\int\varepsilon_X|K_d|\ \ \text{(twice)},\qquad
Ct^{2d+2}\Bigl|\int K_d\Bigr|^2,
\end{align*}
each at most $Ct^{2d+2}(\sigma_t^{-(d+1)}\mu_Xf(t)h_t\Lambda_X)^2$ by (K2), (K3).
With the $n^d/d!$ sets $T$ and the prefactor, class ($\alpha$) contributes at
most
\begin{align*}
Cn^{d+2}\rho_n^{2d+2}t^{2d+2}\sigma_t^{-2d-2}\mu_X^2(f(t)h_t)^2\Lambda_X^2
 &= C\Bigl(\frac{n\rho_n^2t^2}{\sigma_t^2}\Bigr)^d\frac{n^2\rho_n^2t^2h_t^2}{\sigma_t^2}f(t)^2\mu_X^2\Lambda_X^2
\\ &\leq C(C\rho_nt)^d\,\Xi_X ,
\end{align*}
using $n^2\rho_n^2t^2h_t^2/\sigma_t^2 = (n/n'')^2t^2$.  Summed over $d\leq d_0+1$
this is at most $C\rho_nt\,\Xi_X$.

\emph{Class ($\beta$).}  Write
\begin{align*}
a(x,y)\Gamma_{d+1}(x,y) &= \bigl[a(x,y)-a(x,y_t)\bigr]\Gamma_{d+1}(x,y)
+a(x,y_t)\bigl[\Gamma_{d+1}(x,y)-\Gamma_{d+1}(x,y_t)\bigr]
\\&\quad+a(x,y_t)\Gamma_{d+1}(x,y_t).
\end{align*}
For fixed $y\in W_X$, Cauchy--Schwarz in $x$ with ($\Gamma$1), ($\Gamma$2),
($\Gamma$4) and (K1) bounds the $x$-integral of each of the first two terms
against $|K_{d+1}(x)|$ by
\begin{align*}
Ct^{2d+3}\varepsilon_X(y)\sigma_t^{-(d+2)}\mu_X(f(t)h_t)^{1/2};
\end{align*}
integrating
against $|K_d(y)|$ with (K2) gives
\begin{align*}
Ct^{2d+3}\sigma_t^{-(2d+3)}\mu_X^2(f(t)h_t)^{3/2}(\bar\ell_X/t)^\vartheta .
\end{align*}
The third term factorizes, and by (K3), ($\Gamma$4) and (K1),
\begin{align*}
\Bigl|\int K_d(y)\,dy\cdot\int a(x,y_t)\Gamma_{d+1}(x,y_t)K_{d+1}(x)\,dx\Bigr|
\leq C\sigma_t^{-(d+1)}\mu_Xf(t)h_t\Lambda_X\cdot Ct^{2d+3}\sigma_t^{-(d+2)}\mu_X(f(t)h_t)^{1/2}.
\end{align*}  With the $n^d$ sets $T$, the prefactor and the
mirror class, class ($\beta$) contributes at most $\sum_{d\geq1}B_d$, where
\begin{align*}
B_d &:= Cn^{d+2}\rho_n^{2d+3}t^{2d+3}\sigma_t^{-(2d+3)}\mu_X^2(f(t)h_t)^{3/2}\Lambda_X,
\qquad B_{d+1}\leq C\rho_ntB_d,
\\ B_1 &= C\Bigl(\frac n{n''}\Bigr)^{5/2}n^{1/2}\rho_n^{5/2}t^{5/2}\mu_X^2(f(t)h_t)^{3/2}\Lambda_X .
\end{align*}
In the small regime $\mu_X^2f(t)h_t = f(t)\ell_A^2/h_t\leq C\Delta_X$ and
$\Lambda_X = (h_t/t)^\vartheta+\Upsilon_n(t)$; in the large regime and for $O_t$,
$\mu_X = 1$, $\Delta_X\geq cf(t)\bar\ell_X$ and
$(h_t/\bar\ell_X)(\bar\ell_X/t)^\vartheta\leq(h_t/t)^\vartheta$.  In all cases,
using $h_t = (t/(n''\rho_n))^{1/2}$ and $(t/h_t)^\vartheta\leq t\,h_t^{-\vartheta}t_0^{\vartheta-1}$,
\begin{align*}
\frac{B_1}{\Delta_X}\leq Cn^{1/2}\rho_n^{5/2}t^{5/2}(f(t)h_t)^{1/2}
 \Bigl[\Bigl(\frac{h_t}t\Bigr)^\vartheta+\Upsilon_n(t)\Bigr]
 \leq Cv_n(t)\bigl[1+\Upsilon_n(t)\,t\,h_t^{-\vartheta}\bigr]\leq C\bar s_n .
\end{align*}
By \ref{condition-s3}(ii) and (iv), $\bar s_n\to0$.

\emph{Class ($\gamma$).}  By $|a(x,y)a(y,x)|\leq\frac12(a(x,y)^2+a(y,x)^2)$,
symmetry, ($\Gamma$1), ($\Gamma$4), (K1) and (K2), the $d$-th term is at most
\begin{align*}
C\rho_n^{2d+2}t^{2d+4}\sup|K_{d+1}|\int|K_{d+1}|\leq C\rho_n^{2d+2}t^{2d+4}\sigma_t^{-2(d+2)}\mu_X^2f(t)h_t ;
\end{align*}
with the $n^d$ sets $T$ and the prefactor, class ($\gamma$) contributes at most
\begin{align*}
C\Bigl(\frac{n\rho_n^2t^2}{\sigma_t^2}\Bigr)^{d+2}\mu_X^2f(t)h_t\leq C(C\rho_nt)^{d+2}Q_X,
\end{align*}
summed to $C(\rho_nt)^2Q_X$.

\emph{Conclusion.}  Collecting,
\begin{align*}
E\Phi_1^2\leq C\bigl[\Xi_X+\rho_ntQ_X+\rho_nt\Xi_X+\bar s_n\Delta_X\bigr]+\mathfrak e_n,
\end{align*}
and $p'_X\leq C\Delta_X$ gives the bound on $\|Z_X\|_2^2$.  For $X^{(-m)}$,
write its centered part as $(X^{\rm m}-EX^{\rm m})+X^{\rm f}$ with
$X^{\rm m} := \rho_n\sum\chi_ia_i\bar g(\xi_i)$ and $X^{\rm f} := \rho_n\sum\chi_ia_ie_i$,
where $EX^{\rm f} = 0$ since $E[e_i\mid\xi_i] = 0$; the first piece is bounded
by the variance term of the main part ($C\rho_ntQ_X$), the second by the
fluctuation part.
\end{proof}

\paragraph{Step Z4: latent part, fourth moment.}
\emph{Analogous.}  The moment inequality of \citet{boucheron2005moment}, the
decomposition $V\leq4(V_1+V_2)+2(V_3+V_4)$, the reduction of
$\mathrm{Cov}(\Phi_1^2,\Phi_2^2)$ and $\mathrm{Cov}(\Psi_1,\Psi_2)$ to
$\mathrm{Var}(\Psi^*)$, and the exact Efron--Stein decomposition of
$\Psi^* = \|X\|^2$ into direct, indirect (mean operator, Hilbert--Schmidt,
operator bound) and quartic parts are those of Lemma~\ref{lem-c-Z4}.
\emph{Additional care.}
\begin{enumerate}[(1)]
\item The $L^4$ norms of the profiles grow like $t^{q_4}$ by
\ref{condition-l3}, so every remainder ratio carries a weight
$t^{\cdots}f(t)^{\cdots}$; constraint (iv) of \eqref{eq-w-theta-prime} is
exactly what keeps these weights bounded, through \ref{condition-d3}(iii).
\item The product $P_X := (\rho_nt)^2(tf(t)\mu_X)^2\Xi_X$ replaces
$\rho_n^2\mu_A^2\Xi_A$; its comparison with $\Delta_X^{4\vartheta'}$ is where
constraint (v), that is \ref{condition-s3}(v), enters, and the weight is
again bounded by constraint (iv).
\item The class-($\beta$) term in the Hilbert--Schmidt bound carries
$(\rho_n^3t/(nf(t)^2))^{1/4}$, which tends to zero by
\ref{condition-s3}(iii); on a compact set this factor is $(\rho_n^3/n)^{1/4}$
and needs nothing.
\item The event $\mathscr G$ in the operator bound requires
$nf(t)h_t\geq n^{c_1}$, again \eqref{eq-w-poly}.
\end{enumerate}

\begin{lemma}[Latent fourth moment, weighted]
\label{lem-w-Z4}
$\|Z_X\|_{L^4}\leq C\Delta_X^{\vartheta'}$ for every $X\in\mathscr C_n$ and $n\geq n_0$.
\end{lemma}

\begin{proof}
By \citet[Theorem~2]{boucheron2005moment} with $q = 4$,
$\|Z_X\|_4^2\leq C\|V\|_{L^2}$ with $V := \sum_mE'_m[(Z_X-Z^{(m)}_X)^2]$.  By
\eqref{eq-w-Z-diff} and $E'_m[(\Phi_m-\Phi'_m)^2] = \Phi_m^2+\Psi_m$,
\begin{align*}
V\leq4(V_1+V_2)+2(V_3+V_4),
\end{align*}
with $V_1 := n^{-1}\sum_m\chi_mq_m^2$, $V_2$ its resampled analogue,
$V_3 := n^{-1}\sum_m\Phi_m^2$ and $V_4 := n^{-1}\sum_m\Psi_m$.  By
Lemma~\ref{lem-w-Z2}, $EV\leq C(\Delta_X+\Xi_X)\leq C\Delta_X^{2\vartheta}$, so it
suffices to show $\mathrm{Var}(V_k)\leq C\Delta_X^{4\vartheta'}$ for each $k$.
Throughout,
\begin{align*}
M := \|\Phi_1\|_{L^4}\leq Ct^{q_4}f(t)\mu_X,\qquad \Delta_X\geq cf(t)h_t\mu_X,
\end{align*}
the first by Minkowski's inequality in $L^4(d\xi_1)$, \eqref{eq-w-a-window},
\eqref{eq-w-g-delta} and $E(\sum_i\chi_i\omega^{[q_4+1]}_X(\xi_i))^4\leq C(nf(t)h_t)^4$,
the second by Lemma~\ref{lem-w-cell}(a).

\emph{$V_1,V_2$.}  Resampling $\xi_m$ changes $nV_1$ by
\begin{align*}
\bigl(\chi_mq_m^2-\chi'_mq'^2_m\bigr)+\sum_{i\neq m}\chi_i\bigl(q_i+q^{(m)}_i\bigr)\bigl(q_i-q^{(m)}_i\bigr).
\end{align*}
By Cauchy--Schwarz with weights $|g^{(m)}_i|$, \eqref{eq-w-g-delta},
\eqref{eq-w-a-window} and Lemma~\ref{lem-w-loc},
\begin{align*}
E\Bigl[\Bigl(\sum_{i\neq m}\chi_i(q_i+q^{(m)}_i)(q_i-q^{(m)}_i)\Bigr)^2\Bigr]
&\leq Cn^2\rho_n^2m_X^2(f(t)h_t)\int\|a(x,\cdot)\|^2\omega^{[2]}_X(x)dx
\\ &\leq Ct^2f(t)^2\mu_X^2,
\end{align*}
so $\mathrm{Var}(V_1)\leq Cn^{-1}(p'_X+t^2f(t)^2\mu_X^2)$.  Here
$n^{-1}p'_X\leq\Delta_X^{4\vartheta'}$ by \eqref{eq-w-absorb}, and
\begin{align*}
n^{-1}t^2f(t)^2\mu_X^2\leq Cn^{-1}t^2h_t^{-2}\Delta_X^2 = C\rho_nt\,\Delta_X^2\leq\Delta_X^{4\vartheta'};
\end{align*}
for $O_s$ use instead $n^{-1}(sf(s))^2\leq Cn^{-1}\Delta_{O_s}\leq\Delta_{O_s}^{4\vartheta'}$.
The same holds for $V_2$.

\emph{$V_3$: the diagonal.}
$\mathrm{Var}(V_3)\leq n^{-1}E\Phi_1^4+|\mathrm{Cov}(\Phi_1^2,\Phi_2^2)|$, and
\begin{align*}
n^{-1}M^4\leq Cn^{-1}(t^{q_4}f(t)\mu_X)^4\leq Cn^{-1}(n\rho_n)^{2\vartheta'}t^{4q_4-2\vartheta'}f(t)^{4-4\vartheta'}\Delta_X^{4\vartheta'}
\leq C\rho_n\Delta_X^{4\vartheta'},
\end{align*}
using $\mu_X^4\leq\mu_X^{4\vartheta'}$, $h_t^{-4\vartheta'} = (n''\rho_n/t)^{2\vartheta'}$,
$n^{-1}(n\rho_n)^{2\vartheta'}\leq\rho_n$ and $\gamma_-\geq2q_4$.  For $O_s$,
\begin{align*}
n^{-1}\bigl(s^{q_4}f(s)\bigr)^4\leq Cn^{-1}\Delta_{O_s},\qquad\text{since } s^{4q_4-1}f(s)^3\leq C .
\end{align*}

\emph{$V_3$: reduction to $\mathrm{Var}(\Psi^*)$.}  Let $g^{(1,2)}_i$ be
computed with vertices $1,2$ removed, and let
\begin{align*}
\Phi^*_l := \rho_n\sum_{i\geq3}\chi_ia(\xi_i,\xi_l)g^{(1,2)}_i\quad(l = 1,2),
\qquad\mathcal F := \sigma(\xi_3,\ldots,\xi_n).
\end{align*}  Given $\mathcal F$,
$\Phi^*_1,\Phi^*_2$ are i.i.d., so
\begin{align*}
\mathrm{Cov}(\Phi_1^{*2},\Phi_2^{*2}) &= \mathrm{Var}(\Psi^*),\qquad
\Psi^* := E[\Phi_1^{*2}\mid\mathcal F] = \|X\|^2,
\\
X &:= \rho_n\sum_{i\geq3}\chi_ia(\xi_i,\cdot)g^{(1,2)}_i\in L^2(d\xi).
\end{align*}
The remainder is $r_1 := \Phi_1-\Phi^*_1 = r^{(1)}+r^{(2)}$, with
\begin{align*}
r^{(1)} := \rho_n\chi_2a(\xi_2,\xi_1)g^{(1)}_2,\qquad
r^{(2)} := \rho_n\sum_{i\geq3}\chi_ia(\xi_i,\xi_1)p_{i2}\delta_i,
\end{align*}
$\delta_i$ the slope of $g^{(1)}_i$ in $p_{i2}$, which does not involve
$\xi_1,\xi_2$.  By \eqref{eq-w-g-delta} and Lemma~\ref{lem-w-loc}:
\begin{enumerate}[(r1)]
\item since $|a|\leq\rho_n^{-1}$,
\begin{align*}
E\bigl[a(\xi_2,\xi_1)^4\bigm|\xi_2\bigr]\leq\rho_n^{-2}\|a(\xi_2,\cdot)\|^2,
\qquad\text{so}\qquad
\|r^{(1)}\|_4\leq C\rho_n^{1/2}t^{1/2}m_X(f(t)h_t)^{1/4};
\end{align*}
\item given $(\xi_l)_{l\geq3}$, Minkowski's inequality in $L^4(d\xi_1d\xi_2)$
gives
\begin{align*}
\|r^{(2)}\|\leq\rho_n\sum_i\chi_i|\delta_i|\,\|a(\xi_i,\cdot)\|_{L^4}\,\|p_{i2}\|_{L^4(d\xi_2)},
\end{align*}
with $\|p_{i2}\|_{L^4(d\xi_2)}\leq\min\{\rho_nK_4\mathfrak D(\xi_i)^{q_4},1\}$ by
\ref{condition-l3}; with $N := \sum_i\chi_i\omega^{[2q_4+2]}_X(\xi_i)$ and
$EN^4\leq C(nf(t)h_t)^4$, this gives
$\|r^{(2)}\|_4\leq C\min\{\rho_nt^{2q_4-1},t^{q_4-1}\}\mu_Xf(t)h_t$.
\end{enumerate}
Expanding $\Phi_l^2 = \Phi_l^{*2}+2\Phi^*_lr_l+r_l^2$, H\"older's inequality
and $x^3y+x^2y^2+xy^3\leq C(x^3y+y^4)$ give
\begin{align*}
\bigl|\mathrm{Cov}(\Phi_1^2,\Phi_2^2)-\mathrm{Var}(\Psi^*)\bigr|\leq C\bigl(\|r_1\|_4M^3+\|r_1\|_4^4\bigr).
\end{align*}
Using the first branch of the minimum in (r2) for the term with $M^3$ and the
second for the fourth power, $\mu_X\leq1$, $h_t\leq1$ and
$\Delta_X^{4\vartheta'}\geq c(f(t)h_t\mu_X)^{4\vartheta'}$:
\begin{align*}
\frac{\|r^{(1)}\|_4M^3}{\Delta_X^{4\vartheta'}}&\leq C\rho_n^{1/2}\,t^{3q_4-1/2}f(t)^{13/4-4\vartheta'},&
\frac{\|r^{(2)}\|_4M^3}{\Delta_X^{4\vartheta'}}&\leq Cn^{\frac{4\vartheta'-1}2}\rho_n^{\frac{4\vartheta'+1}2}\,t^{5q_4-1}f(t)^{4-4\vartheta'},\\
\|r^{(1)}\|_4^4&\leq Cn^{-2}\Delta_X,&
\frac{\|r^{(2)}\|_4^4}{\Delta_X^{4\vartheta'}}&\leq C(n\rho_n)^{2\vartheta'-2}\,t^{4q_4}f(t)^{4-4\vartheta'},
\end{align*}
where negative powers of $t$ have been absorbed into $C$ (as $t\geq t_0$).
For the first we used $h_t^{5/4-4\vartheta'}\leq1$ ($\vartheta'\leq5/16$); for
the third, $\rho_n^2t^2m_X^4f(t)h_t\leq C(\rho_nt)^2\sigma_t^{-4}Q_X$; for the
fourth, $(t^{q_4-1}\mu_Xf(t)h_t)^4 = (t^{q_4}f(t)\mu_X)^4\sigma_t^{-4}$ and
$h_t^{-4\vartheta'} = (n''\rho_n/t)^{2\vartheta'}$.  The weights are bounded by
\ref{condition-d3}(iii) and constraint (iv) of \eqref{eq-w-theta-prime}, the
prefactors tend to $0$ by constraint (iii) and $\vartheta'<1$, and
$n^{-2}\Delta_X\leq\Delta_X^{4\vartheta'}$ by \eqref{eq-w-absorb}.  Hence
$|\mathrm{Cov}(\Phi_1^2,\Phi_2^2)-\mathrm{Var}(\Psi^*)| = o(\Delta_X^{4\vartheta'})$
uniformly.

\emph{$V_4$: reduction to $\mathrm{Var}(\Psi^*)$.}  Since $g^{(m)}_i$ does not
involve $\xi_m$, $\Psi_m = E[\Phi_m^2\mid\xi_l,l\neq m]$, so
$\mathrm{Var}(V_4)\leq n^{-1}E\Phi_1^4+|\mathrm{Cov}(\Psi_1,\Psi_2)|$.  By
multilinearity $g^{(1)}_i = g^{(1,2)}_i+p_{i2}\delta_i$ exactly, so
\begin{align*}
X_1 := \rho_n\sum_{i\neq1}\chi_ia(\xi_i,\cdot)g^{(1)}_i = X+\tilde r_1,\qquad \tilde r_1(\xi_1) = r_1,
\end{align*}
and $\Psi_1 = \Psi^*+e_1$ with $|e_1|\leq2\|X\|\|\tilde r_1\|+\|\tilde r_1\|^2$;
likewise for $\Psi_2$.  Hence, using $E\|X\|^4\leq E\Phi_1^{*4}$ and
$E\|\tilde r_1\|_{L^2(d\xi)}^4\leq Er_1^4$ (Jensen in $d\xi$),
\begin{align*}
\bigl|\mathrm{Cov}(\Psi_1,\Psi_2)-\mathrm{Var}(\Psi^*)\bigr|\leq2\|\Psi^*\|_2\|e_1\|_2+\|e_1\|_2^2
\leq C\bigl(M^3\|r_1\|_4+\|r_1\|_4^4\bigr),
\end{align*}
which is the remainder just bounded.

\emph{Exact Efron--Stein decomposition.}  For $m\geq3$ let
$\mathcal F_{-m} := \sigma(\xi_l:l\geq3,l\neq m)$ and
$X^{(-m)} := \rho_n\sum_{i\geq3,i\neq m}\chi_ia(\xi_i,\cdot)g^{(1,2,m)}_i$, which is
$\mathcal F_{-m}$-measurable.  By multilinearity,
$g^{(1,2)}_i = g^{(1,2,m)}_i+p_{im}\delta_i$ with $\delta_i$ now the second
difference with vertices $1,2,m$ removed, so $X = X^{(-m)}+Y(\xi_m)+V(\xi_m)$,
where, given $\mathcal F_{-m}$,
\begin{align*}
Y(y) := \rho_n\chi(y)a(y,\cdot)g_m(y),\qquad
V(y) := \rho_n^2\sum_{i\geq3,i\neq m}\chi_ia(\xi_i,\cdot)\check w(\xi_i,y)\delta_i,
\end{align*}
$g_m(y)$ denoting $g^{(1,2)}_m$ with $\xi_m = y$.  With $R_m := Y(\xi_m)+V(\xi_m)$
and $R'_m$ its resampled copy, the Efron--Stein inequality gives
\begin{align}
\begin{split}
\label{eq-w-Psi-ES}
\mathrm{Var}(\Psi^*)&\leq\sum_{m\geq3}\Bigl(4E\langle X^{(-m)},R_m-R'_m\rangle^2+4E\|R_m\|^4\Bigr),
\\
E\bigl[\langle X^{(-m)},R_m-R'_m\rangle^2\bigm|\mathcal F_{-m}\bigr]
 &\leq4E_y\langle X^{(-m)},Y(y)\rangle^2+4\mathrm{Var}_y\langle X^{(-m)},V(y)\rangle,
 \end{split}
\end{align}
with $y$ a fresh latent position.

\emph{Direct part.}  For $y\notin\mathscr H$ on $\mathscr E_m$,
$|g_m(y)|\leq Cm_X\omega_X(y)$ and $\|a(y,\cdot)\|\omega_X(y)\leq Ct\,\omega^{[2]}_X(y)$,
so $E_y\langle X^{(-m)},Y(y)\rangle^2\leq C\rho_n^2t^2m_X^2f(t)h_t\|X^{(-m)}\|^2$
with a deterministic factor, and by $n\rho_n^2t^2m_X^2 = (n/n'')\rho_nt\mu_X^2$
and Lemma~\ref{lem-w-Z2},
\begin{align*}
\sum_m4E\bigl[E_y\langle X^{(-m)},Y(y)\rangle^2\bigr]
 \leq Cn\rho_n^2t^2m_X^2f(t)h_tE\|X^{(-m)}\|^2\leq C\rho_nt\,Q_X(\Delta_X+\Xi_X).
\end{align*}

\emph{Indirect part: an exact identity.}  Since
$\check w(x,y) = \check{\mathfrak D}(x)+a(x,y)$ and $E_ya(x,y) = 0$,
\begin{align}
\label{eq-w-T}
\mathrm{Var}_y\langle X^{(-m)},V(y)\rangle = \rho_n^4\|\mathcal T_mX^{(-m)}\|^2,
\qquad
\mathcal T_m := \sum_{i\geq3,i\neq m}\chi_i\delta_i\,a(\xi_i,\cdot)\otimes a(\xi_i,\cdot),
\end{align}
where $(u\otimes u)h := \langle u,h\rangle u$.  With $\bar{\mathcal T} := E\mathcal T_m$,
$\tilde\varphi := EX^{(-m)}$ and $\tilde X := X^{(-m)}-\tilde\varphi$,
\begin{align}
\label{eq-w-T-split}
\|\mathcal T_mX^{(-m)}\|^2\leq3\|\bar{\mathcal T}X^{(-m)}\|^2+3\|(\mathcal T_m-\bar{\mathcal T})\tilde\varphi\|^2+3\|(\mathcal T_m-\bar{\mathcal T})\tilde X\|^2 .
\end{align}

\emph{Mean operator.}  We have
\begin{align*}
\bar{\mathcal T} = (n-3)\int_{[0,1]\setminus\mathscr H}\bar\delta(x)\,a(x,\cdot)\otimes a(x,\cdot)\,dx,
\end{align*}
with $\bar\delta = \mathsf K^{(1)}_{X,n-4}(\rho_n\check{\mathfrak D}(\cdot))$ up
to sign.  With $a_t := a(x_t,\cdot)$, split
\begin{align*}
a(x)\otimes a(x) = a_t\otimes a_t+\bigl[a(x)\otimes a(x)-a_t\otimes a_t\bigr].
\end{align*}  The second part has operator norm at most
$Ct^2\varepsilon_X(x)$ on $W_X$ by \eqref{eq-w-a-window}, so by
Lemma~\ref{lem-w-cancel}(b) its integral is at most
$Ct^2\sigma_t^{-2}\mu_Xf(t)h_t(\bar\ell_X/t)^\vartheta$; the first part is
bounded by Lemma~\ref{lem-w-cancel}(c) with $d = 1$.  Since
$nh_t\sigma_t^{-2} = (n/n'')(\rho_n\sigma_t)^{-1}$, both give
\begin{align*}
\|\bar{\mathcal T}\|_{\rm op}&\leq\frac{Ct^2f(t)\mu_X\Lambda_X}{\rho_n\sigma_t}+\mathfrak e_n,
\qquad\text{hence}
\\
n\rho_n^4\|\bar{\mathcal T}\|_{\rm op}^2&\leq C\rho_nt\,f(t)^2t^2\mu_X^2\Lambda_X^2\leq C\rho_nt\,\Xi_X,
\end{align*}
and the first term of \eqref{eq-w-T-split} contributes at most
$C\rho_nt\,\Xi_X(\Delta_X+\Xi_X)+\mathfrak e_n$.

\emph{Hilbert--Schmidt bound.}  Write $a_i := a(\xi_i,\cdot)$,
$\delta_i = \bar\delta(\xi_i)+e^\delta_i$, and
$\mathcal T_m-\bar{\mathcal T} = \mathcal T^{\rm i}+\mathcal T^{\rm f}$ with
\begin{align*}
\mathcal T^{\rm i} := \sum_i\chi_i\bar\delta(\xi_i)a_i\otimes a_i-\bar{\mathcal T},\qquad
\mathcal T^{\rm f} := \sum_i\chi_ie^\delta_ia_i\otimes a_i .
\end{align*}
We show, without cancellation,
\begin{align}
\label{eq-w-T-HS}
n\rho_n^4E\|\mathcal T_m-\bar{\mathcal T}\|_{\rm HS}^2
 \leq C(\rho_nt)^2\bigl[Q_X+(tf(t)\mu_X)^2\bigr]+\mathfrak e_n .
\end{align}
$\mathcal T^{\rm i}$ is a centered i.i.d.\ sum, so by
Lemma~\ref{lem-w-cancel}(a),(b) with $d = 1$,
\begin{align*}
E\|\mathcal T^{\rm i}\|_{\rm HS}^2\leq nE[\chi\bar\delta^2\|a\|^4]\leq Cnm_X^2\sigma_t^{-2}t^4f(t)h_t,
\end{align*}
which times $n\rho_n^4$ is $C(\rho_nt)^2Q_X$.  Since
$\langle u\otimes u,v\otimes v\rangle_{\rm HS} = \langle u,v\rangle^2$,
\begin{align*}
E\|\mathcal T^{\rm f}\|_{\rm HS}^2 = \sum_iE\bigl[\chi_i(e^\delta_i)^2\|a_i\|^4\bigr]+\sum_{i\neq k}E\bigl[\chi_i\chi_ke^\delta_ie^\delta_kC(\xi_i,\xi_k)^2\bigr];
\end{align*}
the diagonal is bounded like $\mathcal T^{\rm i}$, because
$E[(e^\delta_i)^2\mid\xi_i]\leq E[\delta_i^2\mid\xi_i]\leq Cm_X^2\sigma_t^{-2}\omega^{[2]}_X(\xi_i)^2+\mathfrak e_n$.
For the off-diagonal sum, $\delta_i$ is multilinear in the $p_{il}$,
$l\notin\{1,2,m,i\}$, so Lemma~\ref{lem-multilinear} gives
\begin{align*}
e^\delta_i = \sum_{S\neq\emptyset}\rho_n^{|S|}\prod_{j\in S}a(\xi_i,\xi_j)\tilde{\mathsf K}_S(\xi_i),
\qquad \tilde{\mathsf K}_S = \pm\mathsf K^{(|S|+1)}_{X,n-4-|S|}(\rho_n\check{\mathfrak D}(\cdot)),
\end{align*}
and by Lemma~\ref{lem-w-cancel}(a),(b), for $j = 1,2$,
\begin{align*}
|\tilde{\mathsf K}_S|\leq C\sigma_t^{-(|S|+2)}\mu_X\omega^{[|S|+2]}_X+\mathfrak e_n,\qquad
\int|\tilde{\mathsf K}_S|^j\leq C\sigma_t^{-j(|S|+2)}\mu_X^jf(t)h_t .
\end{align*}  The
off-diagonal sum is thus the fluctuation sum of Lemma~\ref{lem-w-Z2} with
$C(\xi_2,\xi_3)$ replaced by $C(\xi_2,\xi_3)^2$ and kernels of one higher
order, with the same truncation at $d_0$, $n^2$ pairs $(i,k)$ and the outer
factor $n\rho_n^4$.  Using ($\Gamma$1), ($\Gamma$4) and the integral bounds
above, the three classes contribute at most
\begin{align*}
(\alpha)\ &\sum_{d\geq1}\frac{n^{d+3}\rho_n^{2d+4}}{d!}(Ct)^{2d+4}\bigl(\sigma_t^{-(d+2)}\mu_Xf(t)h_t\bigr)^2
 \leq C(\rho_nt)^2(tf(t)\mu_X)^2,\\
(\beta)\ &\sum_{d\geq1}n^{d+3}\rho_n^{2d+5}(Ct)^{2d+5}\sigma_t^{-(2d+5)}\mu_X^2(f(t)h_t)^{3/2}
 \leq C(\rho_nt)^2(tf(t)\mu_X)^2\Bigl(\frac{\rho_n^3t}{nf(t)^2}\Bigr)^{1/4},\\
(\gamma)\ &\sum_{d\geq0}n^{d+3}\rho_n^{2d+6}(Ct)^{2d+6}\sigma_t^{-2(d+3)}\mu_X^2f(t)h_t
 \leq C(\rho_nt)^3Q_X .
\end{align*}
In ($\beta$) the $x$-integral is bounded by Cauchy--Schwarz,
$\int|a(x,y)||\tilde{\mathsf K}(x)|dx\leq Ct\|\tilde{\mathsf K}\|_{L^2}$ for $y$ in the window, and
the $y$-integral by $\int|\tilde{\mathsf K}|$.  Moreover, by \ref{condition-d3}(iii),
\ref{condition-s3}(iii) and $\rho_n\lambda_n^2\to0$,
\begin{align*}
\frac{\rho_n^3t}{nf(t)^2}\leq C\rho_n^2\lambda_n\cdot n^{-1}\rho_n\lambda_n^{2\gamma_+}\to0 .
\end{align*}
This proves \eqref{eq-w-T-HS}.

\emph{Second term.}  $\tilde\varphi = EX^{(-m)}$ is the H\'ajek kernel of
Lemma~\ref{lem-w-kernels} with $n-4$ trials in place of $n''$, so
Proposition~\ref{prop-w-hajek} (whose proof allows $|M-n''|\leq d_1$) gives
$\|\tilde\varphi\|^2\leq C\Xi_X+\mathfrak e_n$.  As $\tilde\varphi$ is
deterministic and $\|\cdot\|_{\rm op}\leq\|\cdot\|_{\rm HS}$, \eqref{eq-w-T-HS}
gives
\begin{align*}
12n\rho_n^4E\|(\mathcal T_m-\bar{\mathcal T})\tilde\varphi\|^2
 \leq C(\rho_nt)^2\Xi_XQ_X+CP_X+\mathfrak e_n,
\qquad P_X := (\rho_nt)^2\bigl(tf(t)\mu_X\bigr)^2\Xi_X .
\end{align*}

\emph{Third term.}  Let $N := \sum_{i\geq3,i\neq m}\chi_i\omega^{[4]}_X(\xi_i)$, a
sum of i.i.d.\ variables in $[0,C]$ with $EN\leq C_2nf(t)h_t$.  Since
$nf(t)h_t\geq n^{c_1}$ by \eqref{eq-w-poly}, Bernstein's inequality shows that
$\mathscr G := \{N\leq2C_2nf(t)h_t\}$ has negligible complement.  By
\eqref{eq-w-g-delta}, $|\delta_i|\|a_i\|^2\leq Cm_X\sigma_t^{-1}t^2\omega^{[4]}_X(\xi_i)+\mathfrak e_n$,
and by Lemma~\ref{lem-w-cancel}(a),(b),
$\|\bar{\mathcal T}\|_{\rm op}\leq Cm_X\sigma_t^{-1}t^2nf(t)h_t$.  Hence on
$\mathscr G$,
\begin{align*}
n\rho_n^4\|\mathcal T_m-\bar{\mathcal T}\|_{\rm op}^2
 \leq Cn\rho_n^4\bigl(m_X\sigma_t^{-1}t^2nf(t)h_t\bigr)^2
 = C\rho_nt\,(tf(t)\mu_X)^2 =: \tau_X\leq C\rho_nt ,
\end{align*}
using $tf(t)\leq C$; off $\mathscr G$, $\|\mathcal T_m\|_{\rm op}\leq n\rho_n^{-2}$
and $\|\tilde X\|\leq2n$, so that part is negligible.  With
Lemma~\ref{lem-w-Z2} for $E\|\tilde X\|^2$,
\begin{align*}
12n\rho_n^4E\|(\mathcal T_m-\bar{\mathcal T})\tilde X\|^2
 \leq12\tau_XE\|\tilde X\|^2+\mathfrak e_n
 \leq C\rho_n\lambda_n(\rho_n\lambda_n+\bar s_n)\Delta_X+CP_X+\mathfrak e_n .
\end{align*}

\emph{The product $P_X$.}  Both of the last two terms contain $P_X$.  By
\eqref{eq-w-chi},
\begin{align*}
\Xi_X\leq Cf(t)^2t^2\mu_X^2\Bigl(\frac{\bar\ell_X}t\Bigr)^{2\vartheta},
\end{align*}
and $\Delta_X\geq cf(t)\ell_A$ for a short cell.  In the small regime
($\mu_A = \ell_A/h_t$, $\bar\ell_A = h_t$), using $\ell_A\leq h_t$ and
$4-4\vartheta'>0$,
\begin{align*}
\frac{P_A}{\Delta_A^{4\vartheta'}}\leq C\rho_n^2t^{6-2\vartheta}f(t)^{4-4\vartheta'}\ell_A^{4-4\vartheta'}h_t^{2\vartheta-4}
 \leq C\rho_n^2t^{6-2\vartheta}f(t)^{4-4\vartheta'}h_t^{-(4\vartheta'-2\vartheta)};
\end{align*}
in the large regime ($\mu_A = 1$, $\bar\ell_A = \ell_A\geq h_t$) the ratio is
$C\rho_n^2t^{6-2\vartheta}f(t)^{4-4\vartheta'}\ell_A^{-(4\vartheta'-2\vartheta)}$, which
obeys the same bound because $4\vartheta'-2\vartheta>0$.  Substituting
$h_t = (t/(n''\rho_n))^{1/2}$,
\begin{align}
\label{eq-w-PX}
\frac{P_A}{\Delta_A^{4\vartheta'}}\leq C\,n^{2\vartheta'-\vartheta}\rho_n^{2+2\vartheta'-\vartheta}\;
 t^{6-\vartheta-2\vartheta'}f(t)^{4-4\vartheta'} ,
\end{align}
which is bounded uniformly in $t\in[t_0,\lambda_n]$: constraint (iv) of
\eqref{eq-w-theta-prime} bounds the weight through \ref{condition-d3}(iii),
and constraint (v) gives
\begin{align*}
n^{2\vartheta'-\vartheta}\rho_n^{2+2\vartheta'-\vartheta} = n^{(2\vartheta'-\vartheta)(1-b)-2b}\leq1 .
\end{align*}
This is the only place where \ref{condition-s3}(v) is used.  For one-sided
cells, $P_{O_s} = (\rho_ns)^2(sf(s))^4\leq C\Delta_{O_s}^{4\vartheta'}$ directly.
The remaining terms are small:
$(\rho_nt)^2\Xi_XQ_X\leq C(\rho_n\lambda_n)^2\Delta_X^{1+2\vartheta}\leq C\Delta_X^{4\vartheta'}$
by $1+2\vartheta\geq4\vartheta'$, and
$\rho_n\lambda_n(\rho_n\lambda_n+\bar s_n)\Delta_X\leq C\Delta_X^{4\vartheta'}$
by \eqref{eq-w-absorb}.

\emph{Quartic part.}  First,
$E\|Y(\xi_m)\|^4\leq\rho_n^4E[\chi\|a(\xi_m,\cdot)\|^4g_m^4]\leq C\rho_n^4t^4m_X^4f(t)h_t$.
Next,
\begin{align*}
\|V(y)\|\leq C\rho_n^2tm_X\sigma_t^{-1}S(y),\qquad
S(y) := \sum_i\chi_i\check w(\xi_i,y)\omega^{[3]}_X(\xi_i).
\end{align*}
Given $y$, $S(y)$ is a sum of i.i.d.\ nonnegative terms $Y_i$, so by
Rosenthal's inequality,
$E[S(y)^4\mid y]\leq C[(\sum_iE[Y_i\mid y])^4+\sum_iE[Y_i^4\mid y]]$.  By
Minkowski's inequality in $L^4(dy)$ and \ref{condition-l3},
\begin{align*}
\Bigl\|\sum_iE[Y_i\mid y]\Bigr\|_{L^4(dy)}\leq n\int\|w(x,\cdot)\|_{L^4}\omega^{[3]}_X(x)dx\leq Cnt^{q_4}f(t)h_t,
\end{align*}
and since $\check w\leq\rho_n^{-1}$,
$E\sum_iY_i^4\leq n\rho_n^{-2}\iint\check w(x,y)^2\omega^{[3]}_X(x)^4dxdy\leq Cn\rho_n^{-2}t^2f(t)h_t$.
Hence, with $\sigma_t^2 = n''\rho_nt$, $\mu_X^4f(t)h_t\leq Q_X\leq C\Delta_X$ and
$t^{1+4q_4}f(t)^2\leq Ct$ (by $\gamma_-\geq2q_4$),
\begin{align*}
\sum_m4E\|R_m\|^4
 &\leq C\Bigl(n\rho_n^4t^4m_X^4f(t)h_t+n\rho_n^8t^4m_X^4\sigma_t^{-4}\bigl[n^4t^{4q_4}(f(t)h_t)^4+n\rho_n^{-2}t^2f(t)h_t\bigr]\Bigr)
\\ &\leq C\rho_n\lambda_n\bigl(n^{-1}+\rho_n^2\bigr)\Delta_X .
\end{align*}

\emph{Conclusion.}  Collecting the direct part, the three indirect terms and
the quartic part,
\begin{align*}
\mathrm{Var}(\Psi^*)\leq C\rho_n\lambda_n(\Delta_X+\Xi_X)^2+C\rho_n\lambda_n\Delta_X+C\Delta_X^{4\vartheta'}+\mathfrak e_n\leq C\Delta_X^{4\vartheta'},
\end{align*}
because $(\Delta_X+\Xi_X)^2\leq C\Delta_X^{4\vartheta}\leq C\Delta_X^{4\vartheta'}$ and
$\rho_n\lambda_n\Delta_X\leq\Delta_X^{4\vartheta'}$ by \eqref{eq-w-absorb}.
Therefore $\|V\|_2\leq EV+\sum_k\mathrm{Var}(V_k)^{1/2}\leq C\Delta_X^{2\vartheta'}$,
and $\|Z_X\|_{L^4}\leq C\Delta_X^{\vartheta'}$.
\end{proof}

\paragraph{Step A: assembly.}
\emph{Analogous.}  The single-cell bound combines Lemmas~\ref{lem-w-edge},
\ref{lem-w-Z2} and \ref{lem-w-Z4} exactly as Lemma~\ref{lem-c-single}, and the
product condition of Proposition~\ref{prop-grid} follows by Cauchy--Schwarz.
\emph{Additional care.}
\begin{enumerate}[(1)]
\item The term $n^{-1}\Xi_X$ in the fourth moment of the edge noise must be
absorbed into $\Delta_X^{4\vartheta'}$; on a compact set this used
$h_t\asymp(n\rho_n)^{-1/2}$ only, while here the weights produce a power of
$t$ and $f(t)$ that constraint (ii) of \eqref{eq-w-theta-prime} keeps bounded
(Lemma~\ref{lem-w-absorb}).
\item Long cells are assembled from one-sided cells through
\eqref{eq-w-long-split} (Corollary~\ref{cor-w-long}).
\item The proof of the theorem must verify condition~(iii) of
Proposition~\ref{prop-grid}, which is where \ref{condition-s3}(i) enters, and
continuity of the limit at the right endpoint $u = 1$ of the reparametrized
index set.
\end{enumerate}

\begin{lemma}[Absorption]
\label{lem-w-absorb}
For $n\geq n_0$ and every $X\in\mathscr C_n$,
$n^{-1}\Delta_X\leq\Delta_X^{4\vartheta'}$ and $n^{-1}\Xi_X\leq C\Delta_X^{4\vartheta'}$.
\end{lemma}

\begin{proof}
The first is \eqref{eq-w-absorb}.  For the second, with scale $s$: for $O_s$,
$n^{-1}\Xi_{O_s}\leq\Delta_{O_s}^2\leq\Delta_{O_s}^{4\vartheta'}$.  For a short cell,
by \eqref{eq-w-chi}, $\Xi_A\leq Cf(s)^2s^{2-2\vartheta}\mu_A^2\bar\ell_A^{2\vartheta}$, and
$\Delta_A\geq cf(s)\ell_A$.  In the small regime, using
$\mu_A^2 = (\ell_A/h_s)^2\leq(\ell_A/h_s)^{4\vartheta'}$ (as $4\vartheta'\leq2$),
\begin{align*}
\frac{n^{-1}\Xi_A}{\Delta_A^{4\vartheta'}}
 \leq Cn^{-1}f(s)^{2-4\vartheta'}s^{2-2\vartheta}h_s^{2\vartheta-4\vartheta'}
 = Cn^{-1}(n''\rho_n)^{2\vartheta'-\vartheta}f(s)^{2-4\vartheta'}s^{2-\vartheta-2\vartheta'},
\end{align*}
where $n^{-1}(n\rho_n)^{2\vartheta'-\vartheta}\leq n^{-1/2}$ since $2\vartheta'-\vartheta\leq1/2$,
and
\begin{align*}
f(s)^{2-4\vartheta'}s^{2-\vartheta-2\vartheta'}\leq Cs^{2-\vartheta-2\vartheta'-\gamma_-(2-4\vartheta')}\leq C
\end{align*}
by \ref{condition-d3}(iii) and constraint (ii) of \eqref{eq-w-theta-prime}.
In the large regime $\mu_A = 1$ and $\bar\ell_A = \ell_A\geq h_s$; since
$2\vartheta-4\vartheta'<0$, $\ell_A^{2\vartheta-4\vartheta'}\leq h_s^{2\vartheta-4\vartheta'}$ and
the same bound results.
\end{proof}

\begin{lemma}[Single-cell bound]
\label{lem-w-single}
$\|\Delta'_X\|_{L^4}\leq C\Delta_X^{\vartheta'}$ for every $X\in\mathscr C_n$ and $n\geq n_0$.
\end{lemma}

\begin{proof}
$\|\Delta'_X\|_4\leq\|B_X\|_4+\|Z_X\|_4$.  By Lemmas~\ref{lem-w-edge},
\ref{lem-w-Z2} and \ref{lem-w-absorb},
\begin{align*}
EB_X^4\leq C\bigl(\Delta_X^2+n^{-1}\Delta_X+n^{-1}\Xi_X\bigr)+\mathfrak e_n\leq C\Delta_X^{4\vartheta'},
\end{align*}
since $\Delta_X\leq1$, $4\vartheta'\leq2$ and $\mathfrak e_n\leq\Delta_X^{4\vartheta'}$;
and $\|Z_X\|_4\leq C\Delta_X^{\vartheta'}$ by Lemma~\ref{lem-w-Z4}.
\end{proof}

\begin{corollary}[Long cells]
\label{cor-w-long}
For every long cell $A = A_{rt}$ with $t_0\leq r<t\leq\lambda_n$ and
$n\geq n_0$, $\|\Delta'_A\|_{L^4}\leq C\Delta_A^{\vartheta'}$.
\end{corollary}

\begin{proof}
By \eqref{eq-w-long-split}, Minkowski's inequality and Lemma~\ref{lem-w-single},
\begin{align*}
\|\Delta'_A\|_4\leq\|\Delta'_{O_r}\|_4+\|\Delta'_{O_t}\|_4\leq C\bigl(\Delta_{O_r}^{\vartheta'}+\Delta_{O_t}^{\vartheta'}\bigr),
\end{align*}
and $\Delta_{O_s}\leq C\Delta_A$ for $s\in\{r,t\}$ by Lemma~\ref{lem-w-cell}(b).
\end{proof}

\begin{proposition}[Fourth-moment increment bound, weighted]
\label{prop-w-fourth}
There is $C<\infty$ such that for all lattice points
$t_0\leq r<s<t\leq\lambda_n$ and $n\geq n_0$,
\begin{align}
\label{eq-w-fourth}
E\Bigl[\bigl(\mathbb G'_n(s)-\mathbb G'_n(r)\bigr)^2\bigl(\mathbb G'_n(t)-\mathbb G'_n(s)\bigr)^2\Bigr]
 \leq C\bigl(F(t)-F(r)\bigr)^{4\vartheta'},\qquad4\vartheta'>1 .
\end{align}
\end{proposition}

\begin{proof}
With $A := A_{rs}$ and $B := A_{st}$, each short or long, Cauchy--Schwarz,
Lemma~\ref{lem-w-single} and Corollary~\ref{cor-w-long} give
\begin{align*}
E\bigl[(\Delta'_A)^2(\Delta'_B)^2\bigr]\leq\|\Delta'_A\|_4^2\|\Delta'_B\|_4^2\leq C(\Delta_A\Delta_B)^{2\vartheta'}\leq C(\Delta_A+\Delta_B)^{4\vartheta'}.
\end{align*}
\end{proof}

\begin{remark}[Why single cells suffice]
\label{rem-w-single}
For the i.i.d.\ empirical process the product form of
Proposition~\ref{prop-grid}(ii) is essential, because the fourth moment of a
single increment contains a term $n^{-1}\Delta$, which is not bounded by
$\Delta^{2\nu}$, $2\nu>1$, when $\Delta\ll n^{-1}$.  Here such cells do not
occur: on the lattice and below the cutoff $\Delta_X\geq n^{-K_0}$, and
constraint (i) of \eqref{eq-w-theta-prime} gives $n^{-1}\Delta_X\leq\Delta_X^{4\vartheta'}$.
This is what allows the bounds to be stated for single cells, which removes
all cross-cell quantities and makes long cells routine through
\eqref{eq-w-long-split}.
\end{remark}

\begin{proof}[Proof of Theorem~\ref{thm-tightness-ring}]
By Lemma~\ref{lem-c-hub-reduction} (Step~H) it suffices to prove
$\mathbb G'_n\rightsquigarrow\mathbb G$ in $\ell^\infty([t_0,\infty))$.  We verify
Proposition~\ref{prop-grid} for $\mathbb H_n := \mathbb G'_n$ with
$T_n := \lambda_n$ and $\nu := 2\vartheta'>1/2$ (Step~L).

\emph{Condition (i).}  $\mathbb G'_n\stackrel{fdd}{\rightsquigarrow}\mathbb G$ by
Lemma~\ref{lem-c-hub-reduction}.  For continuity of the limit, let
$t_0\leq s<t<\infty$.  For $n$ large $t\leq\lambda_n$, and
$\mathbb G'_n(t)-\mathbb G'_n(s) = \Delta'_{A_n}$ with $A_n$ the cell between
the lattice points $\lambda_{s'}$ and $\lambda_{t'}$ (or $0$ if $s' = t'$).
By Lemma~\ref{lem-w-single} and Corollary~\ref{cor-w-long},
$E|\Delta'_{A_n}|^4\leq C\Delta_{A_n}^{4\vartheta'}$, and $\Delta_{A_n}\to F(t)-F(s)$ by
continuity of $F$.  Since $\Delta'_{A_n}\rightsquigarrow\mathbb G(t)-\mathbb G(s)$,
Fatou's lemma (via the Skorokhod representation) gives
\begin{align*}
E|\mathbb G(t)-\mathbb G(s)|^4\leq C\bigl(F(t)-F(s)\bigr)^{4\vartheta'}.
\end{align*}
In the parametrization $u = F(t)$, the Kolmogorov--Chentsov theorem gives a
version of $\mathbb G$ that is H\"older continuous on $[F(t_0),1)$ with a
single random constant, hence uniformly continuous there; it extends
continuously to $u = 1$, with $\mathbb G(\infty) := \lim_{t\to\infty}\mathbb G(t)$.

\emph{Condition (ii)} is Proposition~\ref{prop-w-fourth}, and the mesh
condition is Lemma~\ref{lem-w-cell}(c).

\emph{Condition (iii).}  Since $\hat F$ and $F^\star$ are nondecreasing,
\begin{align*}
\sup_{t\geq\lambda_n}\bigl|\mathring{\mathbb G}_n(t)-\mathring{\mathbb G}_n(\lambda_n)\bigr|
\leq\sqrt n\bigl(1-\hat F(\lambda_n)\bigr)+\sqrt n\bigl(1-F^\star(\lambda_n)\bigr),
\end{align*}
whose expectation is
\begin{align*}
2\sqrt nP(D_1>\lambda'_n)\leq2\sqrt nP(\mathfrak D>\lambda_n/2)+\mathfrak e_n\to0
\end{align*}
by the binomial Chernoff bound and \ref{condition-s3}(i); Markov's
inequality gives (iii) for $\mathring{\mathbb G}_n$, and
Lemma~\ref{lem-c-hub-reduction} transfers it to $\mathbb G'_n$.

\emph{Endpoints.}  In the proof of Proposition~\ref{prop-grid} the
right-endpoint condition of \citet[Theorem~13.3]{billingsley1999convergence}
at $u = 1$ reads $\mathbb H_n(1)-\mathbb H_n(1-\delta)\to0$ in probability as
$\delta\downarrow0$, uniformly in $n$.  Here
$\mathbb H_n(1) = \mathbb G'_n(\infty) = n^{-1/2}\sum_i(\chi_i-E\chi_1)$, and
\begin{align*}
\bigl|\mathbb H_n(1)-\mathbb H_n(F(\lambda_n))\bigr|\leq\sup_{t\geq\lambda_n}\bigl|\mathbb G'_n(t)-\mathbb G'_n(\lambda_n)\bigr| = o_P(1)
\end{align*}
by (iii), while the increments below $F(\lambda_n)$ are controlled by (ii) and
Chebyshev's inequality.  The left-endpoint condition at $u = F(t_0)$ follows
from (ii) in the same way.  Hence $\mathbb H_n\rightsquigarrow\mathbb H$ in
$D[F(t_0),1]$, and as the limit is continuous, in $\ell^\infty([t_0,\infty))$.
\end{proof}

\section{Proof of Lemma~\ref{lem-feasibility}}
\label{app-feasibility}

Throughout, $f\asymp\lambda^{-\gamma}$ on $J_\infty$, so that
\ref{condition-d3}(iii) holds with $\gamma_+ = \gamma_- = \gamma$, and, integrating
the two bounds on $f$,
\begin{align}
\label{eq-feas-tail}
c\,\lambda^{1-\gamma}\leq1-F(\lambda)\leq C\lambda^{1-\gamma}\qquad(\lambda\geq t_0^-).
\end{align}
With $\rho_n = n^{-b}$ and $\lambda_n = n^a$, every condition in
\ref{condition-s1} and \ref{condition-s3} asks that a power of $n$ tend to zero,
that is, that an exponent, affine in $(a,b)$, be negative.  We compute the
exponents one at a time.

\emph{\ref{condition-s1}.}  Since $b<1$, $n\rho_n/\log n = n^{1-b}/\log n\to\infty$,
and $n\rho_n^{2+\beta} = n^{1-b(2+\beta)}\to0$ if and only if $b>1/(2+\beta)$.

\emph{\ref{condition-s3}(i).}  By \eqref{eq-feas-tail},
$\sqrt n\,P(\mathfrak D>\lambda_n/2)\asymp n^{1/2-a(\gamma-1)}$, which tends to
zero if and only if $a>a_0 = 1/(2(\gamma-1))$.

\emph{\ref{condition-s3}(ii).}  By \eqref{eq-chi-explicit} and
$t\,h_t^{-\vartheta} = t^{1-\vartheta/2}(n''\rho_n)^{\vartheta/2}$,
\begin{align*}
\Upsilon_n(t)\,t\,h_t^{-\vartheta}\leq C\,t^{\,5/2+\beta-\vartheta/2}\rho_n^{1+\beta}(n\rho_n)^{(1+\vartheta)/2}.
\end{align*}
The exponent of $t$ is positive, so the supremum over $t\leq\lambda_n$ is
attained at $t = \lambda_n$, where the bound equals $Cn^{E(a,b)}$ with
\begin{align*}
E(a,b) := a\Bigl(\frac52+\beta-\frac\vartheta2\Bigr)+\frac{1+\vartheta}2
 -b\Bigl(1+\beta+\frac{1+\vartheta}2\Bigr).
\end{align*}
Thus the bound tends to zero if and only if $b>b(a)$, where
\begin{align}
\label{eq-feas-ba}
b(a) := \frac{\tfrac{1+\vartheta}2+a\bigl(\tfrac52+\beta-\tfrac\vartheta2\bigr)}
 {1+\beta+\tfrac{1+\vartheta}2}
\end{align}
is continuous and strictly increasing in $a$, with $b(a_0) = b_*$.

\emph{\ref{condition-s3}(iii) and $a<b/2$.}  $n^{-1}\rho_n\lambda_n^{2\gamma} = n^{-1-b+2a\gamma}$,
which tends to zero if and only if $a<(1+b)/(2\gamma)$.  At $a = a_0$ this reads
$\gamma<(1+b)(\gamma-1)$, that is $b>1/(\gamma-1)$; and $a_0<b/2$ also reads
$b>1/(\gamma-1)$.  Both hold whenever $b>1/(2+\beta)$, because
$\gamma-1>2+\beta$.

\emph{\ref{condition-s3}(iv).}  With $f(t)\asymp t^{-\gamma}$,
\begin{align*}
v_n(t)\asymp n^{(1-2\vartheta)/4}\rho_n^{(9-2\vartheta)/4}\,t^{\,e},\qquad
e := \frac{11-2\vartheta}4-\frac\gamma2 .
\end{align*}
If $e\leq0$, the supremum is attained at $t = t_0$, and (iv) reads
$1-2\vartheta<b(9-2\vartheta)$, which is implied by \ref{condition-s3}(v) since
$(1-2\vartheta)/(9-2\vartheta)\leq(1-2\vartheta)/(5-2\vartheta)$.  If $e>0$, the
supremum is attained at $t = \lambda_n$, and at $a = a_0$ (iv) reads
\begin{align*}
b>b_4 := \frac{(1-2\vartheta)+4a_0e}{9-2\vartheta},\qquad
4a_0e = \frac{11-2\vartheta-2\gamma}{2(\gamma-1)} .
\end{align*}
We claim $b_4<1/(2+\beta)$, so that (iv) is implied by \ref{condition-s1}.
Indeed, $\gamma>3$ gives $4a_0e<(5-2\vartheta)/4$, whence
\begin{align*}
b_4<\frac{(1-2\vartheta)+(5-2\vartheta)/4}{9-2\vartheta} = \frac{9-10\vartheta}{4(9-2\vartheta)}
<\frac2{9-2\vartheta},
\end{align*}
the last inequality because $\vartheta>1/10$.  On the other hand $e>0$ means
$\gamma<(11-2\vartheta)/2$, and with $\gamma>3+\beta$ this gives
$\beta<(5-2\vartheta)/2$, that is $1/(2+\beta)>2/(9-2\vartheta)$.

\emph{\ref{condition-s3}(v)} is the condition $b>(1-2\vartheta)/(5-2\vartheta)$ itself.

\begin{proof}[Proof of Lemma~\ref{lem-feasibility}]
Let $b\in(b_{\min},1)$.  By the computations above, \ref{condition-s1} and
\ref{condition-s3}(v) hold, and at $a = a_0$ the conditions (ii), (iii), (iv) and
$a<b/2$ hold with strict inequalities between the relevant exponents: (ii)
because $b>b_* = b(a_0)$, (iii) and $a<b/2$ because $b>1/(2+\beta)>1/(\gamma-1)$,
and (iv) because $b>1/(2+\beta)>b_4$ (or by (v) if $e\leq0$).  The exponents
are continuous in $a$, so these strict inequalities persist for
$a\in[a_0,a_0+\eta)$ for some $\eta>0$, and for $a\in(a_0,a_0+\eta)$ condition
(i) holds as well.  Hence \ref{condition-s3} holds with $\lambda_n = n^a$ for
every such $a$.

Conversely, let $b\leq b_*$.  By the lower bound in \eqref{eq-feas-tail},
\ref{condition-s3}(i) fails unless $a>a_0$.  For $a>a_0$, $b(a)>b(a_0) = b_*\geq b$,
so $E(a,b)>0$ and the bound on \ref{condition-s3}(ii) supplied by
\eqref{eq-chi-explicit} diverges.
\end{proof}

\begin{remark}[Which condition binds]
\label{rem-feas-binding}
The proof shows that \ref{condition-s3}(iii) and (iv) are never binding for
power-law densities: they are implied by \ref{condition-s1} once the cutoff is
taken just above $a_0$.  The binding constraint is the tension between the
cutoff, which must be large enough for the tail beyond it to be negligible
((i), $a>a_0$), and the clipping ratio, which must be small up to the cutoff
((ii), $b>b(a)$).  For example, with $\vartheta = 1/2$, $\beta = 0.1$ and
$\gamma = 3.2$, $a_0 = 0.227$ and $b_* = 0.694$, while $1/(2+\beta) = 0.476$.
\end{remark}

\section{Proof of Proposition~\ref{prop-centering-tail}}
\label{app-centering-tail}

The proof follows that of Lemma~\ref{lem-centering}, with the H\"older
condition replaced by its relative form \ref{condition-d3}(ii), and with the
range $t>\lambda_n$ treated separately.  Let
\begin{align*}
c_* := \min\Bigl\{\frac12,\ 1-\frac{t_0^-}{t_0}\Bigr\},
\qquad\text{so that}\qquad (1-c_*)t\geq t_0^-\quad\text{for all }t\geq t_0 .
\end{align*}
By \ref{condition-d3}(i) and (iii), $f$ is bounded on $J_\infty$, and since
$\gamma_->3$,
\begin{align}
\label{eq-sd-weights}
\bar f := \sup_{J_\infty}f<\infty,\qquad
\sup_{t\geq t_0}t\,f(t)<\infty,\qquad
\sup_{t\geq t_0}f(t)\,t^{(1-\theta_f)/2}<\infty .
\end{align}
Write $L_n := (4\rho_n)^{-1}$; since $\rho_n\lambda_n\to0$, $\lambda_n\leq L_n$
for $n$ large.

\emph{Step 0: beyond the cutoff.}  Since $F^\star$ and $F$ are nondecreasing and
bounded by one, for $t\geq\lambda_n$,
\begin{align*}
\bigl|F^\star(t)-F(t)\bigr|\leq\bigl(1-F^\star(\lambda_n)\bigr)+\bigl(1-F(\lambda_n)\bigr).
\end{align*}
First, $1-F(\lambda_n)\leq P(\mathfrak D>\lambda_n/2)$.  Second,
$1-F^\star(\lambda_n) = P(D_1>\lambda_n')$.  On $\{\check{\mathfrak D}(\xi_1)\leq\lambda_n/2\}$
the conditional mean of $D_1$ is at most $n'\rho_n\lambda_n/2$, while
$\lambda_n'\geq n'\rho_n\lambda_n-1$, so the multiplicative Chernoff bound gives
$P(D_1>\lambda_n'\mid\xi_1)\leq e^{-cn\rho_n\lambda_n}$ there; and
$\check{\mathfrak D}\leq\mathfrak D$.  Hence
\begin{align*}
\sqrt n\sup_{t\geq\lambda_n}\bigl|F^\star(t)-F(t)\bigr|
 \leq2\sqrt n\,P(\mathfrak D>\lambda_n/2)+\sqrt n\,e^{-cn\rho_nt_0}\;\rightarrow\;0
\end{align*}
by \ref{condition-s3}(i) and $n\rho_n/\log n\to\infty$.  It remains to treat
$t\in[t_0,\lambda_n]$.

\emph{Step 1: reduction to the unclipped binomial.}  The two bounds of the
proof of Lemma~\ref{lem-centering} that remove the clipping on the good set
$\mathscr G_n$, \eqref{eq-bias-clip}, and the complement of the good set do
not depend on $t$.  They give, uniformly in $t\in[t_0,\lambda_n]$,
\begin{align*}
\sqrt n\bigl[F^\star(t)-F(t)\bigr]
 = \sqrt n\int_0^{L_n}K_t(\lambda)\,f(\lambda)\,d\lambda
 + O\bigl(n\rho_n^{\,3/2+\beta}\bigr) + o(1),
\end{align*}
with $K_t(\lambda) := \mathcal B_{n'}(t';\rho_n\lambda)-\mathbbm 1(\lambda\leq t)$
as in \eqref{eq-bias-kernel}; here we used $t\leq\lambda_n\leq L_n$.  Let
$U_t := [(1-c_*)t,(1+c_*)t]\subseteq J_\infty$, and split
\begin{align*}
\int_0^{L_n}K_tf\,d\lambda
 = \underbrace{f(t)\int_0^{L_n}K_t\,d\lambda}_{\mathrm{(A)}}
 +\underbrace{\int_{U_t}K_t\bigl[f-f(t)\bigr]d\lambda}_{\mathrm{(B)}}
 +\underbrace{\int_{[0,L_n]\setminus U_t}K_t\bigl[f-f(t)\bigr]d\lambda}_{\mathrm{(C)}} .
\end{align*}

\emph{Step 2: the lattice term (A).}  The computation leading to
\eqref{eq-bias-lattice} uses only $t\leq L_n$ and $t'\leq n'p/2$ for
$p\geq\frac14$, and the latter holds uniformly in $t\leq\lambda_n$ because
$t'\leq n'\rho_n\lambda_n = o(n')$.  It gives
\begin{align*}
\Bigl|\int_0^{L_n}K_t(\lambda)\,d\lambda\Bigr|
 \leq\frac{\max\{1,\rho_nt\}}{(n'+1)\rho_n}+O\bigl(\rho_n^{-1}e^{-n'/32}\bigr),
\qquad\text{so}\qquad
\sqrt n\,|\mathrm{(A)}|\leq\frac{C\bar f}{\sqrt n\,\rho_n}+o(1),
\end{align*}
uniformly in $t\in[t_0,\lambda_n]$.

\emph{Step 3: the smoothing bias (B).}  Let $X_\lambda\sim\mathrm{Bin}(n',\rho_n\lambda)$.
As in the proof of \eqref{eq-bernstein-K}, for $|\lambda-t|\leq c_*t\leq t/2$ the
variance of $X_\lambda$ is at most $\frac32n'\rho_nt$ and the relevant deviation
is at most $\frac12n'\rho_nt$, so Bernstein's inequality gives
\begin{align*}
|K_t(\lambda)|\leq C\exp\bigl(-\zeta(\lambda-t)^2/h_t^2\bigr),
\end{align*}
with $C,\zeta$ not depending on $t$.  By \ref{condition-d3}(ii) at the base
point $t$, $|f(\lambda)-f(t)|\leq K_ff(t)(|\lambda-t|/t)^{\theta_f}$ on $U_t$.
Hence, by \eqref{eq-sd-weights},
\begin{align*}
\sqrt n\,|\mathrm{(B)}|
 &\leq CK_f\sqrt n\,f(t)\,t^{-\theta_f}\int_{\mathbb R}|u|^{\theta_f}e^{-\zeta u^2/h_t^2}\,du
 = C'\sqrt n\,f(t)\,t^{-\theta_f}h_t^{1+\theta_f}
\\ &= C'\sqrt n\,f(t)\,t^{(1-\theta_f)/2}\,(n''\rho_n)^{-(1+\theta_f)/2}
 \leq C\sqrt n\,(n\rho_n)^{-(1+\theta_f)/2},
\end{align*}
uniformly in $t\in[t_0,\lambda_n]$.  This is the only place where the tail
regularity enters: the weight $f(t)t^{(1-\theta_f)/2}$ is bounded on
$[t_0,\infty)$, so the bound is no worse than at $t = t_0$.

\emph{Step 4: away from $t$ (C).}  For $\lambda<(1-c_*)t$, $K_t(\lambda) = -P(X_\lambda>t')$,
the deviation $t'-n'\rho_n\lambda$ is at least $n'\rho_nc_*t-1$ and the variance
is at most $n'\rho_nt$; for $\lambda>(1+c_*)t$, $K_t(\lambda) = P(X_\lambda\leq t')$,
the deviation $n'\rho_n\lambda-t'$ is at least $n'\rho_n\lambda\,c_*/(1+c_*)$ and
the variance is at most $n'\rho_n\lambda$.  In both cases Bernstein's
inequality gives $|K_t(\lambda)|\leq2e^{-cn\rho_nt_0}$.  Since
$\int_0^{L_n}|f(\lambda)-f(t)|\,d\lambda\leq1+L_n\bar f$,
\begin{align*}
\sqrt n\,|\mathrm{(C)}|\leq2\sqrt n\,e^{-cn\rho_nt_0}\bigl(1+\bar f/(4\rho_n)\bigr)\;\rightarrow\;0,
\end{align*}
uniformly in $t$, by $n\rho_n/\log n\to\infty$.  Note that (C) involves $f$
only through its integrability, so no regularity of $f$ below $t_0^-$ is
needed.

\emph{Conclusion.}  Collecting Steps 0--4,
\begin{align*}
\sup_{t\geq t_0}\sqrt n\,\bigl|F^\star(t)-F(t)\bigr|
 \leq C\Bigl[n\rho_n^{\,3/2+\beta}+\frac1{\sqrt n\,\rho_n}
 +\sqrt n\,(n\rho_n)^{-(1+\theta_f)/2}\Bigr]+o(1).
\end{align*}
As in the proof of Lemma~\ref{lem-centering}, the middle term is dominated by
the last, and the right side tends to zero under
\eqref{eq-centering-conditions} with $\theta = \theta_f$.
\qed

\end{document}